\documentclass[11pt]{article}
\usepackage{preambule}

\usepackage{notations}
\graphicspath{{./}}
\pgfplotsset{table/search path={figures},}
\input{./figures.tex}

\usepackage{xcolor}

\newboolean{afficherGraphes}
\setboolean{afficherGraphes}{false}
\title{Wave propagation in the frequency regime in one-dimensional quasiperiodic media - Limiting absorption principle}
\author{Pierre Amenoagbadji$^{(1)}$, Sonia Fliss$^{(2)}$, Patrick Joly$^{(2)}$\\
\small{
${(1)}$ APAM, Columbia University, United States;}\\
\small{${(2)}$ POEMS, CNRS, Inria, ENSTA, Institut Polytechnique de Paris,91120 Palaiseau, France.}\\
\small{Corresponding author: \texttt{sonia.fliss@ensta.fr}.}
}
\begin{document}
% % Titre du document
%
%\tableofcontents

\maketitle
\begin{abstract}
We study the one-dimensional Helmholtz equation with (possibly perturbed) quasiperiodic coefficients. Quasiperiodic functions are the restriction of higher dimensional periodic functions along a certain (irrational) direction. In classical settings, for real-valued frequencies, this equation is generally not well-posed: existence of solutions in $L^2$ is not guaranteed and uniqueness in $L^\infty$ may fail. This is a well-known difficulty of Helmholtz equations, but it has never been addressed in the quasiperiodic case. We tackle this issue by using the limiting absorption principle, which consists in adding some imaginary part (also called absorption) to the frequency in order to make the equation well-posed in $L^2$, and then defining the physically relevant solution by making the absorption tend to zero. In previous work, we introduced a definition of the solution of the equation with absorption based on Dirichlet-to-Neumann (DtN) boundary conditions. This approach offers two key advantages: it facilitates the limiting process and has a direct numerical counterpart. In this work, we first explain why the DtN boundary conditions have to be replaced by Robin-to-Robin boundary conditions to make the absorption go to zero. We then prove, under technical assumptions on the frequency, that the limiting absorption principle holds and we propose a numerical method to compute the physical solution.
\end{abstract}
\section{Presentation} \label{sec_position}

\subsection{The model problem and the assumptions on the coefficients}\label{sec_QD}
We want to study and solve numerically the Helmholtz equation
\begin{equation}
\label{eq:whole_line_problem}\tag{$\mathscr{P}$}
\displaystyle
- \frac{d}{d \xvi} \Big( \mu \, \frac{d u}{d \xvi} \Big) - \rho \, \omega^2 \, u = f \quad \textnormal{in} \quad \R,
\end{equation}
under the assumptions:
\begin{itemize}
	\item the coefficients $\mu$ and $\rho$ have positive upper and lower bounds:
	\begin{equation}\label{eq:coef_ellipt}
		\displaystyle
		\spexists \mu_\pm, \rho_\pm, \quad \spforall x \in \R, \qquad 0 < \mu_- \leq \mu(\xvi) \leq \mu_+ \quad \textnormal{and} \quad 0 < \rho_- \leq \rho(\xvi) \leq \rho_+,
	\end{equation}
	and either they coincide with continuous quasiperiodic functions of order $2$ (see below for the definition), 
    denoted respectively $\mu_\cutvec$ and $\rho_\cutvec$, or they are local perturbations of these quasiperiodic functions 
	\begin{equation}\label{eq:coef_support}
		\spexists a^\gauche < a^\droite, \quad \mu\equiv\mu_\cutvec \quad \text{and}\quad \rho\equiv\rho_\cutvec\quad \text{on}\;\R\setminus(a^\gauche,a^\droite);
	\end{equation}
	\item the source term $f$ belongs to $L^2(\R)$ and is assumed to have a compact support, which can be supposed to be in $(a^\gauche,a^\droite)$ without any loss of generality:
	\begin{equation}
		\label{eq:source_terme}
		 \supp f \subset (a^\gauche,a^\droite).
	\end{equation}
\end{itemize}
Here and in what follows, the superscript \textquote{$\gauche$}, \resp \textquote{$\droite$}, stands for \textquote{\emph{left}}, \resp \textquote{\emph{right}}.

\begin{figure}[H]
	\centering
	\includegraphics[page=1]{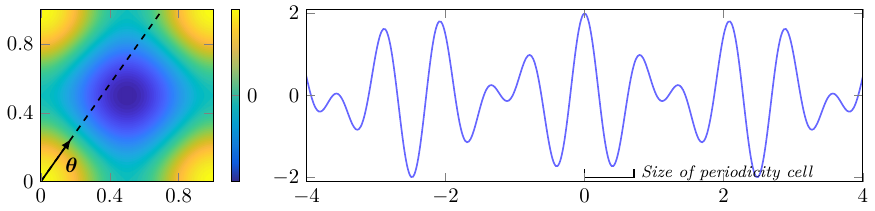}
	\caption{Function $F: (\yvi_1, \yvi_2) \mapsto \cos 2\pi \yvi_1 + \cos 2\pi \yvi_2$ in its periodicity cell $(0,1)^2$ (left), and whose trace along $\cutvec = (1, \sqrt{2})$ leads to a quasiperiodic function (right). \label{fig:example_quasiperiodic_function}}
\end{figure}

\vspspe
The functions $\mu_\cutvec$ and $\rho_\cutvec$ are supposed to be quasiperiodic of order $2$, meaning that there exist functions $\mu_p,\rho_p: \R^2 \to \R$ which are continuous and $1$--periodic with respect to each variable, and a vector $\cutvec\coloneqq (\cuti_1,\cuti_2) \in \R^2$ such that
\[
	\spforall x \in \R, \quad\mu_\cutvec=\mu_p(\xvi\, \cutvec) \quad \textnormal{and} \quad
	\rho_\cutvec=\rho_p(\xvi\, \cutvec).
\]
We refer to \cite{bohr2018almost,besicovitch1954almost,levitan1982almost} for complete presentations of the theory of quasiperiodic functions. A geometrical interpretation of this definition is to see the one-dimensional function $\mu_\theta$ (\resp $\rho_\theta$) as the trace of a $2$-dimensional function $\mu_p$ (\resp $\rho_p$) along the line passing through $(0, 0)$ and parallel to the vector $\cutvec$. This is illustrated in Figure \ref{fig:example_quasiperiodic_function} for $\cutvec = (1, \sqrt{2})$.

\vspspe
Finally, in order to simplify the presentation, we choose $a^\gauche$ and $a^\droite$ such that
\begin{equation}
	\label{eq:hyp_alar}
	a^\gauche\, \cuti_2 \in \Z \quad \textnormal{and} \quad a^\droite\, \cuti_2 \in \Z.
\end{equation}
The interest of this assumption lies in the fact that
\begin{equation}\label{eq:cell_lr}
	\mu_p|_{\vts (0, 1) \times (a^\droite\, \cuti_2, a^\droite\, \cuti_2 + 1)} = \mu_p|_{(0, 1)^2} = \mu_p|_{\vts (0, 1) \times (a^\gauche\, \cuti_2, a^\gauche\, \cuti_2 + 1)},
\end{equation}
and similarly for $\rho_p$.

\noindent
If the ratio $\delta \coloneqq  \cuti_1/\cuti_2$ is rational, then $\mu_\cutvec$ and $\rho_\cutvec$ are $q/\cuti_2$--periodic, with $q$ being the denominator of $\delta$. One can use for instance \cite{fliss2009analyse,fliss2021dirichlet} to solve the problem in this periodic context. In this paper, we shall assume the opposite, that is,
\begin{equation}\label{eq:delta}
	\delta\coloneqq \frac{\theta_1}{\theta_2} \in \R\setminus\Q,
\end{equation}
so that $\mu_\cutvec$ and $\rho_\cutvec$ are not periodic in general.

\begin{rmk}
{Note that the method described in this paper works also for periodic coefficients, \emph{i.e.} when $\delta\equiv p/q\in\Q$ (except Proposition \ref{thm_spectre}). This could be of particular interest if the period of the medium is very large, \emph{i.e.} when $q$ is very large. Indeed, the computational cost of the method developed in \cite{joly2006exact} increases with the period whereas the computational cost of the method described in this paper is the same whatever the period is.}
\end{rmk} \noindent
As $\delta$ is irrational, Kronecker's approximation theorem \cite[Theorem 442]{hardy1979introduction} ensures that the functions $\mu_p$ and $\rho_p$ are entirely determined by their restrictions on the line $\R\, \cutvec$. In particular, $\mu_p$ (\resp $\rho_p$) has the same lower and upper bounds as $\mu$ (\resp $\rho$).

\begin{thm}[Kronecker's theorem]\label{thm:kronecker}
	If $\cuti_1 / \cuti_2 \in \R \setminus \Q$, then the set $\cutvec\, \R + \Z^2$ is dense in $\R^2$.
\end{thm}
\noindent Thus, if $\cuti_1 / \cuti_2 \in \R \setminus \Q$, and if $F$ is a continuous and 1-periodic function of $\R^2$ satisfying $F(\cutvec\, \R) = 0$, then it follows from Theorem \ref{thm:kronecker} that $F = 0$. 

\begin{figure}[ht!]
	\includegraphics[page=2]{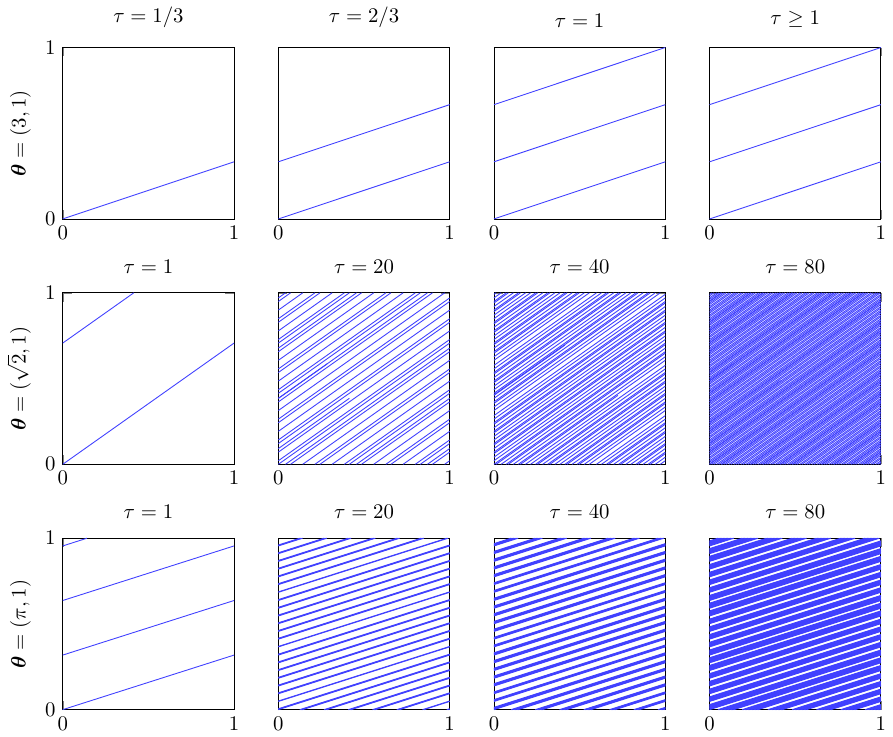}
	\caption{Representation of the set $\big\{(\{\cuti_1\, \xvi\}, \{\cuti_2\, x\}),\ x \in (0, \tau)\big\}$ in $(0, 1)^2$ for different values of $\tau$, when $\cuti_1/\cuti_2 \in \Q$ (first row), and when $\cuti_1/\cuti_2 \in \R \setminus \Q$ (second row for $\cutvec = (\sqrt{2}, 1)$ and third row for $\cutvec = (\pi, 1)$).\label{fig:fibrage}}
\end{figure}

\vspace{1\baselineskip}\noindent 
% Let $\{\xvi\}$ denote the Theorem \ref{thm:kronecker} implies that the broken line $\big\{(x\,\cuti_1[1],x\,\cuti_2[1]),\;x\in\R\big\}$ is dense in the unit cell $(0,1)^2$. To illustrate this, Figure \ref{fig:fibrage} represents the set $\big\{(x\,\cuti_1[1],x\,\cuti_2[1]),\;x\in(0,M)\big\}$ in the unit cell for different values of $M$, when (\emph{1}) $\theta_1/\theta_2 \in \Q$ (see the first row), and when (\emph{2}) $\theta_1/\theta_2 \in \R \setminus \Q$ (see the second row for $\cutvec = (\sqrt{2}, 1)$ and the third one for $\cutvec = (\pi, 1)$). For $M$ large enough, in the first case, this set is reduced to a \textit{finite} union of segments, whereas in the second case, it seems to fill the unit cell without ever passing through the same positions.
%\textcolor{blue}{
Let $\partent{\xvi}$  and $\{\xvi\}$ denote respectively the integer part and the fractional part of $\xvi\in\R$, given by $\partent{\xvi}\in\N\cap(\xvi-1,\xvi]$ and $\{\xvi\} \coloneqq  \xvi - \partent{\xvi}$. For $\ord = 2$, defining 
\begin{equation*}
	\label{QPtheo:eq:def_broken_line}
	\spforall I\subset \R,\quad \mathscr{D}_{I} \coloneqq  \big\{(\{\cuti_1\, \xvi\},\ \{\cuti_2\, \xvi\}),\quad  \xvi \in I\big\},
\end{equation*}
 Theorem \ref{thm:kronecker} implies that the broken line $\mathscr{D}_{\R}$ is dense in the unit cell $(0,1)^2$. To illustrate this result, Figure \ref{fig:fibrage} represents the set $\mathscr{D}_{(0, \tau)}$ in the unit cell for different values of $\tau > 0$, when $(1)$ $\theta_1/\theta_2 \in \Q$ (see the first row), and when $(2)$ $\theta_1/\theta_2 \in \R \setminus \Q$ (see the second row for $\cutvec = (1, \sqrt{2})$, and the third row for $\cutvec = (1, \pi)$). For $\tau$ large enough, in the first case, this set is reduced to a \textit{finite} union of segments, whereas in the second case, it seems to fill the unit cell without ever passing through the same location.
%} 

% As $\tau$ increases, in the first case, $\mathscr{D}_{(0, \tau)}$ is reduced to a \textit{finite} union of segments, contrary to the second case where it seems to fill the unit cell without ever passing through the same positions. It is also interesting to see that for $\cutvec = (1, \sqrt{2})$, the unit cell is in some sense filled uniformly, whereas for $\cutvec = (1, 9\, C_{10})$, the filling pattern is much closer to that of a rational case. The difference of behaviour seems to suggest that some irrationals (here $C_{10}$) may be \textquote{closer} to rationals than others ($\sqrt{2}$). This question is given further consideration in Section \ref{QPtheo:sec:irrationality_measure}. 

\noindent 
Another important notion is the \emph{Liouville-Roth irrationality measure}, which is an indicator on \textquote{how far the irrational $\delta$ is close to rational numbers}: 
\begin{equation}\label{eq:meas_irrat}
	\spforall \delta \in \R\setminus\Q, \; m(\delta) \coloneqq  \sup\left\{\nu>0, \;\spexists (p_n,q_n)\in\Z \times \N^*,\ (p_n,q_n)\rightarrow \infty, \ \left|\delta-\frac{p_n}{q_n}\right|\leq \frac{1}{q_n^\nu} \right\}.
\end{equation}
By Dirichlet's theorem, we have 
\[
	\spforall\delta\in\R\setminus\Q, \quad m(\delta) \geq 2.
\]
By the Thue-Siegel-Roth's theorem, the algebraic irrationals (for instance $\sqrt{2}$) have a measure of irrationality equal to $2$, which means that they are irrationals that, in some sense, are the furthest from the rationals. On the other hand, numbers whose measure of irrationality is infinite are called \emph{Liouville numbers}. 
%The set of Liouville numbers %
% \begin{equation}\label{eq:Liouville}
% \{\vts\delta\ /\ m(\delta) = +\infty\vts\}
% \end{equation}
%will play an important role in our study. 
%We note from Kintchine-Groshev's theorem (see the note \cite{hussain2014note} and references therein) that this set has Lebesgue measure $0$. 
Note that if $\delta$ is not a Liouville number, then $\delta$ satisfies the so-called Diophantine condition
\begin{equation}\label{eq:notLiouville}
	\spforall \nu > m(\delta), \quad \spexists c \equiv c(\delta, \nu) > 0, \quad  \spforall (p,q) \in \Z \times \N^*, \quad \left|\delta - \frac{p}{q}\right| > \frac{c}{q^\nu}.
\end{equation}

\subsection{Theoretical difficulties and limiting absorption principle}\label{sec_theory}
Equation \eqref{eq:whole_line_problem} is encountered when solving the linear wave equation with a time-harmonic source term $f(x)\, e^{-\icplx \omega t}$, and when one is looking for a time-harmonic solution $u(x)\, e^{-\icplx \omega t}$. For real frequencies $\omega$, the well-posedness of this problem is unclear. This is linked to the spectrum of the non-negative self-adjoint operators
\begin{equation}
	\label{eq:op_pb}
	D(\mathcal{A}) \coloneqq  \Big\{u\in H^1(\R) \ \big/\ \mu \, \frac{d u}{d \xvi}\in H^1(\R)\Big\} \quad \textnormal{and} \quad \mathcal{A}u \coloneqq  -\frac{1}{\rho}\frac{d}{d \xvi} \Big( \mu \, \frac{d u}{d \xvi} \Big),% \quad \spforall u \in D(\mathcal{A}),
\end{equation}
and 
\begin{equation}
	\label{eq:op_quasiper}
	D(\mathcal{A}_\cutvec) \coloneqq  \Big\{u\in H^1(\R) \ \big/\ \mu_\cutvec \, \frac{d u}{d \xvi}\in H^1(\R)\Big\} \quad \textnormal{and} \quad \mathcal{A}_\cutvec u \coloneqq  -\frac{1}{\rho_\cutvec}\frac{d}{d \xvi} \Big( \mu_\cutvec \, \frac{d u}{d \xvi} \Big).% \quad \spforall u \in D(\mathcal{A}_\cutvec).
\end{equation}
By Weyl's theorem \cite{reed1978iv}, we know that the spectrum of $\mathcal{A}$ denoted $\sigma(\mathcal{A})$ is the union of its discrete spectrum denoted $\sigma_d(\mathcal{A})$ and the spectrum of $\mathcal{A}_\cutvec$ denoted $\sigma(\mathcal{A}_\cutvec)$:
\begin{equation}
	\label{eq:spect_A}
	\sigma(\mathcal{A})=\sigma_d(\mathcal{A})\, \cup\, \sigma(\mathcal{A}_\cutvec).
\end{equation}
If $\omega^2\notin \sigma(\mathcal{A})$, then $\mathcal{A}-\omega^2$ is invertible and for any $f\in L^2(\R)$, there exists a unique solution $u\in D(\mathcal{A})$ to \eqref{eq:whole_line_problem}. 

\vspace{1\baselineskip}\noindent 
{When $\omega^2\in \sigma_d(\mathcal{A})$, one cannot expect uniqueness of a solution to \eqref{eq:whole_line_problem}: if $u$ is solution, any element of the kernel of $\mathcal{A}-\omega^2$ added to $u$ is also solution. Actually, it is well known  that the time-harmonic regime does not make sense in general for these frequencies (associated to the so called trapped modes of the problem). }

\vspace{1\baselineskip}\noindent 
{Let us now focus on the case where $\omega^2\in \sigma(\mathcal{A}_\cutvec)$, that is in our point of view, the most interesting but also tricky situation. To explain the difficulties, let us focus in this paragraph on the periodic case, \emph{i.e.} the case where $\delta\in\Q$.  We know from \cite{fliss2009analyse,fliss2013dirichlet} and references therein, that the problem cannot be well-posed in $H^1(\R)$.
In this case, on one hand, the physical solution $u$, when it exists, does not belong to $H^1(\R)$ due to wave propagation phenomena and a lack of decay at infinity. On the other hand, uniqueness of a solution in $H^1_{\textit{loc}}(\R)$ does not hold in general since within this framework, one cannot make the difference between the so-called outgoing wave (propagating towards infinity), the ingoing wave (propagating from infinity), or any linear combination of them. To select the physical and outgoing wave and restore the uniqueness, one has to add a so-called radiation or outgoing condition.}  To obtain such condition in practice, one uses the \emph{limiting absorption principle}, which consists in 
\begin{itemize}
	\item adding some absorption denoted $\veps$, that is replacing $\omega^2$ by $\omega^2+\icplx\vts \veps$ in \eqref{eq:whole_line_problem};% (see Remark \ref{rmk:assumptions_absorption});
	\item solving the problem with absorption, well-posed in $H^1(\R)$:
	\begin{equation} 
		\label{eq:whole_line_problem_eps}\tag{$\mathscr{P}_\veps$}
		\displaystyle
		- \frac{d }{d\xvi}\Big( \mu \,\frac{d u_\veps}{d \xvi} \Big) - \rho \, (\omega^2 + \icplx\vts \veps) \, u_\veps = f \quad \textnormal{in} \quad \R;
		\end{equation}
	%	problem which has a unique solution $u_\veps\in H^1(\R)$.
		\item studying the limit of the solution $u_\veps$ as the absorption $\veps$ tends to $0$.
\end{itemize}

% \begin{rmk}\label{rmk:assumptions_absorption}
% 	Unlike the other chapters, we do not assume that $\Imag \omega > 0$, but rather that $\Imag(\omega^2) > 0$. However, since we are interested in the case where $\omega \geq 0$, these two assumptions are equivalent.
% \end{rmk}

% \vspspe
\noindent
The limiting absorption principle is a classical approach to study time-harmonic wave propagation problems in unbounded domains; see for instance \cite{agmon1975spectral, eidus1986limiting, wilcox1966wave}. More recently, it has been successfully applied to locally perturbed periodic media \cite{joly2006exact, hoang2011limiting, kirsch2018radiation, radosz2015new}. 
It is worth mentioning that these works implicitly rely on the purely absolutely continuous nature of the spectrum of second-order differential operators with periodic coefficients. On the other hand, the spectrum of $\mathcal{A}_\cutvec$ is in general more intricate \cite{last2007exotic,simon1982almost}. In fact, the spectrum has in general an absolutely continuous part (as for periodic media) \cite{dinaburg1975one,chierchia1987absolutely,russmann1980one,sarnak1982spectral,eliasson1992floquet,pastur1992spectra}, but it can also have a pure point part \cite{sarnak1982spectral,frohlich1990localization,bjerklov2021coexistence} and a singular continuous part \cite{pearson1978singular,avron1983almost,del1994operators,del1996operators,jitomirskaya1994operators,simon1995operators,simon1996operators,moser1981example,fedotov2002anderson}. Moreover, in the aforementioned works, it is highlighted that some parts of the spectrum could even be a Cantor set, that is a closed set with no isolated points and with dense complement. While the spectral properties of discrete or continuous Schrödinger operators with quasiperiodic potentials has been widely studied since the 1980s, the limiting absorption principle remains, to the best of our knowledge, an open question.

 \vspspe
%  \sfcomment{The objective of this paper is to prove that the limiting absorption principle holds provided some assumptions on the coefficients according to the irrationality measure of $\delta$ (see Assumption \ref{assumption_Dioph}) and some indirect assumptions (see Assumptions \ref{assumption} and \ref{ass:energy_flux}) are satisfied. These assumptions should be linked to the frequency and the nature of the spectrum, but for now, this link is an open question. We also propose a method to characterize and compute the physical solution. This is, as far as we know, the first result on the limiting absorption principle for quasiperiodic media.}
{Our objective in this paper is twofold:
\begin{enumerate} 
	\item contributing to the analysis of the limiting absorption process for $(\mathscr{P}_\varepsilon)$ via a constructive approach;
	\item proposing a numerical method for computing the corresponding physical solution based on the approach of the previous point.
	\end{enumerate}
	For achieving this, we shall be led to make two types of assumptions (that cannot be exposed in detail at this early stage of the paper).
	\begin{itemize}
		\item [(i)] The first assumption is technical and relates the smoothness of the 2D functions $\rho_p$ and $\mu_p$ to the irrationality measure of $\delta$: this is Assumption \ref{assumption_Dioph} of Section \ref{sec:Riccati_eps}.
		\item [(ii)] The second assumption requires the frequency $\omega$ to be a{\it regular frequency} (see Definition \ref{def:regfreq}). This notion is related to a spectral problem depending on $\omega$ and the absorption parameter $\varepsilon$ (see Section \ref{sec:spectrum_propagation_operator}) and more precisely its fundamental eigenpair (defined in Theorem \ref{thm_spectre}): a frequency is defined as regular if the fundamental eigenpair admits a limit when $\eps \rightarrow 0$
		(as defined more precisely in Assumption \ref{assumption}) and this limit satisfies a non degeneracy property (Assumption \ref{ass:energy_flux}).
	\end{itemize}
	\begin{rmk} The assumptions evoked in the point (ii) are not explicit at all and a natural direction for future work would consist in getting a more tractable description and more explicit properties about the set of regular frequencies. We conjecture that if a subinterval of the spectrum of $\mathcal{A}_\cutvec$ is purely absolutely continuous with a smooth enough spectral density, then its elements are regular frequencies (see Section \ref{sec:conclusions}).
		\end{rmk}}

 \vspspe
We first recall briefly in the next subsection the analysis we have performed in \cite{amenoagbadji2023wave} for the sake of completeness.% \cite{amenoagbadji2023wave} for the problem with absorption. % \sfnote{{Citer les travaux de B. Davies et L. Morini "Super band gaps and periodic approximants of generalised Fibonacci tilings"?}}

\subsection{A brief recap of the DtN method for the absorbing case} \label{sec_DTN}
The method we use in \cite{amenoagbadji2023wave} consists in restricting the problem \eqref{eq:whole_line_problem_eps} to the bounded interval
\[
	I^0\coloneqq (a^\gauche,a^\droite)
\]
by constructing transparent boundary conditions of Dirichlet-to-Neumann (DtN) type, defined thanks to problems set on the unbounded intervals 
\[
	I^\gauche \coloneqq  (-\infty, a^\gauche) \quad \textnormal{and} \quad I^\droite\coloneqq (a^\droite,\infty).	
\]
More precisely, let us introduce the so-called DtN coefficients $\lambdadir{\idomlap}_{\veps}$ defined for $\idomlap \in \{\gauche, \droite\}$ by
\begin{equation}
	\displaystyle
	\lambdadir{\idomlap}_{\veps} \coloneqq  \nu^\idomlap\,\Big[\mu\, \frac{d \udir{\idomlap}_\veps}{d \xvi}\Big](a^\idomlap),\quad\text{with}\quad \nu^\gauche \coloneqq  1\quad \text{and}\quad \nu^\droite \coloneqq  -1, \label{eq:DtN_coefficients}
\end{equation}
where for $\idomlap \in \{\gauche, \droite\}$, $\udir{\idomlap}_\veps$ is the unique solution in $H^1(I^\idomlap)$ of
\begin{equation}\label{eq:halfline_pb_eps}%\tag{$\mathcal{P}_{\veps,D}^j$}
	\left|
		\begin{array}{r@{\ =\ }l}
			\displaystyle - \frac{d}{dx} \Big(\mu_\cutvec \,  \frac{d \udir{\idomlap}_{\veps}}{dx}\Big)-  \rho_\cutvec \, (\omega^2 + \icplx \veps) \, \udir{\idomlap}_{\veps} & 0 \quad \mbox{in}\ \ I^\idomlap
			\ret
			\udir{\idomlap}_{\veps}(a^\idomlap) & 1.
		\end{array}
	\right.
\end{equation}
Then one can show that if $u_\veps$ is the unique solution of \eqref{eq:whole_line_problem_eps}, its restriction to $I^0$ (denoted by $\udir{\interieur}_\veps$) is the unique solution of the problem
\begin{equation}
	\left|
		\begin{array}{r@{\ =\ }l}
			\displaystyle - \frac{d}{dx} \Big(\mu_\cutvec \,  \frac{d \udir{\interieur}_\veps}{dx}\Big)-  \rho_\cutvec \, (\omega^2 + \icplx \veps) \, \udir{\interieur}_\veps & f \quad \mbox{in}\ \ I^0
			\rett
			\multicolumn{2}{c}{\displaystyle \nu^\idomlap\,\Big[\mu\, \frac{d \udir{\interieur}_\veps}{d \xvi}\Big](a^\idomlap) = \lambdadir{\idomlap}_\veps  \, \udir{\interieur}_\veps(a^\idomlap), \quad \idomlap \in \{\gauche, \droite\}.}
		\end{array}
	\right.
\label{eq:interior_problem_eps}%\tag{$\mathcal{P}_\veps^0$}
\end{equation}
Moreover, by linearity, the solution $u_\veps$ of \eqref{eq:whole_line_problem_eps} is given by
\begin{equation}\label{eq:solution_of_whole_line_problem_eps}
\spforall \xvi \in \R, \quad u_\veps(\xvi) = 
\left\{
	\begin{array}{l@{\quad}l}
		\displaystyle \udir{\interieur}_{\veps}(\xvi), & \xvi \in I^0,
		\ret
		\displaystyle \udir{\interieur}_{\veps}(a^\idomlap)\; \udir{\idomlap}_{\veps}(\xvi), & \xvi \in I^\idomlap,\quad \idomlap \in \{\gauche, \droite\}.
	\end{array}
\right.
\end{equation}
In order to solve \eqref{eq:halfline_pb_eps} for $\idomlap \in \{\gauche, \droite\}$ and deduce the associated DtN coefficients $\lambdadir{\idomlap}_\veps$, we use the quasiperiodic nature of $\mu_\cutvec$ and $\rho_\cutvec$. We explain how to solve \eqref{eq:halfline_pb_eps} for $\idomlap = \droite$, the case of $\idomlap = \gauche$ being completely similar, and we suppose for simplicity that $a^\droite = 0$.

\vspspe
{To solve the half-line problem \eqref{eq:halfline_pb_eps}, we use the so-called \textit{lifting approach} which has been used for the homogenization of quasicrystals and Penrose tilings, and for numerical purposes \cite{kozlov1979averaging,rodriguez2008computation,bouchitte2010homogenization,wellander2019homogenization,ferreira2021homogenization,wang2025convergence} as well as for the analysis of some boundary layers in presence of periodic halfspaces \cite{gerard2012homogenization,gerard2011homogenization} or periodic structures separated by an interface \cite{blanc2015local}.} More precisely, as the coefficients $\mu_\cutvec$ and $\rho_\cutvec$ in \eqref{eq:halfline_pb_eps} are by definition traces of two--dimensional functions along the half-line $\cutvec\, \R_+$, it is natural to seek $\udir{\droite}_{\veps}$ as the trace along the same half-line of a two--dimensional function $\Udir{\droite}_{\veps}:\R^2\mapsto\C $ that would be characterized as the solution of a two--dimensional PDE with periodic coefficients. Thus, if $\udir{\droite}_{\veps}(\xvi)=\Udir{\droite}_{\veps}(\xvi\, \cutvec)$ holds, according to the chain rule:
\begin{equation}
	\label{eq:chain_rule}
	\displaystyle
	\spforall \xvi \in \R, \quad \frac{d}{d \xvi}[ \Udir{\droite}_{\veps}(\xvi\, \cutvec) ] = ( \Dt{} \Udir{\droite}_{\veps} ) (\xvi\, \cutvec), \quad \text{with} \quad \Dt{} \coloneqq  \cutvec \cdot \nabla = \cuti_1\, \frac{\partial}{\partial \yvi_1}+\cuti_2\, \frac{\partial}{\partial \yvi_2},
\end{equation}
the function $\Udir{\droite}_{\veps}$ must solve 
\begin{equation}
	\label{eq:half-space_problem}
	\displaystyle - \Dt{} \big( \mu_p \, \Dt{} \Udir{\droite}_{\veps} \big) - \rho_p \, (\omega^2+\icplx \veps) \, \Udir{\droite}_{\veps} = 0, \quad \textnormal{for} \quad \yvi_2 > 0,
\end{equation}
where we recall that the coefficients $\mu_p,\, \rho_p: \R^2 \to \R$ are continuous and $1$--periodic with respect to each variable. In addition, the boundary condition in \eqref{eq:halfline_pb_eps} can be lifted onto the Dirichlet boundary condition
\begin{equation}
	\label{eq:half-space_problem_Dir}
	\displaystyle \Udir{\droite}_{\veps} = {\varphi}, \quad \textnormal{on} \quad \yvi_2 = 0,
\end{equation}
\noindent
where the data ${\varphi}: \R \to \C$ could be chosen continuous and must satisfy ${\varphi}(0) = 1$, by consistency with the fact that $\udir{\droite}_{\veps}(0) = 1$. Furthermore, to exploit the periodicity of the coefficients $\mu_p$ and $\rho_p$ with respect to the transverse variable $\yvi_1$, we could choose
\begin{equation}
	\label{eq:varphi_per}
	\text{${\varphi}$ $1$--periodic,}
\end{equation}
so that it is natural to impose that
\begin{equation}
	\label{eq:U_theta_periodic}
	\text{$\Udir{\droite}_{\veps}(\varphi)$ is $1$--periodic with respect to the transverse variable $\yvi_1$.}
\end{equation}
% Let us define the half-space and its 
We consider the half-space $\Omega^\droite \coloneqq  \R \times \R_+$ and the half-cylinder $\Omega^\droite_\periodic\coloneqq (0,1)\times\R_+$ in the following. Let us introduce also the sets%, for $a \in \{0, 1\}$ and for any integer $i\in\llbracket 1, \ord\rrbracket $,
\begin{equation}\label{interfaces}
	% \spforall a>0, \quad \Sigmaper^{a} \coloneqq  \{\yv \in \Omegaper^\droite\ /\ \yvi_2 = a\}\quad \text{and}\quad \spforall\ell\in\{0,1\},\quad \Gamma_{\ell} = \{\yv \in \Omega,\ \yvi_1 =\ell \}.
	\spforall a \geq 0, \quad \Sigma^{\droite, a} \coloneqq  \R \times \{a\} \quad \textnormal{and} \quad \Sigmaper^{\droite, a} \coloneqq  (0, 1) \times \{a\}.% \quad \text{and} \quad \spforall\ell\in\{0,1\},\quad \Gamma^\droite_{\ell} = \{\yv \in \Omega^\droite,\ \yvi_1 =\ell \}.
\end{equation}
For any open set $\mathcal{O} \subset \R^\ord$, the directional Sobolev space
\begin{equation}
		H^1_\cutvec(\mathcal{O}) \coloneqq  \big\{U \in L^2(\mathcal{O})\ /\ \Dt{} U \in L^2(\mathcal{O}) \big\},
\end{equation}
is a Hilbert space equipped with the scalar product
\[
(U, V)_{H^1_\cutvec(\mathcal{O})} \coloneqq  \int_{\mathcal{O}} \Big(\Dt{} U\, \Dt{} \overline{V} + U\, \overline{V}\Big).
\]
Let $\|\cdot\|_{H^1_\cutvec(\mathcal{O})}$ denote the induced norm. 

\vspspe
We can then introduce, for any boundary data $\varphi \in L^2(\Sigmaper^{\droite, 0})$, the unique solution $\Udir{\droite}_{\veps}(\varphi)$ in 
$H^1_{\cutvec}(\Omegaper^\droite)$ of the half-guide problem
\begin{equation}
	\label{eq:half_guide_problem}
	\left|
	\begin{array}{r@{\ }c@{\ }ll}
	\displaystyle - \Dt{} \big( \mu_p \, \Dt{} \Udir{\droite}_{\veps} \big) - \rho_p \, (\omega^2+\icplx \veps) \, \Udir{\droite}_{\veps} &=& 0, \quad \textnormal{in}\quad \Omegaper^\droite,
	\ret
	\displaystyle \Udir{\droite}_{\veps}\Big|_{\Sigmaper^{\droite, 0}} &=& \varphi,
	\ret
	\multicolumn{4}{c}{\displaystyle \Udir{\droite}_{\veps}\Big|_{\yvi_1 = 0}= \Udir{\droite}_{\veps}\Big|_{\yvi_1 = 1}\quad\text{and} \quad \mu_p\, \Dt{} \Udir{\droite}_{\veps}\Big|_{\yvi_1 = 0}= \mu_p\, \Dt{} \Udir{\droite}_{\veps}\Big|_{\yvi_1 = 1}}.
	\end{array}
	\right.
\end{equation}
Note that the definition of the traces of functions in directional Sobolev spaces is really subtle, see \cite{amenoagbadji2023wave} for more details. In this paper, we use directly the results of \cite{amenoagbadji2023wave}, which implies in particular that the trace and the directional normal trace on $\Sigmaper^a,\;a>0$ of a periodic function in $ H^1_{\cutvec}(\Omegaper^\droite)$ are in $L^2(\Sigmaper^a)$. To simplify the notation in the rest of the paper, we will write 
\begin{equation}\label{eq:per_cond}
	V\Big|_{\yvi_1=0}= V\Big|_{\yvi_1=1}\;\;\text{and}\;\; \mu_p\, \Dt{} V\Big|_{\yvi_1=0}= \mu_p\, \Dt{} V\Big|_{\yvi_1=1} \quad \Longleftrightarrow \quad V\;\text{is periodic w.r.t. }\yvi_1
\end{equation}
The advantage of \eqref{eq:half_guide_problem} compared to 
\eqref{eq:halfline_pb_eps} lies in its periodic nature which allows us to exploit tools that are well suited for periodic waveguides, such as the method developed in \cite{joly1992some,fliss2009analyse} for the classical Helmholtz equation $- \nabla \cdot (\mu_p\, \nabla U) - \rho_p\, \omega^2\, U = 0$.

\vspspe
In what follows, $\cellper^{\droite, m}$ is the periodicity cell defined for every $m \in \N$ by
\begin{equation}\label{cellules} 
\cellper^{\droite, 0} = (0, 1)^2 \quad \textnormal{and} \quad \cellper^{\droite, m} = \cellper^{\droite, 0} + m\, \vv{\ev}_2, \quad \textnormal{so that} \quad \overline{\Omegaper^\droite} = \bigcup_{m \in \N} \overline{\cellper^{\droite, m}}.
\end{equation}
By periodicity, each cell $\cellper^{\droite, m}$ can be identified to $\cellper^{\droite, 0}$. Similarly, each interface $\Sigmaper^{\droite, m}$ can be identified to $\Sigmaper^{\droite, 0}$. 
%The cells and interfaces are represented in Figure ....

\vspspe
The solution $\Udir{\droite}_{\veps}(\varphi)$ of \eqref{eq:half_guide_problem} has a particular structure that we now highlight. Let $\Pdir{, \droite}_{\veps} \in \mathscr{L}\big(L^2(\Sigmaper^{\droite, 0})\big)$ be the operator
\begin{equation}
	\label{eq:definition_propagation_operator_D}
	\spforall \varphi \in L^2(\Sigmaper^{\droite, 0}), \quad \Pdir{, \droite}_{\veps}\varphi \coloneqq  \Udir{\droite}_{\veps}(\varphi)\Big|_{\Sigmaper^{\droite, 1}},
\end{equation}
where $L^2(\Sigmaper^{\droite, 1})$ and $L^2(\Sigmaper^{\droite, 0})$ have been identified to each other in an obvious manner. This identification will be used systematically in what follows, even if not mentioned. Since the coefficients $\mu_p$ and $\rho_p$ are $1$-- periodic w.r.t. $y_2$, we can show by induction that for any $\varphi$ in $L^2(\Sigmaper^{\droite, 0})$, we have
	\begin{equation}
		\label{eq:periodic_structure_solution_D}
		\spforall n \in \N, \quad \aeforall \yv \in \Omegaper^\droite, \quad \left[\Udir{\droite}_{\veps}(\varphi)\right](\yv + n\, \vv{\ev}_2) = \left[\Udir{\droite}_{\veps}\big((\Pdir{, \droite}_{\veps})^n \varphi\big)\right](\yv).
	\end{equation}
This allows to show, since $\Udir{\droite}_{\veps}(\varphi)\in H^1_\cutvec(\Omegaper^\droite)$, that the spectral radius of $\Pdir{, \droite}_{\veps}$, namely $\rho(\Pdir{, \droite}_{\veps})$ is strictly less than one. We can rewrite \eqref{eq:periodic_structure_solution_D} thanks to local cell solutions defined as follows: %
for all $\varphi\in L^2(\Sigmaper^{\droite, 0})$, for $\ell\in\{0,1\}$, let $\Edir{, \droite}_{\veps, \ell}(\varphi)\in H^1_{\cutvec}(\mathcal{C}^\periodic_0) $ satisfy
\begin{equation}
	\label{eq:local_cell_problem}
	\left|
	\begin{array}{r@{\ }c@{\ }l@{\quad}l}
	\displaystyle - \Dt{} \big( \mu_p \, \Dt{} \Edir{, \droite}_{\veps, \ell} \big) - \rho_p \, (\omega^2+\icplx \veps) \, \Edir{, \droite}_{\veps, \ell} &=& 0, \quad \textnormal{in} & \cellper^{\droite, 0},
	\\[8pt]
	% \displaystyle \restr{E^j}{\yvi_i = 0} &=&\multicolumn{2}{@{}l}{\restr{E^j}{\yvi_i = 1}, \spforall i = 1,\dots,\ord-1,}
	% \\[4pt]
	\Edir{, \droite}_{\veps, \ell}\;\text{is periodic w.r.t. }\yvi_1,
	\end{array}
	\right.
\end{equation}
with the boundary conditions
\begin{equation}
	\label{eq:local_cell_BC}
	\left|
	\begin{array}{c@{\quad}c@{\quad}c}
		\Edir{, \droite}_{\veps, 0}\Big|_{\Sigmaper^{\droite, 0}} = \varphi &\textnormal{and}& \Edir{, \droite}_{\veps, 0}\Big|_{\Sigmaper^{\droite, 1}} = 0,
		\\[8pt]
		\Edir{, \droite}_{\veps, 1}\Big|_{\Sigmaper^{\droite, 0}} = 0 &\textnormal{and}& \Edir{, \droite}_{\veps, 1}\Big|_{\Sigmaper^{\droite, 1}} = \varphi.
	\end{array}
	\right.
\end{equation}
Hence, by \eqref{eq:periodic_structure_solution_D}, if the propagation operator $\Pdir{, \droite}_{\veps}$ is known, by linearity and uniqueness of the solution, the solution of the half-guide problem can be represented cell by cell as follows:
\begin{equation}
	\displaystyle
	\spforall n \in \N, \quad \left[\Udir{\droite}_{\veps}(\varphi)\right](\cdot + n\, \vv{\ev}_2)|_{\cellper^{\droite, 0}} = \Edir{, \droite}_{\veps, 0}\big((\Pdir{, \droite}_{\veps})^{n} \varphi\big) + \Edir{, \droite}_{\veps, 1}\big((\Pdir{, \droite}_{\veps}\big)^{n+1} \varphi).
	\label{eq:UfromEis}
\end{equation} 
Moreover, the characterization of the operator $\Pdir{, \droite}_{\veps}$ comes from the continuity of the directional derivative of $\Udir{\droite}_{\veps}(\varphi)$ across the interface $\Sigmaper^{\droite, 1}$. If we introduce the \emph{local DtN operators} $\Tdir{, \droite}_{\veps, \ell k} \in \mathscr{L}(L^2(\Sigmaper^{\droite, 0}))$, defined for $\ell, k = 0, 1$ by
\begin{equation}
	\label{eq:DtNloc}
	\spforall \varphi \in L^2(\Sigmaper^{\droite, 0}), \quad \Tdir{, \droite}_{\veps, \ell k} \varphi \coloneqq (-1)^{k+1}\, \cuti_2\, \left[\mu_p\, \Dt{} \Edir{, \droite}_{\veps, \ell}(\varphi)\right]\Big|_{\Sigmaper^k},
\end{equation}
we can prove that the propagation operator $\Pdir{, \droite}_{\veps}$ is the unique solution of the constrained Riccati equation
\begin{equation}
	\label{eq:Riccati}
	\left|
	\begin{array}{l}
		\textnormal{\textit{Find $P \in \mathscr{L}(L^2(\Sigmaper^{\droite, 0}))$ such that $\varrho(P) < 1$ and}}
		\ret
		\multicolumn{1}{c}{\displaystyle \Tdir{, \droite}_{\veps, 10}P^2 + (\Tdir{, \droite}_{\veps, 00} + \Tdir{, \droite}_{\veps, 11})\, P + \Tdir{, \droite}_{\veps, 01} = 0.}
	\end{array}
	\right.
\end{equation}
This characterizes the solution $\Udir{\droite}_{\veps}(\varphi)$ and also the solution $\udir{\droite}_{\veps}$ since we have shown in \cite{amenoagbadji2023wave} that if $\varphi$ is continuous in a neighborhood of $0$ and $\varphi(0)=1$, we have
\[
	\spforall \xvi\in\R_+,\quad \udir{\droite}_{\veps}(\xvi)=\left[ \Udir{\droite}_{\veps}(\varphi)\right](\cuti_1\,\xvi,\cuti_2\,\xvi),
	\]
	where we have identified $\Udir{\droite}_{\veps}(\varphi)$ and its periodic extension with respect to $y_1$.
We can finally deduce the DtN coefficient computed thanks to the DtN operator
	\begin{equation}
		t^{\dirbc, \droite}_{\veps} \coloneqq \frac{1}{\cuti_2}\; \left[\Tdir{, \droite}_{\veps} \varphi\right](0) \quad \text{where} \quad  \Tdir{, \droite}_{\veps} \coloneqq \Tdir{, \droite}_{\veps, 00}+ \Tdir{, \droite}_{\veps, 10}\,\Pdir{, \droite}_{\veps}.\label{eq:DtNcoeff_from_DtNoperator}
	\end{equation}

\subsection{Objectives of the work and outline}
% The method we have performed in Chapter \ref{ch:H1Dabs} consists in restricting the problem \eqref{eq:whole_line_problem_eps} to the bounded interval $(I^\interieur)$ %
% % \[
% % 	I^0\coloneqq (a^\gauche,a^\droite)
% % \]
% by constructing transparent boundary conditions of Dirichlet-to-Neumann (DtN) type, featuring DtN coefficients defined thanks to problems set on the half-lines $(-\infty, a^\gauche)$ and $(a^\droite, +\infty)$. %
% % \[
% % 	I^l\coloneqq (-\infty,a^\gauche)\quad\text{and}\quad I^r\coloneqq (a^\droite,\infty).	
% % \]
% To solve these problems, we have used the lifting approach, which consists in introducing 2D periodic half-guide problems that can be solved using the solutions of local cell problems, as well as a propagation operator, which satisfies a constrained Riccati equation.
% 
% \vspspe
Our goal is to make the absorption $\veps$ tend to $0$ in the computation of the DtN coefficients, in order to propose a limit problem set in the bounded interval $I^\interieur$. If the DtN coefficients have a limit then the limit problem set in $I^\interieur$ will be Fredholm of index 0 in the sense of the associated operator: uniqueness then implies existence. Moreover, we can construct a limit solution on the whole line and show that it is indeed the limit of $u_\veps$. This is exactly what is done in \cite{joly2006exact,fliss2009analyse,fliss2021dirichlet} for the Helmholtz equation with locally perturbed periodic coefficients in dimension $1$ or in waveguides.

\vspspe
However, this approach does not apply directly. Indeed, we show in Appendix \ref{sec_DtN} that in general, the solutions of the local Dirichlet cell problems (\ref{eq:local_cell_problem}, \ref{eq:local_cell_BC}) (which are involved in the construction of the DtN coefficients) do not have a limit when $\veps$ goes to $0$, the limit Dirichlet cell problems being not well-posed, when $\omega^2$ lies in a semi-infinite interval of $\R_+$. This drawback is actually artificial since it is directly linked to the fact that we want to construct DtN coefficients. Constructing so-called \emph{Robin-to-Robin} (RtR) coefficients instead allows to circumvent this difficulty, as we show in Section \ref{sec:RtR_eps}. However, the construction of the associated Robin half-line solution is more involved than the one described previously for the Dirichlet half-line solutions, even if the underlying ideas are really similar. In Sections \ref{sec:RtR_eps} and \ref{sec:Riccati_eps}, we explain how to solve the Robin half-line problem in presence of absorption $(\veps > 0)$. In Section \ref{sec:prelim_LAP}, we establish the links between the spectra of some differential operators that appear in the study. This corresponds to a preliminary section whose results will be used afterwards. Finally, under the assumptions evoked in Section \ref{sec_theory}, we are able to study the limit of the Robin half-line solutions and the associated RtR coefficients in Section \ref{sec:LAP}. Section \ref{H1Dlap:sec:results} provides some numerical results to illustrate the method.

\paragraph*{Notation used throughout the paper} 
In what follows, 
\begin{enumerate}
	\item We denote $\N\coloneqq \{0,1,2,\ldots\}$ the set of natural numbers and $\N^*\coloneqq \N\setminus\{0\}$ the set of positive ones.
	\item For all $p, q \in \N,\ p < q$, we set $\llbracket p, q\rrbracket \coloneqq  \{j \in \N\ /\ p\leq j\leq q\}$.
	\item {For all $a\in\C$, $r>0$, $\mathcal{C}(a,r)\coloneqq \{z\in\C,\; |z-a|=r\}$.}
 \item For $i \in\llbracket 1,2\rrbracket$, we denote by $\vv{\ev}_i$ the $i$-th unit vector from the canonical basis of $\R^2$. For any elements $\yv = (\yvi_1, \yvi_2)$ and $\zv = (z_1,z_2)$ in $\R^2$, the Euclidean inner product of $\yv$ and $\zv$ is denoted $\yv \cdot \zv \coloneqq  \yvi_1\,z_1 + \yvi_2\,z_2$, and the associated norm is $|\yv| \coloneqq  \sqrt{\yv \cdot \yv}$.
 \item We introduce $\mathscr{C}^0_{\textit{per}}(\R^2)$ as the space of continuous functions $F: \R^2 \to \R$ that are $1$-periodic with respect to each variable, and $\mathscr{C}^\infty_0(\mathcal{O})$ as the space of smooth functions that are compactly supported in $\mathcal{O} \subset \R^2$.
 \item In any domain $D \subset \R^2$ of the form $D = \R \times I$, where $I = (\alpha, \beta), (\alpha, + \infty)$ or $( -\infty, \beta)$, the notion of periodicity only makes sense with respect to $y_1$, the corresponding periodicity cell will be systematically labeled with the $\periodic$-index 
 $$
 D_\periodic \coloneqq  [0,1] \times I.
 $$
 Accordingly, we shall introduce the functional spaces of $y_1$-periodic functions
 \[
  \begin{array}{|lll}
			\mathscr{C}^0_{\textit{per}}(D) \coloneqq  \big\{ v \in \mathscr{C}^0(D) \ / \ v(\yv + \vv{\ev}_1) = v(\yv) \ \ \spforall \yv \in D \big\} 
			\ret
			L^2_{\textit{per}}(D) \coloneqq  \big\{ v \in 	L^2_{loc}(D) \ / \ v(\yv + \vv{\ev}_1) = v(\yv)\ \ \mbox{ for } \aeforall \yv \in D \big\} 
		\end{array} 
 \]
 As $L^2_{\textit{per}}(D)$ is trivially isomorphic to $L^2(D_\periodic)$ via the restriction operator from $D$ to $D_\periodic$, it is naturally equipped 
 %with an Hilbert space structure 
 with the $L^2(D_\periodic)$-inner product.
 \item In the same way, for any line $\Gamma \coloneqq  \R \times \{ b \}$, with $b \in \R$, we shall set $\Gamma_\periodic \coloneqq  \; (0,1) \times \{ b \}$ and
 $$
 L^2_{\textit{per}}(\Gamma) \coloneqq  \big\{ v \in 	L^2_{loc}(\Gamma) \ / \ v(\yv + \vv{\ev}_1) = v(\yv), \mbox{ for } \aeforall \yv \in \Gamma \big\},
 $$
 which we equip with an Hilbert space structure using the $L^2(\Gamma_\periodic)$ inner product. We will often make implicitly the trivial identifications $\Gamma \equiv \R$ and $L^2_{\textit{per}}(\Gamma)\equiv L^2_{\textit{per}}(\R) $.
\item The spectral radius of a bounded linear operator $X$ is denoted $\rho(X)$.
\end{enumerate}

\section{A Robin-to-Robin (RtR) approach for the absorbing case}\label{sec:RtR_eps}

\subsection{The RtR transparent boundary conditions}\label{sub:RtR_eps_intro}
As an introductory material, we need to introduce two {\it outgoing Robin} differential operators for 1D functions $u$, associated with a given impedance $z > 0$, namely 
\begin{equation} \label{defRobinoperators_out} 
	R^\idomlap_+ u \coloneqq  \nu^\idomlap\vts \mu_\cutvec \, \frac{du}{dx} - \icplx\vts  z\vts  u \quad \spforall \idomlap \in \{\gauche, \droite\}, \quad \textnormal{with} \quad \nu^\gauche \coloneqq  1 \quad \textnormal{and} \quad \nu^\droite \coloneqq  -1.
	% \quad \textnormal{and}\quad R^\droite_+ u \coloneqq  - \mu \, \frac{du}{dx} - \icplx\vts  z\vts  u.
	\end{equation}
% In the sequel, these operators  $R^+$ will be called the {\it left-going} Robin operator, because in involves the derivative in the $-x$ right direction and 
% oppositely $R^-$ will be called the {\it right-going} Robin operator.
{Here we recall that $\nu^\gauche$, \resp $\nu^\droite$,  corresponds to the outgoing normal with respect to $I^\gauche$, \resp $I^\droite$.} 
This allows us to introduce two problems on the half-lines $I^\gauche$ and $I^\droite$ (which replace the Dirichlet half-line problems \eqref{eq:halfline_pb_eps}), defining two Robin half-line solutions $u^\gauche_\veps$ and $u^\droite_\veps$ respectively: for $\idomlap \in \{\gauche, \droite\}$, \textit{find $u_{\veps}^\idomlap \in H^1(I^\idomlap)$ such that }
\begin{equation}  \label{halfline_problems}
	\left|
	\begin{array}{r@{\ =\ }l}
		\displaystyle - \frac{d}{dx} \Big(\mu_\cutvec \,  \frac{d u_{\veps}^\idomlap}{dx}\Big)-  \rho_\cutvec \, (\omega^2 + \icplx \veps) \, u_{\veps}^\idomlap & 0 \quad \mbox{in}\ \ I^\idomlap
		\ret
		R^\idomlap_+ u_{\veps}^\idomlap(a^\idomlap) & 1.
	\end{array}
	\right.
	\tag{$\mathscr{P}^\idomlap_\veps$}
\end{equation}
Note that for the boundary condition in these half-line  problems, we impose the {\it outgoing Robin trace} with respect to $I^\idomlap$.
Imposing the outgoing trace is essential for ensuring the well-posedness of Problem \eqref{halfline_problems} in $H^1(I^\idomlap)$. Without going into details, let us simply point out that the variational formulation of \eqref{halfline_problems} involves the sesquilinear form 
\[
\displaystyle
b^\idomlap_\veps(u,v) \coloneqq  \int_{I^\idomlap} \, \Big(\mu_\cutvec \,  \frac{d u}{dx}\, \overline{\frac{d v}{dx}} - \rho_\cutvec \, (\omega^2 + \icplx  \veps)\, u \, \overline{v} \Big) \, dx  -  \icplx \vts z \vts u(0)\,  \overline{v(0)},
\]
the key point being that the imaginary part of the corresponding quadratic form has a sign 
\begin{equation} \label{Imag} 
-\Imag b^\idomlap_\veps(u,u) \coloneqq  \veps \int_{I^\idomlap}  \rho_\cutvec \, |u|^2 \, dx  + z \, |u(0)|^2 \geq 0.
\end{equation}
One interpretation is that the impedance $z>0$ brings an absorption term whose role is similar to the one due to $\veps > 0$. The following result can then be deduced without any difficulty.
\begin{prop}\label{exist_1D}
	Let $\idomlap \in \{\gauche, \droite\}$ and $\veps>0$. Problem \eqref{halfline_problems} has a unique solution $u_{\veps}^\idomlap\in H^1(I^\idomlap)$. %
	% Moreover, one has the regularity result 
	% \begin{equation} \label{regularite1D} gauche
	% \spforall m \geq 1, \quad 	(\rho_\cutvec, \mu_\cutvec) \in C^{m-1}(\R) \times C^m(\R), \quad u_\veps^\idomlap\in C^{m+1}(I^j)
	% \end{equation}\sfnote{Je ne suis pas sure que ce resultat de regularite serve}
\end{prop}

\vspspe 
Let us now define for $\idomlap \in \{\gauche, \droite\}$ two {\it ingoing Robin} differential operators for 1D functions $u$ 
\begin{equation} \label{defRobinoperators_in} 
	R^\idomlap_- u \coloneqq  - \nu^\idomlap\vts \mu \, \frac{du}{dx} - \icplx\vts  z\vts  u.
	% R^\gauche_- u \coloneqq  - \mu \, \frac{du}{dx} - \icplx\vts  z\vts  u \quad \textnormal{and}\quad R^\droite_- u \coloneqq  \mu \, \frac{du}{dx} - \icplx\vts  z\vts  u.
\end{equation}
From the solutions $u_\veps^\idomlap$, $\idomlap \in \{\gauche, \droite\}$ we can introduce the {\it Robin-to-Robin} (RtR) coefficients $\RtRcoef^{\idomlap}_\veps$, $\idomlap \in \{\gauche, \droite\}$ defined by
\begin{equation}  \label{Robin_coeff}
\RtRcoef^{\idomlap}_\veps \coloneqq  R^\idomlap_- u_\veps^\idomlap\, (a^\idomlap). 
\end{equation}
This definition consists in evaluating the {\it ingoing Robin trace} of the solution of a half-line problem where we have imposed the {\it outgoing Robin trace}. %
From $\RtRcoef^\idomlap_\veps$, the restriction $u^{\interieur}_\veps$ of the solution $u_\veps$ of \eqref{eq:whole_line_problem_eps} in $(I^\interieur)$ can be characterized using two transparent boundary conditions at $x = a^\gauche$ and $a^\droite$, called Robin-to-Robin (RtR) conditions (that replace the DtN transparent boundary conditions in \eqref{eq:interior_problem_eps}):
\begin{equation}\label{eq:int_pb_eps_RtR}
	\displaystyle
	\left|
	\begin{array}{r@{\ =\ }l}
		\displaystyle - \frac{d}{dx} \Big(\mu_\cutvec \,  \frac{d u_\veps^{\interieur}}{dx}\Big)-  \rho_\cutvec \, (\omega^2 + \icplx \veps) \, u_\veps^{\interieur} & 0 \quad \mbox{in}\ \ I^\interieur
		\ret
		\multicolumn{2}{c}{\displaystyle R^\gauche_- u^{\interieur}_{\veps}\, (a^\gauche) = \RtRcoef^\gauche_\veps  \, R^\gauche_+ u^{\interieur}_{\veps}\,(a^\gauche),}
		\ret
		\multicolumn{2}{c}{\displaystyle R^\droite_- u^{\interieur}_{\veps}\, (a^\droite) = \RtRcoef^\droite_\veps  \, R^\droite_+ u^{\interieur}_{\veps}\,(a^\droite).}
	\end{array}
	\right.
	\tag{$\mathscr{P}^{\texttt{int}}_\veps$}
\end{equation}
Then, $u_\veps$ can be characterized as
\begin{equation}\label{eq:sol_eps}
	\displaystyle
	\aeforall \xvi \in \R, \quad u(\xvi) = \left\{
		\begin{array}{c@{\quad}l}
			u^\interieur_\veps (\xvi), & \xvi \in I^\interieur,
			\\[6pt]
			R^\idomlap_+ u^\interieur_\veps\, (a^\idomlap)\; u^\idomlap_\veps (\xvi), & \xvi \in I^\idomlap,\; \idomlap \in \{\gauche, \droite\}.
		\end{array}
	\right.
\end{equation}

\begin{rmk}
	Let $\idomlap \in \{\gauche, \droite\}$ and $\veps>0$. By definition \eqref{defRobinoperators_out} and \eqref{defRobinoperators_in} of the ingoing and outgoing Robin operators, one can rewrite the transparent boundary conditions in \eqref{eq:int_pb_eps_RtR} in the following form:% with $\nu^\gauche = 1$ and $\nu^\droite = -1$: 
	\[
	 \Big[ -\nu^\idomlap \, \mu \, \frac{du_{\veps}}{dx} -  \, \icplx  z  u_{\veps} \Big] (a^\idomlap) = \RtRcoef^{\idomlap}_\veps  \, \Big[ \nu^\idomlap \mu \,  \frac{du_{\veps}}{dx} - \icplx z  u_{\veps} \Big] (a^\idomlap)
	\quad \Longleftrightarrow \quad  \Big[\nu^\idomlap \, \mu \, \frac{du_{\veps}}{dx} \Big](a^\idomlap) = \lambdadir{\idomlap}_{\veps}  \, u_{\veps} (a^\idomlap)
	\]
	where $\lambdadir{\idomlap}_{\veps}$ are the DtN coefficients defined in \eqref{eq:DtN_coefficients}, and are then related to the RtR ones by 
	\begin{equation}  \label{RtR-DtN}
		\lambdadir{\idomlap}_{\veps}  = \icplx z \;  \Big( \frac{\RtRcoef^{\idomlap}_\veps-1 }{\RtRcoef^{\idomlap}_\veps+1}\Big).
	\end{equation}
At this stage, as far as $\veps > 0$, there is no clear interest of passing from Dirichlet half-line problems to Robin half-line problems and 
from DtN transparent boundary conditions to RtR ones. The interest will appear more clearly when passing to the limit $\veps \rightarrow 0$, in particular in the method for computing $\RtRcoef^{\idomlap}_\veps$ that we shall develop in the next section. In particular, we shall not suffer any longer from the problems explained in Appendix \ref{sec_DtN}.
\end{rmk}

\subsection{The Robin half-line problems and their 2D lifting}\label{sub:Robin2D_eps}
Since the half-line problems satisfied by $u^\gauche_\veps$ and $u^\droite_\veps$ are quite similar to each other, \emph{we first restrict ourselves to the half-line problem satisfied by} $u^\droite_\veps$. We observe first that the quasiperiodic half-line problem \eqref{halfline_problems} with coefficients $(\mu_\cutvec, \rho_\cutvec)$ belongs to a family of similar quasiperiodic half-line problems parametrized by $s \in \R$, where the quasiperiodic coefficients $(\mu_{s, \cutvec}, \rho_{s, \cutvec})$ are, for any $s \in [0,1]$, the traces of the 2D functions $\rho_p$ and $\mu_p$ along the line $\R\, \cutvec+s\,  \vv{\ev}_1$:
\begin{equation} \label{defrhosmus}
\spforall \xvi \in \R, \quad \mu_{s, \cutvec}(x)\coloneqq  \mu_p(s + \cuti_1\vts \xvi, \cuti_2\vts \xvi) \quad \textnormal{and} \quad \rho_{s, \cutvec}(x)\coloneqq  \rho_p(s + \cuti_1\vts \xvi, \cuti_2\vts \xvi),
\end{equation}
The family of quasiperiodic half-line problems are given by
\begin{equation}  \label{halfline_sproblems}
	\left|
	\begin{array}{r@{\ =\ }l}
		\displaystyle - \frac{d}{dx} \Big(\mu_{s, \cutvec} \,  \frac{d u_{s,\veps}^\droite}{dx}\Big) - \rho_{s, \cutvec} \, (\omega^2 + \icplx \veps) \, u_{s,\veps}^\droite & 0 \quad \mbox{in}\ \ \R_+,
		\ret
		R^\droite_{+, s} u_{s,\veps}^\droite(0) & 1,
	\end{array}
	\right.
	\tag{$\mathscr{P}^\droite_{\veps, s}$}
\end{equation}
where $R^\droite_{\pm, s}$ are  1D Robin differential operators for 1D functions:
\begin{equation} \label{defsRobinoperators} 
	R_{\pm,s}^\idomlap u \coloneqq  \pm \, \nu^\idomlap \mu_{s, \cutvec} \, \frac{du}{dx} - \icplx \vts z \vts u,\quad \idomlap\in\{\gauche,\droite\}.
\end{equation}
{Similarly to Proposition \ref{exist_1D}, we can show without any difficulty the following well-posedness result.}
\begin{prop}\label{exist_1D_s}
	Let $s\in\R$ and $\veps>0$. Problem \eqref{halfline_sproblems} has a unique solution $u_{s,\veps}^\droite\in H^1(\R_+)$. %
	% Moreover, one has the regularity result 
	% \begin{equation} \label{regularite1D} gauche
	% \spforall m \geq 1, \quad 	(\rho_\cutvec, \mu_\cutvec) \in C^{m-1}(\R) \times C^m(\R), \quad u_\veps^\idomlap\in C^{m+1}(I^j)
	% \end{equation}\sfnote{Je ne suis pas sure que ce resultat de regularite serve}
\end{prop}
\noindent Note that for all $s$, $u_{s,\veps}^\droite$ is defined on $\R_+$ whereas $u^\droite_\veps$ is defined on $[a^\droite,+\infty)$. Actually these functions are linked, by Propositions \ref{exist_1D} and \ref{exist_1D_s}, as follows: 
\begin{equation}\label{eq:halflinesol_rel}
	\spforall \xvi \in \R_+, \quad u^\droite_\veps (\xvi + a^\droite) = u^\droite_{s,\veps}(\xvi) \quad \textnormal{for}\ s = \cuti_1\, a^\droite.
\end{equation}
Indeed, let us note that %  $\mu_\cutvec (\xvi + a^\droite) = \mu_p(\cuti_1\, \xvi + \cuti_1\, \xvi, \cuti_2\, \xvi + \cuti_2\, a^\droite) = \mu_p(\cuti_1\, \xvi + \cuti_1\, \xvi, \cuti_2\, \xvi)$ because $a^\droite\, \cuti_2 \in \Z$ \eqref{eq:hyp}.
\begin{align*}
	\displaystyle
	\spforall \xvi \in \R_+, \quad \mu_\cutvec (\xvi + a^\droite) &\coloneqq  \mu_p\big(\cuti_1\, (\xvi + a^\droite), \cuti_2\, (\xvi + a^\droite)\big) = \mu_p\big(\cuti_1\, (\xvi + a^\droite), \cuti_2\, \xvi\big)
\end{align*}
because $a^\droite\, \cuti_2 \in \Z$ (see \eqref{eq:hyp_alar}),  so that $\mu_\cutvec(\cdot + a^\droite)$ corresponds to $\mu_{s, \cutvec}$ for $s = \cuti_1\, a^\droite$. The same goes for $\rho_\cutvec(\cdot + a^\droite)$ which corresponds to $\rho_{s, \cutvec}$ for $s = \cuti_1\, a^\droite$. 
~\\\\
The periodicity of $s \mapsto (\mu_{s, \cutvec}, \rho_{s, \cutvec})$ (which follows from the periodicity of $(\mu_p, \rho_p)$ with respect to $\yvi_1$) implies that
\begin{equation} \label{period1d} 
	\spforall s \in \R, \quad u_{s+1,\veps}^\droite = u_{s,\veps}^\droite.
\end{equation}
Moreover, from the uniform continuity of $s \mapsto (\mu_{s, \cutvec}, \rho_{s, \cutvec})$ (which results from the periodicity and the continuity of $(\mu_p, \rho_p)$) and the well-posedness of \eqref{halfline_sproblems}, it follows that $s\mapsto u_{s,\veps}^\droite$ defines a uniformly continuous function from $\R$ to $H^1(\R_+)$. The proof is the same as the one of \cite[Proposition 3.17]{amenoagbadji2023wave}. More generally, using the same arguments as in \cite[Proposition 3.17]{amenoagbadji2023wave}, one shows that
 \begin{equation}\label{reg_halflinesol}
	\mu_p, \rho_p\in \mathscr{C}^m(\R^2) \quad \Longrightarrow \quad s\mapsto u_{s,\veps}^\droite\in \mathscr{C}^m(\R, H^1(\R_+)).
\end{equation}
{As explained in Section \ref{sec_DTN}, the solutions of the half-line problems \eqref{halfline_problems}, or more generally \eqref{halfline_sproblems} can be lifted, \emph{i.e.} their solutions can be seen as traces along a line of the solution of a 2D problem posed in half-spaces.} More precisely, let us define the half-space and its boundary:
\begin{equation} \label{defHalfspaces} 
	\Omega^\droite \coloneqq  \R \times \R_+ \quad \textnormal{and} \quad \Sigma^\droite \coloneqq  \partial \Omega^\droite = \R \times \{0\}.
\end{equation}
%In the sequel, the two 1D boundaries $\Sigma_0^\pm = \partial \Omega^\pm_a \equiv \mathbb{R}$ will be parametrized by $s \in \R$ as follows 
%\begin{equation} \label{parametrization}
%\Sigma_0^\pm   = \big\{ (s, \pm \, a \, \theta_2) , s \in \R\}
%	\end{equation}
In the sequel, we shall play with the change of variables in $\R^2$
\begin{equation} \label{change_of_variable}
	(s,x)  \mapsto (y_1, y_2)= (s + x \, \theta_1, x \, \theta_2) \quad \Longleftrightarrow \quad (y_1,y_2) \mapsto (s,x) = ( y_1-y_2 \, \theta_1/\theta_2  , y_2 / \theta_2)
\end{equation} 
which corresponds to seeing the plane $\R^2$ with the following fibered structure
\begin{equation*} % \label{fibering}
\R^2 \coloneqq  \bigcup_{s \in \R} (s\, \vv{\ev}_1 + \R\, \cutvec). % \widetilde {\cal L}_{s, \cutvec}, \quad \mbox{with} \quad \widetilde {\cal L}_{s, \cutvec}\coloneqq  \big\{ (s + x \, \theta_1, x \, \theta_2),x \in \R \big\}.
\end{equation*} 

% \begin{figure}
% 	\centering
% 	\begin{tikzpicture}[thick, scale]
% \end{figure}

\begin{figure}
	\centering
	\includegraphics[page=5]{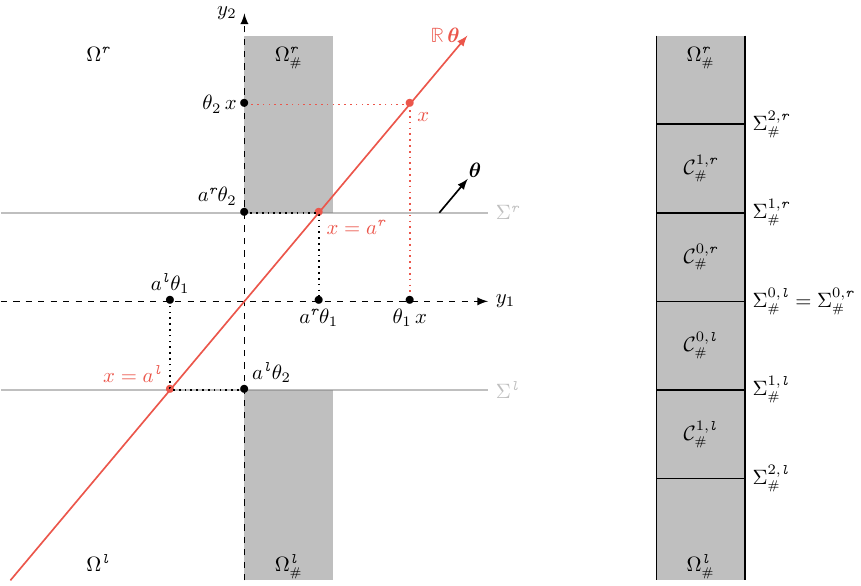}
	\caption{Notation introduced in Section \ref{sub:Robin2D_eps}\label{fig:illust_domains}}
\end{figure}

\noindent
Note that the change of variables defined by \eqref{change_of_variable} maps ${(s, \xvi) \in \R \times \R_+}$ into $\Omega^\droite$. %
Thus, through this change of variables, we can build, from the half-line solutions $u_{s,\veps}^\droite$ and from a boundary function $\varphi \in L^2_{\textit{loc}} (\Sigma^\droite)\equiv L^2_{\textit{loc}} (\R)$, a halfspace function $U_\veps^\droite: \Omega^\droite \rightarrow \C$ via the so-called lifting formula
\begin{equation} \label{defUeps}
\mbox{a.e. } (s,x) \in \R \times \R_+, \quad 	\big[ U^\droite_\veps(\varphi) \big]  (s + x \, \theta_1, x \, \theta_2) \coloneqq  	\varphi(s) \, u_{s,\veps}^\droite (x).
\end{equation} 
Roughly speaking, $U^\droite_\veps(\varphi)$ can be seen as a \textquote{weighted and oblique concatenation} of the 1D solutions $u_{s,\veps}^\droite$. 
Then, using the chain rule,
\begin{equation*}
	\displaystyle
	\frac{d}{dx} \big[ V(s + \theta_1\, \xvi, \theta_2\, \xvi) \big] = \Dt{} V(s + \theta_1\, \xvi, \theta_2\, \xvi),
\end{equation*}
it is easy, by adapting the arguments developed in \cite{amenoagbadji2023wave}, to characterize $U^\droite_\veps (\varphi)$ as the solution of  a 2D boundary value problem in $\Omega^\droite$  associated with the partial differential equation 
\begin{equation} \label{pde}
 - \Dt{} \big( \mu_p \, \Dt{}U^\droite_\veps  \big) - \rho_p \,(\omega^2 + \icplx \veps) \, U^\droite_\veps = 0 \quad \textnormal{in}\ \Omega^\droite.
	\end{equation}
In order to associate boundary conditions to \eqref{pde}, we need to introduce directional Robin operators as 2D extensions of the 1D differential operators $R^\droite_\pm,\;\idomlap \in \{\gauche, \droite\}$ defined in \eqref{defRobinoperators_out}  and \eqref{defRobinoperators_in}. We define the operators $\mathcal{R}^\droite_\pm$ for 2D functions $U$ as 
\begin{equation} \label{defdirectRobinoperators} 
	\mathcal{R}^\droite_\pm U  \coloneqq  \pm\, \nu^\droite \mu_p\,  \Dt{}  U -  \icplx\vts  z\vts  U.
\end{equation}
The operator $\mathcal{R}^\droite_+$ is the outgoing Robin operator whereas $\mathcal{R}^\droite_-$ is the ingoing one. Then from the Robin boundary condition in \eqref{halfline_sproblems}, we deduce the Robin boundary condition for $U^\droite_\veps$
\begin{equation}
	\label{eq:half-space_problem_RtR}
	\displaystyle \mathcal{R}^\droite_+ U^\droite_\veps  = {\varphi}, \quad \textnormal{on} \quad \yvi_2 = 0.
\end{equation} 
If $\varphi$ is 1--periodic, then from the periodicity of $s \mapsto u_{s,\veps}^\droite$ (see \eqref{period1d}), the function $U^\droite_\veps(\varphi)$  is 1--periodic with respect to $y_1$ and therefore belongs to the space 
\begin{equation} \label{defH1thetaper} 
	H^1_{\cutvec,per}(\Omega^\droite) = \big\{ U \in L^2_ {per}(\Omega^\droite) \; / \; \Dt{} U \in L^2_ {per}(\Omega^\droite)\}.
\end{equation}
Modulo an adaptation of the arguments developed in \cite[Proposition 3.18]{amenoagbadji2023wave}, see also Remark \ref{rem_thm}, one easily proves the following result.
\begin{prop}\label{thm-wellposedness}
	Let $\veps>0$. For any $\varphi \in L^2_{per}(\R)$, the function $U^\droite_\veps \coloneqq  U^\droite_\veps(\varphi)$ is the unique solution of the well-posed 2D problem: \textit{Find $U^\droite_\veps \in H^1_{\cutvec,per}(\Omega^\droite)$ such that}
	\begin{equation} \label{RobinHalfspaceproblems}
		\left|
		 \begin{array}{r@{\ =\ }l} 
				\displaystyle - \Dt{} \big( \mu_p \, \Dt{} U^\droite_\veps \big) - \rho_p \, (\omega^2+ \icplx \veps) \, U^\droite_\veps & 0 \quad \textnormal{in }\,\Omega^\droite, 
				\ret
				\displaystyle \mathcal{R}^\droite_+ \,  U^\droite_\veps\Big|_{\Sigma^\droite} & \varphi,
				\ret
				\multicolumn{2}{c}{\textnormal{$U^\droite_{\veps}$ is periodic w.r.t. $\yvi_1$.}}
		\end{array}
		\right.
	\end{equation} 
\end{prop}	
\begin{rmk} \label{rem_thm}
	In \cite[Section 3.2]{amenoagbadji2023wave}, we develop in detail the functional analytic tools to give a rigorous sense to Problem \eqref{RobinHalfspaceproblems}. In particular, functions in  $H^1_{\cutvec,loc}(\Omega^\droite)$ have well-defined traces on $\Sigma^\droite \coloneqq  \partial \Omega^\droite$ as elements of $L^2_{loc}(\Sigma^\droite)$. Note also that the PDE in \eqref{RobinHalfspaceproblems} implies that $\mu_p \, \Dt{} U^\droite_\veps \in H^1_{\cutvec,loc}(\Omega^\droite)$. Thus, the trace of $\mu_p \Dt{} U^\droite_\veps$ on $\Sigma^\droite$ is well-defined in 
	$L^2_{loc}(\Sigma^\droite)$, so that a proper sense can be given to the boundary condition in \eqref{RobinHalfspaceproblems}. 
\end{rmk} 
\noindent Using the definition \eqref{defUeps} of $U^\droite_\veps (\varphi)$, the uniform continuity of the map $s \mapsto u_{s, \veps}^\droite$ and \eqref{eq:halflinesol_rel}, we deduce the following result.
\begin{cor}\label{coro-wellposedness} Let $\veps>0$. If $\varphi$ is continuous in the neighborhood of $\yvi_1 = \cuti_1\, a^\droite$ and satisfies $\varphi(\cuti_1\, a^\droite) = 1$, then
	\begin{equation} \label{cutformulas}  
		\aeforall \xvi \in \R_+, \quad u_{\veps}^\droite (x + a^\droite) = u^\droite_{a^\droite\, \cuti_1, \veps} (\xvi) = U^\droite_\veps (\theta_1\vts (\xvi + a^\droite), \theta_2\vts \xvi),
	\end{equation}
	where $U^\droite_\veps = U^\droite_\veps (\varphi)$.
\end{cor} 
\noindent From Corollary \ref{coro-wellposedness}, we deduce a framework for finding the solution $U_\veps^\droite$ of the half-line problem \eqref{halfline_problems}:
\begin{itemize} 
	\item Choose a $1$--periodic function $\varphi$ which is continuous around $\cuti_1\, a^\droite$, such that $\varphi(\cuti_1\, a^\droite) = 1$ ($\varphi = 1$ for instance);
	\item Compute the corresponding solution $U^\droite_\veps(\varphi)$ of \eqref{RobinHalfspaceproblems};
	\item Deduce the function $u_{\veps}^\droite$ from the cut formula \eqref{cutformulas}.
\end{itemize}
\begin{rmk} Based on the above algorithm, one could wonder whether it is necessary to consider other functions than  $\varphi = 1$. Nevertheless, to derive the properties and the expression of $U^\droite_\veps(\varphi)$, it is necessary to consider general functions $\varphi \in L^2_{\textit{per}}(\R)$.
\end{rmk}
%From the numerical point of view, the new question is: how can one solve the 2D problem \eqref{RobinHalfspaceproblems}? 
At first sight, passing from \eqref{halfline_problems} to \eqref{RobinHalfspaceproblems} has added some difficulties: \emph{(1)} we went from a 1D problem to a 2D problem, and \emph{(2)} we went from a PDE involving an operator with an elliptic principal part to a PDE involving the operator $ - \Dt{} \big( \mu_p \, \Dt{}  \big) $ which is elliptically degenerate.
However, the great gain from \eqref{RobinHalfspaceproblems} is its periodic nature which enables us to use tools that are well-suited for periodic media.

\vspspe
A first trivial but important observation is that since $U^\droite_\veps(\varphi)$ is periodic in the $\vv{\ev}_1$-direction, it suffices to compute its restriction to the cylinder $\Omega^\droite_\periodic\coloneqq \{(y_1,y_2)\in \Omega^\droite\ / \ 0< y_1 <1\}$ (with $\Sigmaper^{\droite, 0} \coloneqq  (0, 1) \times \{0\}$) by solving the problem: \textit{Find $U^\droite_{\veps}   \in  H^1_\cutvec(\Omega_\periodic)$ such that} 
\begin{equation} \label{RobinHalfguideproblems}
	\left|
	 \begin{array}{r@{\ =\ }l} 
			\displaystyle - \Dt{} \big( \mu_p \, \Dt{} U^\droite_{\veps} \big) - \rho_p \, (\omega^2+ \icplx \veps) \, U^\droite_{\veps} & 0 \quad \textnormal{in }\,\Omegaper^\droite, 
			\ret
			\displaystyle \mathcal{R}^\droite_+ \,  U^\droite_{\veps}\Big|_{\Sigmaper^{\droite, 0}} & \varphi,
			\ret
			\multicolumn{2}{c}{\textnormal{$U^\droite_{\veps}$ is periodic w.r.t. $\yvi_1$.}}
	\end{array}
	\right.
\end{equation}
In the sequel, we will play alternatively with the two equivalent versions \eqref{RobinHalfspaceproblems} and \eqref{RobinHalfguideproblems} of the same half-space problem. The form \eqref{RobinHalfspaceproblems} will be used for analytical results, while \eqref{RobinHalfguideproblems} will be used for numerics. We explain in the next sections how to solve \eqref{RobinHalfguideproblems} and compute $U^\droite_{\veps}$. Because of the RtR conditions, we have to adapt the method described in \cite{amenoagbadji2023wave} that was recalled briefly in Section \ref{sec_DTN}.

\subsection{The propagation and scattering operators}
Let us first recall the definitions (\ref{interfaces}, \ref{cellules}) of the lines $\Sigma^{\droite, m} \coloneqq  \R \times \{m\}$, of the cells $\cellper^{\droite, m} \coloneqq  (0,1) \times (m, m+1)$, and of the interfaces $ \Sigmaper^{\droite, m+1}$ between the cells $\cellper^{\droite, m}$ and $\cellper^{\droite, m+1}$. These domains are shown in Figure \ref{fig:illust_domains}. By periodicity, each cell $\cellper^{\droite, m}$ can be identified to $\cellper^{\droite, 0} \equiv \cellper^{\droite}$ and each interface $ \Sigmaper^{\droite, m}$ to $ \Sigmaper^{\droite, 0} \equiv \Sigmaper$.

\vspspe
Proposition \ref{thm-wellposedness} defines a linear operator% for  $\idomlap \in \{\gauche, \droite\}$
$$
\varphi  \in L^2_{per}(\R) \mapsto U^\droite_\veps(\varphi) \in  H^1_{\cutvec,per}(\Omega^\droite) \mbox{ solution of } \eqref{RobinHalfspaceproblems} ,
$$ 
from which we introduce two operators $(\mathcal{P}^\droite_\veps, \mathcal{S}^\droite_\veps) \in \mathscr{L}\big(L^2_{per}(\R)\big)^2$ defined by 
\begin{equation} \label{defPS}
	\mathcal{P}^\droite_\veps \varphi \coloneqq  \mathcal{R}_+^\droite \, U^\droite_\veps(\varphi)\Big|_{\Sigma^{\droite, 1}} \quad \textnormal{and} \quad	\mathcal{S}^\droite_\veps\varphi\coloneqq  {\cal R}_-^\droite \, U^\droite_\veps(\varphi)\Big|_{\Sigma^{\droite, 1}}, 
\end{equation}
where the two interfaces $\Sigma^{\droite, 0}$ and $\Sigma^{\droite, 1}$  are identified to $\R$ in a trivial manner.

\vspspe
These operators play an important role in the structure of the solution. In fact, using the same arguments as for the proof of \cite[Proposition 4.1]{amenoagbadji2023wave}, one deduces in particular from the periodicity of the medium in the $y_2$ direction, the following result.
\begin{prop}\label{prop:property_PS}
	For any $\varphi\in L^2_{per}(\R)$, the Robin traces of $U^\droite_\veps(\varphi)$ on the interfaces ${\Sigma^{\droite, m}}$, $\spforall m \geq 1$, satisfy 
	\begin{equation}
		\displaystyle
		\mathcal{R}_+^\droite \,U^\droite_\veps(\varphi)\Big|_{\Sigma^{\droite, m}}= (\mathcal{P}^\droite_\veps)^{m}  \varphi \quad \textnormal{and} \quad \mathcal{R}_-^\droite \,U^\droite_\veps(\varphi)\Big|_{\Sigma^{\droite, m}}= \mathcal{S}^\droite_\veps\, (\mathcal{P}^\droite_\veps)^{m-1}\, \varphi.
\end{equation}
	\end{prop}
\noindent The operator $\mathcal{P}^\droite_\veps$ is called a \emph{Robin propagation} operator because it expresses how the outgoing Robin trace, that is the trace of $\mathcal{R}_+^\droite \, U^\droite_\veps$, propagates from one interface $\Sigma^{\droite, m}$ to the next one $\Sigma^{\droite, m+1}$.
The operator $\mathcal{S}^\droite_\veps$ transforms an outgoing Robin trace on $\Sigma^{\droite, 0}$ into an ingoing Robin trace on $\Sigma^{\droite, 1}$. We call it a \emph{scattering operator} by analogy with some terminology used in scattering theory.

\vspspe
As it can be expected, the fibered structure of the problem \eqref{RobinHalfspaceproblems} provides a particular structure to the operators $(\mathcal{P}^\droite_\veps,\mathcal{S}^\droite_\veps)$. 
{
\vspspe Let us first introduce the functions
\vspace*{-0.6cm}
% \begin{equation} \label{defsymbols} 
% \spforall s \in \R, \quad \aeforall x \in \R_+ , \quad 
%  \begin{array}{|r@{\ \coloneqq \ }l@{\quad}r@{\ \coloneqq \ }l}
% 	{\bf p}^\droite_\veps (s,x) &  \big(R^\droite_{+,s} \, u_{s,\veps}^\droite \big) (x), & \boldsymbol{p}^\droite_\veps (s) & {\bf p}^\droite_\veps \big(s, 1/\cuti_2\big),
% 	\ret
% 	{\bf s}^\droite_\veps (s,x) &  \big(R^\droite_{-,s} \, u_{s,\veps}^\droite \big) (x), & \boldsymbol{s}^\droite_\veps (s) & {\bf s}^\droite_\veps \big(s, 1/\cuti_2 \big),
% \end{array} 
% \end{equation}
\begin{equation} \label{defsymbols} 
\spforall s \in \R, \quad 
 \begin{array}{|r@{\ \coloneqq \ }l}
	 \boldsymbol{p}^\droite_\veps (s) & \big(R^\droite_{+,s} \, u_{s,\veps}^\droite \big) (1/\cuti_2),
	\ret
	 \boldsymbol{s}^\droite_\veps (s) & \big(R^\droite_{-,s} \, u_{s,\veps}^\droite \big) (1/\cuti_2),
\end{array} 
\end{equation}
where $R^\droite_{\pm,s}$ are the 1D outgoing and ingoing Robin operators, defined for all $s \in \R$ in \eqref{defsRobinoperators}.}
\begin{prop}\label{prop:propshift}
 {
	Let $\veps > 0$. The functions $\big(\boldsymbol{p}^\droite_\veps(s), 	\boldsymbol{s}^\droite_\veps(s) \big)$ 
	are continuous and $1$-periodic with respect to $s$. Moreover, we have
	\[
		\boldsymbol{p}^\droite_\veps(s)\neq 0\quad\spforall s \in\R. %\quad \textnormal{and} \quad \boldsymbol{p}^\droite_\veps(s)\neq 0 \quad \spforall s \in \R.
	\]
	Finally, if $\mu_p, \rho_p \in \mathscr{C}^m (\R^2)$, then
	\begin{equation}\label{reg_p_eps}
		s \mapsto  \boldsymbol{p}^\droite_{\veps}(s),\ s \mapsto  \boldsymbol{s}^\droite_{\veps}(s) \in \mathscr{C}^m(\R; \C).
	\end{equation}}
\end{prop}
\begin{dem}
The continuity properties follow from the uniform continuity of $s \mapsto  u_{s,\veps}^\droite$ and the fact that $u_{s,\veps}^\droite \in H^1(\R_+)\subset \mathscr{C}^0(\R_+)$ for all $s$. The property \eqref{reg_p_eps} is linked to \eqref{reg_halflinesol}.

\vspspe
{
It remains to show that $\boldsymbol{p}^\droite_\veps(s) \neq 0$ for any $s\in \R$. %(this implies directly that $\boldsymbol{p}^\droite_\veps(s)\neq 0$ for any $s\in\R$). 
Assume by contradiction that, for some $s_0 \in \R$, we have $\boldsymbol{p}^\droite_\veps(s_0)= 0$. Since $\boldsymbol{p}^\droite_\veps(s_0)\coloneqq  R^\droite_{+,s_0} u_{s_0,\veps}^\droite (1/\cuti_2)$ by definition \eqref{defsymbols}, where $u_{s_0,\veps}^\droite$ is the solution of \eqref{halfline_sproblems} for $s = s_0$, we have in particular
\begin{equation}  \label{halfline_s0problem}
	\left|
	\begin{array}{r@{\ =\ }l}
		\displaystyle - \frac{d}{dx} \Big(\mu_{s_0, \cutvec} \,  \frac{d u_{s_0,\veps}^\droite}{dx}\Big)-  \rho_{s_0, \cutvec} \, (\omega^2 + \icplx \veps) \, u_{s_0,\veps}^\droite & 0 \quad \mbox{in}\ \ (1/\cuti_2, +\infty),
		\ret
		R^\droite_{+, s} u_{s_0,\veps}^\droite(1/\cuti_2) & 0,
	\end{array}
	\right.
\end{equation}
By well-posedness of the above boundary problem for $z>0$ (see also Proposition \ref{exist_1D_s} and \eqref{halfline_sproblems} for a similar problem set in $\R^+$), we deduce that $u_{s_0,\veps}^\droite = 0$ in $(1/\cuti_2, +\infty)$. % 
Since $u_{s_0,\veps}^\droite$ satisfies the PDE in \eqref{halfline_s0problem} in all $(1/\cuti_2, +\infty)$, by Cauchy's uniqueness result $u_{s_0,\veps}^\droite = 0$ in $\R_+$, which contradicts the fact that $R^\droite_{+,s_0} u_{s_0,\veps}^\droite (0) = 1$.}
\end{dem}
\begin{rmk}
	Contrary to $\boldsymbol{p}^\droite_\veps$, the function $\boldsymbol{s}^\droite_\veps$ can vanish at certain points. For instance, if $\mu_p = \rho_p = 1$ and if the impedance $z$ is allowed to be a complex number with $z \coloneqq  \sqrt{\omega^2 + \icplx\vts \veps}$ and $\Real z > 0$, then it can be computed without difficulty that $\boldsymbol{s}^\droite_\veps = 0$.
\end{rmk}

\noindent The propagation operator and the scattering operator have a particular structure, highlighted in the next result.
\begin{prop}\label{thm_weighted-shift} 
	We have for all $\varphi\in L^2_{per}(\R)$,
	\begin{equation*}
		\spforall s \in \R, \quad \left\{
		\begin{array}{l@{\ =\ }l}
			\mathcal{P}^\droite_\veps \varphi(s) & \boldsymbol{p}^\droite_\veps (s-\delta) \, \varphi(s-\delta), 
			\ret
			\mathcal{S}^\droite_\veps \varphi (s) & \boldsymbol{s}^\droite_\veps (s-\delta) \, \varphi(s-\delta),
		\end{array}
		\right.
	\end{equation*}
	where $\mathcal{P}^\droite_\veps$ and $\mathcal{S}^\droite_\veps$ are defined in \eqref{defPS}, $\boldsymbol{p}^\droite_\veps$ and $\boldsymbol{s}^\droite_\veps$ in \eqref{defsymbols} and $\delta$ defined in \eqref{eq:delta} plays the role of a shift parameter. Moreover, $\mathcal{P}^\droite_\veps$ is bijective.
\end{prop}

\begin{dem} 
	We give the proof of the first identity, the proof of the second one being identical. It follows from the fibered structure \eqref{defUeps} of $U^\droite_\veps$ that
	\begin{equation} \label{fibered-structure-R} 
		\mathcal{R}_+^\droite U^\droite_\veps (s + \cuti_1\, \xvi, \cuti_2\, \xvi) = 	\big(R^\droite_{+,s} \, u_{s,\veps}^\droite \big) (x) \, \varphi(s)  ,\quad s\in\R, \ x \in \R_+.
	\end{equation}	
	To evaluate the trace of $R_+^\droite U^\droite_\veps$ on $\Sigma^{\droite, 1}$, it suffices to take $x = 1 / \cuti_2$ in the above expression, which gives
	\begin{equation} \label{fibered-structure} 
		\mathcal{R}_+^\droite U^\droite_\veps (s + \delta, 1) = \big(R^\droite_{+,s} \, u_{s,\veps}^\droite \big) (1 / \cuti_2) \, \varphi(s), \quad s \in \R,
	\end{equation}	
	that is to say, by definitions \eqref{defPS} and \eqref{defsymbols} of $\mathcal{P}^\droite_\veps$ and $\boldsymbol{p}^\droite_\veps$,
	\begin{equation} \label{fibered-structur2} 
		[\mathcal{P}^\droite_\veps \, \varphi](s + \delta) = \boldsymbol{p}^\droite_\veps (s) \, \varphi(s), \quad s \in \R,
	\end{equation}	
	which is nothing but the first identity after changing $s \mapsto s + \delta$.% The proof of the second identity in \eqref{weightedshiftoperators} is identical after setting 
	
	\vspspe
	The bijectivity of $\mathcal{P}^\droite_\veps$ is a direct consequence of Proposition \ref{prop:propshift} which states that $\boldsymbol{p}^\droite_\veps$ never vanishes.
\end{dem} 

\vspspe
Operators that have the form of $\mathcal{P}^\droite_\veps$ and $\mathcal{S}^\droite_\veps$ are known in the literature (see \cite{antonevich2012linear} for instance) as \emph{weighted shift operators}. The function
$\boldsymbol{p}^\droite_\veps$ (\resp $\boldsymbol{s}^\droite_\veps$) will be called the {\it symbol} of  $\mathcal{P}^\droite_\veps$ (\resp $\mathcal{S}^\droite_\veps$).

\subsection{Local cell problems and reconstruction formula}
We introduce two problems set in $\cellper^{\droite, 0}\!\! \coloneqq  (0, 1)^2$, and given for a boundary data $\varphi \in L^2_{per}(\R)$: \textit{Find $E_\veps^{\droite, 0} (\varphi) \in H^1_{\cutvec} (\cellper^{\droite, 0})$ such that}
\begin{equation} \label{cellpb0}
	\left|
	 \begin{array}{r@{\ =\ }l} 
			\displaystyle - \Dt{} \big( \mu_p \, \Dt{} E_\veps^{\droite, 0} (\varphi) \big) - \rho_p \, (\omega^2+ \icplx \veps) \, E_\veps^{\droite, 0} (\varphi) & 0 \quad \textnormal{in }\,\cellper^{\droite, 0}, 
			\ret
			\multicolumn{2}{c}{\displaystyle \mathcal{R}^\droite_+ \,  E_\veps^{\droite, 0} (\varphi)\Big|_{\Sigmaper^{\droite, 0}} = \varphi \quad \textnormal{and} \quad \mathcal{R}^\droite_- \,  E_\veps^{\droite, 0} (\varphi)\Big|_{\Sigmaper^{\droite, 1}} = 0,}
			\ret
			\multicolumn{2}{c}{\textnormal{$E_\veps^{\droite, 0} (\varphi)$ is periodic w.r.t. $\yvi_1$;}}
	\end{array}
	\right.
\end{equation}
and \textit{Find $E_\veps^{\droite, 1} (\varphi) \in H^1_{\cutvec} (\cellper^{\droite, 0})$ such that}
\begin{equation} \label{cellpb1}
	\left|
	 \begin{array}{r@{\ =\ }l} 
			\displaystyle - \Dt{} \big( \mu_p \, \Dt{} E_\veps^{\droite, 1} (\varphi) \big) - \rho_p \, (\omega^2+ \icplx \veps) \, E_\veps^{\droite, 1} (\varphi) & 0 \quad \textnormal{in }\,\cellper^{\droite, 0}, 
			\ret
			\multicolumn{2}{c}{\displaystyle \mathcal{R}^\droite_+ \,  E_\veps^{\droite, 1} (\varphi)\Big|_{\Sigmaper^{\droite, 0}} = 0 \quad \textnormal{and} \quad \mathcal{R}^\droite_- \,  E_\veps^{\droite, 1} (\varphi)\Big|_{\Sigmaper^{\droite, 1}} = \varphi,}
			\ret
			\multicolumn{2}{c}{\textnormal{$E_\veps^{\droite, 1} (\varphi)$ is periodic w.r.t. $\yvi_1$.}}
	\end{array}
	\right.
\end{equation}
Using Lax-Milgram lemma, one easily obtains the following result.
\begin{prop}
	Let $\veps>0$. The problems \eqref{cellpb0} and \eqref{cellpb1} are well posed. 
\end{prop}
	
\vspspe	
In the sequel, for any $\varphi \in L^2_{per}(\R)$, \textit{we identify} $E_\veps^{\droite, 0} (\varphi)$ \textit{and} $E_\veps^{\droite, 1} (\varphi)$ \textit{with their periodic extensions with respect to $y_1$}. 
	
\vspspe	
By linearity, if $U \in H^1_{\theta} (\cellper^{\droite, 0})$ is periodic with respect to $\yvi_1$ and satisfies 
\begin{equation}\label{hom_pde} 
	- \Dt{} \big( \mu_p \, \Dt{}U  \big) - \rho_p \; (\omega^2 + \icplx \veps) \; U = 0 \quad \mbox{in } \cellper^{\droite, 0},
\end{equation}
then
\begin{equation} \label{reconstruction}
	\mathcal{R}_+^\droite U\big|_{\Sigmaper^{\droite, 0}} = \varphi \quad \mbox{and} \quad \mathcal{R}_-^\droite U\big|_{\Sigmaper^{\droite, 1}} = \psi\quad \Longrightarrow \quad U =  E_\veps^{\droite, 0}  (\varphi) + E_\veps^{\droite, 1}  (\psi).
\end{equation}
From this remark, by definition \eqref{defPS} of $(\mathcal{P}^\droite_\veps, \mathcal{S}^\droite_\veps)$ and by Proposition \ref{prop:property_PS}, it is clear that the solution $U^\droite_\veps(\varphi)$ of \eqref{RobinHalfspaceproblems} can be reconstructed cell by cell using $E_\veps^{\droite, 0} (\varphi)$ and $E_\veps^{\droite, 1} (\varphi)$. 
\begin{prop}\label{thm_reconstruction}
	For any $\varphi \in L^2_{per}(\R)$, the solution $U^\droite_\veps (\varphi)$ of \eqref{RobinHalfspaceproblems} is given by
	\begin{equation*}
	\aeforall \yv \in \cellper^{\droite, 0}, \quad \big[U^\droite_\veps (\varphi)\big](\yv + n\, \vv{\ev}_2) = \big[ E_\veps^{\droite, 0} ((\mathcal{P}^\droite_\veps)^n\,\varphi)  + E_\veps^{\droite, 1} (\mathcal{S}^\droite_\veps\, (\mathcal{P}^\droite_\veps)^{n}\varphi)\big] (\yv), \quad \spforall n \in \N,
	\end{equation*}
	where the operators $\mathcal{P}^\droite_\veps$ and $\mathcal{S}^\droite_\veps$ are defined in \eqref{defPS}.
\end{prop}

\vspspe
Similarly to $U_\veps^\droite$, the solutions $(E^{\droite, 0}_\veps, E^{\droite, 1}_\veps)$ of Problems \eqref{cellpb0} and \eqref{cellpb1} benefit from a quasi-1D structure. Indeed, let us introduce the two families $(e_{s, \veps}^{\droite, 0})_{s \in \R}$ and $(e_{s, \veps}^{\droite, 1})_{s \in \R}$ of functions in $H^1(0, 1/\cuti_2)$ which satisfy the 1D problems 
\begin{equation}  \label{segment_problems0} 
	\left|
	\begin{array}{lll}
	\displaystyle - \, \frac{d}{dx} \Big(\mu_{s, \cutvec} \,  \frac{d e_{s, \veps}^{\droite, 0}}{dx}\Big)-  \rho_{s, \cutvec} \, (\omega^2 + \icplx \veps) \, e_{s, \veps}^{\droite, 0} = 0 \quad \mbox{in}\ \ (0, 1/\cuti_2), 
	\ret
	R^\droite_{+,s} e_{s, \veps}^{\droite, 0}(0) = 1 \quad \textnormal{and} \quad  R^\droite_{-,s} e_{s, \veps}^{\droite, 0}(1/\cuti_2) = 0,
	\end{array}	
	\right. 
\end{equation}
and
\begin{equation}  \label{segment_problems1}
	\left| \begin{array}{lll}
\displaystyle 	- \, \frac{d}{dx} \Big(\mu_{s, \cutvec} \,  \frac{d e_{s, \veps}^{\droite, 1}}{dx}\Big)-  \rho_{s, \cutvec} \, (\omega^2 + \icplx \veps) \,  e_{s, \veps}^{\droite, 1} = 0 \quad \mbox{in}\ \ (0, 1/\cuti_2), \ret
R^\droite_{+,s} e_{s, \veps}^{\droite, 1}(0) = 0 \quad \textnormal{and} \quad  R^\droite_{-,s} e_{s, \veps}^{\droite, 1}(1/\cuti_2) = 1
\end{array}	\right.
\end{equation}
We then have the next property, which shows the fibered structure of the cell problems.
\begin{prop}\label{thm_fibrage_cells}
	Let $\veps> 0$. The solutions $E_\veps^{\droite, 0} (\varphi)$ and $E_\veps^{\droite, 1} (\varphi)$ of the local cell problems \eqref{cellpb0} and  \eqref{cellpb1} and the families $(e_{s, \veps}^{\droite, 0})_{s \in \R}$ and $(e_{s, \veps}^{\droite, 1})_{s \in \R}$ defined by \eqref{segment_problems0} and  \eqref{segment_problems1} are related by:
	\begin{equation} \label{fibered-structure-cell} 
		\left\{ 
			\begin{array}{l@{\quad}l}
				E_\veps^{\droite, 0} (\varphi) (s + x \, \theta_1, x \, \theta_2) = 	\varphi(s) \, e_{s, \veps}^{\droite, 0} (x), &  x \in \; (0, 1/\cuti_2), 
			\; s \in \R , 
			\ret
			E_\veps^{\droite, 1} (\varphi)(s + x \, \theta_1, x \, \theta_2) = 	\varphi(s + \delta) \, e_{s, \veps}^{\droite, 1} (x), &  x \in \; (0, 1/\cuti_2), 
				\; s \in \R. 
	\end{array} \right. 
\end{equation}
where $s\mapsto e_{s, \veps}^{\droite, j},\;j\in\{0,1\}$ is $1$-periodic.
\end{prop} 
\noindent 
Hence, the evaluation of $E_\veps^{\droite, 0} (\varphi)$ and $E_\veps^{\droite, 1} (\varphi)$ simply consist in solving a family of 1D-problems posed on a segment, namely (\ref{segment_problems0}) and  (\ref{segment_problems1}). %
To conclude this section, we provide a 1D equivalent of the cell by cell reconstruction formula given in Proposition \ref{thm_reconstruction} in terms of the half-line solutions $u^\droite_{s, \veps} \in H^1(I^\droite)$.
\begin{cor}
	Given $\varphi \in \mathscr{C}^0_{\textit{per}}(\R)$ such that $\varphi(s) \neq 0$ for any $s \in \R$, define
	\begin{equation*}
		\displaystyle
		\spforall n \in \N, \quad \varphi^\veps_n \coloneqq  (\mathcal{P}^\droite_\veps)^n\, \varphi \quad \textnormal{and} \quad \psi^\veps_n \coloneqq  \mathcal{S}^\droite_\veps\, (\mathcal{P}^\droite_\veps)^n\, \varphi.
	\end{equation*}
	Then for any $s \in \R$, the solution $u^\droite_{s, \veps}$ of the half-line problem \eqref{halfline_sproblems} is given for $\xvi \in (0, 1/\cuti_2)$ by
	\begin{equation}\label{eq:cell_by_cell_1D}
		\displaystyle
		\varphi(s)\, u^\droite_{s, \veps}(\xvi + n/\cuti_2) = \varphi^\veps_n (s + n \delta)\, e^{\droite, 0}_{s + n\delta, \veps} (\xvi) + \psi^\veps_n (s + (n + 1) \delta)\, e^{\droite, 1}_{s + n\delta, \veps} (\xvi).
	\end{equation}
\end{cor}

\begin{dem}
	If $\xvi \in (0, 1/\cuti_2)$, then $\yv \coloneqq  (s + \cuti_1\, \xvi, \cuti_2\, \xvi)$ belongs to $\R \times (0, 1)$, so that by Proposition \ref{thm_reconstruction},
	\begin{equation*}
		\displaystyle
		\big[U^\droite_\veps (\varphi)\big](s + \cuti_1\, \xvi,\; \cuti_2\, \xvi + n) = \big[ E_\veps^{\droite, 0} ((\mathcal{P}^\droite_\veps)^n\,\varphi) + E_\veps^{\droite, 1} (\mathcal{S}^\droite_\veps\, (\mathcal{P}^\droite_\veps)^{n}\varphi)\big] (s + \cuti_1\, \xvi, \cuti_2\, \xvi),
	\end{equation*}
	where $E_\veps^{\droite, 0} ((\mathcal{P}^\droite_\veps)^n\,\varphi)$, and $E_\veps^{\droite, 1} (\mathcal{S}^\droite_\veps\, (\mathcal{P}^\droite_\veps)^{n}\varphi)$ have been identified with their periodic extensions with respect to $\yvi_1$. By using the fibered structure of $U^\droite_\veps$ \eqref{defUeps} and $E_\veps^{\droite, 0}, E_\veps^{\droite, 1}$ (Proposition \ref{thm_fibrage_cells}), the above equality becomes
	\begin{equation*}
		\displaystyle
		\varphi(s - n\delta)\, u^\droite_{ s - n\delta,\veps} (\xvi + n/\cuti_2) = [(\mathcal{P}^\droite_\veps)^n\,\varphi] (s)\, e^{\droite, 0}_{s, \veps} (\xvi) + [\mathcal{S}^\droite_\veps\, (\mathcal{P}^\droite_\veps)^{n}\, \varphi] (s + \delta)\, e^{\droite, 1}_{s, \veps} (\xvi),
	\end{equation*}
	which corresponds to the desired equality, up to a translation $s \mapsto s + n\delta$.
\end{dem}

\noindent
The only missing
step for evaluating  $U_\veps^\droite (\varphi)$ via the reconstruction formula given in Proposition \ref{thm_reconstruction} is the determination of the operators $(\mathcal{P}^\droite_\veps, \mathcal{S}^\droite_\veps)$. This is the object of the next section.

\subsection{Local RtR operators and  Riccati system}
The cell by cell reconstruction formula given in Proposition \ref{thm_reconstruction} only ensures that $U_{\veps}^\droite (\varphi)$ satisfies the partial differential equation \eqref{pde} in each cell $\cellper^{\droite, n}$. For \eqref{pde} to be satisfied in the whole half-space $\Omega^\droite$ (or equivalently in the whole half-guide $\Omegaper^\droite$), we have to ensure that \eqref{pde} is satisfied across each interface $\Sigmaper^{\droite, n}$, which corresponds to the continuity of $U_\veps^\droite (\varphi)$ and $\mu_p\, \Dt{} U_\veps^\droite (\varphi)$
across $\Sigmaper^{\droite, n}$, or equivalently the continuity of the two Robin traces $\mathcal{R}_+^\droite U_\veps^\droite (\varphi)$ and $\mathcal{R}_-^\droite  U_\veps^\droite (\varphi)$. It is by writing these two continuity conditions that we will obtain two equations for $(\mathcal{P}^\droite_\veps,\mathcal{S}^\droite_\veps)$.

\vspspe
To this end, we introduce for $\veps> 0$ four local \emph{Robin-to-Robin} (RtR) operators $\mathcal{T}^{\droite, \ell k}_\veps \in \mathscr{L}\big( L^2_{\textit{per}} (\R))$ defined for $\ell, k \in \{0, 1\}$ by
\begin{equation} \label{defRtR}
	\begin{array}{r@{\ \coloneqq \ }l@{\qquad \textnormal{and} \qquad}r@{\ \coloneqq \ }l}	
	\mathcal{T}^{\droite, 00}_\veps\, \varphi & \mathcal{R}_-^\droite E_\veps^{\droite, 0} (\varphi)\Big|_{\Sigma^{\droite, 0}} &  	\mathcal{T}^{\droite, 01}_\veps\, \varphi & \mathcal{R}_+^\droite E_\veps^{\droite, 0} (\varphi)\Big|_{\Sigma^{\droite, 1}},
	\ret
	\mathcal{T}^{\droite, 10}_\veps\, \varphi & \mathcal{R}_-^\droite E_\veps^{\droite, 1} (\varphi)\Big|_{\Sigma^{\droite, 0}} &  	\mathcal{T}^{\droite, 11}_\veps\, \varphi & \mathcal{R}_+^\droite E_\veps^{\droite, 1} (\varphi)\Big|_{\Sigma^{\droite, 1}}. 	
		\end{array} 
\end{equation}

\begin{rmk} \label{rempm} As the understanding of the manipulations of $\pm$ for the Robin differential operators is not immediate, it may be helpful to remember that in the cell problems we impose the {\it outgoing} Robin trace (outgoing with respect to the periodicity cell) while with the operators $\mathcal{T}^{\droite, \ell k}_\veps$, we evaluate the {\it ingoing} Robin trace (ingoing with respect to the periodicity cell).
\end{rmk} 

\vspspe
Proceeding as in the proof of Proposition \ref{thm_weighted-shift}, and using the fibered structure of the solutions $E^{\droite, 0}_\veps$ and $E^{\droite, 1}_\veps$ (see Theorem \ref{thm_fibrage_cells}), one easily proves the following result.
\begin{prop}\label{thm_weighted-shift2}
	Let $\veps> 0$. The operators $\mathcal{T}^{\droite, \ell k}_\veps
	$, $\ell,k\in\{0,1\}$, are weighted shift operators defined for any $\varphi \in L^2_{\textit{per}}(\R)$ by
	\begin{equation} \label{RtRweighted}
		\spforall s\in\R, \quad \left\{
		\begin{array}{r@{\ =\ }l@{ \quad }r@{\ =\ }l}	
			\mathcal{T}^{\droite, 00}_\veps\, \varphi(s) &	\boldsymbol{t}^{\droite, 00}_{\veps}(s) \, \varphi(s), & \quad \mathcal{T}^{\droite, 01}_\veps\, \varphi(s) &	\boldsymbol{t}^{\droite, 01}_{\veps}(s - \delta)\, \varphi(s - \delta), 
			\rett
			\mathcal{T}^{\droite, 10}_\veps\, \varphi(s) & \boldsymbol{t}^{\droite, 10}_{\veps}(s) \, \varphi(s + \delta), &  \quad \mathcal{T}^{\droite, 11}_\veps\, \varphi(s) &	\boldsymbol{t}^{\droite, 11}_{\veps}(s - \delta) \, \varphi(s),
		\end{array} 
		\right. 
	\end{equation}
	where the RtR symbols $\boldsymbol{t}^{\droite, k\ell}_{\veps}(s)$ are defined from $e_{s, \veps}^{\droite, 0}$ and $e_{s, \veps}^{\droite, 1}$ (see \eqref{segment_problems0} and \eqref{segment_problems1}) by:
	\begin{equation} \label{defRtRcoeff}
		\spforall s\in\R, \quad\left\{ \begin{array}{r@{\ =\ }l@{ \qquad }r@{\ =\ }l}	
			\boldsymbol{t}^{\droite, 00}_{\veps}(s) & R_{-,s}^\droite e_{s, \veps}^{\droite, 0}(0), &%
			\boldsymbol{t}^{\droite, 01}_{\veps}(s) & R_{+,s}^\droite e_{s, \veps}^{\droite, 0}(1/\theta_2),
			\ret
			\boldsymbol{t}^{\droite, 10}_{\veps}(s) & R_{-,s}^\droite e_{s, \veps}^{\droite, 1}(0), &%
			\boldsymbol{t}^{\droite, 11}_{\veps}(s) & R_{+,s}^\droite e_{s, \veps}^{\droite, 1}(1/\theta_2).
		\end{array} \right. 
	\end{equation}
\end{prop}

\vspspe
If $B \in \mathscr{L}\big(L_{per}^2(\R)\big)$, its transpose $\transp{B} \in \mathscr{L}\big(L_{per}^2(\R)\big)$ is defined by
\begin{eqnarray}
\spforall (\varphi, \psi)	\in L_{per}^2(\R)^2, \quad \int_{0}^1 \transp{B} \varphi(s) \, \psi(s) \, ds  =  \int_{0}^1  \varphi(s) \, B \psi(s) \, ds,
\end{eqnarray}
and is related to its adjoint $B^*$ by conjugation: $\transp{B} \varphi = \overline{B^* \overline{\varphi} }$. Moreover, $B$ and $\transp{B}$ have the same spectrum.

\vspspe
The operators $\mathcal{T}^{\droite, \ell k}_\veps$ are related to one another by transposition properties that reflect the formal symmetry of the operator $\Dt{} \big( \mu_p \Dt{}\big)$, which is closely related to the well-known reciprocity principle in wave propagation.
\begin{lem}\label{lem_transposition_properties}
	Let $\veps> 0$. We have
	\begin{equation}\label{transposition_properties}
		\transp{\mathcal{T}^{\droite, 00}_\veps} = {\mathcal{T}^{\droite, 00}_\veps}, \quad \transp{\mathcal{T}^{\droite, 11}_\veps} = \mathcal{T}^{\droite, 11}_\veps, \quad 
\transp{\mathcal{T}^{\droite, 01}_\veps} = \mathcal{T}^{\droite, 10}_\veps.
\end{equation}
\end{lem} 
\begin{dem} If $U, V \in H^1_{\cutvec}(\cellper^{\droite, 0})$ are periodic with respect to $y_1$ and satisfy the PDE \eqref{hom_pde} in $\cellper^{\droite, 0}$, then using a Green's formula and the periodicity properties, it comes
\[\displaystyle
	\int_{\Sigmaper^{\droite, 0}} \big( \mu_p\,V\, \Dt{} U \,  - \mu_p\,U\, \Dt{} V \, \big) =	\int_{\Sigmaper^{\droite, 1}}  \big( \mu_p\,V\, \Dt{} U \, - \mu_p\, U\, \Dt{} V \,  \big) .
\] 
As $2 \nu^\droite \mu_p\, \Dt{} U = \mathcal{R}_+^\droite U- \mathcal{R}_-^\droite U$ and $ -2\vts\icplx\vts z U =  \mathcal{R}_+^\droite U+ \mathcal{R}_-^\droite U$, and the same replacing $U$ by $V$, the above rewrites (we omit the details):
%	$$
%	\int_{\Sigma_{0,\periodic}^\pm}  \big(R_\theta^\pm U \,  R_\theta^\mp V  - R_\theta^\pm V \,  R_\theta^\mp U  \big)= 	\int_{\Sigma_{1,\periodic}^\pm}  \big(R_\theta^\pm U \,  R_\theta^\mp V  - R_\theta^\pm V \,  R_\theta^\mp U  \big)
%			$$
			$$
			\int_{\Sigmaper^{\droite, 0}}  \big(\mathcal{R}_+^\droite U \,  \mathcal{R}_-^\droite V  - \mathcal{R}_-^\droite U \, \mathcal{R}_+^\droite V  \big)= 	\int_{\Sigmaper^{\droite, 1}}  \big(\mathcal{R}_+^\droite U \,  \mathcal{R}_-^\droite V  - \mathcal{R}_-^\droite U \, \mathcal{R}_+^\droite V  \big).
			$$
	With $U= E_\veps^{\droite, 0} (\varphi)$ and $V= E_\veps^{\droite, 1} (\psi)$, the above identity gives, identifying $\Sigmaper^{\droite, 0}$ and $\Sigmaper^{\droite, 1}$,
	$$
	\int_{\Sigmaper^{\droite, 0}} (\mathcal{T}^{\droite, 10}_\veps \psi) \, \varphi  =  \int_{\Sigmaper^{\droite, 0}} (\mathcal{T}^{\droite, 01}_\veps \varphi) \, \psi , \quad \mbox{which yields } \transp{\mathcal{T}^{\droite, 10}_\veps} = \mathcal{T}^{\droite, 01}_\veps.
	$$
	To obtain $\transp{\mathcal{T}^{\droite, 00}_\veps} = {\mathcal{T}^{\droite, 00}_\veps}$, it suffices to take $U= E_\veps^{\droite, 0} (\varphi)$ and $V= E_\veps^{\droite, 0} (\psi)$ and to get $\transp{\mathcal{T}^{\droite, 11}_\veps} = {\mathcal{T}^{\droite, 11}_\veps}$, it suffices to take $U= E_\veps^{\droite, 1} (\varphi)$ and $V= E_\veps^{\droite, 1} (\psi)$.
\end{dem}

\vspspe
Using the cell by cell expression (see Proposition \eqref{thm_reconstruction}) of the solution of \eqref{RobinHalfspaceproblems} in terms of the cell solutions defined in \eqref{cellpb0} and \eqref{cellpb1}, we can show the following result.
\begin{prop} Let $\veps>0$. The pair $(\mathcal{P}^\droite_\veps,\mathcal{S}^\droite_\veps)$ is the unique solution of the problem
	\begin{equation}  \label{Riccatiplus}
		%\left| \begin{array}{l}
			\mbox{Find } (P,S) \in \mathscr{L}(L^2_{per}(\R))^2 \mbox{ such that } \varrho(P) < 1 \mbox{ and }
			%\ret
			\quad \left\{ \begin{array}{llll}
				P=\mathcal{T}^{\droite, 01}_\veps + \mathcal{T}^{\droite, 11}_\veps S  ,
				\ret
				S   = \mathcal{T}^{\droite, 00}_\veps P+ \mathcal{T}^{\droite, 10}_\veps S \, P.
			\end{array}  \right.
		%\end{array}\right.
	\end{equation}
\end{prop}
\begin{dem} We only sketch the proof since it is similar to the one of \cite[Proposition 4.4]{amenoagbadji2023wave}.
\noindent Let us first show that the pair $(\mathcal{P}^\droite_\veps,\mathcal{S}^\droite_\veps)$ is solution of \eqref{Riccatiplus}. Let $\varphi \in L^2_{per}(\R)$ and $U_\veps^\droite(\varphi)$ be the unique solution of \eqref{RobinHalfspaceproblems}. The cell by cell expression of $U^\droite_\veps$ given in Proposition \eqref{thm_reconstruction}, taken for $n=0$ and $n = 1$ gives:
\begin{equation} \label{rec1}
	\begin{array}{c}
 [U_{\veps}^\droite(\varphi)](\yv) = [E_\veps^{\droite, 0} (\varphi) + E_\veps^{\droite, 1} \big( \mathcal{S}^\droite_\veps\, \varphi \big)](\yv)\quad  \yv\in \cellper^{\droite, 0},
\ret  [U_{\veps}^\droite(\varphi)](\yv) = [E_\veps^{\droite, 0} (\mathcal{P}^\droite_\veps\, \varphi) + E_\veps^{\droite, 1} \big( \mathcal{S}^\droite_\veps \mathcal{P}^\droite_\veps\,  \varphi \big)] (\yv - \vv{\ev}_2)\quad \yv\in \cellper^{\droite,1}.
	\end{array}
\end{equation}
By writing that $\mathcal{R}_+^\droite U_{\veps}^\droite(\varphi)$ and $\mathcal{R}_-^\droite U_{\veps}^\droite(\varphi)$ are continuous across the interface $\Sigmaper^{\droite, 1}$, we obtain that $(\mathcal{P}^\droite_\veps,\mathcal{S}^\droite_\veps)$ satisfies respectively the first and second equations of \eqref{Riccatiplus}.
\\\\
\noindent The fact that $\varrho(\mathcal{P}^\droite_\veps) < 1$ is a consequence of the exponential decay of the solution in the $y_2$ direction  (with a decay rate that is proportional to $\veps$). We refer to the proof of \cite[Proposition 4.1]{amenoagbadji2023wave} for more details.
\\\\ 
 \noindent For the uniqueness result, we show that if $(P,S)$ is a solution of \eqref{Riccatiplus}, then for any $\varphi \in L^2_{per}(\R)$, the function $U(\varphi)$ defined cell by cell in the half-guide $\Omegaper^\droite$ by (we mimic the formula given in Proposition \ref{thm_reconstruction})
	$$
	\spforall \yv \in \cellper^{\droite, n}, \quad \big[U(\varphi)\big](\yv) = \big[ E_\veps^{\droite, 0} (P^n\varphi)  + E_\veps^{\droite, 1} (S P^n\varphi)\big] (\yv -n \vv{\ev}_2)
	$$
	and extended by periodicity in the $y_1$-direction, belongs to $H^1_{\cutvec, \textit{per}}(\Omega)$ {because $\varrho(P)<1$} and is solution of the half-space problem \eqref{RobinHalfspaceproblems}. From the well-posedness of \eqref{RobinHalfspaceproblems}, we deduce the equality $U(\varphi) = U_\veps^\pm(\varphi)$, which implies that $P \varphi = \mathcal{P}^\droite_\veps \varphi$ and  $S \varphi = \mathcal{S}^\droite_\veps \varphi$. 
	\end{dem}

\begin{rmk} 
	We call \eqref{Riccatiplus} a \textbf{Riccati system} because if one formally eliminates $\mathcal{S}^\droite_\veps$ in the system, one obtains a quadratic equation in $\mathcal{P}^\droite_\veps$, as in the Dirichlet case. Note however that eliminating $\mathcal{S}^\droite_\veps$ involves inverting $\mathcal{T}^{\droite, 11}_\veps$ which is not necessarily invertible. In fact, if $\mu_p = \rho_p = 1$ for instance, and if $z$ is allowed to be complex with $z \coloneqq  \sqrt{\omega^2 + \icplx\vts \veps}$ with the convention that $\Real \sqrt{\phantom{2}} > 0$, then it can be computed without difficulty that $\mathcal{T}^{\droite, 11}_\veps = 0$. 
	
	\vspspe	In the proof we only used the continuity conditions across $\Sigmaper^{\droite, 1}$ to derive the Riccati system \eqref{Riccatiplus}. However, it is easy to check that if \eqref{Riccatiplus} is satisfied, the same continuity conditions are satisfied across any $\Sigmaper^{\droite,n}$.
\end{rmk}

\subsection{RtR half-space operator and transparent boundary condition} \label{sec_RtR}
In link with the RtR  coefficient $\RtRcoef_\veps^\droite$ defined in  \eqref{Robin_coeff}, we introduce for $\veps>0$ the RtR half-space operator $\RtRop_\veps^\droite \in \mathscr{L}\big(L^2_{per}(\R)\big) $ defined by
\begin{equation} \label{def_Lambda}
	\RtRop_\veps^\droite \varphi \coloneqq  {\cal R}^{\droite}_- U_\veps^\droite (\varphi)|_{\Sigma^{\droite, 0}}.
	\end{equation}
From the cell by cell reconstruction formula given in Proposition \ref{thm_reconstruction} and the definition \eqref{defRtR} of the RtR operators, we deduce 
\begin{equation} \label{form_Lambda}
	\RtRop_\veps^\droite =\mathcal{T}^{\droite, 00}_\veps + \mathcal{T}^{\droite, 10}_\veps \mathcal{S}^\droite_\veps.
\end{equation}
Using the same arguments as in the proof of Lemma \ref{lem_transposition_properties}, we can show that
\begin{equation} \label{transpose_Lambda}
	\transp{\RtRop_\veps^\droite} = \RtRop_\veps^\droite.
\end{equation}
Moreover, using a Green's formula in \eqref{RobinHalfguideproblems}, we can show that there exists $c_\eps>0$  such that 
\[
	\spforall \varphi\in L^2_{\textit{per}}(\R),\quad
		\Imag\int_{\Sigmaper^{\droite, 0}} \mu_p\, \Dt{} U^\droite_\veps (\varphi)\; \overline{U^\droite_\veps (\varphi)} \geq c_\eps\, \|\varphi\|^2_{L^2(\Sigmaper^{\droite, 0})}
\]
which implies, using the expressions $-2 \mu_p\, \Dt{}U=({\cal R}^\droite_{+}-{\cal R}^\droite_{-})\, U$ and $-2\vts \icplx\vts z\, U=({\cal R}^\droite_{+}+{\cal R}^\droite_{-})\,  U$, as well as the definition \eqref{defRtR} of the RtR operators, that
\begin{equation}\label{eq:flux_eps}
	 \spforall \varphi\in L^2_{\textit{per}}(\R),\quad\frac{1}{4z}\Real\, \int_{\Sigmaper^{\droite, 0}} (I-\RtRop_\veps^\droite)\, \varphi\; \overline{(I+\RtRop_\veps^\droite)\, \varphi} \geq c_\eps\, \|\varphi\|^2_{L^2(\Sigmaper^{\droite, 0})}.
\end{equation}
From the fibered structure \eqref{defUeps} of $U^\droite_\veps$, it is immediate to see that $\RtRop_\veps^\droite$ is a multiplication operator 
\begin{equation} \label{prop_Lambda}
	\big(\RtRop_\veps^\droite \varphi\big) (s) = \RtRcoef_\veps^\droite(s)\, \varphi(s) \quad \textnormal{with} \quad  \RtRcoef_\veps^\droite(s) \coloneqq  R_{-,s}^\droite  u^\droite_{s,\veps}(0),
\end{equation}
and where $u^\droite_{s,\veps}$ is defined for any $s \in \R$ by the half-line problem \eqref{halfline_sproblems}. 
In particular, similarly to \eqref{cutformulas}, the RtR coefficient $\RtRcoef_\veps^\droite$ can be computed thanks to
\begin{equation} \label{def_lambda_1D}
\RtRcoef_\veps^\droite  = \RtRcoef_\veps^\droite(\cuti_1\, a^\droite) =	\big(\RtRop_\veps^\droite \varphi\big) (\cuti_1\, a^\droite) \; \mbox{ for any } \varphi \in \mathscr{C}^0_{per}(\R) \mbox{ satisfying } \varphi(\cuti_1\, a^\droite) = 1,
\end{equation}
that can be rewritten, using Propositions \ref{thm_weighted-shift} and \ref{RtRweighted}, as
\begin{equation} \label{def_lambda_1D_bis}
\RtRcoef_\veps^\droite  = [\boldsymbol{t}^{\droite, 00}_{\veps}+\boldsymbol{t}^{\droite, 10}_{\veps}\,\boldsymbol{s}^\droite_\veps](\cuti_1\, a^\droite) .
\end{equation}

\section{Analysis of the Riccati system}\label{sec:Riccati_eps}
In this section, we provide some spectral results that will be used to compute the propagation operator and the scattering operator. Section \ref{sec:spectrum_propagation_operator} is devoted to describing the spectral structure of the propagation operator in terms of its eigenpairs, {provided that $\delta$ is not a Liouville number and if $(\rho_p,\mu_p)$ are smooth enough (the required smoothness depending on the measure of irrationality of $\delta$)}. In particular, our analysis will rely heavily on the notion of \textbf{fundamental eigenpair} which we shall introduce below. We then reformulate the Riccati system \eqref{Riccatiplus} in terms of these eigenpairs in Section \ref{sec:spectrum_Riccati}. 

\subsection{Spectral description of the propagation operators}\label{sec:spectrum_propagation_operator}
The spectral theory of weighted shift operators has been thoroughly studied and detailed in \cite{antonevich2012linear}. In particular, it is shown that for all $\delta\in\R\setminus\Q$, the spectrum of $\mathcal{P}^\droite_\veps$ is given by
\begin{equation}\label{eq:spectre_Peps} \sigma(\mathcal{P}^\droite_\veps)= \mathcal{C}(0,\varrho(\mathcal{P}^\droite_\veps))\quad\text{with}\quad\varrho(\mathcal{P}^\droite_\veps)=\exp\left(\int_0^1\log|  \boldsymbol{p}^\droite_\veps (s)|\,ds\right).
\end{equation}
Let us now describe more precisely the spectral properties of $\mathcal{P}^\droite_\veps$ in the case where the ratio $\delta \coloneqq  \cuti_1/\cuti_2$ is not a Liouville number, so that its irrationality measure $m(\delta)$ (defined by \eqref{eq:meas_irrat}) is finite and the Diophantine condition \eqref{eq:notLiouville} holds. 
% namely
% \begin{equation} \label{diophantine}
% 	\spforall \nu > m(\delta), \quad \spexists c \equiv c(\delta, \nu) > 0\ \ \big/\ \ \spforall (p,q) \in \Z \times \N^*, \quad \left|\delta - \frac{p}{q}\right| > \frac{c}{q^\nu}.
% 	% m(\delta) < + \infty \quad \mbox{(with $m(\delta)$ the measure of irrationality defined by \eqref{eq:meas_irrat})}
% \end{equation}
\\\\
\noindent We rely on the notion of \emph{winding number}, also called index; see \cite[Chapter 7]{beardon1979complex}.%\sfnote{Rajouter une ref}

\begin{defi}\label{def_windingnumber}
	 Let ${\cal N}_{\textit{per}}$ be the open subset of $\mathscr{C}_{\textit{per}}(\R; \C)$ defined by
	\begin{equation} \label{defN01}
		{\cal N}_{\textit{per}} \coloneqq  \big\{ f \in \mathscr{C}_{\textit{per}}(\R; \C) \mbox{ such that } f(s) \neq 0, \spforall s \in \R \big \}.
	\end{equation}
	Given $f \in {\cal N}_{\textit{per}}$, there exists a function $s\mapsto \argument f(s) \in \mathscr{C}^0(\R;\R)$, where $\argument f(s)$ denotes the argument of $f(s)$, such that 
	\begin{equation} \label{argument}
		\spforall s \in [0,1], \quad f(s) = |f(s)| \, \euler^{\icplx \argument f(s)} \quad 
	\end{equation} 
	and one can define the {\it winding number} ${\bf w}(f)$ of $f$ as 
	\begin{equation} \label{windingnumber}
		{\bf w}(f) = \big(\argument f(1) -\argument f(0)\big)/2\pi \in \Z.
	\end{equation} 
	Moreover, the map $f \in {\cal N}_{\textit{per}} \mapsto \argument f|_{[0,1]} \in \mathscr{C}^0(0,1)$ is continuous. 
\end{defi}

\noindent
In more geometrical (but less precise) terms, ${\bf w}(f)$ is the number of turns around the origin made by the (closed) trajectory in the complex plane described by $f(s)$ when $s$ varies from $0$ to $1$ (hence its name). The turn counts positively (respectively negatively) when the trajectory travels counterclockwise (respectively clockwise) around the origin. It is easy to show the following result.

\begin{cor}\label{coro_WN0} 
	For any $f \in {\cal N}_{\textit{per}}$, one has the equivalences 
	$$
			{\bf w}(f)= 0 \quad \Longleftrightarrow \quad \argument f \in \mathscr{C}_{\textit{per}}^0(\R; \R) \quad \Longleftrightarrow \quad f(s) = \euler^{g(s)}, \quad g \in \mathscr{C}_{\textit{per}}(\R; \C).
	$$
\end{cor}

\begin{lem}\label{lem_WN}
	% Let $\veps > 0$ and $\mathbf{p}^\droite_\veps (x, s)$ be the function defined by \eqref{defsymbols}, then 
	% $$
	% \spforall x \in \R_+, \quad s\mapsto\mathbf{p}^\droite_\veps(x, s) \in {\cal N}_{\textit{per}} \quad \textnormal{and} \quad {\bf w} (\mathbf{p}^\droite_\veps (x,\cdot) ) = 0.
	% $$
	{Let $\veps > 0$ and $s\mapsto \boldsymbol{p}^\droite_\veps(s)$ be the function defined by \eqref{defsymbols}, then \[\boldsymbol{p}^\droite_\veps \in {\cal N}_{\textit{per}}\mbox{ and } {\bf w}\big( \boldsymbol{p}^\droite_\veps ) = 0.\]}
\end{lem}
\begin{dem} 
	{From its definition \eqref{defsymbols} and Proposition \ref{prop:propshift}, it is clear that the function $\boldsymbol{p}^\droite_\veps$ belongs to ${\cal N}_{\textit{per}}$ and thus has a well-defined winding number ${\bf w}\big(\boldsymbol{p}^\droite_\veps\big)$.
    \\\\ Let us introduce the function
    \[ 
    \spforall s \in \R, \quad \aeforall x \in \R_+ , \quad 
{\bf p}^\droite_\veps (s,x) \coloneqq  \big(R^\droite_{+,s} \, u_{s,\veps}^\droite \big) (x).
    \]
    Note that $\boldsymbol{p}^\droite_\veps(s)\equiv {\bf p}^\droite_\veps (s,1/\cuti_2)$. Using a similar argument to that used in the proof of Proposition \ref{prop:propshift}, one shows easily that
   $\mathbf{p}^\droite_\veps$ is continuous with respect to $s$ and $x$. We deduce from Proposition \ref{def_windingnumber} that $x\mapsto \textnormal{Arg}\,\mathbf{p}^\droite_\veps(x,\cdot)$ is continuous. As a consequence, $\xvi \mapsto {\bf w}\big(\mathbf{p}^\droite_\veps(x,\cdot))$ is a continuous and integer-valued function, which means that it is constant. In particular, ${\bf w}\big(\mathbf{p}^\droite_\veps(x,\cdot)) = {\bf w}\big(\mathbf{p}^\droite_\veps(0,\cdot))$ for any $\xvi \in \R_+$. But since $\mathbf{p}^\droite_\veps(0, s) = 1$ for any $s \in \R$, it follows that ${\bf w}\big(\mathbf{p}^\droite_\veps(0,\cdot)) = 0$ which implies ${\bf w}\big(\boldsymbol{p}^\droite_\veps\big)=0$.}
\end{dem}

\vspspe
We can now state the main result of this section which requires the following assumption.
\begin{assumption}[label={assumption_Dioph}]%{[Diophantine condition]} 
	The parameter $\delta$ has a finite irrationality measure $m(\delta)$ (see \eqref{eq:meas_irrat}), and the coefficients $(\mu_p, \rho_p)$ have the regularity
	\begin{equation}\label{regularity}
		\mu_p, \rho_p \in \mathscr{C}^m (\R^2) \quad \textnormal{for some $m \in \N$ such that $m > m(\delta)$ .}
	\end{equation}
\end{assumption}
\begin{thm}\label{thm_spectre} Let $\veps>0$. Under Assumption \ref{assumption_Dioph} there exists a unique pair $(\lambda^\droite_{0, \veps} , \varphi^\droite_{0,\veps}) \in \C \times {\cal N}_{\textit{per}}$ with $|\lambda^\droite_{0, \veps}| < 1$ and ${\bf w} \big(\varphi^\droite_{0,\veps}\big) = 0$ such that
	\begin{equation} \label{fundamental_eigenfunction}
	 \mathcal{P}^\droite_\veps \, \varphi^\droite_{0,\veps} = \lambda^\droite_{0, \veps} \; \varphi^\droite_{0,\veps} \quad \textnormal{and} \quad \varphi^\droite_{0,\veps}(0)= 1.
		\end{equation}
	If we define $\varphi^\droite_{k,\veps} \in {\cal N}_{\textit{per}}$ for any $k \in \Z$ by $\varphi^\droite_{k,\veps}(s) \coloneqq  \euler^{2 \icplx \pi k s} \, \varphi^\droite_{0,\veps}(s)$, then
		\begin{equation} \label{eigenfunctions}
	\mathcal{P}^\droite_\veps \; \varphi^\droite_{k,\veps} = \lambda^\droite_{k,\veps} \; \varphi^\droite_{k,\veps} \quad \mbox{with} \quad \lambda^\droite_{k,\veps} \coloneqq  \euler^{-2 \icplx \pi k \delta} \, \lambda^\droite_{0, \veps}
	\end{equation}
	The point spectrum of $\mathcal{P}^\droite_\veps$ is the countable set
	\begin{equation} \label{point_spectrum}
		\sigma_p\big(\mathcal{P}^\droite_\veps) = \{ \lambda^\droite_{k,\veps},\ k\in \Z \}, 
	\end{equation}
	where the eigenvalues $\lambda^\droite_{k,\veps}$ are all simple. Moreover, the point spectrum of $\mathcal{P}^\droite_\veps$ is a dense subset of the  spectrum $\sigma\big(\mathcal{P}^\droite_\veps) $ of $\mathcal{P}^\droite_\veps$.
	%circle ${\cal C}(0, r_\veps)$ and its spectrum $\sigma\big(\mathcal{P}^\droite_\veps) $ is ${\cal C}(0, r_\veps)$.

	\vspspe
	Finally, the eigenfunctions $\big\{\varphi^\droite_{k,\veps},\ k \in \Z\}$ form an orthogonal basis of $L^2_{\textit{per}}(\R)$ for the inner product 
	\begin{equation} \label{def_innerproduct}
		\big( \varphi, \psi \big)_{\boldsymbol{\rho}^\droite_{\veps} } \coloneqq  \int_0^1 \varphi(s) \, \overline \psi(s) \, \boldsymbol{\rho}^\droite_{\veps} (s) \; ds \quad \mbox{ where } \quad \boldsymbol{\rho}^\droite_{\veps} (s) \coloneqq  |\varphi^\droite_{0,\veps}(s)|^{-2},
	\end{equation}
	where $\boldsymbol{\rho}^\droite_{\veps} $ is bounded from below and above since  $\varphi^\droite_{0,\veps} \in {\cal N}_{\textit{per}}$.
\end{thm} 

\begin{dem} The proof that we are going to give is a slight adaptation of \cite{antonevich2012linear}.
	
\vspspe
\textbf{Step 1: A preliminary observation about the structure of the spectrum of} $\mathcal{P}^\droite_\veps$.~\\ 
Let $J_k, k \in \Z$ be the multiplication operator in $L^2_{\textit{per}}(\R)$ by the 1-periodic function $\euler^{2 \icplx \pi k s}$,
$$
\spforall \varphi \in L^2_{\textit{per}}(\R), \quad \big(J_k \, \varphi\big)(s) \coloneqq  \euler^{2 \icplx \pi k s} \, \varphi(s),
$$
which is unitary in $L^2_{\textit{per}}(\R)$ with $J_k^ {-1} = J_{-k}$. Since $\mathcal{P}^\droite_\veps$ is a weighted-shift operator (see Theorem \ref{thm_weighted-shift}), one checks immediately that 
\begin{equation} \label{commutation}
\spforall k \in \Z, \quad J_{-k} \, \mathcal{P}^\droite_\veps \, J_k = \euler^{2\icplx \pi k \delta} \; \mathcal{P}^\droite_\veps.
\end{equation}
%In particular that the point spectrum $\sigma_p((\mathcal{P}^\droite_\veps)^+)$ is invariant by the multiplication by $\euler^{2 \icplx \pi \, k \, \delta}$.
%$$
%\euler^{2 \icplx \pi \, k \, \delta} \, \sigma_p((\mathcal{P}^\droite_\veps)^+) = \sigma_p((\mathcal{P}^\droite_\veps)^+)
%$$
Let us now make an assumption which will be shown at the step 2 of the proof:
\begin{equation} \label{AssumptionFEV}
	\mbox {The operator $\mathcal{P}^\droite_\veps$ admits an eigenfunction $\varphi^\droite_{0,\veps} \in {\cal N}_{\textit{per}}$ such that ${\bf w} (\varphi^\droite_{0,\veps}) = 0$.}
\end{equation}
Let $\lambda^\droite_{0, \veps} \neq 0$ be the associated eigenvalue, with $r_{\veps} \coloneqq  |\lambda^\droite_{0, \veps}| > 0$ since 
$\mathcal{P}^\droite_\veps$ is bijective (see Theorem \ref{thm_weighted-shift}) and $r_\veps=\varrho(\mathcal{P}^\droite_\veps)$ by \eqref{eq:spectre_Peps}. From the commutation property \eqref{commutation}, we deduce easily \eqref{eigenfunctions}.
As $\delta \in \R \setminus \Q$, the eigenvalues $\lambda^\droite_{k,\veps}$ are distinct and $\big\{ \lambda^\droite_{k,\veps}, k \in \Z\}$ is dense in the circle ${\cal C}\big(0, r_{\veps} \big)$, which coincides with the whole spectrum of $\mathcal{P}^\droite_\veps$ by \eqref{eq:spectre_Peps}. \\[12pt]
%Since $\varphi^+_{0,\veps} \in {\cal N}_{\textit{per}}$, the weight function defined by
%\begin{equation} \label{def_weight}
%\boldsymbol{\rho}^\droite_0^\veps (s) \coloneqq  |\varphi^+_{0,\veps}|^{-2}
%\end{equation}
%is bounded from below and above, so that the weighted inner product
%\begin{equation} \label{def_weight}
% \big( \varphi, \psi \big)_{\boldsymbol{\rho}^\droite_0^\veps} \coloneqq  \int_0^1 \varphi(s) \, \overline \psi(s) \, \boldsymbol{\rho}^\droite_0^\veps (s) \; ds
%\end{equation}
%defines a norm in $L^ 2_{\textit{per}}(\R)$, equivalent to the norm . 
Moreover, from Fourier series theory in $L^ 2(0,1)$, we deduce that 
$
\big\{ \varphi^\droite_{k,\veps},\ k \in \Z \big\}$ is an orthonormal basis of $L^2_{\textit{per}}(\R)$ for the inner product \eqref{def_innerproduct},
in which $\mathcal{P}^\droite_\veps$ is {\it diagonal} in the sense that
\begin{equation} \label{diag}
\spforall \varphi \in L^2_{\textit{per}}(\R), \quad \varphi = \sum_{k \in \Z} \big( \varphi, \varphi^\droite_{k,\veps} \big)_{\boldsymbol{\rho}^\droite_{\veps}} \; \varphi^\droite_{k,\veps} \quad \Longrightarrow \quad \mathcal{P}^\droite_\veps \varphi = \sum_{k \in \Z} \lambda^\droite_{k,\veps} \, \big( \varphi, \varphi^\droite_{k,\veps} \big)_{\boldsymbol{\rho}^\droite_{\veps}} \, \varphi^\droite_{k,\veps}.
\end{equation}
From \eqref{diag}, we deduce that the point spectrum of $\mathcal{P}^\droite_\veps$ is given by \eqref{point_spectrum}.
%  and that its whole spectrum is the circle ${\cal C}\big(0, r_{\veps})$. 
 Moreover, one sees that 
each eigenvalue is simple: $\matrixkernel \big(\mathcal{P}^\droite_\veps - \lambda^\droite_{k,\veps} \big) = \vect \, (\varphi^\droite_{k,\veps} )$.

\vspspe
\textbf{Step 2: Existence of the fundamental eigenfunction} $\varphi^\droite_{0,\veps}$.~\\
Since, according to Lemma \ref{lem_WN}, $\boldsymbol{p}^\droite_\veps \in {\cal N}_{\textit{per}}\mbox{ and } {\bf w}\big( \boldsymbol{p}^\droite_\veps) = 0$, we know from Corollary \ref{coro_WN0}, that% (up to a change of variable in order to simplify the expressions later on)
\begin{equation} \label{expo_p}
	\boldsymbol{p}^\droite_\veps (s) = \euler^{ \, g_\veps(s)} \quad \mbox{for some $g_\veps \in \mathscr{C}^0_{\textit{per}}(\R, \C)$}.
	\end{equation}
Similarly, since we look for $\varphi^\droite_{0,\veps} \in {\cal N}_{\textit{per}}\mbox{ with } {\bf w}\big( \varphi^\droite_{0,\veps}) = 0$ we can look for 
	\begin{equation} \label{expo_phi}
		\varphi^\droite_{k,\veps} (s) = \euler^{ \, v_\veps(s)} \quad \mbox{where} \quad v_\veps \in \mathscr{C}^0_{\textit{per}}(\R, \C)
	\end{equation} 
so that the existence of $\varphi^\droite_{0,\veps}$ is equivalent to the one of $v_\veps$. Substituting \eqref{expo_p} and \eqref{expo_phi} into the eigenvalue equation, that is to say, 
\begin{equation}\label{eq:link_p_phi0}
\boldsymbol{p}^\droite_\veps(s-\delta)\; 	\varphi^\droite_{k,\veps} (s-\delta) = \lambda^\droite_{0, \veps} \; \varphi^\droite_{0,\veps} (s)
\end{equation}
naturally suggests to look for $v_\veps(s)$ as a solution of the following difference equation 
\begin{equation} \label{diff_eq}
	g_\veps(s - \delta) + v_\veps(s - \delta) = \log \lambda^\droite_{0, \veps} + v_\veps(s).
\end{equation} 
Integrating \eqref{diff_eq} between $0$ and $1$, we obtain, by periodicity of $v_\veps$ (and $g_\veps$),
\begin{equation*}
\log \lambda^\droite_{0, \veps} = \int_0^1 g_\veps(s) \; ds, \quad \mbox{thus} \quad \lambda^\droite_{0, \veps} = \exp {\int_0^1 g_\veps(s) \; ds}.
\end{equation*}
Multiplying \eqref{diff_eq} by $\euler^{-2\icplx\pi k s}$, with $k \neq 0$, and integrating the result over $[0,1]$, we obtain for the $k$-th Fourier coefficient $v_{k,\veps}$ of $v_\veps$ ($g_{k,\veps}$ being the one of $g_{\veps}$) the following equation 
\begin{equation} \label{diff_eq_fourier}
\spforall k \in \Z, \quad v_{k,\veps} \; \big( \, 1 - \euler^{2\icplx \pi k\delta} \, \big) = 	-g_{k,\veps} \quad \Longrightarrow \quad v_{k,\veps} = \frac{g_{k,\veps}}{\euler^{2\icplx \pi k\delta} - 1},
\end{equation} 
where the division by $\euler^{2\icplx \pi k\delta}-1$ is possible because $\delta \in \R \setminus \Q$. To conclude, it remains to investigate the convergence, in $\mathscr{C}^0_{\textit{per}}(0,1)$, of the series 
\begin{equation} \label{fourier_series}
\sum_{k \neq 0} \frac{g_{k,\veps}}{\euler^{2\icplx \pi k\delta} - 1} \; \euler^{2\icplx \pi k s}.
\end{equation} 
Note that this convergence is not obvious since, by density of $\big\{ \euler^{2\icplx \pi k\delta }, k \in \Z \big\}$ in the unit circle, $\euler^{2\icplx \pi k\delta} - 1 $ can be arbitrarily small for infinitely many $k$'s. %
This is where the regularity of the symbol $\boldsymbol{p}^\droite_\veps $, through that of the coefficients $(\mu_p, \rho_p)$ will come into play to bound $g_{k,\veps}$, as well as the measure of irrationality $m(\delta)$ of $\delta$ to get a lower bound for $|\euler^{ 2\icplx \pi k\delta} - 1|$.

\vspspe
First, due to the regularity assumption, we know by \eqref{reg_p_eps} that $\boldsymbol{p}^\droite_\veps \in \mathscr{C}^m(\R, \C).\mbox{ Thus } g_\veps \in \mathscr{C}^m(\R, \C)$,
%Let us assume that, for some $m \in \N^*$ to be fixed later, 
%$$
%\boldsymbol{p}^\droite_\veps^+ \in \mathscr{C}^m(\R, \C), \quad \mbox{ thus } \quad \boldsymbol{g}_\veps^+ \in \mathscr{C}^m(\R, \C)
%$$
and we have the following upper bound for the coefficients $g_{k,\veps}$
\begin{equation} \label{upper_bound}
\spforall k \neq 0, \quad |g_{k,\veps}| \leq \frac{C_m( g_\veps )}{|k|^m} , \quad C_m(g_\veps) \coloneqq  \sup_{s \in [0,1]}|g_{\veps}^{(m)}(s)|.
\end{equation} 
To get a lower bound for $|\euler^{2\icplx \pi k\delta} - 1|$, we use the concavity of the sine function in $[0,\pi/2]$ to obtain 
\begin{equation} \label{classic} 
\spforall \theta \in \; [- \pi, \pi]\, , \quad \big| \euler^{i\theta} -1 \big| = 2\,|\sin (\theta/2)| \geq 2\, |\theta|/\pi
\end{equation} 
Next, for any $k \in \Z$ such that $k\vts \delta \notin \Z$, there exists $n \equiv n(k) \in \Z$ such that $|k \delta - n| < 1/2$. Thus,
$$
\big|\euler^{2\icplx \pi k\delta}-1 \big| = \big|\euler^{2\icplx \pi (k\delta-n)}-1 \big| \geq 4 \, |k \, \delta - n|,
$$
where the last inequality follows from \eqref{classic}, since $\theta \coloneqq  2\pi (k\delta-n) \in \;]\!-\! \pi, \pi[ \,.$ However, thanks to the Diophantine condition \eqref{eq:notLiouville}, we have for any $\nu > m(\delta)$
\begin{equation} \label{lower_bound} 
|k \, \delta - n| = |k| \, |\delta - n/k| \geq \frac{C_\nu}{|k|^{\nu-1}} \quad \mbox{which yields } \quad \big|\euler^{2\icplx \pi k\delta}-1 \big| \geq \frac{\widetilde{C}_\nu}{|k|^{\nu-1}}\; .
\end{equation} 
Finally, by combining \eqref{upper_bound} and \eqref{lower_bound}, we obtain the following uniform (in $s$) bound for the general term of the series 
\eqref{fourier_series}
\begin{equation} \label{general_term}
\Big|\frac{g_{k,\veps}}{\euler^{2\icplx \pi k\delta} - 1}\Big| \leq \frac{C_m( g_\veps )}{\widetilde{C}_\nu} \, \frac{1}{|k|^{1 + (m - \nu)}},
\end{equation} 
which ensures the (normal) convergence of the series as soon as $m > m(\delta)$ (since one can choose $\nu$ such that $m(\delta) < \nu < m$) so that
\begin{equation} \label{fourier_series2}
v_\veps(s) \coloneqq  	\sum_{k \neq 0} \frac{g_{k,\veps}}{\euler^{2\icplx \pi k\delta} - 1} \; \euler^{2\icplx \pi k s} \in \mathscr{C}^0_{\textit{per}}(\R,\C).
\end{equation} 
We then take the exponential of $v_\veps$, see \eqref{expo_phi}, to construct $\varphi^\droite_{0,\veps}$.

\vspspe
Note that, in the above construction of $v_\veps$, the mean value $v_{0,\veps}$ of $v_\veps$ (an additive constant for $v_\veps$ which becomes a multiplicative constant for $\varphi^\droite_{0,\veps}$) is completely free. Since $\varphi^\droite_{0,\veps}$ never vanishes, we can thus build $\varphi^\droite_{0,\veps}$
%as
% $$
% \varphi^\droite_{0,\veps} (s) = \euler^{\, \big[\, v_\veps(s) - v_\veps(0) \, \big]},
% $$
so that $\varphi^\droite_{0,\veps} (0) = 1$ by construction. This concludes the proof.
\end{dem}

% \vspace{0\baselineskip}\noindent 
\begin{rmk}\label{rmk:eig_fund}
	Note that since ${\bf w} ( \varphi^\droite_{k,\veps}) = k$, only $\varphi^\droite_{0,\veps}$ has a winding number $0$, reason why it can be easily distinguished from the other eigenfunctions: we call it the \textbf{fundamental eigenfunction} of $\mathcal{P}^\droite_\veps$. By extension, $\lambda^\droite_{0, \veps}$ and $(\varphi^\droite_{0,\veps}, \lambda^\droite_{0, \veps})$ are respectively the \textbf{fundamental eigenvalue} and the \textbf{fundamental eigenpair} of $\mathcal{P}^\droite_\veps$.
%
%	We have chosen in this theorem to normalize $\varphi^\droite_{0,\veps}$ by $\varphi^\droite_{0,\veps}(0)=1$. But since the eigenfunctions are defined up to an additive constant, we could use any other normalization such as $\|\varphi^\droite_{0,\veps}\|_{L^2}=1$ or $\|\varphi^\droite_{0,\veps}\|_{H^1}=1$ or even $\|\varphi^\droite_{0,\veps}\|_{H^2}=1$ if $(\rho_p,\mu_p)$ and then $\varphi^\droite_{0,\veps}$ are regular enough. This is what is done in Section \ref{sec:LAP}.
\end{rmk}

\noindent
From the proof of Proposition \ref{thm_spectre}, one can show the following corollary (we omit the details).
\begin{cor}\label{cor_spectre}
	Assume that $\delta$ is not a Liouville number and that the coefficients $(\mu_p, \rho_p)$ have the regularity given in \eqref{regularity} with $m>m(\delta)+k$ for a given $k\in\N$. % 
	Then the fundamental eigenfunction is such that $\varphi^\droite_{0,\veps}\in \mathscr{C}^k_{\textit{per}}(\R,\C)$.
\end{cor}

\vspace{1\baselineskip}	\noindent 
To each eigenfunction $\varphi^\droite_{k,\veps}$ of $\mathcal{P}^\droite_\veps$, we can associate a {\it 2D mode}, which is the solution of the half-space problem	\eqref{RobinHalfspaceproblems} taking this eigenfunction as the boundary data, namely
\begin{equation} \label{2Dmode_eps}
U^\droite_{k, \varepsilon} \coloneqq  U_{\varepsilon}^\droite (\varphi^\droite_{k,\veps}).
\end{equation}
% \varphi(s)\, u^\droite_{s, \veps}(\xvi + n/\cuti_2) = \varphi^\veps_n (s + n \delta)\, e^0_{s + n\delta, \veps} (\xvi) + \psi^\veps_n (s + (n + 1) \delta)\, e^1_{s + n\delta, \veps} (\xvi).
Similarly to the cell by cell expression \eqref{eq:cell_by_cell_1D} of $u^\droite_{s, \veps}$, one has the following result, which we shall use in Section \ref{sec:propagative_freq}.
%
% With the help of the 1D reconstruction formula \eqref{eq:cell_by_cell_1D}, one can makes explicit the fibered structure of $U_{k, \varepsilon}^\droite$ which we shall use in Section \ref{}.
\begin{lem}\label{rmk:recons_fund_sol}
	We introduce for $k \in \Z, s \in \R$
	\begin{equation*}
		\aeforall x\in \R_+, \quad u^\droite_{s,k,\veps}(x) \coloneqq  U^\droite_{k, \varepsilon} (s + \cuti_1\, \xvi, \cuti_2\, \xvi).
	\end{equation*}
	Note that $u^\droite_{s,k,\veps}$ is well-defined and belongs to $H^1(\R_+)$, because $\varphi^\droite_{k, \veps}$ is continuous. One has the following formula, with $\psi^\droite_{k,\veps} \coloneqq  \mathcal{S}^\droite_\veps\, \varphi^\droite_{k,\veps}$: for all $x\in (0,1/\cuti_2)$, for all $n\in\N$
	\begin{equation}\label{eq:recons_1dk}
		u^\droite_{s, k, \veps}(\xvi + n/\cuti_2) = (\lambda^\droite_{k,\veps})^n \, \big[\varphi^\droite_{k, \veps} (s + n \delta)\, e^{\droite, 0}_{s + n\delta, \veps} (\xvi) + \psi^\droite_{k, \veps} (s + (n + 1) \delta)\, e^{\droite, 1}_{s + n\delta, \veps} (\xvi) \big].
	\end{equation} 
	
\end{lem}
\begin{dem}
	It suffices to use \eqref{eq:cell_by_cell_1D} for $\varphi = \varphi^\droite_{k,\veps}$ and remark that 
	\begin{equation*}
		(\mathcal{P}^\droite_\veps)^n\, \varphi^\droite_{k,\veps} = (\lambda^\droite_{k,\veps})^n \, \varphi^\droite_{k,\veps} \quad \mbox{and} \quad \mathcal{S}^\droite_\veps\, (\mathcal{P}^\droite_\veps)^n\, \varphi^\droite_{k,\veps} = (\lambda^\droite_{k,\veps})^n \, \psi^\droite_{k,\veps}.
	\end{equation*}

	\vspace{-1.5\baselineskip} \noindent
\end{dem}

% \begin{rmk}[Quasi-Floquet modes]
% 	The 2D fields $U^\droite_{k,\veps}$ can be seen as the equivalent of (evanescent) Floquet modes associated with the classical Helmholtz equation in a periodic media. By analogy, the 1D functions $u^\droite_{s,k,\veps}$ can be called the \textbf{quasi-Floquet modes} associated with the quasiperiodic medium $(\mu_{s, \cutvec}, \rho_{s, \cutvec})$.\footnote{\sfcomment{Je ne suis pas sure que cette remarque soit pertinente ici puisqu en general on ne definit les modes de floquet que pour le cas sans dissipation et souvent ce ne sont que les modes propagatifs.  Je propose de reporter cette remarque dans la partie sans dissipation. Si on enleve cette remarque, on pourrait tout enlever a partir de "To each eigenfunction" ... a la fin de la remarque. Et on ne mettrait que \eqref{eq:recons_1d0} eventuellement qui est utile pour la reconstruction du probleme de demi-ligne.}}
% 	\end{rmk}
\begin{rmk}\label{rmk:Spsikeps}
		Since the operator $\mathcal{S}^\droite_\veps$ is a weighted shift operator (see Proposition \ref{thm_weighted-shift}), one deduces from $\varphi^\droite_{k,\veps}(s) \coloneqq  \euler^{2 \icplx \pi k s} \, \varphi^\droite_{0,\veps}(s)$ (see Theorem \ref{thm_spectre}) that $\psi^\droite_{k,\veps} \coloneqq  \mathcal{S}^\droite_\veps\, \varphi^\droite_{k,\veps}$ can be computed directly from $\varphi^\droite_{0,\veps}$ as follows 
		\[ \psi^\droite_{k,\veps}(s) = \euler^{2\icplx\pi k (s-\delta)}\, [\mathcal{S}^\droite_\veps\, \varphi^\droite_{0,\veps}](s). \]
	\end{rmk}
\noindent A direct application of Lemma \ref{rmk:recons_fund_sol} is that the solution $u^\droite_{\veps}$ of the half-line problem \eqref{halfline_problems} for $\idomlap=\droite$ (which satisfies \eqref{cutformulas}) is given by: for all $x\in (0,1/\cuti_2)$, for all $n\in\N$
	\begin{multline}\label{eq:recons_1d0}
		u^\droite_{\veps}(a^\droite+\xvi + n/\cuti_2) = (\lambda^\droite_{0,\veps})^n\,[\varphi^\droite_{0, \veps} (a^\droite\cuti_1)]^{-1} \, \\ \left[\varphi^\droite_{0, \veps} (a^\droite\cuti_1 + n \delta)\, e^{\droite, 0}_{a^\droite\cuti_1 + n \delta, \veps} (\xvi) + \psi^\droite_{0, \veps} (a^\droite\cuti_1 +  (n + 1) \delta)\, e^{\droite, 1}_{a^\droite\cuti_1 + n\delta, \veps} (\xvi) \right],
	\end{multline}  and the RtR coefficient $\RtRcoef_\veps^\droite$ can be computed thanks to
\begin{equation} \label{def_lambda_1D_dir}
\RtRcoef_\veps^\droite  =	 \boldsymbol{t}^{\droite, 00}_{\veps}(a^\droite\cuti_1)+[\varphi^\droite_{0, \veps} (a^\droite\cuti_1)]^{-1}\,\psi^\droite_{0, \veps} (a^\droite\cuti_1+\delta)\, \boldsymbol{t}^{\droite, 10}_{\veps}(a^\droite\cuti_1).
\end{equation} 
\subsection{Spectral resolution of the Riccati system}\label{sec:spectrum_Riccati}
Let us suppose in this section that Assumption \ref{assumption_Dioph} is satisfied. We know that the operator $\mathcal{P}^\droite_\veps$ is entirely determined by its fundamental eigenpair $(\lambda^\droite_{0, \veps}, \varphi^\droite_{0,\veps})$. Our goal now is to find this eigenpair from the Riccati system \eqref{Riccatiplus}.

\vspspe
Let us introduce the two operators $\mathscr{M}^\droite_\veps,\mathscr{N}^\droite_\veps\in \mathscr{L}\big(L^2_{\textit{per}}(\R) \times L^2_{\textit{per}}(\R) \big)$, defined as
\begin{equation}\label{eq:Riccati_op}
\mathscr{M}^\droite_\veps \coloneqq  \begin{pmatrix} \mathcal{T}^{\droite, 01}_\veps & \mathcal{T}^{\droite, 11}_\veps \\[2pt] 0 & I \end{pmatrix} \quad \textnormal{and} \quad \mathscr{N}^\droite_\veps \coloneqq  \begin{pmatrix} I & 0 \\[2pt] \mathcal{T}^{\droite, 00}_\veps & \mathcal{T}^{\droite, 10}_\veps \end{pmatrix},
\end{equation}
where the local RtR operators $\mathcal{T}^{\droite, \ell k}_\veps$ are defined by \eqref{defRtR}.
It is easy to see that for any solution $(P,S)$ of the Riccati system \eqref{Riccatiplus}, if $P\varphi = \lambda \, \varphi$ with $\varphi\in L^2_{\textit{per}}(\R)\setminus\{0\}$ then 
\[\mathscr{M}^\droite_\veps \,  \begin{pmatrix}  \varphi \\[2pt] S \vts \varphi \end{pmatrix} = \lambda \, \mathscr{N}^\droite_\veps \, \begin{pmatrix}  \varphi \\[4pt] S \vts \varphi \end{pmatrix}. \]
It is then natural to introduce the augmented (and generalized) eigenvalue problem 
\begin{equation} \label{augmented_EVP}
	\mbox{Find } \lambda \in \C \mbox{ and } (\varphi, \psi) \in L^2_{\textit{per}}(\R)^2\setminus\{0\}\} \mbox{ such that }\quad  \mathscr{M}^\droite_\veps \, \begin{pmatrix} \varphi \\ \psi \end{pmatrix} = \lambda\, \mathscr{N}^\droite_\veps \, \begin{pmatrix} \varphi \\ \psi \end{pmatrix}
\end{equation} 
and the corresponding point spectrum, called the \textbf{Riccati point spectrum}, defined as
\begin{equation} \label{eq:Ric_spec}
	\sigma^\droite_{R, p, \veps} \coloneqq  \big\{ \lambda \in \C \ /\ 0 \in \sigma_p ( \mathscr{M}^\droite_\veps - \lambda \, \mathscr{N}^\droite_\veps ) \}. 
\end{equation}
We have noticed above that for any solution $(P,S)$ of the Riccati system \eqref{Riccatiplus}, one has the inclusion
	\begin{equation}\label{inclusion}
\sigma_p\big( P\big) \subset \sigma^\droite_{R,p,\veps}.
\end{equation} 
In particular, $\sigma_p \big( \mathcal{P}^\droite_\veps\big) \subset \sigma^\droite_{R,p,\veps}$. {However, the inclusion is strict, since as we shall see later on (see Proposition \ref{thm_spectreaugmented}), one can find $\lambda\in\sigma^\droite_{R,p,\veps}$ such that $|\lambda|>1$.}

\begin{rmk} \label{rmk_size} The Riccati point spectrum is the one that we shall be able to compute numerically because it involves the RtR operators defined thanks to problems set in a bounded domain, namely the cell $\cellper^{\droite, 0} = (0, 1)^2$. After discretization, the RtR operators will be approximated by matrices of dimension $N$, with $N$ being the number of degrees of freedom of the interface $\Sigmaper^{\droite,0}$. One then expects the approximated Riccati spectrum to contain 2N elements while the approximated spectrum of the propagation operator contains only N elements.% In other words one expects $\sigma^\droite_{R,\veps} \setminus \sigma_p\big( \mathcal{P}^\droite_\veps\big)$ to be "of the size of" $\sigma_p\big( \mathcal{P}^\droite_\veps\big)$.
\end{rmk} 
\noindent In order to extract $\sigma_p\big( \mathcal{P}^\droite_\veps\big)$ from $\sigma^\droite_{R, p, \veps}$,we shall use a complete description of the Riccati point spectrum.
This will be based on the factorization of $\mathscr{M}^\droite_\veps - \lambda \, \mathscr{N}^\droite_\veps$, inspired by \cite[Lemma 4.6]{amenoagbadji2023wave} obtained for a similar but simpler Riccati equation.
\begin{prop}\label{thm:facto}
	Given $\veps>0$, let $\mathcal{P}^\droite_\veps, \mathcal{S}^\droite_\veps$ be the operators defined by \eqref{defPS}. For any $\lambda\in\C$, we have
	\begin{equation} \label{factor} 
		\mathscr{M}^\droite_\veps - \lambda \, \mathscr{N}^\droite_\veps = \begin{pmatrix} I & 0 \\[2pt] \lambda  \transp{\mathcal{S}^\droite_\veps} & \lambda \transp{\mathcal{P}^\droite_\veps} - I \end{pmatrix}\; \mathcal{L}^\droite_\veps\;\begin{pmatrix} \mathcal{P}^\droite_\veps - \lambda  & 0 \\[2pt] \mathcal{S}^\droite_\veps & - I \end{pmatrix},\quad \mathcal{L}^\droite_\veps \coloneqq  \begin{pmatrix} I  & - \mathcal{T}^{\droite, 11}_\veps \\[2pt] - \RtRop^\droite_\veps & I
		\end{pmatrix},
	\end{equation}
	where $\mathcal{T}^{\droite, 11}_\veps$ is defined in \eqref{defRtR} and $\RtRop^\droite_\veps$ in \eqref{def_Lambda}. Moreover, $\mathcal{L}^\droite_\veps$ is invertible.
\end{prop}

\begin{dem}
	\textbf{Step 1: Factorization of $\mathscr{M}^\droite_\veps - \lambda \, \mathscr{N}^\droite_\veps$.} Let $\lambda\in\C$; we have 
	\[
	\mathscr{M}^\droite_\veps - \lambda \, \mathscr{N}^\droite_\veps = \begin{pmatrix} \mathcal{T}^{\droite, 01}_\veps- \lambda & \mathcal{T}^{\droite, 11}_\veps \\[2pt] - \lambda \, \mathcal{T}^{\droite, 00}_\veps & I-\lambda \mathcal{T}^{\droite, 10}_\veps \end{pmatrix}.
	\]
	We do not change the above expression by adding the two operators $\mathcal{P}^\droite_\veps-\mathcal{T}^{\droite, 01}_\veps - \mathcal{T}^{\droite, 11}_\veps \mathcal{S}^\droite_\veps$ and $\mathcal{T}^{\droite, 00}_\veps\, \mathcal{P}^\droite_\veps + \mathcal{T}^{\droite, 10}_\veps\, \mathcal{S}^\droite_\veps \, \mathcal{P}^\droite_\veps - \mathcal{S}^\droite_\veps$ respectively to the upper and lower left blocks: indeed, as $(\mathcal{P}^\droite_\veps, \mathcal{S}^\droite_\veps)$ is solution of the Riccati system, these two operators are equal to 0. Doing so, we obtain an expression whose interest is to make appear $\mathcal{P}^\droite_\veps - \lambda $
	\[
	\mathscr{M}^\droite_\veps - \lambda \, \mathscr{N}^\droite_\veps = \begin{pmatrix} \mathcal{P}^\droite_\veps - \lambda - \mathcal{T}^{\droite, 11}_\veps\, \mathcal{S}^\droite_\veps & \mathcal{T}^{\droite, 11}_\veps \\[2pt] \mathcal{T}^{\droite, 00}_\veps (\mathcal{P}^\droite_\veps - \lambda)+\mathcal{T}^{\droite, 10}_\veps\, \mathcal{S}^\droite_\veps\, \mathcal{P}^\droite_\veps - \mathcal{S}^\droite_\veps & I-\lambda \mathcal{T}^{\droite, 10}_\veps \end{pmatrix}. 
	\]
	Thanks to \eqref{form_Lambda}, we note that $\mathcal{T}^{\droite, 00}_\veps (\mathcal{P}^\droite_\veps - \lambda)+\mathcal{T}^{\droite, 10}_\veps \, \mathcal{S}^\droite_\veps\, \mathcal{P}^\droite_\veps - \mathcal{S}^\droite_\veps = \RtRop^\droite_\veps(\mathcal{P}^\droite_\veps - \lambda)+ (\lambda \mathcal{T}^{\droite, 10}_\veps-I)\, \mathcal{S}^\droite_\veps$, so that we deduce a first factorization of $	\mathscr{M}^\droite_\veps - \lambda \, \mathscr{N}^\droite_\veps$, namely
	% \[
	% \mathscr{M}^\droite_\veps - \lambda \, \mathscr{N}^\droite_\veps = \begin{pmatrix} I \times (\mathcal{P}^\droite_\veps - \lambda) - \mathcal{T}^{\droite, 11}_\veps \mathcal{S}^\droite_\veps &I \times 0 - \mathcal{T}^{\droite, 11}_\veps (-I) \\[2pt] \RtRop^\droite_\veps(\mathcal{P}^\droite_\veps - \lambda) + (\lambda \mathcal{T}^{\droite, 10}_\veps-I)\, \mathcal{S}^\droite_\veps & \RtRop^\droite_\veps \times 0 + (\lambda \mathcal{T}^{\droite, 10}_\veps-I) (-I) \end{pmatrix}. 
	% \]
	% We recognize a first factorization of $	\mathscr{M}^\droite_\veps - \lambda \, \mathscr{N}^\droite_\veps$, namely
	\begin{equation}\label{fact1}
		\mathscr{M}^\droite_\veps - \lambda \, \mathscr{N}^\droite_\veps=\mathscr{D}_\veps(\lambda)\;\begin{pmatrix} \mathcal{P}^\droite_\veps - \lambda  & 0 \\[2pt] \mathcal{S}^\droite_\veps & - I \end{pmatrix}\quad\textnormal{with}\quad\mathscr{D}_\veps(\lambda) \coloneqq  \begin{pmatrix} I& - \mathcal{T}^{\droite, 11}_\veps  \\[2pt] \RtRop^\droite_\veps &  \lambda \, \mathcal{T}^{\droite, 10}_\veps -I  \end{pmatrix}.
	\end{equation}
	Let us now factorize the operator $\mathscr{D}_\veps(\lambda)$. We proceed essentially along the same lines as for the first step of the factorization. To exhibit the transpose operators that arise in the first factor of \eqref{factor}, the trick is to transpose the Riccati system \eqref{Riccatiplus}, and use the symmetry properties of the operators $\mathcal{T}^{\droite, \ell k}_\veps$ (Lemma \ref{lem_transposition_properties}), as well as the expression \eqref{form_Lambda} of $\RtRop^\droite_\veps$ and the fact that $\transp{\RtRop^\droite_\veps} = \RtRop^\droite_\veps$ (see \eqref{transpose_Lambda}) to obtain
	\[
		\left\{ \begin{array}{llll} 
		\transp{\mathcal{P}^\droite_\veps}-\mathcal{T}^{\droite, 10}_\veps - \transp{\mathcal{S}^\droite_\veps} \mathcal{T}^{\droite, 11}_\veps =0,\\ [5pt]
		\transp{\mathcal{S}^\droite_\veps} - \transp{\mathcal{P}^\droite_\veps} \RtRop^\droite_\veps=0.
		\end{array} \right.
	\]
	Thus, we do not change anything by adding $\lambda \, \big(\transp{\mathcal{S}^\droite_\veps} - \transp{\mathcal{P}^\droite_\veps} \,\RtRop^\droite_\veps \big)$ and $\lambda \, \big(\transp{\mathcal{P}^\droite_\veps}-\mathcal{T}^{\droite, 10}_\veps + \transp{\mathcal{S}^\droite_\veps}\,\mathcal{T}^{\droite, 11}_\veps\big)$ respectively to the left and right lower blocks of $\mathscr{D}_\veps(\lambda)$ in \eqref{fact1}. This leads to an expression that has the interest to make appear $\lambda \transp{\mathcal{P}^\droite_\veps} - I$, namely
	\[
		\begin{array}{rl}
			\mathscr{D}_\veps(\lambda) &= \begin{pmatrix} I& - \mathcal{T}^{\droite, 11}_\veps  \\[2pt] \lambda \transp{\mathcal{S}^\droite_\veps}- (\lambda \transp{\mathcal{P}^\droite_\veps} - I)\, \RtRop^\droite_\veps& -\lambda \transp{\mathcal{S}^\droite_\veps}\, \mathcal{T}^{\droite, 11}_\veps +(\lambda \transp{\mathcal{P}^\droite_\veps} - I) \end{pmatrix}.
		\end{array}
	\]
	We then recognize, similarly to the first part of the proof, the factorization 
	\[
		\mathscr{D}_\veps(\lambda) =\begin{pmatrix} I & 0 \\[2pt] \lambda \transp{\mathcal{S}^\droite_\veps} & \lambda \transp{\mathcal{P}^\droite_\veps} - I\end{pmatrix}\; \begin{pmatrix} I  & - \mathcal{T}^{\droite, 11}_\veps \\[2pt] - \RtRop^\droite_\veps & I
		\end{pmatrix},
	\]
	which, combined with \eqref{fact1}, leads to the desired factorization of $\mathscr{M}^\droite_\veps - \lambda \, \mathscr{N}^\droite_\veps$.
		
		\vspspe
		\textbf{Step 2: Invertibility of $\mathcal{L}^\droite_\veps$.} We begin by showing that the invertibility of $\mathcal{L}^\droite_\veps$ is equivalent to the well-posedness of a certain problem. Given $(f,g)$ in $ L^2_{\textit{per}}(\R)^2$, we want to construct $(\varphi,\psi) \in L^2_{\textit{per}}(\R)^2$ such that $ \mathcal{L}^\droite_\veps(\varphi,\psi)=(f,g)$, or equivalently
			\[
				\left\{
					\begin{array}{r@{\ =\ }l}
						\varphi-\mathcal{T}^{\droite, 11}_\veps\,\psi&f
						\\[8pt]
						- \RtRop^\droite_\veps\, \varphi + \psi &g
					\end{array}
				\right.
				\quad \Longleftrightarrow \quad 
				\left\{
					\begin{array}{r@{\ =\ }l}
						\varphi- {\cal R}^\droite_- E^{\droite, 1}_\veps (\psi)|_{\Sigmaper^{\droite, 1}}&f
						\\[8pt]
						- \RtRop^\droite_\veps\, \varphi+ {\cal R}^\droite_+\, E^{\droite, 1}_\veps (\psi)|_{\Sigmaper^{\droite, 1}}&g,
					\end{array}
				\right.
			\] 
			where the equivalence follows from the definitions of $E^{\droite, 1}_\veps (\psi)$ and $\mathcal{T}^{\droite, 11}_\veps$. Therefore, since 
			${\cal R}^\droite_+ - {\cal R}^\droite_- = 2 \nu^\droite \mu_p\, \Dt{}$ and ${\cal R}^\droite_+ + {\cal R}^\droite_- = -2\vts \icplx\vts z$, we also have
			\begin{equation*}
				\displaystyle
				\mathcal{L}^\droite_\veps(\varphi,\psi) = (f,g) \quad \Longleftrightarrow \quad \left\{
					\begin{array}{ll}
						(I + \RtRop^\droite_\veps)\varphi + 2\icplx z\, E^{\droite, 1}_\veps (\psi)=f-g & \mbox{on } {\Sigmaper^{\droite, 1}}
						\\[8pt]
						(I - \RtRop^\droite_\veps)\varphi - 2\mu_p\, \Dt{} {E^{\droite, 1}_\veps (\psi)}=f+g & \mbox{on } {\Sigmaper^{\droite, 1}}
					\end{array}
				\right.
			\end{equation*}
			Consequently, if $\mathcal{L}^\droite_\veps(\varphi,\psi)=(f,g)$, setting $U \coloneqq  E_\veps^{1} (\psi)$, then the pair $(\varphi,U) \in {\bf V} \coloneqq L^2_{\textit{per}} (\R)\times H^1_\cutvec (\cellper^{\droite, 0})$ satisfies
			\begin{equation} \label{pbfort}
				\left| \begin{array}{r@{\ =\ }l@{\quad}l} 
							 -\Dt{} \big( \mu_p \, \Dt{} U\big) - \rho_p \, (\omega^2+\icplx\veps) \, U & 0, & \textnormal{in }\cellper^{\droite, 0},	
							 \\[8pt]
							 {\cal R}_+ \, {U} & 0, & \mbox{on } {\Sigmaper^{\droite, 0},}
							 \\[8pt]
							 (I+\RtRop^\droite_\veps)\varphi+2 \, \icplx z\,{U} & f-g,\; & \mbox{on } {\Sigmaper^{\droite, 1},}
							 	\\[8pt]
							 (I-\RtRop^\droite_\veps)\varphi - 2 \vts \mu_p\,{\Dt{}U} & f+g, & \mbox{on } {\Sigmaper^{\droite, 1},}
								\\[8pt]
								\multicolumn{3}{c}{\textnormal{$U$ is periodic w.r.t. $\yvi_1$.}}
				\end{array} \right.
			\end{equation} 
			Conversely, if $(\varphi, U)\in {\bf V}$ satisfies \eqref{pbfort}, then by setting $\psi \coloneqq  {\cal R}^\droite_- U|_{\Sigmaper^{\droite, 0}}$, we obtain $\mathcal{L}^\droite_\veps (\varphi, \psi) = (f, g)$. This proves the equivalence between the invertibility of $\mathcal{L}^\droite_\veps$ and the well-posedness of \eqref{pbfort}. Hence, it remains to prove that \eqref{pbfort} is well-posed.
			
			\vspspe
			Problem \eqref{pbfort} is equivalent to a variational problem involving the closed affine and linear subspaces ${\bf V} (h) \subset {\bf V} $ defined for any $h \in L^2_{\textit{per}}(\R)$ by
			$$
			{\bf V} (h)\coloneqq  \big\{ (\psi,V) \in {\bf V} \; / \; 
				(I+\RtRop^\droite_\veps)\, \psi+2\icplx\vts z\,V = h\mbox{ on }{\Sigmaper^{\droite, 1}}\big\} , \quad {\bf V} _0 \coloneqq  {\bf V} (0).
			$$
			Introducing the continuous sesquilinear forms 
			$$
			 \left\{ \begin{array}{r@{\ \coloneqq \ }l@{\quad}l} 
				a_{\cal C}(U,V) & \displaystyle \int_{\cellper^{\droite, 0}} \big(\mu_p \, \Dt{} U\overline{\Dt{} V} - \rho_p \, (\omega^2+\icplx\veps) \, U \, \overline{V}\big) - \icplx \vts z \vts \int_{\Sigmaper^{\droite, 0}} U \, \overline{V}, & U,V\in H^1_\cutvec(\cellper^{\droite, 0})
				\rett
				a_\Sigma (\varphi, \psi) & \displaystyle -\frac{\icplx}{4z}\int_{\Sigmaper^{\droite, 1}} (I-\RtRop^\droite_\veps)\, \varphi \; \overline{(I+\RtRop^\droite_\veps)\, \psi}, & \varphi, \psi\in L^2_{\textit{per}}(\R),
			\end{array} \right.
			$$
			Problem \eqref{pbfort} is equivalent to: \textit{finding $(\varphi, U) \in {\bf V} (f-g)$ (for the third line of \eqref{pbfort}) such that}
			\begin{equation} \label{pbvar}
				a_{\cal C}(U,V) + a_\Sigma (\varphi, \psi) = \int_{\Sigmaper^{\droite, 0}}(f+g) \, \overline{V}, \quad \spforall (\psi, V) \in {\bf V}(0).
			\end{equation}
			Since $\varepsilon > 0$, the sesquilinear form $a_{\cal C}(\cdot, \cdot)$ is coercive in $H^1_\cutvec(\cellper^{\droite, 0})$ while, thanks to the inequality \eqref{eq:flux_eps} involving $\RtRop^\droite_\veps$, the form $a_{\Sigma}(\cdot, \cdot)$ is coercive in {$L^2_{\textit{per}}(\R)\equiv L^2(\Sigmaper^{\droite, 1})$}. Thus, Lax-Milgram theorem ensures the well-posedness of \eqref{pbfort}, which gives the invertibility of $\mathcal{L}^\droite_\veps$.
\end{dem}

\vspspe
The previous factorization allows to completely link the Riccati point spectrum to the point spectrum of $\mathcal{P}^\droite_\veps$ as follows.
\begin{prop}\label{thm_spectreaugmented} 
	Let $\veps>0$. {Under Assumption \ref{assumption_Dioph},} referring to the notation of Proposition \ref{thm_spectre}, the Riccati point spectrum $\sigma^\droite_{R, p, \veps}$ satisfies
	\begin{equation} \label{decomp_spectre_augmente0}
		\sigma^\droite_{R, p, \veps} = \big\{ \lambda^\droite_{k,\veps},\ k \in \Z \big\} \cup 	\big\{1/{\lambda^\droite_{k,\veps}}, \ k \in \Z \big \}.
	\end{equation}
	Moreover,
	\begin{equation} \label{sous_espaces_propres}
		\matrixkernel \big( \mathscr{M}^\droite_\veps - \lambda^\droite_{k,\veps} \, \mathscr{N}^\droite_\veps \big) = \vect \begin{pmatrix} \varphi^\droite_{k,\veps} \\[4pt]
		\mathcal{S}^\droite_\veps \varphi^\droite_{k,\veps} 
	\end{pmatrix}.
	\end{equation}
\end{prop}
\begin{dem} For proving \eqref{decomp_spectre_augmente0}, we prove successively
	$$
		(a)\ \ \sigma^\droite_{R, p, \veps} \cap \big\{ |\lambda| = 1 \big\} = \varnothing, \quad (b)\ \ \sigma^\droite_{R, p, \veps} \cap \big\{ |\lambda| > 1 \big \} = \big\{ 1/\lambda^\droite_{k,\veps} \big\}, \quad (c)\ \ \sigma^\droite_{R, p, \veps} \cap \big\{ |\lambda| < 1 \big \} = \big\{ \lambda^\droite_{k,\veps}\big\}.
	$$
			We use all along this proof the following obvious observation that, in a Banach space $E$, 
			\begin{equation} \label{observ} 
			\mbox{For } A, B, C\in\mathscr{L}(E),\quad A, B\ \textnormal{are invertible}\quad \Longleftrightarrow \quad \begin{pmatrix} A&0  \\[2pt]C&B \end{pmatrix} \ \textnormal{is invertible,}
			\end{equation} 
			$(a)$ For $|\lambda|=1$, we notice that, as $\sigma_p(\transp{\mathcal{P}^\droite_\veps})=\sigma_p(\mathcal{P}^\droite_\veps)$ is included in the circle $\mathcal{C}(0,\varrho(\mathcal{P}^\droite_\veps))$ {(recall that $\varrho(\mathcal{P}^\droite_\veps)<1$}), the operators $\lambda \transp{\mathcal{P}^\droite_\veps} - I$ and $\mathcal{P}^\droite_\veps - \lambda \, I$ are invertible. Therefore, by \eqref{observ}, the three factors in the factorization \eqref{factor} of $\mathscr{M}^\droite_\veps - \lambda \, \mathscr{N}^\droite_\veps$ are invertible.

			\vspspe
			$(b)$ For $|\lambda|>1$, for the same reason as in $(a)$, $\mathcal{P}^\droite_\veps - \lambda \, I$ is invertible. Thus, by \eqref{observ}, the last two factors in \eqref{factor} are invertible, and only the first factor can have a kernel. In other words, by \eqref{observ} again, we obtain
			$$
			\matrixkernel\big( \mathscr{M}^\droite_\veps - \lambda \, \mathscr{N}^\droite_\veps \big) = \{0\} \quad \Longleftrightarrow \quad \matrixkernel\big( \lambda \transp{\mathcal{P}^\droite_\veps} - I \big) = \{0\},
			$$
			which proves $(b)$.

			\vspspe
			$(c)$ For $|\lambda|<1$, now $\lambda \transp{\mathcal{P}^\droite_\veps} - I$ is invertible so that, by \eqref{observ}, the first two factors in \eqref{factor} are invertible and only the third factor can have a kernel. In other words, 
			$$
			\matrixkernel\big( \mathscr{M}^\droite_\veps - \lambda \, \mathscr{N}^\droite_\veps \big) = \{0\} \quad \Longleftrightarrow \quad  \matrixkernel\big( \mathcal{P}^\droite_\veps - \lambda \, I \big) = \{0\} , \quad \mbox{which proves }(c).
			$$
			Moreover, in that case, the kernel of $\mathscr{M}^\droite_\veps - \lambda \, \mathscr{N}^\droite_\veps $ coincides with the one of the third factor in \eqref{factor}
which means nothing but \eqref{sous_espaces_propres}.
\end{dem}
\begin{rmk}\label{rmk:riccati_Pdroite}
{Under Assumption \ref{assumption_Dioph} and by Theorem \ref{thm_spectre}, the Riccati point spectrum is therefore a dense subset of the union of the two circles
$\mathcal{C}(0,|\lambda^\droite_{0,\veps}|)\;\cup\; \mathcal{C}(0,1/|\lambda^\droite_{0,\veps}|).$
\\\\
% Note also that from \eqref{sous_espaces_propres}, the kernel $\matrixkernel \big( \mathscr{M}^\droite_\veps - \lambda^\droite_{k,\veps} \, \mathscr{N}^\droite_\veps \big)$ is linked to the eigenvector of $\mathcal{P}^\droite_\veps$ associated with the eigenvalue $\lambda^\droite_{k,\veps}$.
% \\\\
We now have all the arguments to find the operator $\mathcal{P}^\droite_\veps$ from the  Riccati point spectrum. It suffices first to compute $\sigma^\droite_{R, p, \veps}$ given by (\ref{augmented_EVP}-\ref{eq:Ric_spec}) and select the unique eigenvalue $\lambda_{0,\eps}$ such that $|\lambda_{0,\eps}|<1$ and
\[ \spexists (\varphi_{0,\eps}, \psi_{0,\eps}) \in L^2_{\textit{per}}(\R)^2\setminus\{0\} \mbox{ such that } \mathscr{M}^\droite_\veps \, \begin{pmatrix} \varphi_{0,\eps} \\ \psi_{0,\eps} \end{pmatrix} = \lambda_{0,\eps}\, \mathscr{N}^\droite_\veps \, \begin{pmatrix} \varphi_{0,\eps} \\ \psi_{0,\eps} \end{pmatrix} \]
satisfying $\varphi_{0,\eps}(0)=1$ and ${\bf w}(\varphi_{0,\eps})=0$. 
The pair $(\lambda_{0,\eps},\varphi_{0,\eps})$ is nothing else but the fundamental eigenpair of $\mathcal{P}^\droite_\veps$. The operator $\mathcal{P}^\droite_\veps$ can then be determined entirely using Theorem \ref{thm_spectre}. 
 Moreover, according to \eqref{sous_espaces_propres}, $\psi_{0,\eps}$ is nothing else but $\mathcal{S}^\droite_\veps\varphi_{0,\eps}$ and the operator $\mathcal{S}^\droite_\veps$ can be computed entirely using Remark \ref{rmk:Spsikeps}: 
\[ \mathcal{S}^\droite_\veps\varphi_{k,\eps} = \euler^{2\icplx\pi k (s-\delta)}\, \psi_{0,\eps}(s). \] }
\end{rmk}
\subsection{Relationship between left and right half-waveguides}\label{sec:link_left_right_half_guides}
For now, we have analyzed the half-waveguide problem, the periodicity cell problems and all the associated solutions and associated operators, for $\idomlap = \droite$ only. Note that up to some slight modifications {(replacing in particular $\delta$ by $-\delta$ in all the formulas where $\delta$ is involved)}, the results above can be extended to the half-line problem set on $I^\gauche$, thus leading to solutions of local cell problems $(E^{\gauche, 0}_\veps, E^{\gauche, 1}_\veps)$, local RtR operators $(\mathcal{T}^{\gauche, 00}_\veps, \mathcal{T}^{\gauche, 10}_\veps, \mathcal{T}^{\gauche, 01}_\veps, \mathcal{T}^{\gauche, 11}_\veps)$, propagation and scattering operators $(\mathcal{P}^\gauche_\veps, \mathcal{S}^\gauche_\veps)$, as well as the operators $(\mathscr{M}^\gauche_\veps, \mathscr{N}^\gauche_\veps)$ associated with the Riccati system. 
{Note in particular that Proposition \ref{thm_weighted-shift} extends as follows
\begin{equation}\label{eq:PlSl}
		\spforall s \in \R, \quad \left\{
		\begin{array}{l@{\ =\ }l}
			\mathcal{P}^\gauche_\veps \varphi(s) & \boldsymbol{p}^\gauche_\veps (s+\delta) \, \varphi(s+\delta), 
			\ret
			\mathcal{S}^\gauche_\veps \varphi (s) & \boldsymbol{s}^\gauche_\veps (s+\delta) \, \varphi(s+\delta).
		\end{array}
		\right.
	\end{equation}
	where 
	\begin{equation} \label{defsymbols_l} 
\spforall s \in \R, \quad 
 \begin{array}{|r@{\ \coloneqq \ }l}
	 \boldsymbol{p}^\gauche_\veps (s) & \big(R^\gauche_{+,s} \, u_{s,\veps}^\gauche \big) (-1/\cuti_2),
	\ret
	 \boldsymbol{s}^\gauche_\veps (s) & \big(R^\gauche_{-,s} \, u_{s,\veps}^\gauche \big) (-1/\cuti_2),
\end{array} 
\end{equation}
where $R^\gauche_{\pm,s}$ are the 1D outgoing and ingoing Robin operators, defined for all $s \in \R$ in \eqref{defsRobinoperators}.
		}
~\\\\
The goal of this section is to exhibit the links between the objects defined for $\idomlap = \gauche$ and $\idomlap = \droite$ provided that \eqref{eq:hyp_alar} is satisfied.
% \begin{equation}%\label{eq:hyp}
% (a^\droite-a^\gauche) \, \theta_2\;\in\;\N.
% \end{equation}
% Note that $a^\droite$ and $a^\gauche$, while satisfying \eqref{eq:coef_support} and \eqref{eq:source_terme} can always be chosen so that \eqref{eq:hyp} holds. 

\vspspe
{By \eqref{eq:hyp_alar} and its direct consequence \eqref{eq:cell_lr}}, $\mu_p$ (\resp $\rho_p$) is the same in the periodicity cells $\cellper^{\droite, 0} \coloneqq  (0, 1)^2$ and $\cellper^{\gauche, 0} \coloneqq  (0, 1) \times (-1, 0)$ in which $(E^{\droite, 0}_\veps, E^{\droite, 1}_\veps)$ and $(E^{\gauche, 0}_\veps, E^{\gauche, 1}_\veps)$ are respectively defined. By identifying these cells, it is then easy to show (see the definition of the cell problems \eqref{cellpb0}, \eqref{cellpb1} and also Figure \ref{fig:illust_domains} for the notations) that
\begin{equation}\label{eq:ident_l_r}
	\begin{array}{l}
	\spforall \varphi \in L^2_{\textit{per}}(\R),\quad E^{\gauche,0}_\veps(\varphi)=E^{\droite,1}_\veps(\varphi),\quad E^{\gauche,1}_\veps(\varphi)=E^{\droite,0}_\veps(\varphi),
	\ret
	{\cal T}_\veps^{\gauche, 00}={\cal T}_\veps^{\droite,11},\quad {\cal T}_\veps^{\gauche, 10}={\cal T}_\veps^{\droite,01},\quad {\cal T}_\veps^{\gauche, 01}={\cal T}_\veps^{\droite,10}, \quad {\cal T}_\veps^{\gauche, 11}={\cal T}_\veps^{\droite,00}.
	\end{array}
\end{equation}
\begin{rmk}
	In practice, we only need to solve the local cell problems defined on $\cellper^{\droite, 0}$. The solutions of the local cell problems defined on $\cellper^{\gauche, 0}$ and the associated local RtR operators can then be deduced thanks to \eqref{eq:ident_l_r}.
\end{rmk}
\noindent We deduce in particular that 
\[
	\mathscr{M}_\veps^\gauche = \begin{pmatrix} {\cal T}^{\droite, 10}_\veps & {\cal T}^{\droite, 00}_\veps \\[2pt] 0 & I \end{pmatrix} \quad \textnormal{and} \quad \mathscr{N}_\veps^\gauche = \begin{pmatrix} I & 0 \\[2pt] {\cal T}^{\droite, 11}_\veps & {\cal T}^{\droite, 01}_\veps \end{pmatrix}.
\]
Therefore, from simple algebraic manipulations, it follows that for any $\lambda\in\C$, $\lambda\neq 0$
\begin{equation}\label{eq:link_l_r}
	\mathscr{M}_\veps^\gauche-\lambda \mathscr{N}_\veps^\gauche = -\frac{1}{\lambda}\mathscr{J}^{-1}\Big(\mathscr{M}_\veps^\droite-\frac{1}\lambda \mathscr{N}_\veps^\droite \Big)\mathscr{J}\quad\textnormal{with}\; \mathscr{J}\coloneqq  \begin{pmatrix} 0 & I \\[2pt] I & 0 \end{pmatrix}.
\end{equation}
These relations enables to show the following important result.
\begin{prop}\label{thm_spectreaugmented_2} 
Let $\veps>0$, under Assumption \ref{assumption_Dioph}, we have
\[
	\sigma_p(\mathcal{P}^\gauche_\veps) = \sigma_p(\mathcal{P}^\droite_\veps)\quad\textnormal{and}\quad\matrixkernel \Big( \mathscr{M}^\idomlap_\veps-\frac{1}{\lambda^{\idomlap'}_{k,\veps}} \mathscr{N}^\idomlap_\veps \Big) = \vect \begin{pmatrix} \mathcal{S}^{\idomlap'}_\veps \varphi^{\idomlap'}_{k,\veps} \\[4pt]
		\varphi^{\idomlap'}_{k,\veps} 
	\end{pmatrix} \quad\textnormal{for}\ \idomlap \neq \idomlap' \in \{\gauche, \droite\}.
\]
\end{prop}			
\begin{dem}
	Let $\lambda^\gauche_\veps\in \sigma_p(\mathcal{P}^\gauche_\veps)$. Then by Proposition \ref{thm_spectre}, there exists $k\in\Z$ such that $\lambda^\gauche_\veps=\lambda^\gauche_{k,\veps}$. Note that $|\lambda^\gauche_{k,\veps}|<1$. Using \eqref{sous_espaces_propres} and \eqref{eq:link_l_r}, we deduce that
	\[
		\matrixkernel \Big( \mathscr{M}^\droite_\veps - \frac{1}{\lambda^\gauche_{k,\veps}} \, \mathscr{N}^\droite_\veps \Big) = \vect \begin{pmatrix} \mathcal{S}^\gauche_\veps \varphi^\gauche_{k,\veps} \\[4pt]
			 \varphi^\gauche_{k,\veps}
			\end{pmatrix} 
	\]
	This implies that $(\lambda^\gauche_{k,\veps})^{-1}\in \sigma^\droite_{R,\veps}$, and therefore by \eqref{decomp_spectre_augmente0}, 
	$$ \spforall k \in \Z, \quad \exists \; k'\in\Z \ /\ \lambda^\gauche_{k,\veps}=\lambda^\droite_{k',\veps}\in \sigma_p(\mathcal{P}^\gauche_\veps).$$
	By inverting the roles of $\gauche$ and $\droite$, one shows that 
	$$ \spforall k \in \Z, \quad \exists \; k'\in\Z \ /\ \lambda^\droite_{k,\veps}=\lambda^\gauche_{k',\veps}\in \sigma_p(\mathcal{P}^\droite_\veps).$$ 

	\vspace{-1.5\baselineskip} \noindent
\end{dem}

\begin{rmk}\label{rmk:riccati_Pgauche}
{
	 In practice to find the operator $\mathcal{P}^\gauche_\veps$, under \eqref{eq:hyp_alar}, we only need to solve \eqref{augmented_EVP} and compute the Riccati point spectrum (see Remark \ref{rmk:riccati_Pdroite} for the computation of $\mathcal{P}^\droite_\veps$). 
	Indeed, according to Proposition \ref{thm_spectreaugmented_2}, it suffices to compute first $\sigma^\droite_{R, p, \veps}$ given by (\ref{augmented_EVP}-\ref{eq:Ric_spec}) and select the unique eigenvalue $\tilde{\lambda}_{0,\eps}$ such that $|\tilde{\lambda}_{0,\eps}|>1$ and
\[ \spexists (\tilde{\varphi}_{0,\eps}, \tilde{\psi}_{0,\eps}) \in L^2_{\textit{per}}(\R)^2\setminus\{0\} \mbox{ such that } \mathscr{M}^\droite_\veps \, \begin{pmatrix} \tilde{\varphi}_{0,\eps} \\ \tilde{\psi}_{0,\eps} \end{pmatrix} = \tilde{\lambda}_{0,\eps}\, \mathscr{N}^\droite_\veps \, \begin{pmatrix}  \tilde{\varphi}_{0,\eps}  \\ \tilde{\psi}_{0,\eps} \end{pmatrix}, \]
satisfying $\tilde{\psi}_{0,\eps}(0)=1$ and ${\bf w}(\tilde{\psi}_{0,\eps})=0$. 
The pair $(1/\tilde{\lambda}_{0,\eps},\tilde{\psi}_{0,\eps})$ is nothing else but the fundamental eigenpair of $\mathcal{P}^\gauche_\veps$. The operator can then be determined entirely using a slight modification of Theorem \ref{thm_spectre} for $\idomlap=\gauche$ (replacing $\delta$ by $-\delta$ in all the formulas): 
if we define $\varphi^\gauche_{k,\veps} \in {\cal N}_{\textit{per}}$ for any $k \in \Z$ by $\varphi^\gauche_{k,\veps}(s) \coloneqq  \euler^{2 \icplx \pi k s} \, \tilde{\psi}_{0,\eps}(s)$, then
		\begin{equation} \label{eigenfunctions_gauche}
	\mathcal{P}^\gauche_\veps \; \varphi^\gauche_{k,\veps} = \lambda^\gauche_{k,\veps} \; \varphi^\gauche_{k,\veps} \quad \mbox{with} \quad \lambda^\gauche_{k,\veps} \coloneqq  \euler^{2 \icplx \pi k \delta}/\tilde{\lambda}_{0,\eps}.
	\end{equation}
Moreover, $\tilde{\varphi}_{0,\eps}=\mathcal{S}^\gauche_\veps\tilde{\psi}_{0,\eps}$ and the operator $\mathcal{S}^\gauche_\veps$ can be computed entirely using a slight adaptation {(again by replacing $\delta$ by $-\delta$)} of Remark \ref{rmk:Spsikeps}: \[ \mathcal{S}^\gauche_\veps\varphi^\gauche_{k,\eps} = \euler^{2\icplx\pi k (s+\delta)}\, \tilde{\varphi}_{0,\eps}(s). \]}
\end{rmk}

\subsection{Solution algorithm in the absorbing case} \label{sec:algo}
Let us still suppose that Assumption \ref{assumption_Dioph} is satisfied.  In order to compute the solution of \eqref{eq:whole_line_problem_eps}, the previous sections provide an algorithm which sums up as follows.
\begin{enumerate} 
	\item Compute for all $s\in[0,1)$ the functions $e_{s, \veps}^{\droite, 0}$ and $e_{s, \veps}^{\droite, 1}$, solutions respectively of \eqref{segment_problems0} and \eqref{segment_problems1} as well as the local RtR symbols $t^{\droite, \ell k}_\varepsilon(s)$  for $\ell,k\in\{0,1\}$ defined in \eqref{defRtRcoeff}. Deduce the cell operators $E^{\droite,0}_\veps$ and $E^{\droite,1}_\veps$ via formulas \eqref{fibered-structure-cell}  as well as the local RtR operators ${\cal T}^{\droite,\ell k}_\veps$  for $\ell,k\in\{0,1\}$ via \eqref{RtRweighted}. Form the operators $\mathscr{M}^\droite_\veps$ and $ \mathscr{N}^\droite_\veps$ defined in \eqref{eq:Riccati_op}.
	\item Compute the cell operators $E^{\gauche,0}_\veps$ and $E^{\gauche,1}_\veps$ as well as the local RtR operators ${\cal T}^{\gauche,\ell k}_\veps$ for $\ell,k\in\{0,1\}$ using \eqref{eq:ident_l_r}. Replacing $\delta$ by $-\delta$ and the superscript $\droite$ by $\gauche$, deduce the functions $e_{s, \veps}^{\gauche, 0}$ and $e_{s, \veps}^{\gauche, 1}$ via formulas \eqref{fibered-structure-cell} applied with $\varphi = 1$ and then the local RtR symbols $t^{\gauche, \ell k}_\varepsilon(s)$ for $\ell,k\in\{0,1\}$ via \eqref{RtRweighted}.
	\item Compute the unique solution $(\lambda_{\varepsilon}, \varphi_{{\varepsilon}}, \psi_{{\varepsilon}}) $ of the augmented spectral problem \eqref{augmented_EVP} such that $|\lambda_{\varepsilon}| < 1$,  $(\varphi_{{\varepsilon}}, \psi_{{\varepsilon}}) \in \mbox{Ker } (\mathscr{M}^\droite_\varepsilon - \lambda_\varepsilon \mathscr{N}^\droite_\varepsilon )$, ${\bf w}(\varphi_{{\varepsilon}}) = 1$ and $\varphi_{{\varepsilon}}(0) = 1$ and set $\lambda_{0,\varepsilon} ^\droite := \lambda_{\varepsilon}$, and $(\varphi_{0,\varepsilon}^\droite, \psi_{0,\varepsilon}^\droite):=(\varphi_{{\varepsilon}},\psi_{{\varepsilon}})$  according to Theorem \ref{thm_spectre} and Proposition \ref{thm_spectreaugmented}.
		\item Compute the unique solution $(\tilde\lambda_{\varepsilon}, \tilde\varphi_{{\varepsilon}}, \tilde\psi_{{\varepsilon}}) $ of the augmented spectral problem \eqref{augmented_EVP} such that $|\tilde\lambda_{\varepsilon}| > 1$,  $(\tilde\varphi_{{\varepsilon}}, \tilde\psi_{{\varepsilon}}) \in \mbox{Ker } (\mathscr{M}^\droite_\varepsilon - \tilde\lambda_\varepsilon \mathscr{N}^\droite_\varepsilon )$, ${\bf w}(\tilde\psi_{{\varepsilon}}) = 1$ and $\tilde\psi_{{\varepsilon}}(0) = 1$ and set $\lambda_{0,\varepsilon} ^\gauche := 1/\lambda_{\varepsilon}$, $\varphi_{0,\varepsilon}^\gauche :=\tilde\psi_{{\varepsilon}}$ and $\psi_{0,\varepsilon}^\gauche :=\tilde\varphi_{{\varepsilon}}$ according to Theorem \ref{thm_spectre} and Proposition \ref{thm_spectreaugmented_2}.
	\item Compute the RtR coefficient $\RtRcoef^\droite_\varepsilon$ via formula \eqref{def_lambda_1D_dir} : note that this formula uses $t^{\droite,00}_\varepsilon(s)$ and $t^{\droite,10}_\varepsilon(s)$ from Step 1 and the functions $\varphi_{0,\varepsilon}^\droite, \psi_{0,\varepsilon}^\droite$ from Step 3. Similarly, compute the RtR coefficient $\RtRcoef^\gauche_\varepsilon$. 
	\item Solve the interior problem \eqref{eq:int_pb_eps_RtR} to determine the solution $u^\varepsilon$ inside $ I^\interieur$. 
	Reconstruct the solution $u^\varepsilon$  in $ I^\droite$ (\resp $I^\gauche$) via the reconstruction formula \eqref{eq:sol_eps}-\eqref{eq:recons_1d0} (\resp  replacing in \eqref{eq:recons_1d0} $\delta$ by $-\delta$ and the superscript $\droite$ by $\gauche$). Note that this formula uses the eigenvalue 
	$\lambda_{0,\varepsilon} ^\droite$ and the functions $\varphi_{0,\varepsilon} ^\droite, \psi_{0,\varepsilon}^\droite$ from Step 3 (\resp $\lambda_{0,\varepsilon} ^\gauche$, $\varphi_{0,\varepsilon} ^\gauche, \psi_{0,\varepsilon}^\gauche$ from Step 4)  as well as the functions $e_{s, \veps}^{\droite, 0}$ and $e_{s, \veps}^{\droite, 1}$ from Step 1 (\resp $e_{s, \veps}^{\gauche, 0}$ and $e_{s, \veps}^{\gauche, 1}$  from Step 2). 
	\end{enumerate} 
\begin{rmk}
The above algorithm is at the theoretical level. Its effective implementation requires a numerical approximation via a discretization procedure as it has been done in \cite{amenoagbadji2023wave}.
	\end{rmk}
	
%	\item Compute the unique solution $(\lambda_{\varepsilon}, \varphi_{\lambda_{\varepsilon}}, \psi_{\lambda_{\varepsilon}}) $ of the augmented spectral problem (4.29) such that $|\lambda_{\varepsilon}| > 1$,  $(\varphi_{\lambda_{\varepsilon}}, \psi_{\lambda_{\varepsilon}}) \in \mbox{Ker } (\mathscr{M}_\varepsilon - \lambda_\varepsilon \mathscr{N}_\varepsilon )$, ${\bf w}(\varphi_{\lambda_{\varepsilon}}) = 1$ and $\varphi_{\lambda_{\varepsilon}}(0) = 1$ and set $\lambda_{0,\varepsilon} ^l := \lambda_{\varepsilon}$, $\varphi_{0,\varepsilon}^l :=\psi_{\lambda_{\varepsilon}}$ and $\psi_{0,\varepsilon}^r :=\varphi_{\lambda_{\varepsilon}}$ according to Theorem 4.5 and Proposition 4.15. 
{\noindent It is also possible to compute directly for any $\varphi\in L^2_{\textit{per}}(\R)$, the half-guide solution $U^{\idomlap}_\veps(\varphi)$ and not only the half-line solution. After the Steps 1. to 4., one has to add the following steps. 
\begin{enumerate}
	 \item[5-bis.] Compute $\mathcal{P}^\droite_\veps$ and $\mathcal{S}^\droite_\veps$ (\resp $\mathcal{P}^\gauche_\veps$ and $\mathcal{S}^\gauche_\veps$) using Theorem \ref{thm_spectre} and Remark  \ref{rmk:riccati_Pdroite} (\resp Remark \ref{rmk:riccati_Pgauche}).
	 \item[6-bis.] Reconstruct for any $\varphi\in L^2_{\textit{per}}(\R)$, the solution $U^\idomlap_\veps(\varphi)$ cell by cell using Proposition \ref{thm_reconstruction} (replacing the superscript $\droite$ by $\gauche$ when $\idomlap=\gauche$).
\end{enumerate}
}
\begin{rmk}
When Assumption \ref{assumption_Dioph} is not satisfied, one cannot use the previous algorithm that is based on Theorem \ref{thm_spectre}. However, one can still solve the Riccati system \eqref{Riccatiplus} using similar algorithms than the one described in \cite{joly2006exact,fliss2009exact,amenoagbadji2023wave}.
\end{rmk}

\section{Towards limiting absorption: study of the spectrum of \\quasiperiodic  operators}\label{sec:prelim_LAP}
As explained in Section \ref{sec_QD}, the problem \eqref{eq:whole_line_problem} can be solved in $H^1(\R)$ only when $\omega^2$ is not in $\sigma(\mathcal{A})$, the spectrum of the positive self-adjoint operator ${\cal A}$, see \eqref{eq:op_pb}. On the other hand, the question of possible limiting absorption, in a space other than $H^1(\R)$, will be studied for values of $\omega^2$ outside the discrete spectrum $\sigma_d(\mathcal{A})$: these values of $\omega^2$, that are related to the existence of possible trapped modes, will be excluded from this study. Thus, according to the decomposition \eqref{eq:spect_A} of $\sigma(\mathcal{A})$, we shall look at the limiting absorption principle for $\omega^2$ in the spectrum of the (purely quasiperiodic) operator ${\cal A}_\cutvec$ defined by \eqref{eq:op_quasiper}. But before doing so, we note that our study also involves the spectra of other differential operators, such as the 2D periodic differential operator associated with $\mathcal{A}_\cutvec$ (namely the operator $\mathcal{A}_p$ defined in \eqref{eq:op_per}), and its equivalent on the half-guide $\Omega^\droite$ (with Robin boundary conditions; see \eqref{eq:op_per_j} below). The goal of the present section is to exhibit for these spectra some equalities or inclusions that will be useful in Section \ref{sec:LAP}.

% or the half-line (\resp half-guide) equivalent of $\mathcal{A}_\cutvec$ (\resp $\mathcal{A}_p$) (with Robin boundary conditions; see (\ref{eq:op_quasipers_j}, \ref{eq:op_per_j}) below). The goal of the present section is to show that these spectra coincide, so that there is no need to distinguish between them.

% \vspspe
%\subsection{The case of the whole quasiperiodic real line}
%We are here interested in the spectrum of the self-adjoint operator ${\cal A}_\cut$. 
\subsection{Quasiperiodic operators on the whole line}
In the spirit of the lifting method presented in our study, a fruitful point of view is to see ${\cal A}_\cutvec$ as one particular element of a family of self-adjoint operators ${\cal A}_{s, \cutvec}, s \in \R$ with 
\begin{equation}
	\label{eq:op_quasipers}
	D(\mathcal{A}_{s, \cutvec}) \coloneqq  \Big\{u\in H^1(\R) \ \big/\ \mu_{s, \cutvec} \, \frac{d u}{d \xvi}\in H^1(\R)\Big\} \quad \textnormal{and} \quad \mathcal{A}_{s, \cutvec}\, u \coloneqq  -\frac{1}{\rho_{s, \cutvec}}\frac{d}{d \xvi} \Big( \mu_{s, \cutvec} \, \frac{d u}{d \xvi} \Big),% \quad \spforall u \in D(\mathcal{A}_{s, \cutvec}),
\end{equation}
where we recall that the functions $(\mu_{s, \cutvec}, \rho_{s, \cutvec})$ have been defined in \eqref{defrhosmus} as the traces of the 2D periodic functions $(\mu_p,\rho_p)$ along the line $\cutvec\, \R + s\, \vv{\ev}_1$. Obviously, $s \mapsto {\cal A}_{s, \cutvec}$ is 1--periodic and according to its definition \eqref{eq:op_quasiper}, ${\cal A}_\cutvec$ coincides with $\mathcal{A}_{s, \cutvec}$ for $s = 0$.
% \begin{equation}
% \label{eq:op_quasiper0}
% \mathcal{A}_{\cutvec} = 	\mathcal{A}_{0, \cutvec}.
% \end{equation}
We can build from the family $\mathcal{A}_{s, \cutvec}, s \in [0,1)$, a self-adjoint operator in $L^2_{\textit{per}}\big(\R, L^2(\R)\big)$, as follows
\begin{equation}
\label{eq:op_quasiper2Dtheta}
\mathcal{A}^{\textnormal{\tiny 2D}}_{\cutvec} = \int^\oplus_{[0,1]}	\mathcal{A}_{s, \cutvec} \, ds
\end{equation}
meaning that (see \cite[Section XIII.16]{reed1978iv} for instance for direct integrals of operators), 
\begin{equation}
\label{eq:op_quasiper2Dtheta2}
U \in D(\mathcal{A}^{\textnormal{\tiny 2D}}_{\cutvec}) \quad \Longleftrightarrow \quad \spforall s\in[0,1], \; U(s,\cdot) \in D(\mathcal{A}_{s, \cutvec}) \; \textnormal{and} \;  [\mathcal{A}^{\textnormal{\tiny 2D}}_{\cutvec}\, U] (s,\cdot)  = \mathcal{A}_{s, \cutvec} [U(s,\cdot)].
\end{equation} 
It is known, see \cite[Theorem XIII.85]{reed1978iv}, that the spectra of $	\mathcal{A}_{\cutvec} $ and the $\mathcal{A}_{s, \cutvec}$ are related by 
\begin{equation}\label{floquet}
\sigma\big(\mathcal{A}_{\cutvec}^{\textnormal{2D}}\big)=\bigcup_{s \in [0,1]} \sigma\big(\mathcal{A}_{s, \cutvec}\big).
\end{equation}
Moreover, the operator $\mathcal{A}_{\cutvec}^{\textnormal{2D}}$ is closely linked to the 2D self-adjoint operator $\mathcal{A}_p$ corresponding to the augmented half-waveguide problems, that is,% (for $\Omega \coloneqq  (0, 1) \times \R$),
\begin{equation}
  \label{eq:op_per}
   D(\mathcal{A}_p) \coloneqq  \big\{U \in H^1_{\cutvec,\, \textit{per}}(\R^2)\ \big/\ \mu_p \, \Dt{} U\in H^1_{\cutvec,\, \textit{per}}(\R^2) \big\}, \quad \mathcal{A}_{p}\, U=-\frac{1}{\rho_p}\, \Dt{} \big( \mu_p \, \Dt{}U\big),
\end{equation}
where we recall that $H^1_{\cutvec,\, \textit{per}}(\R^2) \coloneqq  \{U \in L^2_{\textit{per}}(\R^2)\ /\ \Dt{} U \in L^2_{\textit{per}}(\R^2)\}$. More precisely, the chain rule \eqref{eq:chain_rule} allows to show that $\mathcal{A}_{\cutvec}^{\textnormal{2D}}$ and ${\cal A}_p$ are equivalent in $L^2_{\textit{per}}(\R^2)$:
$$
\mathcal{A}_{\cutvec}^{\textnormal{2D}} = \mathcal{U}_\cutvec \, \mathcal{A}_p \, \mathcal{U}_\cutvec^{-1},
$$
where $\mathcal{U}_\cutvec$ is the invertible operator from $L^2_{\textit{per}}(\R^2)$ into $L^2_{\textit{per}}(\R, L^2_{\textit{loc}}(\R))$ defined by 
\begin{equation*}% \label{defStheta}
 \spforall U \in L^2_{\textit{per}}(\R^2), \quad (\mathcal{U}_\cutvec\, U) (s,x) = U (s + \theta_1 \, x, \theta_2 \, x).
\end{equation*} 
As a consequence, we have 
\begin{equation}\label{floquet2}
  \sigma \big(\mathcal{A}_p\big)=\bigcup_{s \in [0,1]} \sigma\big(\mathcal{A}_{s, \cutvec}\big).
\end{equation}
The spectra $\sigma\big(\mathcal{A}_{s, \cutvec}\big)$ all coincide if the ratio $\delta$ is irrational, as the next result shows.
\begin{prop}\label{thm:spectrunqp_s}
  As soon as $\delta \in \R \setminus \Q$, one has
  \begin{equation} \label{eq_spectres_s} 
    \spforall s \in \R, \quad \sigma(\mathcal{A}_{s, \cutvec}) = \sigma({\cal A}_\cutvec).% \qquad\textnormal{and} \qquad (b) \ \  	\spforall \fbvar\in[-\pi,\pi[,\quad \sigma(\mathcal{A}_{p})=\sigma(\mathcal{A}_p (\fbvar)).
  \end{equation} 
  In particular, $\sigma({\cal A}_\cutvec)= \sigma(\mathcal{A}_{p})$.
\end{prop}
\begin{dem}
  The key idea is to show that $s \mapsto \sigma(\mathcal{A}_{s, \cutvec})$ is a continuous function (in the sense of the Hausdorff distance) which is both $1$ and $\delta$--periodic. It is then the irrationality of $\delta$ that allows to conclude. The proof is detailed in Appendix \ref{app:spectrum}, see Proposition \ref{prop:spec_1D_2D_ligne_indt_s}.
\end{dem}

\vspspe
We shall now highlight another fibered structure for $\mathcal{A}_p$ (similar to the structure \eqref{eq:op_quasiper2Dtheta} for $\mathcal{A}^{\textnormal{\tiny 2D}}_{\cutvec}$) by means of the Floquet-Bloch transform in the $y_2$--direction. To begin, we introduce the space $L^2_{\textit{per}, 2}(\R^2)$ of locally square-integrable functions that are periodic with respect to both $\yvi_1$ \textbf{and} $\yvi_2$ (this is to be distinguished from $L^2_{\textit{per}}(\R^2)$ that only imposes periodicity in the $y_1$--direction), as well as $H^1_{\cutvec,\, \textit{per},2}(\R^2)$, the space defined by:
\begin{equation*}
  \left\{
    \begin{array}{r@{\ \coloneqq \ }l}
      L^2_{\textit{per}, 2}(\R^2) & \big\{ V \in L^2_{loc}(\R^2) \ / \ V(y_1 +1, y_2) = V(y_1 , y_2+1) = V(y_1, y_2) \big\},
      \ret
      H^1_{\cutvec,\, \textit{per}, 2}(\R^2) & \big\{U \in L^2_{\textit{per}, 2}(\R^2)\ /\ \Dt{} U \in L^2_{\textit{per}, 2}(\R^2)\big\}.
    \end{array}
  \right.
\end{equation*}
The Floquet-Bloch transform in the $\yvi_2$--direction is defined by
\begin{equation*}
  \aeforall (\yv, \fbvar) \in \R^2 \times (-\pi, \pi), \quad \mathcal{F}\, U (\yv, \fbvar) \coloneqq  \sqrt{\frac{1}{2\pi}}\; \sum_{n \in \Z} U(\yv + n\, \vv{\ev}_2)\; \euler^{-\icplx\vts \fbvar\vts (\yvi_2 + n)},
  %
  % {\cal F}: U(y_1, y_2) \mapsto U(\fbvar; y_1, y_2), \quad \fbvar \in [ \, -\pi, \pi ),  
\end{equation*}
It is well-known \cite{kuchment1993floquet} that $\mathcal{F}$ defines a unitary map from $L^2_{\textit{per}}(\R^2)$ into $L^2(-\pi, \pi; L^2_{\textit{per}, 2}(\R^2))$, and from $H^1_{\textit{per}}(\R^2)$ to $L^2(-\pi, \pi; H^1_{\cutvec,\, \textit{per}, 2}(\R^2))$. Consider the self-adjoint operators defined for $\fbvar \in [-\pi, \pi]$ as
 \begin{equation}
 \label{eq:op_per_k}
 \left\{ \begin{array}{r@{\ \coloneqq \ }l}
 D(\mathcal{A}_p (\fbvar)) & \big\{ U \in H^1_{\cutvec,\, \textit{per},2}(\R^2)\ /\ \mu_p \, \Dt{} U\in H^1_{\cutvec,\, \textit{per},2}(\R^2) \big\}, 
 \ret
 \mathcal{A}_p (\fbvar)\, U & \displaystyle -\frac{1}{\rho_p}\, (\Dt{} + \icplx\vts \fbvar\vts \theta_2)\; \mu_p \; (\Dt{} + \icplx\vts \fbvar\vts \theta_2 )\, U,
  \end{array} \right.
\end{equation}
From the properties of the Floquet-Bloch transform \cite{kuchment1993floquet}, it follows that the self-adjoint operator given by $\widehat {\mathcal{A}}_p \coloneqq  {\cal F}	\mathcal{A}_p \, {\cal F}^{-1}$ can be expressed as a direct integral
\begin{equation}
 \label{decompAp2}
\widehat {\mathcal{A}}_p = \int^\oplus_{[-\pi,\pi]}	\mathcal{A}_p (\fbvar) \, d\fbvar,
\end{equation}
meaning similarly to \eqref{eq:op_quasiper2Dtheta2} that
\begin{equation}
 \label{eq:op_quasiper2Dtheta3}
 \widehat U \in D\big(	\widehat {\mathcal{A}}_p \big) \quad \Longleftrightarrow \quad \spforall \fbvar, \quad U(\fbvar,\cdot) \in D(\mathcal{A}_p (\fbvar))\quad\text{and} \quad (\widehat {\mathcal{A}}_p\, \widehat U)(\fbvar,\cdot)  = \mathcal{A}_p (\fbvar)\, [\vts \widehat U(\fbvar,\cdot)].
\end{equation} 
As $\widehat {\mathcal{A}}_p$ and ${\mathcal{A}}_p$ are unitary equivalent, we deduce another decomposition of $\sigma\big(\mathcal{A}_p\big)$
\begin{equation}\label{floquet3}
 \sigma\big(\mathcal{A}_p\big)=\bigcup_{\fbvar \in [-\pi, \pi]} \sigma\big(\mathcal{A}_p (\fbvar)\big).
\end{equation}
Note that contrary to elliptic operators, $\mathcal{A}_{p}(\fbvar)$ does not have compact resolvent, and thus its spectrum is not necessarily discrete. Furthermore, similarly to Proposition \ref{thm:spectrunqp_s}, the next result shows that the spectrum of $\mathcal{A}_{p}(\fbvar)$ is surprisingly independent of $\fbvar$ if $\delta$ is irrational. 
% The two decompositions \eqref{floquet2} and \eqref{floquet3} are valid for any real $\delta$, but take a particular interest when $\delta \in \R \setminus \Q$ because all the spectra appearing in  \eqref{floquet2} and \eqref{floquet3} surprisingly coïncide.
\begin{prop}\label{thm:spectrunqp_k}
As soon as $\delta \in \R \setminus \Q$, one has
\begin{equation} \label{eq_spectres_k} 
  \spforall \fbvar\in[-\pi,\pi[,\quad \sigma(\mathcal{A}_{p}) = \sigma(\mathcal{A}_p (\fbvar)).
\end{equation}
\end{prop}
\begin{dem}
The proof, which is detailed in Appendix \ref{app:spectrum} (Proposition \ref{prop:spec_1D_2D_ligne_indt_xi}), relies on showing that $\fbvar \mapsto \sigma(\mathcal{A}_p (\fbvar))$ is a continuous function (in the sense of the Hausdorff distance) which is both $2\pi$ and $2\pi \delta$--periodic. One then concludes using the irrationality of $\delta$.
\end{dem}

\subsection{Quasiperiodic operator on a half-line}
The RtR approach developed in Section \ref{sec:RtR_eps} in the absorbing case involves the operators $\mathcal{A}^\droite_{s, \cutvec}$, $s \in \R$, associated with the half-line problems \eqref{halfline_sproblems}, namely
\begin{equation}
  \label{eq:op_quasipers_j}
  \left\{
    \begin{array}{r@{\ \coloneqq \ }l}
      D(\mathcal{A}^\droite_{s, \cutvec}) & \displaystyle \Big\{u\in H^1(\R_+) \ \big/\ \mu_{s, \cutvec} \, \frac{d u}{d \xvi}\in H^1(\R_+)  \quad \textnormal{and} \quad R^\droite_{+, s}\, u = 0 \Big\} 
      \ret
      \mathcal{A}^\droite_{s, \cutvec}\, u & \displaystyle -\frac{1}{\rho_{s, \cutvec}}\frac{d}{d \xvi} \Big( \mu_{s, \cutvec} \, \frac{d u}{d \xvi} \Big),
    \end{array}
  \right.
\end{equation}
and the operator $\mathcal{A}^\droite_p$ associated with the half-guide problem \eqref{RobinHalfspaceproblems}, namely
\begin{equation}
  \label{eq:op_per_j}
  \left\{
    \begin{array}{r@{\ \coloneqq \ }l}
      D(\mathcal{A}^\droite_p) & \displaystyle\big\{U \in H^1_{\cutvec,\, \textit{per}}(\Omega^\droite)\ \big/\ \mu_p \, \Dt{} U\in H^1_{\cutvec,\, \textit{per}}(\Omega^\droite) \quad \textnormal{and} \quad \mathcal{R}^\droite_+ U = 0\big\} 
      \ret
      \mathcal{A}^\droite_{p}\, U & \displaystyle -\frac{1}{\rho_p}\, \Dt{} \big( \mu_p \, \Dt{}U\big).
    \end{array}
  \right.
\end{equation}
One defines similarly the operators $\mathcal{A}^\gauche_{s, \cutvec}$ and $\mathcal{A}^\gauche_p$. Note that the results of the previous section do not hold a priori, in particular because $\mathcal{A}^\idomlap_{s, \cutvec}$ and $\mathcal{A}^\idomlap_{p}$ are not self-adjoint. However, the next result holds.
  
\begin{prop}\label{prop:spectrum_halves}
  Let $\idomlap \in \{\gauche, \droite\}$. One has 
 \[\spforall\omega\in\R,\quad \omega^2\notin \sigma(\mathcal{A}_{0, \cutvec})\quad\Rightarrow \quad \omega^2\notin \sigma(\mathcal{A}^\idomlap_{s, \cutvec})\;\spforall s\in\R\quad \text{and} \quad \omega^2\notin \sigma(\mathcal{A}^\idomlap_p). \]
\end{prop}
\begin{dem}
	{ 
%   Since $\sigma(\mathcal{A}^\idomlap_{s, \cutvec}) = \sigma(\mathcal{A}^\idomlap_p)$ for $\idomlap \in \{\gauche, \droite, \varnothing\}$, it is sufficient to prove the first part of the proposition, that is, $\sigma(\mathcal{A}^\idomlap_{s, \cutvec}) \subset \sigma(\mathcal{A}_{s, \cutvec})$. 
  For clarity, we fix $\idomlap = \droite$, and we focus on $(\mathcal{A}_{\cutvec}, \mathcal{A}^\droite_{\cutvec}):= (\mathcal{A}_{s, \cutvec}, \mathcal{A}^\droite_{s, \cutvec})$ for $s = 0$. The proof for any $s\in\R$ is very similar. {Similarly, we focus on $\mathcal{A}^\droite_p$. In what follows, $\omega\in\R$ is such that $\omega^2 \not\in \sigma(\mathcal{A}_{\cutvec})$.} }
  
  {\paragraph*{Step 1}
  We first prove that $\omega^2 \not\in \sigma(\mathcal{A}^\droite_{\cutvec})$, or equivalently, that the problem
  % Following the ideas of \cite{hoang2011limiting} (also used in \cite{kirsch2018radiation} and \cite{fliss2021dirichlet}, and which simplify greatly for our 1D operators), we proceed by contradiction: 
  \begin{equation}\label{eq:halfline_problem_f}
    \displaystyle
    \left|
      \begin{array}{r@{\ =\ }l}
        \displaystyle - \frac{d}{dx} \Big(\mu_\cutvec \,  \frac{d u^\droite}{dx}\Big) - \rho_\cutvec \, \omega^2 \, u^\droite & f \quad \mbox{in}\ \ \R_+,
        \ret
        R^\droite_{+} u^\droite(0) & 0,
      \end{array}
    \right.
  \end{equation}
  admits a unique solution $u^\droite \in H^1(\R_+)$ for any $f \in L^2(\R_+)$.}

  \vspspe
  \textbf{\textit{Uniqueness.}} Let $\widetilde{u}^\droite \in H^1(\R_+)$ satisfy \eqref{eq:halfline_problem_f} with $f = 0$. Multiplying the ODE by $\overline{\widetilde{u}^\droite}$, integrating over $\R_+$, and using the Robin boundary condition leads to
  \[
    \displaystyle
    \int_{\R_+} \, \Big(\mu_\cutvec \, \Big|\frac{d \widetilde{u}^\droite}{dx} \Big|^2 - \rho_\cutvec \, \omega^2\, |\widetilde{u}^\droite|^2 \Big) \, dx  -  \icplx \vts z \vts |\widetilde{u}^\droite(0)|^2 = 0.
  \]
  Taking the imaginary part of this expression leads to $\widetilde{u}^\droite(0) = 0$, and then to $(\widetilde{u}^\droite)'(0) = 0$ due to the Robin condition. It then follows from Cauchy uniqueness theorem that $\widetilde{u}^\droite = 0$.%   the difference $\widetilde{u}^\droite:= \widetilde{u}^\droite_1 - u^\droite_2$ satisfies the ODE in \eqref{eq:halfline_problem_f} with . B

  \vspspe
  \textbf{\textit{Existence.}} Given $f \in L^2(\R_+)$, we consider its extension by $0$ on $\R$, still denoted by $f \in L^2(\R)$. We begin by seeking $u^\droite$ as%under the form
  \begin{equation*}
    u^\droite = u_f|_{\R_+} - (R^\droite_+ u_f (0))\, u^\droite_1,
  \end{equation*}
  where $u_f$ is the unique solution of the well-posed ($\omega^2 \not\in \sigma(\mathcal{A}_\cutvec)$) problem 
  \begin{equation*}
    - \frac{d}{dx} \Big(\mu_\cutvec \,  \frac{d u_f }{dx}\Big) - \rho_\cutvec \, \omega^2 \, u_f  = f \quad \mbox{in}\ \ \R,
  \end{equation*}
  and where $u^\droite_1$ is a correction term that, for $u^\droite$ to be solution of \eqref{eq:halfline_problem_f}, should satisfy
  \begin{equation}\label{eq:halfline_problem_0}
    \displaystyle
    \left|
      \begin{array}{r@{\ =\ }l}
        \displaystyle - \frac{d}{dx} \Big(\mu_\cutvec \,  \frac{d u^\droite_1}{dx}\Big) - \rho_\cutvec \, \omega^2 \, u^\droite_1 & 0 \quad \mbox{in}\ \ \R_+,
        \ret
        R^\droite_{+} u^\droite_1(0) & 1.
      \end{array}
    \right.
  \end{equation} 
  In what follows, we construct $u^\droite_1$ following the ideas of \cite[Theorem 4.1, Lemma 4.2]{hoang2011limiting} (also used in \cite{kirsch2018radiation} and \cite{fliss2021dirichlet}, which simplify greatly for our 1D operators). More precisely, we look for $u^\droite_1$ under the form 
  \begin{equation}\label{eq:decomp_usdroite_us1_us2}
    u^\droite_1 = \alpha\, (u^\droite_2 + u_g|_{\R_+}),
  \end{equation}
 where $u^\droite_2$ is solution of the coercive problem
  \begin{equation}\label{eq:halfline_problem_coercive_dem}
    \displaystyle
    \left|
      \begin{array}{r@{\ =\ }l}
        \displaystyle - \frac{d}{dx} \Big(\mu_\cutvec \,  \frac{d u^\droite_2 }{dx}\Big) - \rho_\cutvec \, (\omega^2 + \icplx) \, u^\droite_2  & 0 \quad \mbox{in}\ \ \R_+,
        \ret
        u^\droite_2 (0) & 1.
      \end{array}
    \right.
  \end{equation}
  where $u_g$ satisfies the same equation as $u_f$ on $\R$, but with $g \in L^2(\R)$ as a source term. 
  Our goal is to find the pair $(\alpha, g) \in \C \times L^2(\R)$ such that the function $u^\droite_1$ given by \eqref{eq:decomp_usdroite_us1_us2} satisfies \eqref{eq:halfline_problem_0}. From \eqref{eq:decomp_usdroite_us1_us2} and from the equation verified by $u_g$, it follows that $g$ must verify $g|_{\R_+} = - \icplx \, \rho_\cutvec \, u^\droite_2$.
%   \begin{equation*}
%     0 = - \frac{d}{dx} \Big(\mu_\cutvec \,  \frac{d u^\droite_1}{dx}\Big) - \rho_\cutvec \, \omega^2 \, u^\droite_1 = \alpha\, ( \kappa \, \rho_\cutvec \, u^\droite_2 + g|_{\R_+}), \quad \textnormal{so that} \quad g|_{\R_+} = - \kappa \, \rho_\cutvec \, u^\droite_2.
%   \end{equation*}
  Aside from this condition, $g$ can be extended arbitrarily on $\R$. We extend $g$ by $0$ on $\R$, and still call the extension $g$. 
  %, by setting
  % \begin{equation*}
  %   f:= \alpha\, f^* \quad \textnormal{with} \quad f^*:= \left\{ 
  %     \begin{array}{c@{\quad}l}
  %       2\, \rho_\cutvec \, \omega^2 \, u^\droite_1 & \textnormal{in $\R_+$}
  %       \\[8pt]
  %       0 & \textnormal{in $\R \setminus \R_+$}.
  %     \end{array}
  %   \right.
  % \end{equation*}
  % Then \eqref{eq:decomp_usdroite_us1_us2} can be rewritten as $u^\droite = \alpha\, (u^\droite_1 + u^*|_{\R_+})$, where $u^\droite_1$ and $u^*:= u_{f^*}$ do not depend on $\alpha$. 
  From the Robin boundary condition in \eqref{eq:halfline_problem_0}, we deduce the equation for $\alpha$
  \begin{equation*}
    \alpha\, \big(R^\droite_{+} u^\droite_2 (0) + R^\droite_{+} u_g (0)\big)  = 1.
  \end{equation*}
  In order to conclude, the last point is to prove that $R^\droite_+ (u^\droite_2 + u_g) (0) \neq 0$, so that one can compute $\alpha$. We prove this by contradiction. If $R^\droite_+ (u^\droite_2 + u_g) (0) = 0$, then $\widetilde{u}^\droite:= u^\droite_2 + u_g|_{\R_+}$ belongs to $H^1(\R_+)$ and satisfies the ODE in \eqref{eq:halfline_problem_0} in $\R_+$, with the condition $R^\droite_{+} \widetilde{u}^\droite (0) = 0$. Therefore, according to Step 1, it follows that $\widetilde{u}^\droite = 0$, or equivalently $u^\droite_2 = -u_g$ in $\R_+$. By multiplying the equality $\overline{u}^\droite_2 = -\overline{u_g}$ by $\rho_\cutvec\, u^\droite_2$, and by integrating over $\R_+$ ($u^\droite_2 \in H^1(\R_+)$ and $u_g \in H^1(\R)$), we obtain
  \begin{align}
    \int_{\R_+} \rho_\cutvec\,  |u^\droite_2|^2 &= -\int_{\R_+} \rho_\cutvec\, u^\droite_2 \; \overline{u_g} = -\icplx\, \int_{\R} g\; \overline{u_g} \quad \textnormal{from the definition of $g$} \nonumber
    \\
    &= -\icplx\, \int_\R \Big(\mu_\cutvec \, \Big| \frac{d u_g }{dx}\Big|^2 - \rho_\cutvec\, \omega^2\, |u_g|^2\Big) \quad \textnormal{from the equation satisfied by $u_g$.} \label{eq:incl_spec_dem}
  \end{align}
 This implies that $u^\droite_2 = 0$ in $\R_+$, which contradicts the condition $u^\droite_2 (0) = 1$.
 
{
    \paragraph*{Step 2.}%
    We now show that $\omega^2 \not\in \sigma(\mathcal{A}^\droite_p)$, namely, that the problem
    \begin{equation} \label{eq:robin_half_guide_with_source_term}
	    \left|
	    \begin{array}{r@{\ =\ }l} 
			\displaystyle - \Dt{} \big( \mu_p \, \Dt{} U^\droite \big) - \rho_p \, (\omega^2+ \icplx \veps) \, U^\droite & F \quad \textnormal{in }\,\Omegaper^\droite, 
			\ret
			\displaystyle \mathcal{R}^\droite_+ \,  U^\droite_{\veps}\Big|_{\Sigmaper^{\droite, 0}} & 0,
			\ret
			\multicolumn{2}{c}{\textnormal{$U^\droite_{\veps}$ is periodic w.r.t. $\yvi_1$.}}
    	\end{array}
    	\right.
    \end{equation}
    admits a unique solution $U^\droite \in  H^1_\cutvec(\Omega_\periodic)$ for any $F \in L^2(\Omega_\periodic)$. To prove the existence, we construct a solution of \eqref{eq:robin_half_guide_with_source_term} using Step $1$. More precisely, using the change of variables \eqref{change_of_variable}, and identifying $F$ with its periodic extension $F \in L^2_{\textit{loc}} (\Omega^\droite)$ with respect to $y_1$, we introduce the function
    \begin{equation*}
        \aeforall (s, x) \in \R \times \R_+, \quad f_s (x):= F(s + x \cuti_1, x \cuti_2).
    \end{equation*}
    Note that $f_s$ is well-defined, and belongs to $L^2(\R_+)$ for almost every $s \in \R$. By Step $1$, there exists $u^\droite_s \in H^1(\R_+)$ such that $(\mathcal{A}^\droite_{s, \cutvec} - \omega^2)\, u^\droite_s = f_s$, see \eqref{eq:halfline_problem_f}. Therefore, the function
    \begin{equation*}
        (y_1, y_2) \mapsto u^\droite_{y_1 - y_2 \cuti_1 / \cuti_2} (y_2 / \cuti_2),
    \end{equation*}
    belongs to $H^1_{\cutvec,per}(\Omega^\droite)$, and satisfies \eqref{eq:robin_half_guide_with_source_term}. It remains to prove uniqueness. Let $U^\droite \in H^1_{\cutvec,per}(\Omega^\droite)$ satisfy \eqref{eq:robin_half_guide_with_source_term} with $F = 0$. Then the function $x \mapsto u^\droite_s (x):= U^\droite (s + x \cuti_1, x \cuti_2) \in H^1(\R_+)$ satisfies $(\mathcal{A}^\droite_{s, \cutvec} - \omega^2)\, u^\droite_s = 0$ for almost any $s \in \R$. Since $\omega^2 \not\in \sigma(\mathcal{A}^\droite_{\cutvec})$, we deduce that $u^\droite_s = 0$ for almost any $s \in \R$, or equivalently, that $U^\droite = 0$. This completes the proof.
}
\end{dem}

\section{Limiting absorption principle}\label{sec:LAP}
We now study the limit when $\veps$ goes to $0$ of the solution $u_\veps$ of \eqref{eq:whole_line_problem_eps}. If the limit exists in a certain sense and if it satisfies \eqref{eq:whole_line_problem}, we say that the limiting absorption principle holds and this limit is the physical solution. Since $u_\veps$ is constructed and characterized by the method described in Sections \ref{sec:RtR_eps} and \ref{sec:Riccati_eps}, our approach is to pass to the limit in each step of the method. When studying objects associated with the half-line problems set on $I^\idomlap$, we shall assume for simplicity that $\idomlap = \droite$.

\vspace{1\baselineskip} 
\noindent We suppose in this section that the assumption \ref{assumption_Dioph} is satisfied.

\subsection{Convergence of the cell solutions and the associated RtR operators}\label{sec:conv_cell_sol_RtR_ops}
The first step is to characterize the limits of the cell solutions $(E^{\droite, 0}_\veps(\varphi), E^{\droite, 1}_\veps(\varphi))$ defined by \eqref{cellpb0} and \eqref{cellpb1} for $\varphi \in L^2_\textnormal{per}(\R)$, and of the associated RtR operators ${\cal T}^{\droite, \ell k}_\veps$ defined in \eqref{defRtR} for $\ell,k\in\{0,1\}$. 
 
\vspspe
Let us introduce the limit cell problems for $\varphi\in L^2_{\textit{per}}(\R)$: \textit{Find $E^{\droite, 0} (\varphi) \in H^1_{\cutvec} (\cellper^{\droite, 0})$ such that}
\begin{equation} \label{cellpb0_sseps}
	\left|
	 \begin{array}{r@{\ =\ }l} 
			\displaystyle - \Dt{} \big( \mu_p \, \Dt{} E^{\droite, 0} (\varphi) \big) - \rho_p \, \omega^2 \, E^{\droite, 0} (\varphi) & 0 \quad \textnormal{in }\,\cellper^{\droite, 0}, 
			\ret
			\multicolumn{2}{c}{\displaystyle \mathcal{R}^\droite_+ \, E^{\droite, 0} (\varphi)|_{\Sigmaper^{\droite, 0}} = \varphi \quad \textnormal{and} \quad \mathcal{R}^\droite_- \, E^{\droite, 0} (\varphi)|_{\Sigmaper^{\droite, 1}} = 0,}
			\ret
			\multicolumn{2}{c}{\textnormal{$E^{\droite, 0} (\varphi)$ is periodic w.r.t. $\yvi_1$;}}
	\end{array}
	\right.
\end{equation}
and: \textit{Find $E^{\droite, 1} (\varphi) \in H^1_{\cutvec} (\cellper^{\droite, 0})$ such that}
\begin{equation} \label{cellpb1_sseps}
	\left|
	 \begin{array}{r@{\ =\ }l} 
			\displaystyle - \Dt{} \big( \mu_p \, \Dt{} E^{\droite, 1} (\varphi) \big) - \rho_p \, \omega^2 \, E^{\droite, 1} (\varphi) & 0 \quad \textnormal{in }\,\cellper^{\droite, 0}, 
			\ret
			\multicolumn{2}{c}{\displaystyle \mathcal{R}^\droite_+ \, E^{\droite, 1} (\varphi)|_{\Sigmaper^{\droite, 0}} = 0 \quad \textnormal{and} \quad \mathcal{R}^\droite_- \, E^{\droite, 1} (\varphi)|_{\Sigmaper^{\droite, 1}} = \varphi,}
			\ret
			\multicolumn{2}{c}{\textnormal{$E^{\droite, 1} (\varphi)$ is periodic w.r.t. $\yvi_1$.}}
	\end{array}
	\right.
\end{equation}
In addition, consider the 1D cell problems defined for any $s \in \R$ by
\begin{equation} \label{segment_problems0_sseps} 
	\left|
	\begin{array}{lll}
	\displaystyle - \, \frac{d}{dx} \Big(\mu_{s, \cutvec} \, \frac{d e_{s}^{\droite, 0}}{dx}\Big)- \rho_{s, \cutvec} \, \omega^2 \, e_{s}^{\droite, 0} = 0 \quad \mbox{in}\ \ (0, 1/\cuti_2), 
	\ret
	R^\droite_{+,s} e_{s}^{\droite, 0}(0) = 1 \quad \textnormal{and} \quad R^\droite_{-,s} e_{s}^{\droite, 0}(1/\cuti_2) = 0,
	\end{array}	
	\right. 
\end{equation}
and
\begin{equation} \label{segment_problems1_sseps}
	\left|
		\begin{array}{lll}
			\displaystyle 	- \, \frac{d}{dx} \Big(\mu_{s, \cutvec} \, \frac{d e_{s}^{\droite, 1}}{dx}\Big)- \rho_{s, \cutvec} \, \omega^2 \, e_{s}^{\droite, 1} = 0 \quad \mbox{in}\ \ (0, 1/\cuti_2), 
			\ret
			R^\droite_{+,s} e_{s}^{\droite, 1}(0) = 0 \quad \textnormal{and} \quad R^\droite_{-,s} e_{s}^{\droite, 1}(1/\cuti_2) = 1.
		\end{array}	
	\right.
\end{equation} 
One advantage of the 1D problems \eqref{segment_problems0_sseps} and \eqref{segment_problems1_sseps} is that Fredholm alternative can be applied in $H^1(0, 1/\cuti_2)$. Moreover, thanks to the RtR boundary conditions, the uniqueness of the solutions can be proved easily. As a consequence, Problem \eqref{segment_problems0_sseps} and \eqref{segment_problems1_sseps} are well-posed. We can then deduce the well-posedness of the 2D cell problems \eqref{cellpb0_sseps} and \eqref{cellpb1_sseps} as well as the link between the 2D solutions $(E^{\droite, 0}(\varphi), E^{\droite, 1}(\varphi))$ and the 1D ones $(e_{s}^{\droite, 0}, e_{s}^{\droite, 1})$.
% \sfnote{Cela signifie que nous ne sommes pas capable de montrer la CV des EF pour les problemes de cellule 2D mais seulement pour les problemes de cellule 1D. Cela justifie encore plus l utilisation de la methode quasi-1D. Est ce qu'on le dit quelque part?}
%
%
%
\begin{lem}\label{lem_equivalence_ssabs}
	The problem \eqref{cellpb0_sseps} (\resp \eqref{cellpb1_sseps}) is equivalent to the family of problems \eqref{segment_problems0_sseps} (\resp \eqref{segment_problems1_sseps}) parametrized by $s \in \R$ in the following sense
	\begin{enumerate}[label=(\roman*)]
		\item If $E^{\droite, 0}$, $E^{\droite, 1}$ are the respective solutions of \eqref{cellpb0_sseps} and \eqref{cellpb1_sseps} for $\varphi = 1$, and are identified with their periodic extensions with respect to $\yvi_1$, then by setting for $\ell=0,1$ and $s\in\R$
			\begin{equation} \label{defcellsol_ssabs}
		e_s^{\droite, \ell}(x) = E^{\droite, \ell} (s + \cuti_1\, \xvi, \cuti_2\, \xvi), \quad x \in (0, 1/\theta_2),
		\end{equation}
		the functions $e_s^{\droite,0}$ and $e_s^{\droite,1}$ are the respective solutions of \eqref{segment_problems0_sseps} and \eqref{segment_problems1_sseps}.
		\item If $e_s^{\droite,0}$ and $e_s^{\droite,1}$ are the respective solutions of \eqref{segment_problems0_sseps} and \eqref{segment_problems1_sseps} for $s \in \R$, then by defining for all $\varphi\in L^2_{\textit{per}}(\R) $, the functions $E^{\droite,0}(\varphi)$ and $E^{\droite,1}(\varphi)$ as follows
		\begin{equation} \label{defcellsol_ssabs2}
			\left\{
			\begin{array}{r@{\ \coloneqq \ }l@{\quad}l}
				E^{\droite,0}(s + \cuti_1\, \xvi, \cuti_2\, \xvi) & \varphi(s)\, e^{\droite,0}_s(x), & x \in (0, 1/\cuti_2),\quad s \in [0,1),
				\\[8pt]
				E^{\droite,1}(s + \cuti_1\, \xvi, \cuti_2\, \xvi) & \varphi(s + \delta)\, e^{\droite,1}_s(x), & x \in (0, 1/\cuti_2),\quad s \in [0,1), 
			\end{array}
			\right.
		\end{equation}
		the functions $E^{\droite,0}(\varphi)$, $E^{\droite,1}(\varphi)$ are the respective solutions of \eqref{cellpb0_sseps} and \eqref{cellpb1_sseps}.
	\end{enumerate}
	In particular,
	\[
		\textnormal{\eqref{cellpb0_sseps} (\resp \eqref{cellpb1_sseps}) is well-posed $\spforall \varphi\in L^2_{\textit{per}}(\R)$} \ \ \Longleftrightarrow\ \ \textnormal{\eqref{segment_problems0_sseps} (\resp \eqref{segment_problems1_sseps}) is well-posed $\spforall s\in\R$}.
		\]
\end{lem}
\noindent To $s \mapsto (e_{s}^{\droite, 0}, e_{s}^{\droite, 1})$ are associated the local RtR coefficients (defined similarly to \eqref{defRtRcoeff})
\begin{equation}\label{defRtRcoeff_0}
	\left\{ \begin{array}{r@{\ =\ }l@{ \qquad }r@{\ =\ }l}	
		\boldsymbol{t}^{\droite, 00}(s) & R_{-,s}^\droite e_{s}^{\droite, 0}(0), &%
		\boldsymbol{t}^{\droite, 01}(s) & R_{+,s}^\droite e_{s}^{\droite, 0}(1/\theta_2),
		\ret
		\boldsymbol{t}^{\droite, 10}(s) & R_{-,s}^\droite e_{s}^{\droite, 1}(0), &%
		\boldsymbol{t}^{\droite, 11}(s) & R_{+,s}^\droite e_{s}^{\droite, 1}(1/\theta_2).
	\end{array} \right. 
\end{equation}
When $\veps \to 0$, we can show convergence estimates for $(e_{s, \veps}^{\droite, \ell}, \boldsymbol{t}^{\droite, \ell k}_\veps)$ towards $(e_{s}^{\droite, \ell}, \boldsymbol{t}^{\droite, \ell k})$.

\begin{prop}\label{prop:lap_cell1d}
	There exists a constant $C > 0$ (independent of $s$) such that for any $s \in \R$,
	\begin{equation}\label{eq:lap_cell1d}
		\left\{
			\begin{array}{r@{\ }l@{\qquad}l}
				\|e^{\droite,0}_s-e^{\droite,0}_{s,\veps}\|_{H^1}+ \|e^{\droite,1}_s-e^{\droite,1}_{s,\veps}\|_{H^1} &\leq C\veps, & (i)
				\ret
				\spforall \ell, k \in \{0, 1\}^2, \quad \|\boldsymbol{t}^{\droite, \ell k}_\veps - \boldsymbol{t}^{\droite, \ell k}\|_\infty &\leq C\veps. & (ii)
			\end{array}
		\right.
	\end{equation}
\end{prop}
	
\begin{dem}
	One begins by writing the problem satisfied by the difference $v \coloneqq  e^{\droite, 0}_s-e^{\droite, 0}_{s,\veps}$:
	\begin{equation}\label{eq:e_veps_minus_e_0}
		\left|
			\begin{array}{lll}
			\displaystyle - \, \frac{d}{dx} \Big(\mu_{s, \cutvec} \, \frac{d v}{dx}\Big)- \rho_{s, \cutvec} \, \omega^2 \, v = - \rho_{s, \cutvec}\, \veps\, e^{\droite, 0}_{s,\veps} \quad \mbox{in}\ \ (0, 1/\cuti_2), 
			\ret
			R^\droite_{+,s} v(0) = R^\droite_{-,s} v(1/\cuti_2) = 0.
			\end{array}	
		\right. 
	\end{equation}
	From the well-posedness of this problem, we deduce the first line $(i)$ of \eqref{eq:lap_cell1d}, with a constant $C(s)$ which depend a priori on $s$. To show that the constant can be chosen independently of $s$, we assume by contradiction that there exists a point $s^*$ and a sequence $s_n$ such that $s_n \to s^*$ and $\|v_n\|_{H^1} \to +\infty$ as $n \to +\infty$ with $v_n \coloneqq  e^{\droite, 0}_{s_n}-e^{\droite, 0}_{s_n,\veps}$. Then,
	\begin{equation*}
		\displaystyle
		\widehat{v}_n \coloneqq  \frac{v_n}{\|v_n\|_{H^1}} \in H^1(0, 1/\cuti_2) \quad \textnormal{with} \quad \| \widehat{v}_n \|_{H^1} = 1.
	\end{equation*}
	On the other hand, $\widehat{v}_n$ satisfies Problem \eqref{eq:e_veps_minus_e_0} with a source term that tends to $0$ as $n \to +\infty$ {since $s\mapsto e^{\droite, 0}_{s,\veps}$ is bounded in $L^2$ independently of $s$.} Therefore, the uniform continuity of the maps $s \mapsto (\mu_{s, \cutvec}, \rho_{s, \cutvec})$ implies that $\widehat{v}_n \to 0$, which contradicts the equality $\| \widehat{v}_n \|_{H^1} = 1$. By extending these arguments to $e^{\droite, 1}_s-e^{\droite, 1}_{s,\veps}$, we obtain \eqref{eq:lap_cell1d}--$(i)$.

	\vspspe
	The second line $(ii)$ of \eqref{eq:lap_cell1d} is a consequence of the first one $(i)$ combined with the continuity of the normal trace (and the boundedness of $s \mapsto \mu_{s, \cutvec}$). 
\end{dem}	
\noindent
By analogy with \eqref{defRtR}, let us now introduce  the limit RtR operators $\mathcal{T}^{\droite, \ell k}$ defined for  $\varphi \in L^2_{\textit{per}} (\R)$ by
	\begin{equation}\label{eq:def_RtR_loc_lim}
		\begin{array}{r@{\ \coloneqq \ }l@{\qquad \textnormal{and} \qquad}r@{\ \coloneqq \ }l}	
			\mathcal{T}^{\droite, 00}\, \varphi & \mathcal{R}_-^\droite E^{\droite, 0} (\varphi)|_{\Sigma^{\droite, 0}} &  	\mathcal{T}^{\droite, 01}\, \varphi & \mathcal{R}_+^\droite E^{\droite, 0} (\varphi)|_{\Sigma^{\droite, 1}},
			\ret
			\mathcal{T}^{\droite, 10}\, \varphi & \mathcal{R}_-^\droite E^{\droite, 1} (\varphi)|_{\Sigma^{\droite, 0}} &  	\mathcal{T}^{\droite, 11}\, \varphi & \mathcal{R}_+^\droite E^{\droite, 1} (\varphi)|_{\Sigma^{\droite, 1}}. 	
		\end{array} 
	\end{equation}
	Similarly to \eqref{RtRweighted} for $\veps > 0$, the operators $\mathcal{T}^{\droite, \ell k}$ are weighted shift operators: for $\varphi \in L^2_{\textit{per}} (\R)$
	\begin{equation}\label{RtRweighted_0}
		 \left\{
		\begin{array}{r@{\ =\ }l@{ \quad }r@{\ =\ }l}	
			\mathcal{T}^{\droite, 00}\, \varphi(s) &	\boldsymbol{t}^{\droite, 00}(s) \, \varphi(s), & \quad \mathcal{T}^{\droite, 01}\, \varphi(s) &	\boldsymbol{t}^{\droite, 01}(s - \delta)\, \varphi(s - \delta), 
			\rett
			\mathcal{T}^{\droite, 10}\, \varphi(s) & \boldsymbol{t}^{\droite, 10}(s) \, \varphi(s + \delta), &  \quad \mathcal{T}^{\droite, 11}\, \varphi(s) &	\boldsymbol{t}^{\droite, 11}(s - \delta) \, \varphi(s).
		\end{array} 
		\right. 
	\end{equation}
We finally deduce, from Proposition \ref{prop:lap_cell1d} and the fibered structure mentioned above, the convergence results for the solutions of the local cell problems \eqref{cellpb0_sseps} and \eqref{cellpb1_sseps}, and the associated RtR operators, 
\begin{prop}\label{prop:lap_cell2d}
	There exists a constant $C > 0$ such that
	\begin{equation*}
		\left\{
			\begin{array}{r@{\ }l@{\qquad}l}
				\|E^{\droite,0}-E^{\droite,0}_{\veps}\|_{H^1(\cellper^{\droite, 0})}+ \|E^{\droite,1}-E^{\droite,1}_{\veps}\|_{H^1_\cutvec(\cellper^{\droite, 0})} &\leq C\veps, & (i)
				\ret
				\spforall \ell, k \in \{0, 1\}^2, \quad \|\mathcal{T}^{\droite, \ell k}_\veps - \mathcal{T}^{\droite, \ell k}\|_{\mathscr{L}(L^2(0, 1))} &\leq C\veps. & (ii)
			\end{array}
		\right.
	\end{equation*}
\end{prop}
% \begin{dem}
% 	The first line $(i)$ follows by combining the fibered expression \eqref{defcellsol_ssabs2} with the first line $(i)$ in \eqref{eq:lap_cell1d}. To obtain the second line $(ii)$, we exploit the observation that 
% 	so that the second line $(ii)$ in the estimate \eqref{eq:lap_cell1d} allows to conclude. 
% \end{dem}
\noindent An immediate consequence of Proposition \ref{prop:lap_cell2d} is that 
\[
	\spexists C>0,\quad\|\mathscr{M}^\droite-\mathscr{M}^\droite_\veps\|_{\mathscr{L}(L^2_\textit{per}(\R))^2}+\|\mathscr{N}^\droite-\mathscr{N}^\droite_\veps\|_{\mathscr{L}(L^2_\textit{per}(\R))^2}\leq C\,\veps,
\]
where $\mathscr{M}^\droite_\veps$ and $\mathscr{N}^\droite_\veps$ are defined for $\veps>0$ in \eqref{eq:Riccati_op} and
\begin{equation}\label{eq:Riccati_op_ssabs}
	\mathscr{M}^\droite \coloneqq  \begin{pmatrix} {\cal T}^{\droite, 01} & {\cal T}^{\droite, 11} \\[2pt] 0 & I \end{pmatrix} \quad \textnormal{and} \quad \mathscr{N}^\droite = \begin{pmatrix} I & 0 \\[2pt] {\cal T}^{\droite, 00} & {\cal T}^{\droite, 10} \end{pmatrix}.
\end{equation}

\subsection{Link between the spectrum of \texorpdfstring{$\mathcal{A}_\cutvec$}{Atheta} and the Riccati spectrum}
By analogy with the Riccati point spectrum $\sigma^\droite_{R, p, \veps}$ defined by \eqref{eq:Ric_spec}, it is natural to introduce the Riccati point spectrum in the absence of absorption:
\begin{equation} \label{eq:Ric_spec_sseps}
	\sigma^\droite_{R, p} \coloneqq  \big\{ \lambda \in \C \ /\ 0 \in \sigma_p ( \mathscr{M}^\droite - \lambda \, \mathscr{N}^\droite) \}. 
\end{equation}
For the main theorem of this section, namely Proposition \ref{prop:link_spec_Riccati_spec_diff_op}, it is also useful to introduce the Riccati spectrum defined by
\begin{equation} \label{eq:Ric_spec_cont_sseps}
	\sigma^\droite_{R} \coloneqq  \big\{ \lambda \in \C \ /\ 0 \in \sigma ( \mathscr{M}^\droite - \lambda \, \mathscr{N}^\droite) \}. 
\end{equation}

\begin{prop}\label{prop:link_spec_Riccati_spec_diff_op}
	{The spectra defined by (\ref{eq:Ric_spec_sseps}, \ref{eq:Ric_spec_cont_sseps}) are related to the one of $\mathcal{A}_\cutvec$ as follows.
	\begin{equation}\label{eq:link_spec_Riccati_spec_diff_op}
		\begin{array}{r@{\;\; \Rightarrow\;\; }l@{\;\;}l}
			\spexists \lambda \in \C, \; |\lambda| = 1, \; \lambda \in \sigma^\droite_{R, p} & \omega^2 \in \sigma(\mathcal{A}_\cutvec) & (i)
			\rett
			\omega\in\R,\;\omega^2 \notin \sigma(\mathcal{A}_\cutvec)  &  \spexists r_0,\;|r_0|<1,\; \overline{\sigma_{R, p}^\droite} = \mathcal{C}\big(0,|r_0|\big)\vts \cup\vts  \mathcal{C}\big(0,1/{|r_0|}\big), & (ii)
			\rett
			\omega^2 \in \sigma(\mathcal{A}_\cutvec) & \spexists \lambda \in \C, \; |\lambda| = 1, \quad \lambda \in \sigma^\droite_{R}.  & (iii)
		\end{array}
	\end{equation}}
\end{prop}

\begin{rmk}
	Note that in \eqref{eq:link_spec_Riccati_spec_diff_op}--$(i)$, $\lambda \in \sigma^\droite_{R, p}\subset \sigma^\droite_R$, whereas in \eqref{eq:link_spec_Riccati_spec_diff_op}--$(iii)$, we only have $\lambda \in \sigma^\droite_{R}$, and nothing ensures that $\lambda \in \sigma^\droite_{R, p}$. 
	% \sfcomment{It is worth recalling that as soon as a discrete version of the cell problems are solved and the associated RtR operators computed, the Riccati spectrum is reduced to its point part.}
\end{rmk}

\begin{dem}
	% We are going to prove \eqref{eq:link_spec_Riccati_spec_diff_op} using the operator $\mathcal{A}_p (\fbvar)$ defined by \eqref{eq:op_per_k} instead of $\mathcal{A}_\cutvec$. This can be done since from Propositions \ref{thm:spectrunqp_s} and \ref{thm:spectrunqp_k}, we have $\sigma(\mathcal{A}_p (\fbvar)) = \sigma(\mathcal{A}_\cutvec)$.
\textbf{Point $(i)$.} Let $\lambda \in \C$ with $|\lambda|=1$ such that $\lambda\in \sigma^\droite_{R, p}$. Then there exists $\fbvar \in [-\pi,\pi)$ such that $\lambda = \euler^{\icplx\vts \fbvar}$ and $(\varphi_\fbvar, \psi_\fbvar) \neq 0$ such that
  \begin{equation}\label{eq:ric_eigk}
     \mathscr{M}^\droite\begin{pmatrix}\varphi_\fbvar\\\psi_\fbvar\end{pmatrix} = \euler^{\icplx\vts \fbvar} \, \mathscr{N}^\droite \begin{pmatrix}\varphi_\fbvar\\\psi_\fbvar\end{pmatrix}.
  \end{equation}
	It can be easily checked, using the definitions (\ref{cellpb0_sseps}, \ref{cellpb1_sseps}) of the cell solutions $E^{\droite, \ell}$ that the function defined by
  \[
    U_\fbvar\coloneqq E^{\droite,0}(\varphi_\fbvar)+E^{\droite,1}(\psi_\fbvar)\ \in H^1_{\cutvec} (\cellper^{\droite, 0})
	\]
  satisfies $- \Dt{} \big( \mu_p \, \Dt{} U_\fbvar  \big) - \rho_p \, \omega^2 \, U_\fbvar = 0$ in $\cellper^{\droite, 0}$ and is $1$--periodic with respect to $\yvi_1$, while the definition 	\eqref{eq:def_RtR_loc_lim} of the operators $\mathcal{T}^{\droite, \ell k}$ and \eqref{eq:ric_eigk} ensure that
	\[
		{\cal R}_+^\droite U_\fbvar |_{\Sigma^{\droite, 1}}=\euler^{\icplx\vts \fbvar}\,{\cal R}_+^\droite U_\fbvar |_{\Sigma^{\droite, 0}} \quad\textnormal{and}\quad {\cal R}_-^\droite U_\fbvar |_{\Sigma^{\droite, 1}}=\euler^{\icplx\vts \fbvar}\,{\cal R}_-^\droite U_\fbvar |_{\Sigma^{\droite, 0}}.
	\]
	{By defining $V_\fbvar \coloneqq  \euler^{-\icplx\vts \fbvar\, y_2}\, U_\fbvar\neq 0$, we have that $V_\fbvar$ is $1$--periodic with respect to $\yvi_1$, satisfies $- (\Dt{} + \icplx\vts \fbvar\, \cuti_2)\, \mu_p \, (\Dt{}  + \icplx\vts \fbvar\, \cuti_2)\, V_\fbvar  - \rho_p \, \omega^2 \, V_\fbvar = 0$ in $\cellper^{\droite, 0}$,   and
	\begin{equation*}
		{\cal R}_+^\droite V_\fbvar |_{\Sigma^{\droite, 1}} = {\cal R}_+^\droite V_\fbvar |_{\Sigma^{\droite, 0}} \quad\textnormal{and}\quad {\cal R}_-^\droite V_\fbvar |_{\Sigma^{\droite, 1}} = {\cal R}_-^\droite V_\fbvar |_{\Sigma^{\droite, 0}},
	\end{equation*}
	so that $V_\fbvar$ is also $1$--periodic with respect to $\yvi_2$. }Consequently, by definition \eqref{eq:op_per_k} of the operator $\mathcal{A}(\fbvar)$, we deduce the following
	\begin{equation*}
		\displaystyle
		V_\fbvar \in D(\mathcal{A}_p(\fbvar)) \setminus \{0\} \quad \textnormal{and} \quad \mathcal{A}_p(\fbvar)\, V_\fbvar = \omega^2\, V_\fbvar,
	\end{equation*}
	which implies that $\omega^2\in \sigma_p(\mathcal{A}_p (\fbvar)) \subset \sigma(\mathcal{A}_p (\fbvar)) = \sigma(\mathcal{A}_\cutvec)$.

	{\paragraph*{Point $(ii)$.} If $\omega\in\R$, $\omega^2\notin \sigma({\cal A}_\cutvec)$, then according to Propositions \ref{thm:spectrunqp_s}  and \ref{prop:spectrum_halves}, $\omega^2$ does not belong to the spectrum of the 2D periodic differential operator $\mathcal{A}^\droite_p$ defined by \eqref{eq:op_per_j}. Therefore, all the problems introduced in Section \ref{sec:RtR_eps} remain well-posed for $\veps = 0$, so that the approach and the results presented in Sections \ref{sec:RtR_eps} and \ref{sec:Riccati_eps} are still valid for $\veps = 0$. In particular Theorem \ref{thm_spectre} and Proposition \ref{thm_spectreaugmented} hold, which yields Point $(ii)$.}

	\paragraph*{Point $(iii)$.} Let us show the converse of the point $(iii)$ in \eqref{eq:link_spec_Riccati_spec_diff_op}, namely
  \[
  	\spforall \lambda \in \C, \quad |\lambda|=1,\quad \mathscr{M}^\droite - \lambda \, \mathscr{N}^\droite \mbox{ is invertible} \quad \Longrightarrow \quad \omega^2\notin\sigma(\mathcal{A}_\cutvec).
	\]
	Thanks to the equality of spectra $\sigma(\mathcal{A}_p (\fbvar)) = \sigma(\mathcal{A}_\cutvec)$ which follows from Propositions \ref{thm:spectrunqp_s} and \ref{thm:spectrunqp_k}, it is sufficient to find $\fbvar \in [-\pi, \pi)$ such that $\mathcal{A}_p(\fbvar) - \omega^2$ is invertible, that is,
	\begin{equation}\label{eq:V_xi}
		\displaystyle
		\spforall f \in  L^2_{\textit{per}}(\R^2), \quad \exists!\, V_\fbvar \in D\big(\mathcal{A}_p(\fbvar)\big) \quad \big/ \quad \left[\mathcal{A}_p(\fbvar)\, - \omega^2\right]\, V_\fbvar = f.
	\end{equation}
	Given $\lambda$ such that $|\lambda| = 1$, we choose $\fbvar\in [-\pi, \pi[$ such that $\lambda = \euler^{\icplx\vts \fbvar}$. Then for the existence of $V_\fbvar$, setting $g \coloneqq  \rho_p\, f$, we start from the unique solution $U_g \in H^1_{\cutvec} (\cellper^{\droite, 0})$ of
	\begin{equation*}
		\left|
		 \begin{array}{r@{\ =\ }l} 
				\displaystyle - \Dt{} \big( \mu_p \, \Dt{} U_g \big) - \rho_p \, \omega^2 \, U_g & g \quad \textnormal{in }\,\cellper^{\droite, 0}, 
				\ret
				\multicolumn{2}{c}{\displaystyle \mathcal{R}^\droite_+ \, U_g|_{\Sigmaper^{\droite, 0}} = 0 \quad \textnormal{and} \quad \mathcal{R}^\droite_- \, U_g|_{\Sigmaper^{\droite, 1}} = 0,}
				\ret
				\multicolumn{2}{c}{\textnormal{$U_g$ is periodic w.r.t. $\yvi_1$.}}
		\end{array}
		\right.
	\end{equation*}
	Note that this problem is well-posed for the same reasons as the local cell problems (\ref{cellpb0_sseps}, \ref{cellpb1_sseps}). Then for any $(\varphi, \psi) \in [L^2_{\textit{per}}(\R^2)]^2$, $U[\varphi,\psi] \coloneqq  E^{\droite, 0} (\varphi) + E^{\droite, 1} (\psi) + U_g \in H^1_{\cutvec} (\cellper^{\droite, 0})$ clearly satisfies $- \Dt{} \big( \mu_p \, \Dt{} U  \big) - \rho_p \, \omega^2 \, U= g$ in $\cellper^{\droite, 0}$. To construct $V_\fbvar$ that solves \eqref{eq:V_xi}, it suffices to find $(\varphi_\fbvar, \psi_\fbvar) \in L^2_{\textit{per}}(\R)^2$ so that $U_\fbvar\coloneqq U[\varphi_\fbvar, \psi_\fbvar]$ satisfies 
	\begin{equation}\label{eq:Uxi_QP}
		\displaystyle
		{\cal R}_+^\droite U_\fbvar |_{\Sigma^{\droite, 1}}=\euler^{\icplx\vts \fbvar}\,{\cal R}_+^\droite U_\fbvar |_{\Sigma^{\droite, 0}} \quad\textnormal{and}\quad {\cal R}_-^\droite U_\fbvar |_{\Sigma^{\droite, 1}}=\euler^{\icplx\vts \fbvar}\,{\cal R}_-^\droite U_\fbvar |_{\Sigma^{\droite, 0}}.
	\end{equation} 
	Indeed, it can be seen easily that the function defined by $V_\fbvar \coloneqq  \euler^{-\icplx\vts \fbvar\, y_2}\, U_\fbvar$ in $\mathcal{C}^{\droite, 0}$ belongs to $D\big(\mathcal{A}_p(\fbvar)\big)$ and satisfies $\mathcal{A}_p(\fbvar)\, V_\fbvar - \omega^2\, V_\fbvar = f.$ To find such $(\varphi_\fbvar, \psi_\fbvar)$, note that, by definition (\ref{eq:Riccati_op_ssabs}, \ref{eq:def_RtR_loc_lim}) of $(\mathscr{M}^\droite, \mathscr{N}^\droite)$ and of the local RtR operators $\mathcal{T}^{\droite, \ell k}$, for any $(\varphi, \psi) \in L^2_{\textit{per}}(\R^2)^2$, 
	\begin{equation}\label{eq:link_riccati_spec_diff_op_dem_2}
		\big(\mathscr{M}^\droite - \euler^{\icplx\vts \fbvar} \, \mathscr{N}^\droite\big)% 
		\begin{pmatrix}\varphi \\[8pt] \psi\end{pmatrix} %
		= %
		\begin{pmatrix}
			\mathcal{R}_+^\droite E^{\droite, 0} (\varphi) |_{\Sigma^{\droite, 1}} + \mathcal{R}_+^\droite E^{\droite, 1} (\psi) |_{\Sigma^{\droite, 1}} - \euler^{\icplx\vts \fbvar}\, \psi
			\\[8pt]		
			\psi - \euler^{\icplx\vts \fbvar}\, \mathcal{R}_-^\droite E^{\droite, 0} (\varphi) |_{\Sigma^{\droite, 0}} - \euler^{\icplx\vts \fbvar}\, \mathcal{R}_-^\droite E^{\droite, 1} (\psi) |_{\Sigma^{\droite, 0}}
		\end{pmatrix}.
	\end{equation}
		Then \eqref{eq:Uxi_QP} is equivalent to 
		\begin{equation}\label{eq:link_riccati_spec_diff_op_dem_1}
		\big(\mathscr{M}^\droite - \euler^{\icplx\vts \fbvar} \, \mathscr{N}^\droite\big)\begin{pmatrix}\varphi_\fbvar\\[8pt]\psi_\fbvar\end{pmatrix} = \begin{pmatrix}-{\cal R}_+^\droite U_g|_{\Sigma^{\droite, 1}}\\[8pt] \euler^{\icplx\vts \fbvar}{\cal R}_-^\droite U_g|_{\Sigma^{\droite, 0}}\end{pmatrix},
	\end{equation}
	which has a unique solution since $\mathscr{M}^\droite - \lambda \, \mathscr{N}^\droite$ is invertible for $\lambda = \euler^{\icplx\vts \fbvar}$.

	 \vspspe
	 Finally, to prove that $V_\fbvar$ is unique, assume that $\mathcal{A}_p(\fbvar)\, V_\fbvar - \omega^2\, V_\fbvar = 0$. Then by defining $\varphi_\fbvar \coloneqq  \mathcal{R}^\droite_+ V_\fbvar|_{\Sigma^{\droite, 0}}$ and $\psi_\fbvar \coloneqq  \mathcal{R}^\droite_- V_\fbvar|_{\Sigma^{\droite, 1}}$, we obtain $(\mathscr{M}^\droite - \lambda\, \mathscr{N}^\droite)\, \transp{(\varphi_\fbvar, \psi_\fbvar)} = 0$, which, thanks to the injectivity of $\mathscr{M}^\droite - \lambda\, \mathscr{N}^\droite$, leads to $\varphi_\fbvar = \psi_\fbvar = 0$. Hence, for almost any $s \in (0, 1)$, the restriction $v_s(\xvi) \coloneqq  V_\fbvar (\cuti_1\, \xvi + s, \cuti_2\, \xvi)$ satisfies the same ODE as $e^{\droite, 0}_s$ and $e^{\droite, 1}_s$, with $R^\droite_{-, s} v(0) = R^\droite_{+, s} v_s (1/\cuti_2) = 0$. From Cauchy uniqueness theorem, we deduce $v_s = 0$ for almost any $s$, which implies that $V_\fbvar = 0$.	
\end{dem}

\begin{rmk} 
	{From the proof of Proposition \ref{prop:link_spec_Riccati_spec_diff_op}, we can deduce easily that
\begin{equation}\label{eq:link_spec_Riccati_spec_diff_op_bis}
		\begin{array}{r@{\quad \Longrightarrow\quad }l@{\qquad}l}
			\spexists \lambda \in \C, \quad |\lambda| = 1, \quad \lambda \in \sigma^\droite_{R, p} & \spforall k\,\in\R,\quad \lambda\, \euler^{2\icplx\pi k \delta} \in \sigma^\droite_{R, p}.
		\end{array}
	\end{equation}
	Indeed, in the proof of \eqref{eq:link_spec_Riccati_spec_diff_op}-(i), we have constructed when $\lambda=\euler^{\icplx\fbvar}\in \sigma^\droite_{R, p}$, $\fbvar\in[-\pi,\pi[$ an eigenvector $V_\fbvar$ of  $\mathcal{A}_p (\fbvar)$ associated with the eigenvalue $\omega^2$. By using Lemma \ref{lem:unitary_equivalence2}, one easily shows that for all $k\in\R$, $\euler^{2\icplx\pi k y_1}V_\fbvar$ is an eigenvector of $\mathcal{A}_p (\fbvar+2\pi k\delta)$ associated with the eigenvalue $\omega^2$. By reversing the proof of \eqref{eq:link_spec_Riccati_spec_diff_op}-(i), one can then construct a non vanishing element of Ker$(\mathscr{M}^\droite - \lambda\, \euler^{2\icplx\pi k \delta} \, \mathscr{N}^\droite)$.}
\end{rmk}

\subsection{Case of evanescent frequencies: \texorpdfstring{$\omega^2\notin\sigma({\cal A}_\cutvec)$}{omega2 not in sigma(Atheta)}}
 If $\omega^2\notin \sigma({\cal A}_\cutvec)$, then as explained in the proof of Proposition \ref{prop:link_spec_Riccati_spec_diff_op}, 
 %according to Propositions \ref{thm:spectrunqp_s}  and \ref{prop:spectrum_halves}, $\omega^2$ does not belong to the spectrum of the 2D periodic differential operator $\mathcal{A}^\droite_p$ defined by \eqref{eq:op_per_j}. Therefore, all the problems introduced in Section \ref{sec:RtR_eps} (in particular the 2D periodic half-space problem \eqref{RobinHalfspaceproblems} and the Riccati system \eqref{Riccatiplus}) remain well-posed for $\veps = 0$, so that 
the approach and the results presented in Sections \ref{sec:RtR_eps} and \ref{sec:Riccati_eps} are still valid for $\veps = 0$. In consequence, %
% It can also be shown that for $\veps = 0$, the solution of the half-line problem \eqref{halfline_problems} decays exponentially at infinity (with a decaying rate which is related to the distance of $\omega^2$ to $\sigma({\cal A}_\cutvec)$). This is the reason why we refer to $\omega$ as an \textbf{evanescent frequency}.
%
% \vspspe
we can use exactly the same algorithm as the one described in Section \ref{sec:algo}, until Step 5 included. For the last step of the algorithm (namely Step 6), we study the well-posedness of the interior problem in Section \ref{sec:lap_intpb}. 
% Let us note from \eqref{decomp_spectre_augmente0} and \eqref{eq:link_spec_Riccati_spec_diff_op}--$(i)$ that
% \begin{equation} \label{eq:omeganotinsigma_riccati}
% \omega^2 \notin \sigma(\mathcal{A}_\cutvec)\quad \Longrightarrow \quad  \spforall \idomlap \in \{\gauche, \droite\}, \quad \left\{%
% 	\begin{array}{r@{\ }l}
% 		\displaystyle \sigma_{R, p}^\idomlap &= \big\{\lambda_k^\idomlap,\ k\in\Z\big\} \cup \big\{1/{\lambda_k^\idomlap},\ k\in\Z\big\},
% 		\ret
% 		\displaystyle \overline{\sigma_{R, p}^\idomlap} &= \sigma^\idomlap_R = \mathcal{C}\big(0,|\lambda^\idomlap_0|\big)\vts \cup\vts  \mathcal{C}\big(0,1/{|\lambda^\idomlap_0|}\big),
% 		\ret
% 		|\lambda^\idomlap_0| &< 1,
% 	\end{array}
% \right.
% \end{equation}
% where $\sigma_R^\idomlap$ is the Riccati spectrum defined in \eqref{eq:Ric_spec} with $\veps = 0$ and $\lambda_k^\idomlap$ is defined in Proposition \ref{thm_spectre} with $\veps=0$. These objects have been defined for $\idomlap = \droite$, but their definition extends to $\idomlap = \gauche$.

\vspspe
It can also be shown that for $\veps = 0$, the solution of the half-line problem \eqref{halfline_problems} decays exponentially at infinity in the sense of the next result, with a rate of decay which is linked to the modulus of the fundamental eigenvalue $\lambda^\idomlap_0$. This is the reason why we refer to $\omega$ as an \textbf{evanescent frequency}.

\begin{prop}\label{prop:exp_decay_evanescent}
	Let $U^\droite$ and $u^\droite_{\cutvec}$ be the respective solutions of the half-space problem \eqref{RobinHalfspaceproblems} and the half-line problem \eqref{halfline_problems} for $\veps = 0$. Then there exists a constant $C > 0$ such that for any $\varphi \in L^2_{\textit{per}}(\R)$,
	\begin{equation}\label{eq:exp_decay_evanescent}
		\forall n\in\N,\quad \left\{
			\begin{array}{r@{\ }l@{\qquad}l}
				\displaystyle
				\|U^\droite (\varphi) (\cdot + n\, \vv{\ev}_2)\|_{H^1_\cutvec (\cellper^{\droite, 0})} &\leq C\, |\lambda^\droite_0|^n\, \|\varphi\|_{L^2(0, 1)}, &(i)
				\ret
				\displaystyle \|u^\droite_\cutvec (\cdot + a^\droite + n/\cuti_2)\|_{H^1(0, 1/\cuti_2)} &\leq C\, |\lambda^\droite_0|^n. &(ii)
			\end{array}
		\right.
	\end{equation}
\end{prop}

\begin{dem}
	The cell by cell expression given in Proposition \ref{thm_reconstruction}  taken for $\veps = 0$ becomes
	\begin{equation*}
		\aeforall \yv \in \cellper^{\droite, 0}, \quad \big[U^\droite (\varphi)\big](\yv + n\, \vv{\ev}_2) = \big[ E^{\droite, 0} ((\mathcal{P}^\droite)^n\,\varphi)  + E^{\droite, 1} (\mathcal{S}^\droite\, (\mathcal{P}^\droite)^{n}\varphi)\big] (\yv), \quad \spforall n \in \N.
	\end{equation*}
	Since the spectral radius of $\mathcal{P}^\droite$ is $|\lambda^\droite_0|$, Gelfand's formula $|\lambda^\droite_0| = \lim_{n \to +\infty} \|(\mathcal{P}^\droite)^n\|^{1/n}$ and the well-posedness of the cell problems imply the point $(i)$ in \eqref{eq:exp_decay_evanescent}. The point $(ii)$ follows similarly using the 1D version \eqref{eq:cell_by_cell_1D} of the cell by cell expression of $U^\droite$ above.
\end{dem}

\subsection{Case of propagative frequencies \texorpdfstring{$\omega^2\in\sigma({\cal A}_\cutvec)$}{omega2 in sigma(Atheta)}}\label{sec:propagative_freq}
Let us now focus on the case $\omega^2\in\sigma(\mathcal{A}_\cutvec$). Under two assumptions that will be explicited in the next sections, we prove the existence of a limit for the half-guide solution $U^\idomlap_\veps$ and the half-line solution $u^\idomlap_{s, \veps}$.

\subsubsection{Convergence assumption on the fundamental eigenpair}\label{sec:conv_assumption_fund_eigenpair}
In order to study the limit when $\veps$ goes to $0$ of the half-guide solution $U^\idomlap_\veps (\varphi)$ (and subsequently of the associated half-line solutions $u^\idomlap_{s, \veps}$), it can be seen from the cell by cell reconstruction formula given in Proposition \ref{thm_reconstruction} that since we have convergence of $(E^{\droite, 0}_\veps, E^{\droite, 1}_\veps)$ by Proposition \ref{prop:lap_cell2d}, the only missing result is the convergence of the propagation operator $\mathcal{P}_\veps$ and the scattering operator $\mathcal{S}_\veps$ defined in \eqref{defPS}.

\vspspe
According to Theorem \ref{thm_spectre}, a condition for such a convergence is the existence of the limit of the fundamental pairs $(\lambda^\droite_{0,\veps},\varphi^\droite_{0,\veps})$ when $\veps$ goes to $0$. This observation motivates the following assumption.
\begin{assumption}[label={assumption}]%{Convergence of fundamental eigenpair}{} 
	Let $\idomlap \in \{\gauche, \droite\}$. At least up to a subsequence of values of $\veps$ converging to 0, the sequence $(\lambda^\idomlap_{0,\veps},\varphi^\idomlap_{0,\veps}, \psi^\idomlap_{0,\veps})$, with $\psi^\idomlap_{0,\veps}\coloneqq  \mathcal{S}_\veps\,\varphi^\idomlap_{0,\veps}$ converges in $\C \times L^\infty(\R)^2$, that is,
	\begin{equation} \label{convergencelambda}
			\spexists \lambda_0^\idomlap \in \C, \quad \lim_{\veps\to 0}\lambda^\idomlap_{0,\veps}=\lambda^\idomlap_0,
	\end{equation} 
	\begin{equation} \label{convergencephipsi}
		\spexists (\varphi_0^\idomlap, \psi_0^\idomlap) \in \mathscr{C}^0_{\textit{per}}(\R), \quad \lim_{\veps\to 0} \|\varphi_{0,\veps}^\idomlap - \varphi_0^\idomlap\|_{L^\infty(\R)}=\lim_{\veps\to 0}\|\psi_{0,\veps}^\idomlap- \psi_0^\idomlap\|_{L^
			\infty(\R)}=0.
	\end{equation} 
\end{assumption}

\vspspe
Of course, since $|\lambda^\idomlap_{0,\veps}| < 1$ and $\varphi_{0,\veps}^\idomlap (0) = 1$, one deduces that
\begin{equation*}
	|\lambda^\idomlap_{0}| \leq 1 \quad \textnormal{and} \quad \varphi_0^\idomlap (0) = 1, \quad \textnormal{so that} \quad \varphi_0^\idomlap \neq 0.
\end{equation*}
Also, since $(\varphi^\idomlap_{0,\veps}, \psi^\idomlap_{0,\veps})$ is in the kernel of $\mathscr{M}^\idomlap_\veps - \lambda^\idomlap_{0, \veps}\, \mathscr{N}^\idomlap_\veps$ according to \eqref{sous_espaces_propres}, we obtain the following as $\veps \to 0$, using the limit \eqref{eq:Riccati_op_ssabs} of $(\mathscr{M}^\idomlap_\veps, \mathscr{N}^\idomlap_\veps)$:
\begin{equation}\label{eq:Riccati_fundeig}
	\mathscr{M}^\idomlap   \begin{pmatrix}   \varphi^\idomlap_0  \\[2pt]  \psi^\idomlap_0  \end{pmatrix} = \lambda^\idomlap_0 \, \mathscr{N}^\idomlap  \begin{pmatrix}   \varphi^\idomlap_0 \\[4pt]  \psi^\idomlap_0  \end{pmatrix}.
\end{equation}
Consequently, $\lambda_0^\idomlap$ belongs to the Riccati point spectrum $\sigma^\idomlap_{R, p}$.

\begin{rmk} \label{rem_convP} 
	Considering the expression \eqref{diag} for $\mathcal{P}^\idomlap_\veps$ with respect to its fundamental eigenpair, one may object that Assumption \ref{assumption} is almost equivalent to assuming the existence of a limit $\mathcal{P}^\idomlap$ for the propagation operator $\mathcal{P}^\idomlap_\veps$, one difficulty being in which sense the convergence holds. This is in fact not completely immediate but we shall show that, with an additional assumption involving the notion of energy flux, there is convergence in $\mathscr{L}(L^2_{\textit{per}}(\R))$ (i.e. in operator norm).
\end{rmk}
\begin{rmk} \label{rem_convP_2} 
The sequence $(\lambda^\idomlap_{0, \veps})_\veps$ being bounded, the convergence \eqref{convergencelambda} is obvious. The convergence result \eqref{convergencephipsi} on the other hand is less obvious to establish since we are a priori missing compactness for $(\varphi^\idomlap_{0,\veps}, \psi^\idomlap_{0,\veps})$ in $L^\infty(\R)^2$.
Such compactness would be restaured if we normalized $(\varphi^\idomlap_{0,\veps}, \psi^\idomlap_{0,\veps})$ in a space that is compactly embedded in $L^\infty(\R)^2$, for instance in Sobolev spaces (which is possible with sufficiently smooth coefficients $(\mu_p, \rho_p)$). However, the property  $\varphi^\idomlap_{0}(0)=1$ would be lost a priori and we would no longer have the guarantee that the limit $(\varphi^\idomlap_{0}, \psi^\idomlap_{0})$ would not be identically $0$.
\end{rmk}

\noindent{Following \eqref{eigenfunctions} and Remark \ref{rmk:Spsikeps} for $\idomlap=\droite$ and replacing $\delta$ by $-\delta$ for $\idomlap=\gauche$,}{ we define 
\begin{equation} \label{eigenbasis} 
	\spforall k \in \Z, \quad \varphi^\idomlap_{k} \coloneqq  \euler^{2 \icplx\vts \pi\vts k s} \, \varphi^\idomlap_0 \quad \textnormal{and}\quad  \psi^\idomlap_{k} \coloneqq  \euler^{2 \icplx\vts \pi\vts k (s+\nu^\idomlap\,\delta)} \, \psi^\idomlap_0
\end{equation} }
with $\nu^\gauche = -\nu^\droite = 1$, so that, from \eqref{convergencephipsi},
\begin{equation} \label{convergenceeigenbasis} 
	\spforall k \in \Z, \quad \lim_{\veps \to 0}\|\varphi^\idomlap_{k,\veps}- \varphi^\idomlap_k\|_{L^ \infty(\R)} =  \lim_{\veps \to 0}\|\psi^\idomlap_{k,\veps}- \psi^\idomlap_k\|_{L^ \infty(\R)}=0. 
\end{equation} 
Moreover, using the same argument as for proving \eqref{eq:Riccati_fundeig}, we deduce
\begin{equation}\label{eq:Riccati_eig}
	\mathscr{M}^\idomlap \,   \begin{pmatrix}   \varphi^\idomlap_k  \\[2pt]  \psi^\idomlap_k  \end{pmatrix} = \lambda^\idomlap_k \, \mathscr{N}^\idomlap\,  \begin{pmatrix}   \varphi^\idomlap_k \\[4pt]  \psi^\idomlap_k  \end{pmatrix} , \quad \mbox{where} \quad \lambda^\idomlap_{k} \coloneqq  \euler^{2 \icplx\vts \pi\vts \nu^\idomlap k \delta} \, \lambda_{0}^\idomlap
\end{equation}
 This identity justifies calling $(\varphi^\idomlap_k ,\psi^\idomlap_k)$ a \emph{Riccati eigenmode} associated with the eigenvalue $\lambda^\idomlap_{k}$. 

\subsubsection{Notion of flux density}\label{sub:fluxdensity}
In order to ensure that the limiting absorption principle holds, a natural assumption concerns the energy flux of the propagative modes, see for instance \cite{joly2006exact,fliss2009analyse,nazarov2014umov,fliss2016solutions,kirsch2018radiation, fliss2021dirichlet} for periodic problems. Indeed, the energy flux is a physical criterion which appears naturally when studying the limiting absorption principle, and which allows determining if a propagative mode is an outgoing or an ingoing mode. If the energy flux of one of the propagative mode vanishes, the selection between outgoing and ingoing modes is not possible. Note that for certain situations, see for instance \cite{fliss2016solutions} which deals with perfectly periodic waveguides, when the energy flux of one of the propagative mode vanishes, the limiting absorption principle does not hold and it can be shown than $(u_\veps)_\veps$ diverges when $\veps$ tends to $0$. Hence, it is natural to suppose that the energy flux of the propagative modes does not vanish. This is the assumption we make in this study but only for the fundamental eigenfunction. We will see that under this assumption, we are able to show that the limiting absorption holds for our problem.

\vspspe
For simplicity, let $\idomlap = \droite$. To define a good notion of flux, we return to the estimate \eqref{eq:flux_eps} which, using the expression $\RtRop^\droite_\veps = \mathcal{T}^{\droite, 00}_\veps + \mathcal{T}^{\droite, 10}_\veps\, \mathcal{S}_\veps$, becomes
\begin{equation*}
	\spforall \varphi\in L^2_{\textit{per}}(\R) \setminus \{0\},\quad\frac{1}{4z}\Real\, \int_{\Sigmaper^{\droite, 0}} (I - \mathcal{T}^{\droite, 00}_\veps - \mathcal{T}^{\droite, 10}_\veps\, \mathcal{S}_\veps)\, \varphi\; \overline{(I+\mathcal{T}^{\droite, 00}_\veps + \mathcal{T}^{\droite, 10}_\veps\, \mathcal{S}_\veps)\, \varphi} > 0.
\end{equation*}
Using the expression \eqref{RtRweighted} of $(\mathcal{T}^{\droite, 00}_\veps, \mathcal{T}^{\droite, 10}_\veps)$ as weighted shift operators, that is to say $\mathcal{T}^{\droite, 00}_\veps\, \varphi(s) =	\boldsymbol{t}^{\droite, 00}_{\veps}(s) \, \varphi(s)$ and $\mathcal{T}^{\droite, 10}_\veps\, \varphi(s) = \boldsymbol{t}^{\droite, 10}_{\veps}(s) \, \varphi(s + \delta)$, the inequality above simply rewrites
\begin{equation}\label{eq:motiv_flux}
	\int_0^1 Q^\droite_{s, \veps} (\varphi, \mathcal{S}_\veps\, \varphi)\; ds > 0,
\end{equation}
where, for any $s \in \R$, the \emph{flux density} $Q^\droite_{s, \veps} (\varphi, \psi)$ is defined for $\varphi, \psi \in L^2_{\textit{per}}(\R)$ by
\begin{equation*}
	\displaystyle
	Q^\droite_{s, \veps} (\varphi, \psi) \coloneqq  \Real \Big[ \big(1 - \boldsymbol{t}^{\droite, 00}_{\veps}(s)\big) \, \varphi(s) - \boldsymbol{t}^{\droite, 10}_{\veps}(s) \, \psi(s + \delta) \Big]\, \Big[ \overline{\big(1 + \boldsymbol{t}^{\droite, 00}_{\veps}(s)\big) \, \varphi(s) + \boldsymbol{t}^{\droite, 10}_{\veps}(s) \, \psi(s + \delta)} \Big].
\end{equation*}
Since the functions $(\boldsymbol{t}^{\droite, 00}_{\veps}, \boldsymbol{t}^{\droite, 10}_{\veps})$ have a limit as $\veps \to 0$ according to Proposition \ref{prop:lap_cell1d}, it is natural to introduce the limit flux density
\begin{equation}\label{eq:def_limit_flux_density}
	\begin{array}{l}\displaystyle
	Q^\droite_{s} (\varphi, \psi) \coloneqq \\[5pt] \dsp  \Real \Big[ \big(1 - \boldsymbol{t}^{\droite, 00}(s)\big) \, \varphi(s) - \boldsymbol{t}^{\droite, 10}(s) \, \psi(s + \delta) \Big]\, \Big[ \overline{\big(1 + \boldsymbol{t}^{\droite, 00}(s)\big) \, \varphi(s) + \boldsymbol{t}^{\droite, 10}(s) \, \psi(s + \delta)} \Big].\end{array}
\end{equation}
as well as the sesquilinear form $(\Phi_1, \Phi_2) \mapsto q^\droite_s (\Phi_1, \Phi_2)$, defined for $\Phi_1 \coloneqq  (\varphi_1, \psi_1)$, $\Phi_2 \coloneqq  (\varphi_2, \psi_2)$ as
\begin{equation*}
	\displaystyle
	q^\droite_s(\Phi_1, \Phi_2) \coloneqq  \frac{1}{2} \big( q^{\droite, *}_s(\Phi_1, \Phi_2) + \overline{q^{\droite, *}_s(\Phi_2, \Phi_1)} \big),
\end{equation*}
where we have defined 
\begin{equation*}
	q^{\droite, *}_s(\Phi_1, \Phi_2) \coloneqq  \Big[ \big(1 - \boldsymbol{t}^{\droite, 00}(s)\big) \, \varphi_1(s) - \boldsymbol{t}^{\droite, 10}(s) \, \psi_1(s + \delta) \Big]\, \Big[ \overline{\big(1 + \boldsymbol{t}^{\droite, 00}(s)\big) \, \varphi_2(s) + \boldsymbol{t}^{\droite, 10}(s) \, \psi_2(s + \delta)} \Big].
\end{equation*}
By construction, $q^\droite_s (\cdot, \cdot)$ and $Q^\droite_s(\cdot)$ are related by
\begin{equation*}
	\spforall \Phi = (\varphi, \psi) \in L^2_{\textit{per}}(\R)^2, \quad Q^\droite_s(\varphi, \psi) = q^\droite_s (\Phi, \Phi).
\end{equation*}
For the sequel, it is useful to reformulate $q^\droite_s$ in terms of the function $w^\droite_s (\Phi) \in H^1(0, 1/\cuti_2)$ defined for any $\Phi = (\varphi, \psi) \in L^2_{\textit{per}}(\R)^2$ by
\begin{equation}\label{eq:def_ws_Phi}
	\displaystyle
	w^\droite_s(\Phi) \coloneqq  \varphi(s)\, e^{\droite, 0}_s + \psi(s + \delta)\, e^{\droite, 1}_s.
\end{equation}
This is the object of the next lemma.

\begin{lem}\label{lem:link_q_s_w_s}
	The sesquilinear form $q^\droite_s (\cdot, \cdot)$ is given for $\Phi_1, \Phi_2 \in L^2_{\textit{per}}(\R)^2$ by
	\begin{equation}\label{eq:link_q_s_w_s}
		\displaystyle
		q^\droite_s(\Phi_1, \Phi_2) = 2\vvts \icplx\vvts z\nu^\droite\, \Big(\mu_{s, \cutvec}\, \frac{d w^\droite_s(\Phi_1)}{d\xvi}\, \overline{w^\droite_s(\Phi_2)} - \mu_{s, \cutvec}\, w^\droite_s(\Phi_1)\, \overline{\frac{d w^\droite_s(\Phi_2)}{d\xvi}}   \Big) (0).%, \ \ \textnormal{with} \ \ w_{s, 1} \coloneqq  w^\droite_s(\Phi_1), \  w_{s, 2} \coloneqq  w^\droite_s(\Phi_2).
	\end{equation}
	% with
\end{lem}

\begin{dem}
	This result is a pure matter of calculations. Let $w_{s, 1} \coloneqq  w^\droite_s(\Phi_1)$ and $w_{s, 2} \coloneqq  w^\droite_s(\Phi_2)$. From the definition of $R^\droite_{\pm, s}$, one has
	\begin{equation*}
		 2\nu^\droite\mu_{s, \cutvec}\, \frac{d w_{s, 1}}{d\xvi} = (R^\droite_{+, s} - R^\droite_{-, s})\, w_{s, 1} \quad \textnormal{and} \quad -2\vts \icplx\vts z\, w_{s, 2} = (R^\droite_{+, s} + R^\droite_{-, s})\, w_{s, 2}.
	\end{equation*}
	Therefore, by definition of $w_{s, j} \coloneqq  w^\droite_s(\Phi_j)$ and of the RtR coefficients $\boldsymbol{t}^{\droite, \ell k}(s)$ (see \eqref{defRtRcoeff_0}), we obtain
	\begin{equation*}
		\left\{
			\begin{array}{r@{\ =\ }l}
				-2\vts \icplx\vts z\, w_{s, 1} & \big(1 + \boldsymbol{t}^{\droite, 00}(s)\big) \, \varphi_1(s) + \boldsymbol{t}^{\droite, 10}(s) \, \psi_1(s + \delta),
				\ret
				\displaystyle 2\nu^\droite \mu_{s, \cutvec}\, \frac{d w_{s, 2}}{d\xvi} & \big(1 - \boldsymbol{t}^{\droite, 00}(s)\big) \, \varphi_2(s) - \boldsymbol{t}^{\droite, 10}(s) \, \psi_2(s + \delta)
			\end{array}
		\right.
	\end{equation*}
	so that 
	\[\displaystyle
		 -4\vvts \icplx\vvts z\nu^\droite\Big(\mu_{s, \cutvec}\, w_{s, 1}\, \overline{\frac{d w_{s, 2}}{d\xvi}}\Big) (0) =  \overline{q^{\droite, *}_s(\Phi_2, \Phi_1)},\]\\
	and similarly,
		 \[4\vvts \icplx\vvts z\nu^\droite\Big(\mu_{s, \cutvec}\, \frac{d w_{s, 1}}{d\xvi}\, \overline{w_{s, 2}}\Big) (0) =  q^{\droite, *}_s(\Phi_1, \Phi_2).\]
	Taking the difference between these two last two equalities then leads to \eqref{eq:link_q_s_w_s}.
\end{dem}

\vspspe
For the Riccati eigenmodes $(\varphi^\droite_k, \psi^\droite_k)$ defined in Section \ref{sec:conv_assumption_fund_eigenpair}, the flux density has a particular form which we exhibit in the next result.

\begin{prop}\label{prop:flux_density_Riccati_eigs}
	Let $\Phi_k \coloneqq  (\varphi^\droite_k, \psi^\droite_k)$ for any $k \in \Z$. Since $\delta$ is irrational, one has 
	\begin{equation}
		\displaystyle
			\begin{array}{|l@{\quad \Longrightarrow \quad}l}
				|\lambda^\droite_0| < 1 & \spforall k, \ell \in \Z, \quad q^\droite_s (\Phi_k, \Phi_\ell) = 0,
				\ret
				|\lambda^\droite_0| = 1 & \spforall k, \ell \in \Z, \quad \spexists \alpha_{k, \ell} \in \C\setminus\{0\}, \quad q^\droite_s (\Phi_k, \Phi_\ell) = \alpha_{k, \ell}\, \euler^{2 \icplx\vts \pi\vts (k - \ell)\, s}.
			\end{array}
	\end{equation}
	%
	% $q^\droite_s (\Phi_k, \Phi_\ell) = 0$ if $|\lambda_0| < 1$. If $|\lambda| = 1$, then
	% \begin{equation}
	% 	\spforall k, \ell \in \Z, \quad \spexists \alpha_{k, \ell} \in \C, \quad q^\droite_s (\Phi_k, \Phi_\ell) = \alpha_{k, \ell}\, \euler^{2 \icplx\vts \pi\vts (k - \ell)\, s}.
	% \end{equation}
\end{prop}

\begin{dem}
	\textbf{Step 1.}\quad We begin by proving that
	\begin{equation}\label{eq:flux_density_riccati_dem0}
		\spforall k, \ell \in \Z, \quad q^\droite_{s + \delta} (\Phi_k, \Phi_\ell) = |\lambda_0^\droite|^2\, \euler^{2 \icplx\vts \pi\vts (k - \ell)\, \delta} q^\droite_s (\Phi_k, \Phi_\ell).
	\end{equation}
	To do so, the idea is to introduce for $s \in \R$ and for any $k \in \Z$ the limit as $\veps \to 0$ of the quasi-Floquet mode $u^\droite_{s, k, \veps}$ defined in Lemma \ref{rmk:recons_fund_sol}, that is, for $\xvi \in (0, 1/\cuti_2)$ and for any $n \in \N$,
	\begin{equation}\label{eq:recons_1dk_0}
		% u_{s,k,\veps}(x)= (\lambda^\droite_{k,\veps})^n \, \Big(\varphi^\droite_{k,\veps}(s+n\delta) \, \euler^{0}_{s+n\delta,\veps}\big(x-n /\theta_2\big) + \psi^\droite_{k,\veps}\big(s+(n+1)\delta\big)\, \euler^{1}_{s+n\delta,\veps}\big(x-n /\theta_2\big) \Big)
		u^\droite_{s, k}(\xvi + n/\cuti_2) \coloneqq  (\lambda^\droite_{k})^n \, \big[\varphi^\droite_{k} (s + n \delta)\, e^{\droite, 0}_{s + n\delta} (\xvi) + \psi^\droite_{k} (s + (n + 1) \delta)\, e^{\droite, 1}_{s + n\delta} (\xvi) \big].
	\end{equation} 
	From the ODE satisfied by $(e^{\droite, 0}_s, e^{\droite, 1}_s)$ and the property \eqref{eq:Riccati_eig} of the Riccati eigenmodes $(\varphi^\droite_k, \psi^\droite_k)$, it follows that $u^\droite_{s, k}$ satisfies the homogeneous differential equation
	\begin{equation}\label{eq:homogeneous_QP_line_pb}
		\displaystyle 
		- \frac{d}{d \xvi} \Big( \mu_{s, \cutvec} \, \frac{d u}{d \xvi} \Big) - \rho_{s, \cutvec} \, \omega^2 \, u = 0 \quad \textnormal{in} \quad \R_+.
	\end{equation}
	since $\omega^2 \in \R$, another solution of \eqref{eq:homogeneous_QP_line_pb} is $\overline{u^\droite_{s, k}}$.
	Now let $k, \ell \in \Z$. Since $u^\droite_{s, k}$ and $\overline{u^\droite_{s, \ell}}$ are both solutions of \eqref{eq:homogeneous_QP_line_pb}, the Wronskian
	\begin{equation*}
		\displaystyle
		\big[ W(u^\droite_{s, k}, \overline{u^\droite_{s, \ell}}) \big] (\xvi) \coloneqq  \Big(\mu_{s, \cutvec}\, u^\droite_{s, k}\, \overline{\frac{d u^\droite_{s, \ell}}{d\xvi}} - \mu_{s, \cutvec}\, \frac{d u^\droite_{s, k}}{d\xvi}\, \overline{u^\droite_{s, \ell}} \Big) (\xvi)
	\end{equation*}
	is constant with respect to $\xvi$, and in particular
	\begin{equation}\label{eq:flux_density_riccati_dem1}
		\displaystyle
		\big[ W(u^\droite_{s, k}, \overline{u^\droite_{s, \ell}}) \big] (0) = \big[ W(u^\droite_{s, k}, \overline{u^\droite_{s, \ell}}) \big] (1/\cuti_2).
	\end{equation}
	However, writing \eqref{eq:recons_1dk_0} for $n = 0$ and for $(k, \ell)$ shows directly that
	\begin{equation*}
		u^\droite_{s, k} = w^\droite_s (\Phi_k) \quad \textnormal{and} \quad u^\droite_{s, \ell} = w^\droite_s (\Phi_\ell) \quad \textnormal{in}\ (0, 1/\cuti_2),
	\end{equation*}
	with $\Phi_k \coloneqq  (\varphi_k, \psi_k)$, $\Phi_\ell \coloneqq  (\varphi_\ell, \psi_\ell)$, and where $w^\droite_s(\cdot)$ is defined by \eqref{eq:def_ws_Phi}. Thus by Lemma \ref{lem:link_q_s_w_s}, we deduce that $\big[ W(u^\droite_{s, k}, \overline{u^\droite_{s, \ell}}) \big] (0) = -q^\droite_s (\Phi_k, \Phi_\ell) / (2\vvts \icplx\vvts z\nu^\droite)$.

	\vspspe
	On the other hand, writing \eqref{eq:recons_1dk_0} for $n = 1$ leads directly to
	\begin{equation*}
		u^\droite_{s, k}(\cdot + 1/\cuti_2) = \lambda^\droite_k\, w^\droite_{s + \delta} (\Phi_k) \quad \textnormal{and} \quad u^\droite_{s, \ell}(\cdot + 1/\cuti_2) = \lambda^\droite_\ell\, w^\droite_{s + \delta} (\Phi_\ell) \quad \textnormal{in}\ (0, 1/\cuti_2).
	\end{equation*}
	Therefore, $\big[ W(u^\droite_{s, k}, \overline{u^\droite_{s, \ell}}) \big] (1/\cuti_2) = -\lambda^\droite_k\, \overline{\lambda^\droite_\ell}\; q^\droite_{s + \delta} (\Phi_k, \Phi_\ell) / (2\vvts \icplx\vvts z\nu^\droite)$. Then \eqref{eq:flux_density_riccati_dem0} follows from \eqref{eq:flux_density_riccati_dem1}.

	\vspspe
	\textbf{Step 2.}\quad \eqref{eq:flux_density_riccati_dem0} corresponds to an eigenvalue problem associated with the shift operator%, $s \mapsto q^\droite_s (\Phi_k, \Phi_\ell)$ is an eigenfunction of the shift operator 
	\begin{equation*}
		\displaystyle
		\tau_\delta: \varphi \in L^2_{\textit{per}}(\R) \mapsto \varphi(\cdot + \delta) \in L^2_{\textit{per}}(\R).
	\end{equation*}
	Using a Fourier series expansion (in a similar but much simpler way as for Proposition \ref{thm_spectre}), one shows that if $\delta$ is irrational, then the eigenvalues of $\tau_\delta$ are simple, and of the form $\euler^{2\icplx\vts \pi\vts k\vts \delta}$ for $k \in \Z$, with associated eigenfunctions $s \mapsto \euler^{2\icplx\vts \pi\vts k\vts s}$. Therefore, if $|\lambda^\droite_0| < 1$, then $q^\droite_s (\Phi_k, \Phi_\ell) = 0$. On the other hand, if $|\lambda^\droite_0| = 1$, then $s \mapsto q^\droite_s (\Phi_k, \Phi_\ell)$ is nothing but an eigenfunction of $\tau_\delta$	corresponding to the eigenvalue $\euler^{2 \icplx\vts \pi\vts (k - \ell)\, \delta}$, hence the result.
\end{dem}

\vspspe
Choosing $k = \ell = 0$ in Proposition \ref{prop:flux_density_Riccati_eigs} leads directly to the following.
\begin{cor}\label{prop:flux_density_constant_wrt_s}
	As $\delta$ is irrational, the flux density $Q^\droite_s(\varphi^\droite_0, \psi^\droite_0)$ is independent of $s$:
	\begin{equation}
		\displaystyle
		\spexists Q^\droite_0 \geq 0, \quad Q^\droite_s(\varphi^\droite_0, \psi^\droite_0) = Q^\droite_0 \quad \spforall s \in [0, 1].
	\end{equation}
	Moreover, if $|\lambda^\droite_0| < 1$, then $Q^\droite_0 = 0$.
\end{cor}

\begin{rmk}\label{rmk:flux_constant_k}
	Similarly, $Q^\droite_s(\varphi^\droite_k, \psi^\droite_k)$ is independent of $s$ for any $k \in \Z$. Moreover, by letting $Q^\droite_k$ be this constant value, it can be shown, using the definition of $Q^\droite_s(\varphi^\droite_k, \psi^\droite_k)$ that $Q^\droite_k = Q^\droite_0$ for any $k \in \Z$.
\end{rmk}

\vspspe
By extending all the above to the half-line problem on $I^\gauche$, we have a constant flux density $Q^\gauche_0$ associated with the fundamental eigenpair $(\varphi^\gauche_0, \psi^\gauche_0)$. {By positivity of the \eqref{eq:motiv_flux} of the energy flux in the absorbing case, and passing to the limit on this quantity, one obtains
\begin{equation}\label{eq:fluxpos}
		\displaystyle
		Q^\gauche_0 \geq 0.
	\end{equation} Our second assumption concerns the quantity $Q^\idomlap_0$ for $\idomlap \in \{\gauche, \droite\}$. 
	\begin{assumption}[label={ass:energy_flux}]%{Positive flux}{}
	Let $\idomlap \in \{\gauche, \droite\}$. The flux density $Q^\idomlap_0$ (introduced in Corollary \ref{prop:flux_density_constant_wrt_s} for $\idomlap = \droite$) does not vanish, namely
	\begin{equation}
		\displaystyle
		Q^\idomlap_0 \neq 0.
	\end{equation}
\end{assumption}}
%Inspired by the positivity \eqref{eq:motiv_flux} of the energy flux in the absorbing case, our second assumption concerns the quantity $Q^\idomlap_0$ for $\idomlap \in \{\gauche, \droite\}$.

% \begin{assumption}[label={ass:energy_flux}]%{Positive flux}{}
% 	Let $\idomlap \in \{\gauche, \droite\}$. The flux density $Q^\idomlap_0$ (introduced in Corollary \ref{prop:flux_density_constant_wrt_s} for $\idomlap = \droite$) is positive, namely
% 	\begin{equation}
% 		\displaystyle
% 		Q^\droite_0 > 0.
% 	\end{equation}
% \end{assumption}

\noindent
{
\begin{defi}\label{def:regfreq}
A frequency $\omega$ such that Assumptions \ref{assumption} and \ref{ass:energy_flux} are satisfied will be called a \textbf{regular frequency}. 
\end{defi}}
\noindent Note that under Assumption \ref{ass:energy_flux}, we have $|\lambda^\droite_0| = 1$, by Corollary \ref{prop:flux_density_constant_wrt_s}. Another consequence of this assumption is the following.
\begin{prop}\label{prop:non_null_fund_eig}
	{We have
	\begin{equation}
		\displaystyle
		 \spforall s \in \R, \quad |\varphi^\droite_0 (s)|^2 \geq Q^\droite_0,
	\end{equation}
	which implies that, under Assumption \ref{ass:energy_flux},  $\varphi^\droite_0 \in \mathcal{N}_{\textit{per}}$ and $\mathbf{w}(\varphi^\droite_0) = 0$.}
\end{prop} 

\begin{dem}
	Since $(\lambda^\droite_0, (\varphi^\droite_0, \psi^\droite_0))$ is a Riccati eigenpair, we have $\mathcal{T}^{\droite, 00}\, \varphi^\droite_0 + \mathcal{T}^{\droite, 10}\, \psi^\droite_0 = (\lambda^\droite_0)^{-1}\, \psi^\droite_0$, which, using the weighted shift representation \eqref{RtRweighted_0} leads to
	\begin{equation*}
\spforall s\in\R,\quad		\boldsymbol{t}^{\droite, 00}(s)\, \varphi^\droite_0(s) + \boldsymbol{t}^{\droite, 10}(s)\, \psi^\droite_0(s + \delta) = (\lambda^\droite_0)^{-1}\, \psi^\droite_0 (s).
	\end{equation*}
	Substituting this relation in the expression \eqref{eq:def_limit_flux_density} of $Q^\droite_s (\varphi^\droite_0, \psi^\droite_0)$ gives
	\begin{align*}
		\displaystyle
	\spforall s\in\R,\quad		Q^\droite_s (\varphi^\droite_0, \psi^\droite_0) &= \Real \big[ \varphi^\droite_0(s) - (\lambda^\droite_0)^{-1}\, \psi^\droite_0 (s) \big]\, \big[ \overline{\varphi^\droite_0(s) + (\lambda^\droite_0)^{-1}\, \psi^\droite_0 (s)} \big] 
		\\[8pt]
		&= |\varphi^\droite_0(s)|^2 - |\lambda^\droite_0|^{-2}\, |\psi^\droite_0(s)|^2.
	\end{align*}
	{Therefore, we have
	\[ \spforall s\in\R,\quad |\varphi^\droite_0(s)|^2 \geq Q^\droite_s (\varphi^\droite_0, \psi^\droite_0) = Q^\droite_0, \] 
	where we have used Corollary \ref{prop:flux_density_constant_wrt_s}.}
\end{dem}

\vspspe
Proposition \ref{prop:non_null_fund_eig} implies that for a regular frequency, the weight function
\(
	\boldsymbol{\rho}^\droite \coloneqq  |\varphi^\droite_0|^{-2}
\)
is bounded from above and below by positive constants. As a consequence,
\begin{equation}\label{eq:def_weighted_L2_inner_product_norm}
	\displaystyle
	(\varphi, \psi)_{\boldsymbol{\rho}^\droite} \coloneqq  \int_0^1 \varphi(s)\, \overline{\psi(s)}\, \boldsymbol{\rho}^\droite (s)\; ds \quad (\textnormal{with} \quad \|\varphi\|_{\boldsymbol{\rho}^\droite} \coloneqq  \sqrt{(\varphi, \varphi)_{\boldsymbol{\rho}^\droite}})
\end{equation}
defines an inner product in $L^2_{\textit{per}}(\R)$ whose associated norm $\|\cdot\|_{\boldsymbol{\rho}^\droite}$ is equivalent to the usual $L^2$--norm. Moreover, using Fourier series theory in $L^2(0, 1)$, one sees from the definition of $\varphi^\droite_k$ that the family $(\varphi_k^\droite)_{k \in \Z}$ is an orthonormal basis of $L^2_{\textit{per}}(\R)$ for $(\cdot, \cdot)_{\boldsymbol{\rho}^\droite}$. Finally, it follows from Proposition \ref{prop:flux_density_Riccati_eigs} that the Riccati eigenmodes $(\varphi_k^\droite, \psi_k^\droite)$ satisfy a bi-orthogonality property related to the sesquilinear form $q^\droite_s$, and which will be exploited in the sequel.

\begin{prop}\label{prop:biortho}
	Let $\Phi_k \coloneqq  (\varphi^\droite_k, \psi^\droite_k)$ for any $k \in \Z$. Since $\delta$ is irrational, one has
	\begin{equation}\label{eq:biortho}
		\spforall k \neq \ell, \quad \int_0^1 q^\droite_s (\Phi_k, \Phi_\ell)\; ds = 0.
	\end{equation}
\end{prop}

\begin{dem}
	This result follows directly by integrating the equality $q^\droite_s (\Phi_k, \Phi_\ell) = \alpha_{k, \ell}\, \euler^{2 \icplx\vts \pi\vts (k - \ell)\, s}$ (given in Proposition \ref{prop:flux_density_Riccati_eigs}).
\end{dem}

\subsubsection{Convergence of the propagation and scattering operators}
Let $\idomlap = \droite$. Since the operator $\mathcal{P}_\veps^\droite$ is a weighted shift operator (see Proposition \ref{thm_weighted-shift}), and since this weight is linked to the fundamental eigenfunction (see \eqref{eq:link_p_phi0}), it is natural to introduce the symbol
\begin{equation}\label{eq:weight_p}
	\displaystyle
	\spforall s \in \R, \quad \boldsymbol{p}^\droite (s) \coloneqq  \lambda^\droite_0\; \frac{\varphi^\droite_0 (s + \delta)}{\varphi^\droite_0 (s)},
\end{equation}
which is well-defined since $\varphi^\droite_0$ never vanishes according to Proposition \ref{prop:non_null_fund_eig}. We then define the limit propagation operator as follows
\begin{equation}\label{op_P_lim}
	\displaystyle
	\spforall \varphi \in L^2_{\textit{per}}(\R), \quad \spforall s \in \R, \quad \mathcal{P}^\droite\, \varphi (s) \coloneqq  \boldsymbol{p}^\droite (s - \delta)\, \varphi(s - \delta).
\end{equation}
We have by definition of the pair $(\lambda^\droite_k, \varphi^\droite_k)$ that
\begin{equation*}
	\displaystyle
	\spforall k \in \Z, \quad \mathcal{P}^\droite\, \varphi^\droite_k = \lambda^\droite_k\, \varphi^\droite_k \quad \textnormal{and} \quad \mathcal{P}^\droite\, \varphi = \sum_{k \in \Z} \lambda^\droite_k\, (\varphi, \varphi^\droite_k)_{\boldsymbol{\rho}^\droite}\, \varphi^\droite_k.
\end{equation*}
Similarly, by analogy with Proposition \ref{thm_weighted-shift} and Remark \ref{rmk:riccati_Pdroite}, we introduce the weight
\begin{equation}\label{eq:weight_s}
	\spforall s \in \R, \quad \boldsymbol{s}^\droite (s) \coloneqq  \lambda^\droite_0\, \frac{\psi^\droite_0 (s + \delta)}{\varphi^\droite_0 (s)},
\end{equation}
and we also define the limit scattering operator $\mathcal{S}^\droite$ as
\begin{equation} \label{op_S_lim}
	\spforall \varphi \in L^2_{\textit{per}}(\R), \quad \spforall s \in \R, \quad \mathcal{S}^\droite\, \varphi (s) \coloneqq  \boldsymbol{s}^\droite (s - \delta)\, \varphi(s - \delta).
\end{equation}
If $\omega$ is a regular frequency, then the next result can be shown.

\begin{prop}\label{thm:lap_PS}
	If $\omega$ is a regular frequency, then up to a subsequence, we have
	\begin{equation*}
		\lim_{\veps \to 0} \|\mathcal{P}^\droite_\veps - \mathcal{P}^\droite\|_{\mathscr{L}(L^2_{\textit{per}}(\R))} = 0 \quad \textnormal{and} \quad \lim_{\veps \to 0} \|\mathcal{S}^\droite_\veps - \mathcal{S}^\droite\|_{\mathscr{L}(L^2_{\textit{per}}(\R))} = 0.
	\end{equation*}
	Moreover, the pair $(\mathcal{P}^\droite,\mathcal{S}^\droite)$ solves the following Riccati system
	\begin{equation}\label{Riccati_sseps}
		% \left| \begin{array}{l}
		 \mbox{Find } (P,S) \in \mathscr{L}(L^2_{\textit{per}}(\R))^2 \mbox{ such that } 
		% \ret
	\quad \left\{ \begin{array}{llll} 
			P=\mathcal{T}^{\droite, 01} + \mathcal{T}^{\droite, 11} S ,
			\ret
			S   = \mathcal{T}^{\droite, 00} P+ \mathcal{T}^{\droite, 10} S \, P.
			\end{array}  \right.
		% \end{array}\right.
	\end{equation}
\end{prop}

\begin{dem}
	From the definitions of $\mathcal{P}^\droite_\veps$ and $\mathcal{P}^\droite$ as weighted shift operators, we have
	\begin{equation*}
		\spforall \varphi \in L^2_{\textit{per}}(\R), \quad \spforall s \in \R, \quad (\mathcal{P}^\droite_\veps - \mathcal{P}^\droite)\, \varphi (s) = [\boldsymbol{p}^\droite_\veps - \boldsymbol{p}^\droite] (s - \delta)\, \varphi(s - \delta),
	\end{equation*}
	so that $\|\mathcal{P}^\droite_\veps - \mathcal{P}^\droite\|_{\mathscr{L}(L^2_{\textit{per}}(\R))} = \|\boldsymbol{p}^\droite_\veps - \boldsymbol{p}^\droite\|_{\infty}$. But since $\boldsymbol{p}^\droite_\veps$ is linked to the fundamental eigenfunction $\varphi^\droite_{0, \veps}$ (see \eqref{eq:link_p_phi0}) and by definition \eqref{eq:weight_p} of $\boldsymbol{p}^\droite$,
	% \begin{equation*}
	% 	\spforall s \in \R, \quad \boldsymbol{p}^\droite_\veps (s) \coloneqq  \lambda^\droite_{0, \veps}\, \frac{\varphi^\droite_{0, \veps} (s + \delta)}{\varphi^\droite_{0, \veps} (s)} \quad \textnormal{and} \quad \boldsymbol{p}^\droite (s) \coloneqq  \lambda^\droite_0\, \frac{\varphi^\droite_0 (s + \delta)}{\varphi^\droite_0 (s)},
	% \end{equation*}
	and because $\varphi^\droite_0$ is bounded from below by a positive constant, we deduce from Assumption \ref{assumption} that $\boldsymbol{p}^\droite_\veps \to \boldsymbol{p}^\droite$ in $L^\infty (\R)$, which implies that $\mathcal{P}^\droite_\veps \to \mathcal{P}^\droite$ in $\mathscr{L}(L^2_{\textit{per}}(\R))$ as $\veps \to 0$. The convergence of $\mathcal{S}^\droite_\veps$ towards $\mathcal{S}^\droite$ in operator norm is obtained similarly. Finally, since $(\mathcal{P}^\droite_\veps, \mathcal{S}^\droite_\veps)$ is a solution of the Riccati system \eqref{Riccatiplus} and the RtR operators $\mathcal{T}^{\droite, \ell k}_\veps$ converge towards $\mathcal{T}^{\droite, \ell k}$ (Proposition \ref{prop:lap_cell2d}), it follows that $(\mathcal{P}^\droite,\mathcal{S}^\droite)$ is solution of \eqref{Riccati_sseps}.
\end{dem}

\subsubsection{Limiting absorption for the half-line problem}\label{sec:lap_half_line}
We fix $\idomlap = \droite$ for simplicity. Thanks to the limit propagation operator and the limit scattering operator, we can now define for any $\varphi \in L^2_{\textit{per}} (\R)$ the function
\begin{equation}\label{Reconstruction_formula_sseps}
	\displaystyle
	 \yv \in \cellper^{\droite, 0},\;n \in \N, \quad \big[U^\droite (\varphi)\big](\yv + n\, \vv{\ev}_2) = \big[ E^{\droite, 0} ((\mathcal{P}^\droite)^n\,\varphi)  + E^{\droite, 1} (\mathcal{S}^\droite\, (\mathcal{P}^\droite)^{n}\varphi)\big] (\yv), 
\end{equation}
by analogy with the cell by cell expression given in Proposition \ref{thm_reconstruction} of the 2D solution $U^\droite_\veps (\varphi)$ for $\veps > 0$. Thanks to the cell problems \eqref{cellpb0_sseps} and \eqref{cellpb1_sseps} satisfied by $E^{\droite, 0}$ and $E^{\droite, 1}$ respectively, one has that $U^\droite (\varphi)$ satisfies the PDE $- \Dt{} \big( \mu_p \, \Dt{} U^\droite (\varphi) \big) - \rho_p \, \omega^2 \, U^\droite (\varphi)$ in each cell $\cellper^{\droite, n}$ and is $1$--periodic with respect to $\yvi_1$. Furthermore, the Riccati system satisfied by $(\mathcal{P}^\droite,\mathcal{S}^\droite)$ ensures that $U^\droite (\varphi)$ and $\mu_p \, \Dt{} U^\droite (\varphi)$ are continuous across each interface $\Sigmaper^{\droite, n}$. Consequently,
\begin{equation} \label{RobinHalfguideproblems_sseps}
	\left|
	 \begin{array}{r@{\ =\ }l} 
			\displaystyle - \Dt{} \big( \mu_p \, \Dt{} U^\droite (\varphi) \big) - \rho_p \, \omega^2 \, U^\droite (\varphi) & 0 \quad \textnormal{in }\,\Omegaper^\droite, 
			\ret
			\displaystyle \mathcal{R}^\droite_+ \,  U^\droite (\varphi)|_{\Sigmaper^\droite} & \varphi,
			\ret
			\multicolumn{2}{c}{\textnormal{$U^\droite (\varphi)$ is periodic w.r.t. $\yvi_1$.}}
	\end{array}
	\right.
\end{equation} 
In what follows, $U^\droite (\varphi)$ is identified with its periodic extension with respect to $\yvi_1$.
\vspspe
Let us also introduce the RtR operator
\begin{equation}\label{eq:LambdaRtR}
	\RtRop^\droite =\mathcal{T}^{\droite, 00} + \mathcal{T}^{\droite, 10} \mathcal{S}^\droite.
\end{equation}

\begin{prop}\label{thm:lap_ULambda}
	If $\omega$ is a regular frequency, then up to a subsequence extraction, we have  
	\[
		\spforall(\varphi, n) \in L^2_{\textit{per}} (\R) \times \N,\quad \lim_{\veps \to 0} \|U^\droite_\veps (\varphi) - U^\droite (\varphi)\|_{H^1_\cutvec(\cellper^{\droite, n})} = 0  \quad \textnormal{and} \quad \lim_{\veps \to 0} \|\RtRop^\droite_\veps - \RtRop^\droite\|_{\mathscr{L}(L^2_{\textit{per}}(\R))} = 0.
	\]
	Moreover, for any $\varphi \in L^2_{\textit{per}} (\R) \setminus \{0\}$, the energy flux of $U^\droite (\varphi)$ is positive:
	\begin{equation}\label{eq:coercive_flux_2D}
		\displaystyle
		{Q}^\droite(\varphi, \mathcal{S}^\droite\, \varphi) \coloneqq  \Real\, \int_{\Sigmaper^{\droite, 0}} (I - \RtRop^\droite)\, \varphi\; \overline{(I+ \RtRop^\droite)\, \varphi} = Q^\droite_0\; \|\varphi\|^2_{\boldsymbol{\rho}^\droite} > 0
	\end{equation}
	where $Q^\droite_0$ is the positive flux density of the fundamental eigenvalue (see Corollary \ref{prop:flux_density_constant_wrt_s}, \eqref{eq:fluxpos} and Assumption \ref{ass:energy_flux}), and where $\|\cdot\|_{\boldsymbol{\rho}^\droite}$ is the weighted $L^2$--norm defined in \eqref{eq:def_weighted_L2_inner_product_norm}.
\end{prop}

\begin{dem}
	The convergence result follows directly from the convergence of the cell solutions $(E^{\droite, 0}_\veps, E^{\droite, 1}_\veps)$ and the associated local RtR operators $\mathcal{T}^{\droite, \ell k}_\veps$ (Proposition \ref{prop:lap_cell2d}), and from the convergence of the operators $(\mathcal{P}^\droite_\veps,\mathcal{S}^\droite_\veps)$ (Proposition \ref{thm:lap_PS}).

	\vspspe
	In order to derive the estimate \eqref{eq:coercive_flux_2D}, we expand $\varphi \in L^2_{\textit{per}} (\R)$ in the basis $(\varphi^\droite_k)_{k \in \Z}$.
	% and use the Corollary \ref{prop:flux_density_constant_wrt_s} as well as the bi-orthogonality relation in Proposition \ref{prop:biortho}. 
	By setting $\Phi_k \coloneqq  (\varphi^\droite_k, \psi^\droite_k)$ and $\alpha_k \coloneqq  (\varphi, \varphi^\droite_k)_{\boldsymbol{\rho}^\droite}$, the sesquilinearity of $q^\droite$ leads to 
	\begin{align*}
		\mathcal{Q}^\droite(\varphi, \mathcal{S}^\droite\, \varphi) &= \int_0^1 q^\droite_s \Big( \sum_{k \in \Z} \alpha_k\, \Phi_k, \ \sum_{\ell \in \Z} \alpha_\ell\, \Phi_\ell \Big)\; ds = \sum_{k, \ell \in \Z} \alpha_k\, \overline{\alpha_\ell}\, \int_0^1  q^\droite_s(\Phi_k, \Phi_\ell) \; ds
		\\
		&= \sum_{k \in \Z} |\alpha_k|^2\; \int_0^1 q^\droite_s(\Phi_k, \Phi_k) \; ds \quad \textnormal{from the bi-orthogonality relation \eqref{eq:biortho}}
		\\
		&= Q^\droite_0 \, \sum_{k \in \Z} |\alpha_k|^2=Q^\droite_0\; \|\varphi\|^2_{\boldsymbol{\rho}^\droite} \quad \textnormal{from Corollary \ref{prop:flux_density_constant_wrt_s}.}
	\vspace{-1\baselineskip}
	\end{align*}
	%Parseval's theorem ensures that $\sum_{k \in \Z} |\alpha_k|^2 = \|\varphi\|^2_{\boldsymbol{\rho}^\droite}$, thus leading to the desired result.
\end{dem}
%In what follows, we define $s^\droite \coloneqq  a^\droite\, \cuti_1$. 
\noindent Consider $\varphi \in \mathscr{C}^0_{\textit{per}}(\R)$ such that $\varphi(a^\droite\, \cuti_1) = 1$. Substituting the link between the 2D solutions $(E^{\droite, 0}(\varphi), E^{\droite, 1}(\varphi))$ and the 1D functions $(e^{\droite, 0}, e^{\droite, 1})$ in the cell by cell expression \eqref{Reconstruction_formula_sseps} of $U^\droite (\varphi)$, we deduce that $U^\droite (\varphi)$ admits a trace along the line $a^\droite\, \cuti_1\, \vv{\ev}_1 + \R_+\, \cutvec$:
\begin{equation}\label{cutformulas_sseps} 
	\displaystyle 
	\aeforall \xvi \in \R_+, \quad u^\droite (x + a^\droite) \coloneqq  U^\droite (a^\droite\, \cuti_1 + \theta_1\vts \xvi, \theta_2\vts \xvi).
\end{equation}
Moreover, consider the RtR coefficient
\begin{equation}\label{eq:lambda1DRtR}
	\RtRcoef^\droite \coloneqq  (R^\droite_{-, a^\droite\, \cuti_1} u^\droite)(a^\droite) = (\RtRop^\droite\varphi)(a^\droite\, \cuti_1). % \boldsymbol{t}^{\droite, 00}(s^\droite)\, \varphi(s^\droite) + \boldsymbol{t}^{\droite, 10}(s^\droite)\, (\mathcal{S}^\droite\, \varphi)(s^\droite + \delta),% 
\end{equation}
It appears at first sight that $(u^\droite, \RtRcoef^\droite)$ depend on the choice of $\varphi$. To prove that this is not the case, let us note that $u^\droite$ can also be written cell by cell for all $n\in\N$ as
%
%
% we introduce, in the spirit of \eqref{eq:cell_by_cell_1D}, the 1D function defined for any $n \in \N$ as
\begin{equation}\label{cell_by_cell_1D_sseps}
 u^\droite(\cdot + a^\droite + n/\cuti_2) = \varphi_n (a^\droite\, \cuti_1 + n \delta)\, e^{\droite, 0}_{a^\droite\, \cuti_1 + n\delta} + \psi_n (a^\droite\, \cuti_1 + (n + 1) \delta)\, e^{\droite, 1}_{a^\droite\, \cuti_1 + n\delta} \quad \textnormal{in}\ \ (0, 1/\cuti_2),
\end{equation}
with $\varphi_n \coloneqq  (\mathcal{P}^\droite)^n\, \varphi$ and $\psi_n \coloneqq  \mathcal{S}^\droite\, (\mathcal{P}^\droite)^n\, \varphi$, in the spirit of \eqref{eq:cell_by_cell_1D}. %Define also the RtR coefficient
% \begin{equation}\label{eq:lambda1DRtR}
% 	\lambda^\droite_\varphi \coloneqq  (R^\droite_{+, s^\droite} u^\droite_\varphi)(a^\droite) = \boldsymbol{t}^{\droite, 00}(s^\droite)\, \varphi(s^\droite) + \boldsymbol{t}^{\droite, 10}(s^\droite)\, (\mathcal{S}^\droite\, \varphi)(s^\droite + \delta),% (\RtRop^\droite\varphi)(a^\droite\, \theta_1).
% \end{equation}
% by analogy with \eqref{lambda_1D_shifts_eps}. For now, $(u^\droite_\varphi, \lambda^\droite_\varphi)$ depend a priori on $\varphi$. However, %
% one shows, similarly to \eqref{cutformulas} the 1D function
% 
% and the RtR coefficient
% 
% Note that using the cell by cell definition of \eqref{Reconstruction_formula_sseps} of $U^\droite(\varphi)$, $u^\droite$ can also be expressed as
The 1D problems (\ref{segment_problems0_sseps}, \ref{segment_problems1_sseps}) satisfied by $(e^{\droite, 0}_s, e^{\droite, 1}_s)$, and the Riccati system satisfied by $(\mathcal{P}^\droite, \mathcal{S}^\droite)$, imply that $u^\droite \in H^1_{\textit{loc}} (I^\droite)$, and is solution of
\begin{equation}
	\left|
	\begin{array}{r@{\ =\ }l}
		\displaystyle - \frac{d}{dx} \Big(\mu_\cutvec \,  \frac{d u^\droite}{dx}\Big)-  \rho_\cutvec \, \omega^2 \, u^\droite & 0 \quad \mbox{in}\ \ I^\droite
		\ret
		R^\droite_+ u^\droite(a^\droite) & 1.
	\end{array}
	\right.
\end{equation}
Since this problem admits at most one solution because of the Robin boundary condition, we deduce that $(u^\droite, \RtRcoef^\droite)$ does not depend on $\varphi$. The couple can be computed directly from the fundamental eigenpair:
for all $x\in (0,1/\cuti_2)$, for all $n\in\N$
	\begin{multline}\label{eq:recons_1d0_sseps}
		u^\droite(a^\droite+\xvi + n/\cuti_2) = (\lambda^\droite_{0})^n\,[\varphi^\droite_{0} (a^\droite\cuti_1)]^{-1} \, \\ \left[\varphi^\droite_{0} (a^\droite\cuti_1 + n \delta)\, e^{\droite, 0}_{a^\droite\cuti_1 + n \delta} (\xvi) + \psi^\droite_{0} (a^\droite\cuti_1 +  (n + 1) \delta)\, e^{\droite, 1}_{a^\droite\cuti_1 + n\delta} (\xvi) \right],
	\end{multline}  
	and the RtR coefficient $\RtRcoef^\droite$ can be computed thanks to
\begin{equation} \label{def_lambda_1D_dir_sseps}
\RtRcoef^\droite  =	 \boldsymbol{t}^{\droite, 00}(a^\droite\cuti_1)+[\varphi^\droite_{0} (a^\droite\cuti_1)]^{-1}\,\psi^\droite_{0} (a^\droite\cuti_1+\delta)\, \boldsymbol{t}^{\droite, 10}(a^\droite\cuti_1).
\end{equation} 
{Note that, by definition \eqref{eq:lambda1DRtR} of the RtR coefficient, 
% \[ \Big[ -\nu^\droite \, \mu_{\cutvec} \, \frac{d u^\droite}{d\xvi} -  \, \icplx  z  u^\droite \Big] (a^\droite) = \RtRcoef^{\droite}  \, \Big[ \nu^\droite \mu_{\cutvec} \,  \frac{du^\droite}{dx} - \icplx z  u^\droite \Big] (a^\droite)
% \]
a direct computation yields 
\begin{equation}\label{eq:lien_flux_rtr}
	\Imag \Big(\mu_{\cutvec}\, \frac{d u^\droite}{d\xvi}\, \overline{u^\droite}\Big) (a^\droite) = \frac{1}{4z} \Real\,\big[ (1 - \RtRcoef^\droite)\, \overline{(1 + \RtRcoef^\droite)} \big].
\end{equation}}
Furthermore, the next result holds.

\begin{prop}\label{thm:lap_ulambda_1D}
	If $\omega$ is a regular frequency, then up to a subsequence extraction, we have
	\begin{equation}\label{eq:conv_half_line_RtR}
		\spforall n\,\in\N,\quad \lim_{\veps \to 0} \|u^\droite_\veps - u^\droite\|_{H^1(a^\droite + n/\cuti_2+(0,1/\theta_2))} = 0  \quad \textnormal{and} \quad \lim_{\veps \to 0} |\RtRcoef^\droite_\veps - \RtRcoef^\droite| = 0.
	\end{equation}
	Moreover, the energy flux of $u^\droite$ is positive:
	\begin{equation}\label{eq:coercive_flux_1D}
		\displaystyle
		 \Real\,\big[ (1 - \RtRcoef^\droite)\, \overline{(1 + \RtRcoef^\droite)} \big] > 0.
	\end{equation}
\end{prop}

\begin{dem}
	Similarly to the proof of Proposition \ref{thm:lap_ULambda}, the convergence results in \eqref{eq:conv_half_line_RtR} are obtained using the cell by cell expression \eqref{cell_by_cell_1D_sseps} of $u^\droite$, as well as the convergence of $(e^{\droite, 0}_{s, \veps}, e^{\droite, 1}_{s, \veps})$ and the associated RtR coefficients $\boldsymbol{t}^{\droite, \ell k}_\veps$.

	\vspspe
	It remains to prove \eqref{eq:coercive_flux_1D}. Since $u^\droite$ does not depend on the choice of $\varphi$, we consider the expression \eqref{cell_by_cell_1D_sseps} for $\varphi = \alpha\, \varphi^\droite_0$, where $\alpha$ is chosen such that $\varphi(a^\droite\, \cuti_1) = 1$ (we recall that $\varphi^\droite_0(s) \neq 0$ for all $s$). Writing \eqref{cell_by_cell_1D_sseps} for $n = 0$ then shows that $u^\droite$ is simply $\alpha\, w_{a^\droite\, \cuti_1} (\Phi_0)$, where $w_s(\cdot)$ is defined by \eqref{eq:def_ws_Phi}, with $\Phi_0 \coloneqq  (\varphi_0, \mathcal{S}^\droite\, \varphi_0)$. Consequently, according to Lemma \ref{lem:link_q_s_w_s}, we have
	\begin{equation*}
		\displaystyle
		\frac{1}{4z} \Real\,\big[ (1 - \RtRcoef^\droite)\, \overline{(1 + \RtRcoef^\droite)} \big]\underset{\text{by \eqref{eq:lien_flux_rtr}}}{=}\Imag \Big(\mu_{\cutvec}\, \frac{d u^\droite}{d\xvi}\, \overline{u^\droite}\Big) (a^\droite) = \frac{|\alpha|^2}{4 z}\, q^\droite_{a^\droite\, \cuti_1} (\Phi_0, \Phi_0) = \frac{|\alpha|^2}{4 z}\,Q^\droite_0,
	\end{equation*}
	where $Q^\droite_0$ is the constant flux density identified in Corollary \ref{prop:flux_density_constant_wrt_s}. Finally, Assumption \ref{ass:energy_flux} corresponds to the positivity of $Q^\droite_0$, and therefore gives the desired result.
\end{dem}

\subsubsection{The limit Riccati system and spectral characterization of the limit propagation and scattering operators}
The computation of the propagation operator and of the scattering operator, and the subsequent construction of the limit solutions exhibited in Section \ref{sec:lap_half_line} rely on $(\lambda^\droite_0, \varphi^\droite_0, \psi^\droite_0)$. Similarly to the absorbing case (see Section \ref{sec:spectrum_Riccati}), the practical computation of this triple is not direct: it has to be extracted from the larger set of Riccati eigenmodes (which is the one we are able to compute in practice). The goal of this section is therefore to characterize $(\lambda^\droite_0, \varphi^\droite_0, \psi^\droite_0)$ using a criterion which can be computed numerically. 

\vspspe
Our starting point is the following result, which corresponds to Proposition \ref{thm:facto} for $\veps = 0$.% (although the proof is different).
\begin{prop}\label{thm:facto_sseps}
	Given $\veps>0$, let $\mathcal{P}^\droite, \mathcal{S}^\droite$ be the operators defined by \eqref{defPS}. For any $\lambda\in\C$, we have
	\begin{equation} \label{factor_sseps} 
		\mathscr{M}^\droite - \lambda \, \mathscr{N}^\droite = \begin{pmatrix} I & 0 \\[2pt] \lambda  \transp{\mathcal{S}^\droite} & \lambda \transp{\mathcal{P}^\droite} - I \end{pmatrix}\; \mathcal{L}^\droite\;\begin{pmatrix} \mathcal{P}^\droite - \lambda  & 0 \\[2pt] \mathcal{S}^\droite & - I \end{pmatrix},\quad \mathcal{L}^\droite \coloneqq  \begin{pmatrix} I  & - \mathcal{T}^{\droite, 11} \\[2pt] - \RtRop^\droite & I
		\end{pmatrix},
	\end{equation}
	where $\mathcal{T}^{\droite, 11}$ is defined in \eqref{eq:def_RtR_loc_lim} and $\RtRop^\droite$ in \eqref{eq:LambdaRtR}. Moreover, $\mathcal{L}^\droite$ is invertible.
\end{prop} 

\begin{dem}
	All the operators in the factorization formula of Proposition \ref{thm:facto} have been shown to have a limit as $\veps \to 0$. Passing to the limit in this factorization formula then leads directly to \eqref{factor_sseps}.

	\vspspe
	The most delicate part of this result is the invertibility of $\mathcal{L}^\droite$. Similarly to the proof of Proposition \ref{thm:facto}, we construct for some $(f, g) \in L^2_{\textit{per}}(\R)^2$ a vector $(\varphi, \psi) \in L^2_{\textit{per}}(\R)^2$ such that $\mathcal{L}^\droite (\varphi, \psi) = (f, g)$, that is,{
	\begin{equation}\label{eq:invert_Lr}
		\left\{
			\begin{array}{r@{\ =\ }l}
				\varphi-\mathcal{T}^{\droite, 11}\,\psi&f
				\\[8pt]
				- \RtRop^\droite\, \varphi + \psi &g
			\end{array}
		\right.
		\quad \Longleftrightarrow \quad 
		\left\{
			\begin{array}{r@{\ =\ }l}
				\varphi- {\cal R}^\droite_+ E^{\droite, 1} (\psi)|_{\Sigmaper^{\droite, 1}}&f
				\\[8pt]
				- \RtRop^\droite\, \varphi+ {\cal R}^\droite_-\, E^{\droite, 1} (\psi)|_{\Sigmaper^{\droite, 1}}&g
			\end{array}
		\right.
		\end{equation}}
	Similarly to the proof with $\veps > 0$, this problem is equivalent to: \textit{Find $(\varphi, U) \in L^2_{\textit{per}} (\R)\times H^1_\cutvec (\cellper^{\droite, 0})$ such that}
	\begin{equation} \label{pbfort_sseps}
		\left| \begin{array}{r@{\ =\ }l@{\quad}l} 
				-\Dt{} \big( \mu^\droite_p \, \Dt{} U\big) - \rho^\droite_p \, \omega^2 \, U & 0, & \textnormal{in }\cellper^{\droite, 0},	
				\\[8pt]
				{\cal R}_+^\droite \, {U} & 0, & \mbox{on } {\Sigmaper^{\droite, 0},}
				\\[8pt]
				(I+\RtRop^\droite)\varphi+2 \, \icplx z\,{U} & f-g,\; & \mbox{on } {\Sigmaper^{\droite, 1},}
					\\[8pt]
				(I-\RtRop^\droite)\varphi - 2 \vts \mu^\droite_p\,{\Dt{}U} & f+g, & \mbox{on } {\Sigmaper^{\droite, 1},}
				\\[8pt]
				\multicolumn{3}{c}{\textnormal{$U$ is periodic w.r.t. $\yvi_1$.}}
		\end{array} \right.
	\end{equation}
	More precisely, if $\mathcal{L}^\droite (\varphi, \psi) = (f, g)$, then by defining $U \coloneqq  E^1(\psi)$, the pair $(\varphi, U)$ is a solution of \eqref{pbfort_sseps}. Conversely, if $(\varphi, U)$ is a solution of \eqref{pbfort_sseps}, then by setting $\psi \coloneqq  \mathcal{R}^\droite_- U|_{\Sigmaper^{\droite, 0}}$, $(\varphi,\psi)$ satisfies \eqref{eq:invert_Lr}. In other words, proving the invertibility of $\mathcal{L}^\droite$ reduces to showing the well-posedness of the coupled problem \eqref{pbfort_sseps}.

	\vspspe
	Problem \eqref{pbfort_sseps} is not of Fredholm type, but it is equivalent to a family of 1D problems defined for any $s \in \R$ by: \textit{Find $(\varphi_s, u_s) \in \C \times H^1 (0, 1/\cuti_2)$ such that}
	\begin{equation} \label{pbfort_sseps_1D} 
		\left|
		\begin{array}{r@{\ =\ }l}
		\displaystyle - \, \frac{d}{dx} \Big(\mu_{s, \cutvec} \, \frac{d u_s}{dx}\Big)- \rho_{s, \cutvec} \, \omega^2 \, u_s & 0 \quad \mbox{in}\ \ (0, 1/\cuti_2), 
		\\[8pt]
		R^\droite_{+,s} u_s(0) & 0, %\quad \textnormal{and} \quad R^\droite_{-,s} U_s(1/\cuti_2) = 0,
		\\[8pt]
		\multicolumn{2}{l}{(1 + \lambda^\droite_s)\, \varphi_s + 2\vts \icplx\vts z\, u_s (1/\cuti_2) = f_s - g_s,}
		\\[8pt]
		\multicolumn{2}{l}{\displaystyle (1 - \lambda^\droite_s)\, \varphi_s - 2\mu_{s, \cutvec} \, \frac{d u_s}{dx} (1/\cuti_2) = f_s + g_s,}
	\end{array}	
		\right. 
	\end{equation}
	with $f_s \coloneqq  f(s + \delta)\in\C$, $g_s \coloneqq  g(s + \delta)\in\C$, and where $\lambda^\droite_s \coloneqq  (\RtRop^\droite \phi)(s)$, with $\phi = 1$. For any $s \in \R$, Fredholm alternative holds for Problem \eqref{pbfort_sseps_1D}, meaning that uniqueness implies existence and well-posedness of the problem. To prove the uniqueness, it suffices to consider the solution for $f_s = g_s = 0$. Using an integration by parts formula, we obtain that
	% \begin{equation} \label{pbfort_sseps_1D_homogene} 
	% 	\left\{
	% 		\begin{array}{l@{\ }l}
	% 			\displaystyle \int_{0}^{1/\cuti_2} \Big( \mu_{s, \cutvec}\, \Big|\frac{d u_s}{dx}\Big|^2 - \rho_p\, \omega^2\, |u_s|^2 \Big) - \icplx\vts z\, |u_s(0)|^2 - \frac{\icplx}{4z} (1 - \lambda^\droite_s)\, \overline{(1 - \lambda^\droite_s)}\, |\varphi_s|^2 = 0, & (i)
	% 			\rett
	% 			(1 + \lambda^\droite_s)\, \varphi_s + 2\vts \icplx\vts z\, u_s (1/\cuti_2) = 0. & (ii)
	% 		\end{array}
	% 	\right.
	% \end{equation} 
	% \begin{equation} \label{pbfort_sseps_1D_homogene} 
	% 			\displaystyle \int_{0}^{1/\cuti_2} \Big( \mu_{s, \cutvec}\, \Big|\frac{d u_s}{dx}\Big|^2 - \rho_p\, \omega^2\, |u_s|^2 \Big) - \icplx\vts z\, |u_s(0)|^2 - \frac{\icplx}{4z} (1 - \lambda^\droite_s)\, \overline{(1 - \lambda^\droite_s)}\, |\varphi_s|^2 = 0, 
	% \end{equation} 
	% By adapting the arguments of Proposition \ref{thm:lap_ulambda_1D}, one shows that $\Real\,\big[ (1 - \lambda^\droite)\, \overline{(1 + \lambda^\droite)} \big] > 0$. Therefore, by taking the imaginary part in \eqref{pbfort_sseps_1D_homogene}, we obtain that 
	% \begin{equation*}
	% 	\varphi_s = 0 \quad \textnormal{and} \quad U_s(0) = U_s(1/\cuti_2) = 0,
	% \end{equation*}
	% where the second equality follows from \eqref{pbfort_sseps_1D_homogene}--$(ii)$. Furthermore, the last equation in \eqref{pbfort_sseps_1D} implies that $(d U_s/dx)(1/\cuti_2) = 0$. Therefore, from Cauchy uniqueness theorem, $U_s = 0$.
{
	\begin{equation} \label{pbfort_sseps_1D_homogene} 
				\displaystyle \int_{0}^{1/\cuti_2} \Big( \mu_{s, \cutvec}\, \Big|\frac{d u_s}{dx}\Big|^2 - \rho_p\, \omega^2\, |u_s|^2 \Big) - \icplx\vts z\, |u_s(0)|^2 - \frac{\icplx}{4z} (1 - \lambda^\droite_s)\, \overline{(1 - \lambda^\droite_s)}\, |\varphi_s|^2 = 0. 
	\end{equation} 
	By adapting the arguments of Proposition \ref{thm:lap_ulambda_1D}, one shows that $\Real\,\big[ (1 - \lambda^\droite_s)\, \overline{(1 + \lambda^\droite_s)} \big] > 0$. Therefore, by taking the imaginary part in \eqref{pbfort_sseps_1D_homogene}, we obtain that 
	\begin{equation*}
		\varphi_s = 0 \quad \textnormal{and} \quad u_s(0) = 0.
	\end{equation*}
Furthermore, $R^\droite_{+,s} u_s(0)=0$ implies that $(d u_s/dx)(0) = 0$. Therefore, from Cauchy uniqueness theorem, $u_s = 0$.}

	\vspspe
	Since \eqref{pbfort_sseps_1D} is well-posed for any $s \in \R$, it can be shown by contradiction (see for instance the proof of Proposition \ref{prop:lap_cell1d}) that $s \mapsto \varphi_s \in L^2_{\textit{per}}(\R; \C)$ and $s \mapsto u_s \in L^2_{\textit{per}}(\R; H^1 (0, 1/\cuti_2))$. Consequently, the functions
	\begin{equation*}
		\displaystyle
		\aeforall (s, \xvi) \in \R \times (0, 1/\cuti_2), \quad \varphi(s) \coloneqq  \varphi_s \quad \textnormal{and} \quad U(s + \cuti_1\, \xvi, \cuti_2\, \xvi) \coloneqq  u_s(\xvi),
	\end{equation*} 
	are well-defined in $L^2_{\textit{per}}(\R)$ and $H^1_{\cutvec,\, \textit{per}}(\mathcal{C}^{\droite, 0})$ respectively. Moreover, from the properties of $(\varphi_s, u_s)$, one has that $(\varphi, U)$ is the unique solution of \eqref{pbfort_sseps}, and is continuous with respect to $(f, g)$.
\end{dem}

\begin{prop}\label{thm_spectreaugmented_sseps} 
	The Riccati point spectrum $\sigma^\droite_{R, p}$ is given by
	\begin{equation} \label{decomp_spectre_augmente0_sseps}
		\sigma^\droite_{R, p} = \big\{ \lambda^\droite_{k},\ k \in \Z \big\} \cup 	\big\{1/{\lambda^\droite_{k}}, \ k \in \Z \big \}, \quad \textnormal{with} \quad \matrixkernel \big( \mathscr{M}^\droite - \lambda^\droite_{k} \, \mathscr{N}^\droite \big) = \vect \begin{pmatrix} \varphi^\droite_{k} \\[4pt]
			\psi^\droite_{k} 
		\end{pmatrix}.
	\end{equation}
	% Moreover,
	% \begin{equation} \label{sous_espaces_propres_sseps_droite}
	% 	\matrixkernel \big( \mathscr{M}^\droite - \lambda^\droite_{k} \, \mathscr{N}^\droite \big) = \vect \begin{pmatrix} \varphi^\droite_{k} \\[4pt]
	% 	\mathcal{S}^\droite \varphi^\droite_{k} 
	% \end{pmatrix}.
	% \end{equation}
\end{prop}

\begin{dem}
	We begin by proving that $\{ \lambda^\droite_{k},\ k \in \Z \} \cup 	\{1/{\lambda^\droite_{k}}, \ k \in \Z  \} \subset \sigma^\droite_{R, p}$. Given $k \in \Z$, we already have from \eqref{eq:Riccati_eig} that $\lambda^\droite_k$ is a Riccati eigenvalue, with associated eigenfunction ${(\varphi^\droite_k, \mathcal{S}^\droite \varphi^\droite_{k})}$. Furthermore, using the symmetry properties of the operators $\mathcal{T}^{\droite, \ell k}$ (which are the same as in Proposition \ref{transposition_properties}), it can be computed without difficulty that
	\begin{equation*}
		\displaystyle
		\transp{(\mathscr{M}^\droite - \lambda \mathscr{N}^\droite)} = \widetilde{\mathscr{J}}^{-1}\, \Big(\mathscr{M}^\droite - \frac{1}{\lambda} \mathscr{N}^\droite\Big) \, \widetilde{\mathscr{J}} \quad \textnormal{with} \quad \widetilde{\mathscr{J}} \coloneqq  \begin{pmatrix} 0 & -\lambda\, I \\[2pt] I & 0\end{pmatrix}.
	\end{equation*}
	Since $\transp{(\mathscr{M}^\droite - \lambda \mathscr{N}^\droite)}$ and $\mathscr{M}^\droite - \lambda \mathscr{N}^\droite$ have the same eigenvalues, we get for any $\lambda \neq 0$ the equivalence $\lambda \in \sigma^\droite_{R, p} \ \Longleftrightarrow\ 1/\lambda \in \sigma^\droite_{R, p}$. Therefore, $\{1/{\lambda^\droite_{k}}, \ k \in \Z  \} \subset \sigma^\droite_{R, p}$.

	\vspspe
	Conversely, let $\lambda \in \sigma^\droite_{R, p}$ and consider an associated Riccati eigenmode $(\varphi_\lambda, \psi_\lambda) \neq 0$. If $(\lambda, \varphi_\lambda)$ is an eigenpair of $\mathcal{P}^\droite$, then $\lambda \in \{ \lambda^\droite_{k},\ k \in \Z \}$. Otherwise, since $\mathcal{L}^\droite$ is invertible, we have
	\begin{equation*}
		\displaystyle
		\begin{pmatrix} \widetilde{\varphi}_\lambda \\[2pt] \widetilde{\psi}_\lambda \end{pmatrix} \coloneqq  \mathcal{L}^\droite\;\begin{pmatrix} \mathcal{P}^\droite - \lambda  & 0 \\[2pt] \mathcal{S}^\droite & - I \end{pmatrix} \begin{pmatrix} \varphi_\lambda \\[2pt] \psi_\lambda \end{pmatrix} \neq 0 \quad \textnormal{and} \quad \begin{pmatrix} I & 0 \\[2pt] \lambda  \transp{\mathcal{S}^\droite} & \lambda \transp{\mathcal{P}^\droite} - I \end{pmatrix} \begin{pmatrix} \widetilde{\varphi}_\lambda \\[2pt] \widetilde{\psi}_\lambda \end{pmatrix} = 0.
	\end{equation*}
	It then follows that $\widetilde{\varphi}_\lambda = 0$ and $(\lambda \transp{\mathcal{P}^\droite} - I)\, \widetilde{\psi}_\lambda = 0$ with $\widetilde{\psi}_\lambda \neq 0$ because $(\widetilde{\varphi}_\lambda, \widetilde{\psi}_\lambda) \neq 0$. Consequently, $1/\lambda$ is an eigenvalue of $\transp{\mathcal{P}^\droite}$. Since $\mathcal{P}^\droite$ and $\transp{\mathcal{P}^\droite}$ have the same eigenvalues, we get $\lambda \in \{1/{\lambda^\droite_{k}}, \ k \in \Z  \}$.
	% Let us prove that $\lambda \coloneqq  1/\lambda^\droite_k$ also belongs to $\sigma^\droite_{R, p}$ for any $k \in \Z$. Since $\mathcal{P}^\droite$ and $\transp{\mathcal{P}^\droite}$ have the same eigenvalues, there exists $\widetilde{\varphi}_k \neq 0$ such that $(\lambda \transp{\mathcal{P}^\droite} - I)\, \widetilde{\varphi}_k = (\transp{\mathcal{P}^\droite} - \lambda^\droite_k)\, \widetilde{\varphi}_k = 0$. If $\lambda$ is an eigenvalue of $\mathcal{P}^\droite$, then $\lambda \in \sigma^\droite_{R, p}$ according to the above. Otherwise, 
\end{dem}

\vspspe
We recall that $(E^{\gauche, 0}, E^{\gauche, 1})$ are the local cell solutions associated with the half-line problem on $I^\gauche$, and that $(\mathcal{T}^{\gauche, 00}, \mathcal{T}^{\gauche, 10}, \mathcal{T}^{\gauche, 01}, \mathcal{T}^{\gauche, 11})$ are the associated local RtR operators. Similarly to \eqref{eq:ident_l_r}, we have under \eqref{eq:hyp_alar} that 
\begin{equation}\label{eq:ident_l_r_sseps}
	\begin{array}{l}
	\spforall \varphi \in L^2_{\textit{per}}(\R),\quad E^{\gauche,0}(\varphi)=E^{\droite,1}(\varphi),\quad E^{\gauche,1}(\varphi)=E^{\droite,0}(\varphi),
	\ret
	{\cal T}^{\gauche, 00}={\cal T}^{\droite,11},\quad {\cal T}^{\gauche, 10}={\cal T}^{\droite,01},\quad {\cal T}^{\gauche, 01}={\cal T}^{\droite,10}, \quad {\cal T}^{\gauche, 11}={\cal T}^{\droite,00},
	\end{array}
\end{equation}
and the next equality is obtained by making $\veps \to 0$ in \eqref{eq:link_l_r}: for any $\lambda\in\C$, $\lambda\neq 0$
\begin{equation}\label{eq:link_l_r_sseps}
	\mathscr{M}^\gauche-\lambda \mathscr{N}^\gauche = -\frac{1}{\lambda}\mathscr{J}^{-1}\Big(\mathscr{M}^\droite-\frac{1}\lambda \mathscr{N}^\droite \Big)\mathscr{J}\quad\textnormal{with}\; \mathscr{J}\coloneqq  \begin{pmatrix} 0 & I \\[2pt] I & 0 \end{pmatrix}.
\end{equation}
Consequently, we deduce the next result, which is the equivalent of Proposition \ref{thm_spectreaugmented_2}  for $\veps = 0$ (the proof is also the same).
\begin{prop}\label{thm_spectreaugmented_2_sseps} 
	We have
	\[
		\sigma_p(\mathcal{P}^\gauche) = \sigma_p(\mathcal{P}^\droite)\quad\textnormal{and}\quad\matrixkernel \Big( \mathscr{M}^\idomlap-\frac{1}{\lambda^{\idomlap'}_{k}} \mathscr{N}^\idomlap \Big) = \vect \begin{pmatrix} \psi^{\idomlap'}_{k} \\[4pt]
			\varphi^{\idomlap'}_{k} 
		\end{pmatrix} \quad\textnormal{for}\ \idomlap \neq \idomlap' \in \{\gauche, \droite\}.
	\]
\end{prop}			
\noindent We recall from Corollary \ref{prop:flux_density_constant_wrt_s} that if Assumption \ref{ass:energy_flux} holds, then $|\lambda^\droite_0| = 1$. This means according to Proposition \ref{thm_spectreaugmented_sseps} that \textit{all} the Riccati eigenvalues lie on the unit circle. Thus, contrary to the absorbing case, we cannot identify from \textit{all} the Riccati eigenvalues, the eigenvalues of the propagation operator using only their modulus. In order to characterize the fundamental eigenpair, we introduce a criterion using the flux. % If $Q^\droite_s(\varphi, \psi)$ denotes the flux density defined by \eqref{def_limit_flux_density}, then we recall from Corollary \ref{prop:flux_density_constant_wrt_s} and Remark \ref{rmk:flux_constant_k} that $Q^\droite_s (\varphi_k, \psi_k)$ is independent of $s$ and of $k$ (we call this constant value %
\begin{prop}\label{prop:flux_positive_negative}
	{For $\idomlap \in \{\gauche, \droite\}$, let $Q^\idomlap_s(\varphi, \psi)$ denote the flux density (replacing in \eqref{eq:def_limit_flux_density} $\droite$ by $\idomlap$), and let $Q^\idomlap_0$ be the constant flux density (replacing in Corollary \ref{prop:flux_density_constant_wrt_s} $ \droite$ by $\idomlap $). If $\omega^2$ is a regular frequency, then we have} 
	\begin{equation}\label{eq:flux_positive_negative}
		\displaystyle
		\spforall k \in \Z, \quad \spforall s \in \R, \quad Q^\droite_s (\varphi^\droite_k, \psi^\droite_k) = Q^\droite_0 > 0 \quad \textnormal{and} \quad Q^\droite_s (\psi^\gauche_k, \varphi^\gauche_k) = -Q^\gauche_0 < 0.
	\end{equation}
\end{prop}

 \begin{dem}
	{
	The first part of \eqref{eq:flux_positive_negative} is a direct consequence of Corollary \ref{prop:flux_density_constant_wrt_s} and Remark \ref{rmk:flux_constant_k}. For the second part, we begin by using the definitions \eqref{eq:def_limit_flux_density} of $Q^\gauche_s(\varphi, \psi)$ (replacing $\droite$ by $\gauche$). 
	The property ${(\varphi^\gauche_k, \psi^\gauche_k)} \in \matrixkernel \big(\mathscr{M}^\gauche - \lambda^\gauche_k\, \mathscr{N}^\gauche\big)$ on one hand implies $\mathcal{T}^{\gauche, 00}\, \varphi^\gauche_k + \mathcal{T}^{\gauche, 10}\, \psi^\gauche_k = (\lambda^\gauche_k)^{-1}\, \psi^\gauche_k$, so that for all $s$
	\begin{align*}
		Q^\gauche_s (\varphi^\gauche_k, \psi^\gauche_k) &\coloneqq  \Real \big[(I - \mathcal{T}^{\gauche, 00})\, \varphi^\gauche_k - \mathcal{T}^{\gauche, 10}\, \psi^\gauche_k\big]\, \big[\overline{(I + \mathcal{T}^{\gauche, 00})\, \varphi^\gauche_k + \mathcal{T}^{\gauche, 10}\, \psi^\gauche_k}\big](s).
		\\[8pt]
		&= |\varphi^\gauche_k (s)|^2 - |\psi^\gauche_k(s)|^2
	\end{align*}
	where we have used $|\lambda^\gauche_k| = 1$.  Similarly, the property ${(\psi^\gauche_k, \varphi^\gauche_k)} \in \matrixkernel \big(\mathscr{M}^\droite - (\lambda^\gauche_k)^{-1}\, \mathscr{N}^\droite\big)$ from Proposition \ref{thm_spectreaugmented_2_sseps} leads to $\mathcal{T}^{\droite, 00}\, \psi^\gauche_k + \mathcal{T}^{\droite, 10}\, \varphi^\gauche_k = \lambda^\gauche_k\, \varphi^\gauche_k$, so that for all $s$
	\begin{align*}
		Q^\droite_s (\psi^\gauche_k, \varphi^\gauche_k) &\coloneqq  \Real \big[(I - \mathcal{T}^{\droite, 00})\, \psi^\gauche_k - \mathcal{T}^{\droite, 10}\, \varphi^\gauche_k\big]\, \big[\overline{(I + \mathcal{T}^{\droite, 00})\, \psi^\gauche_k + \mathcal{T}^{\droite, 10}\, \varphi^\gauche_k}\big](s).
		\ret
		&= |\psi^\gauche_k (s)|^2 - |\varphi^\gauche_k(s)|^2,
	\end{align*}
	where we have also used $|\lambda^\gauche_k| = 1$. We then deduce $Q^\droite_s (\psi^\gauche_k, \varphi^\gauche_k) = -Q^\gauche_s (\varphi^\gauche_k, \psi^\gauche_k) = -Q^\gauche_0 < 0$ by combining the above two equalities, Corollary \ref{prop:flux_density_constant_wrt_s} and Assumption \ref{ass:energy_flux}.}
\end{dem}

\vspspe
Thanks to Proposition \ref{prop:flux_positive_negative}, we can separate the Riccati eigenmodes.

\begin{prop}\label{thm:sigmapm}
	We have $\sigma^\droite_{R, p} = \sigma^\droite_+ \cup \sigma^\droite_-$, where
	\begin{equation}\label{eq:sigmaRpm}
		\sigma^\droite_\pm \coloneqq  \Big\{ \lambda \in \C\ \  \big/\ \ \spexists (\varphi , \psi) \in \matrixkernel \big(\mathscr{M}^\droite - \lambda\, \mathscr{N}^\droite\big) \quad \textnormal{such that} \quad \pm Q^\droite(\varphi, \psi) > 0  \Big\}.
	\end{equation}
	Moreover,
	\begin{equation*}
		\displaystyle
		\exists\,!\; (\lambda_{0, +}\, (\varphi_{0, +}, \psi_{0, +})) \in \sigma^\droite_+ \times L^2_{\textit{per}}(\R)^2 \quad \textnormal{such that} \quad \varphi_{0, +} (0) = 1 \quad \textnormal{and} \quad \mathbf{w}(\varphi_{0, +}) = 0,
	\end{equation*}
	and
	\begin{equation*}
		\displaystyle
		\exists\,!\; (\lambda_{0, -}\, ( \psi_{0, -}, \varphi_{0, -} )) \in \sigma^\droite_- \times L^2_{\textit{per}}(\R)^2 \quad \textnormal{such that} \quad \varphi_{0, -} (0) = 1 \quad \textnormal{and} \quad \mathbf{w}(\varphi_{0, -}) = 0.
	\end{equation*}
\end{prop}

% \newpage

\subsection{Problem in the bounded interval and definition of the physical solution}\label{sec:lap_intpb}
We can now introduce the following problem in $(a^\gauche,a^\droite)$ using RtR boundary conditions 
\begin{equation}\label{eq:int_pb_RtR}
	\displaystyle
	\left|
	\begin{array}{r@{\ =\ }l}
		\displaystyle - \frac{d}{dx} \Big(\mu_\cutvec \,  \frac{d u^{\interieur}}{dx}\Big)-  \rho_\cutvec \, \omega^2 \, u^{\interieur} & 0 \quad \mbox{in}\ \ I^\interieur
		\ret
		\multicolumn{2}{c}{\displaystyle (R^\gauche_- u^{\interieur})(a^\gauche) = \RtRcoef^\gauche  \, (R^\gauche_+ u^{\interieur})(a^\gauche),}
		\ret
		\multicolumn{2}{c}{\displaystyle (R^\droite_- u^{\interieur})(a^\droite) = \RtRcoef^\droite  \, (R^\droite_+ u^{\interieur})(a^\droite).}
	\end{array}
	\right.
	\tag{$\mathscr{P}^{\texttt{int}}$}
\end{equation}

% \begin{equation}\label{eq:int_pb_RtR}
% 	\begin{array}{|rcl}
% 		\displaystyle - \frac{d}{d \xvi} \Big( \mu \, \frac{d u^0}{d \xvi} \Big) - \rho \, \omega^2 \, u^0 &=& f \quad \textnormal{in} \; I^0,\\[5pt]
% 		\displaystyle R_-^l u^0 &=& \RtRcoef^\gauche \; R_+^l u^0 \quad \mbox{at } x = a^\gauche,\\[5pt]
% 		\displaystyle R_-^r u^0 &=& \RtRcoef^\droite \; R_+^r u^0 \quad \mbox{at } x = a^\droite.
% 	\end{array}
% \end{equation}

\noindent
We show that this problem is well-posed provided that $\omega^2$ is not in the discrete spectrum $\sigma_d({\cal A})$ of $\mathcal{A}$. 
\begin{prop}\label{prop:lap_interior_pb}
	Suppose that either $\omega^2\notin\sigma(\mathcal{A}_\cutvec) \cup \sigma_d({\cal A})$, or $\omega^2\in\sigma(\mathcal{A}_\cutvec)$ and $\omega$ is a regular frequency. Then Problem \eqref{eq:int_pb_RtR} is well-posed in $H^1(I^\interieur)$ and 
	\[
		\lim_{\veps \to  0}\|u^\interieur_\veps-u^\interieur\|_{H^1(I^\interieur)}=0.
		\]
\end{prop}

\begin{dem}
	It is easy to see that Fredholm alternative holds for \eqref{eq:int_pb_RtR}. Therefore, it suffices to show uniqueness to deduce well-posedness. The sesquilinear form associated with \eqref{eq:int_pb_RtR} is given for any $u, v\in H^1(I^\interieur)$ by
	\begin{equation*}
		b^\interieur (u, v) \coloneqq  \int_{I^\interieur} \Big(\mu \,  \frac{d u}{dx}\, \overline{\frac{d v}{dx}} - \rho\, \omega^2\, u\, \overline{v}\Big)\,dx - \frac{\icplx}{4z} \sum_{\idomlap \in \{\gauche, \droite\}} (1 - \RtRcoef^\idomlap)\, R^\idomlap_+\, u (a^\idomlap)\; \overline{(1+\RtRcoef^\idomlap)\, R^\idomlap_+\, v(a^\idomlap)})% + (1-\RtRcoef^\gauche)R^l_+u(a^\gauche)\overline{(1+\RtRcoef^\gauche)R^l_+v(a^\gauche)})\right)
	\end{equation*}
	Suppose that $b^\interieur(u,v) = 0$ for all $v \in H^1(I^\interieur)$, and let us show that $u = 0$.
	
	\vspspe
	\textbf{Case $1$.} \quad %
	If $\omega^2\in\sigma(\mathcal{A}_\cutvec)$, by choosing $v = u$, we obtain $b^\interieur (u, u) = 0$, which implies that
	\[
		0 = \Imag(b^\interieur (u, u))=- \frac{1}{4z} \sum_{\idomlap \in \{\gauche, \droite\}} \Real \big[(1 - \RtRcoef^\idomlap)\, \overline{(1+\RtRcoef^\idomlap)}\big]\; |R^\idomlap_+\, u (a^\idomlap)|^2.
	\]
	Since $\Real \big[(1 - \RtRcoef^\idomlap)\, \overline{(1+\RtRcoef^\idomlap)}\big] > 0$ according to Proposition \ref{thm:lap_ulambda_1D}, it follows that
	\[
		R^\gauche_+u(a^\gauche) =  R^\droite_+u(a^\droite) = 0 \quad \underset{\textnormal{by \eqref{eq:int_pb_RtR}}}{\Longrightarrow} \quad R^\gauche_-u(a^\gauche) =  R^\droite_-u(a^\droite) = 0.
	\]
	This yields in particular $u(a^\droite)=u'(a^\droite)=0$ which, by Cauchy-Lipchitz theorem, gives $u=0$ in $(I^\interieur)$.

	\vspspe
	\textbf{Case $2$.} \quad If $\omega^2\notin \sigma(\mathcal{A}_\cutvec)$, we cannot conclude similarly. Instead, assume that there exists a solution $u^\interieur$ of the homogeneous version of \eqref{eq:int_pb_RtR}. We can then consider the function defined by
	\begin{equation*}
	\displaystyle
	\aeforall \xvi \in \R, \quad u(\xvi) = \left\{
		\begin{array}{c@{\quad}l}
			u^\interieur (\xvi), & \xvi \in I^\interieur,
			\\[6pt]
			R^\idomlap_+ u^\interieur\, (a^\idomlap)\; u^\idomlap (\xvi), & \xvi \in I^\idomlap,\; \idomlap \in \{\gauche, \droite\}.
		\end{array}
	\right.
\end{equation*}
	% \begin{equation*}
	% 	\displaystyle
	% 	\aeforall \xvi \in \R, \quad u(\xvi) \coloneqq  \left\{
	% 		\begin{array}{c@{\quad}l}
	% 			[R^\gauche_+ u^\interieur (a^\gauche)]\; u^\gauche (\xvi), & \xvi < a^\gauche,
	% 			\\[6pt]
	% 			u^\interieur (\xvi), & \xvi \in (I^\interieur),
	% 			\\[6pt]
	% 			[R^\droite_+ u^\interieur (a^\droite)]\; u^\droite (\xvi), & \xvi > a^\droite.
	% 		\end{array}
	% 	\right.
	% \end{equation*}
	Since $\omega^2\notin \sigma(\mathcal{A}_\cutvec)$, $u$ belongs to $H^1(\R)$ and is solution of the homogeneous version of \eqref{eq:whole_line_problem}. In conclusion, $u = 0$ except if $\omega^2\notin\sigma_d({\cal A})$.

	\vspspe
	Let us now show that $u_\veps^\interieur \to u^\interieur$ in $H^1(I^\interieur)$. It suffices to show that $(u^\interieur_\veps)_\veps$ is bounded in $H^1(I^\interieur)$. Indeed, in this case, up to a subsequence extraction, $u^\interieur_\veps\rightharpoonup u^*$ weakly in $H^1(I^\interieur)$ and $u^\interieur_\veps\to u^*$ strongly in $L^2(I^\interieur)$. We can then show that $u^*$ satisfy \eqref{eq:int_pb_RtR}. But since \eqref{eq:int_pb_RtR} is well-posed, $u^*=u^\interieur$ and the whole sequence converges strongly in $H^1(I^\interieur)$. To prove that $(u^\interieur_\veps)_\veps$ is bounded in $H^1(I^\interieur)$, we proceed by contradiction. Let us suppose that $\|u^\interieur_\veps\|_{H^1}\to +\infty$ and introduce $\widetilde{u}^\interieur_\veps\coloneqq u^\interieur_\veps\|u^\interieur_\veps\|^{-1}_{H^1}$. By linearity $\widetilde{u}^\interieur_\veps$ satisfies \eqref{eq:int_pb_eps_RtR} where the source $f$ is replaced by $f \|u^\interieur_\veps\|^{-1}_{H^1}$. The sequence $(\widetilde{u}^\interieur_\veps)_\veps$ is by definition bounded while $f \|u^\interieur_\veps\|^{-1}_{H^1}$ tends to $0$ when $\veps$ tends to 0. It then follows that $\widetilde{u}^\interieur_\veps\to 0$ strongly in $H^1$, which contradicts $\|\widetilde{u}^\interieur_\veps\|_{H^1}=1$.
\end{dem}

\vspspe
At last, we are able to define the physical solution $u$ of \eqref{eq:whole_line_problem} as follows:
\begin{equation}\label{eq:sol}
	\displaystyle
	\aeforall \xvi \in \R, \quad u(\xvi) = \left\{
		\begin{array}{c@{\quad}l}
			u^\interieur (\xvi), & \xvi \in I^\interieur,
			\\[6pt]
			R^\idomlap_+ u^\interieur\, (a^\idomlap)\; u^\idomlap (\xvi), & \xvi \in I^\idomlap,\; \idomlap \in \{\gauche, \droite\},
		\end{array}
	\right.
\end{equation}
where the convergence of $u_\veps$ to $u$ in $H^1_{\textit{loc}}$ results from Propositions \ref{thm:lap_ulambda_1D} and \ref{prop:lap_interior_pb}.

\subsection{Solution algorithm}\label{sec:algo_sseps}
In order to compute the solution of \eqref{eq:whole_line_problem} defined by limiting absorption, the previous sections provide an algorithm which sums up as follows. We suppose here that when $\omega^2\in  \sigma (\mathcal{A}_\cutvec)$, $\omega$ is a regular frequency, or equivalently that Assumptions \eqref{assumption} and \eqref{ass:energy_flux} are satisfied.% and the assumption \eqref{eq:hyp} (which is only a choice on $(a^\gauche,a^\droite)$) hold.

\begin{enumerate}
	\item Compute for all $s\in[0,1)$ the functions $e_{s}^{\droite, 0}$ and $e_{s}^{\droite, 1}$, solutions respectively of \eqref{segment_problems0_sseps} and \eqref{segment_problems1_sseps} as well as the local RtR symbols $t^{\droite, \ell k}(s)$  for $\ell,k\in\{0,1\}$ defined in \eqref{defRtRcoeff_0}. Deduce the cell operators $E^{\droite,0}$ and $E^{\droite,1}$ using Lemma \ref{lem_equivalence_ssabs}  as well as the local RtR operators ${\cal T}^{\droite,\ell k}$  for $\ell,k\in\{0,1\}$ via \eqref{RtRweighted_0}. Form the operators $\mathscr{M}^\droite$ and $ \mathscr{N}^\droite$ defined in \eqref{eq:Riccati_op_ssabs}.
	\item Compute the cell operators $E^{\gauche,0}$ and $E^{\gauche,1}$ as well as the local RtR operators ${\cal T}^{\gauche,\ell k}$ for $\ell,k\in\{0,1\}$ using \eqref{eq:ident_l_r_sseps}. Replacing $\delta$ by $-\delta$ and the superscript $\droite$ by $\gauche$, deduce the functions $e_{s}^{\gauche, 0}$ and $e_{s}^{\gauche, 1}$  using Lemma \ref{lem_equivalence_ssabs}  and then the local RtR symbols $t^{\gauche, \ell k}(s)$ for $\ell,k\in\{0,1\}$ via \eqref{RtRweighted} applied with $\varphi = 1$.
	\item Compute the Riccati point spectrum $\sigma^\droite_{R}$ (see \eqref{eq:Ric_spec_sseps}). For $\lambda \in \sigma^\droite_{R, p}$, denote ${(\varphi_\lambda,\psi_\lambda)} \in \matrixkernel (\mathscr{M}^\droite - \lambda\, \mathscr{N}^\droite)$.
	\begin{enumerate}[label=(\alph*), itemsep=8pt]
    \item \emph{If there is no $\lambda \in \sigma^\droite_{R}$ such that $|\lambda| = 1$, then $\omega^2\notin \sigma (\mathcal{A}_\cutvec)$ according to \eqref{eq:link_spec_Riccati_spec_diff_op}--$(iii)$}. 
    \\[4pt]
		Compute the unique $\lambda \in \sigma^\droite_{R, p}$ such that $|\lambda|<1$, $\varphi_\lambda(0)=1$ and ${\bf w}(\varphi_\lambda)=0$ and set $\lambda^\droite_{0}\coloneqq\lambda$ and $(\varphi^\droite_{0},\psi^\droite_0)\coloneqq (\varphi_\lambda,\psi_\lambda)$ (see Proposition \ref{thm_spectreaugmented_sseps}). Compute the unique $\lambda \in \sigma^\droite_{R, p}$ such that $|\lambda|>1$, $\psi_\lambda(0)=1$ and ${\bf w}(\psi_\lambda) = 0$ and set $\lambda^\gauche_{0}\coloneqq \lambda^{-1}$ and $(\varphi^\gauche_{0},\psi^\gauche_0)\coloneqq (\psi_\lambda,\varphi_\lambda)$ (see Proposition \ref{thm_spectreaugmented_2_sseps}).
  	\item \emph{If there exists $\lambda  \in \sigma^\droite_{R, p}$ such that $|\lambda| = 1$, then $\omega^2\in \sigma (\mathcal{A}_\cutvec)$ according to \eqref{eq:link_spec_Riccati_spec_diff_op}--$(i)$}.
		\\[4pt]
		Decompose $\sigma^\droite_{R,p}=\sigma_+^\droite \cup \sigma_-^\droite$, using Definition \eqref{eq:sigmaRpm}. Compute the unique $\lambda \in \sigma^\droite_{+}$ such that $\varphi_\lambda(0)=1$ and ${\bf w}(\varphi_\lambda)=0$ and set $\lambda^\droite_{0}\coloneqq \lambda$ and $(\varphi^\droite_{0},\psi^\droite_0)\coloneqq (\varphi_\lambda,\psi_\lambda)$ (see Propositions \ref{thm_spectreaugmented_sseps} and \ref{thm:sigmapm}). Compute the unique $\lambda\in \sigma_-^\droite$ such that $\psi_\lambda(0)=1$ and ${\bf w}(\psi_\lambda) = 0$ and set $\lambda^\gauche_{0}\coloneqq \lambda^{-1}$ and $(\varphi^\gauche_{0},\psi^\gauche_0)\coloneqq (\psi_\lambda,\varphi_\lambda)$ (see Propositions \ref{thm_spectreaugmented_2_sseps} and \ref{thm:sigmapm}).
		\item \emph{If there exist both $\lambda \in \sigma^\droite_{R}$ such that $|\lambda| = 1$ and $\lambda'  \in \sigma^\droite_{R, p}$ such that $|\lambda'| < 1$, then $\omega^2\in \sigma (\mathcal{A}_\cutvec)$ \eqref{eq:link_spec_Riccati_spec_diff_op}--$(i)$, whereas $Q^\droite_0 = 0$ (Corollary \ref{prop:flux_density_constant_wrt_s})}.
		In this case, we cannot identify $\lambda^\droite_0$ or $\lambda^\gauche_0$ from $\sigma^\droite_{R, p}$.
	\end{enumerate}
	\item Compute the RtR coefficient $\RtRcoef^\droite$ via formula \eqref{def_lambda_1D_dir_sseps} : note that this formula uses $t^{\droite,00}(s)$ and $t^{\droite,10}(s)$ from Step 1 and the functions $\varphi_{0}^\droite, \psi_{0}^\droite$ from Step 3. Similarly, compute the RtR coefficient $\RtRcoef^\gauche$. 
	\item Solve the interior problem \eqref{eq:int_pb_RtR}  to determine the solution $u$ inside $ I^\interieur$ (except if $\omega^2\in \sigma_d({\cal A})$). 
	Reconstruct the solution $u$  in $ I^\droite$ (\resp $I^\gauche$) via the reconstruction formula \eqref{eq:sol}-\eqref{eq:recons_1d0_sseps} (\resp  replacing in \eqref{eq:recons_1d0} $\delta$ by $-\delta$ and the superscript $\droite$ by $\gauche$). Note that this formula uses the eigenvalue 
	$\lambda_{0} ^\droite$ and the functions $\varphi_{0} ^\droite, \psi_{0}^\droite$ from Step 3 (\resp $\lambda_{0} ^\gauche$, $\varphi_{0} ^\gauche, \psi_{0}^\gauche$ from Step 4)  as well as the functions $e_{s}^{\droite, 0}$ and $e_{s}^{\droite, 1}$ from Step 1 (\resp $e_{s}^{\gauche, 0}$ and $e_{s}^{\gauche, 1}$  from Step 2). 
	\end{enumerate}
It is also possible to compute the half-guide solution $U^\idomlap(\varphi)$ for any $\varphi$. After the steps 1 to 4, one has to add the following steps.
	\begin{enumerate}
			\item[4-bis] Compute $\mathcal{P}^\idomlap$ for $\idomlap = \gauche$ and $\droite$ using \eqref{eq:weight_p} and \eqref{op_P_lim}. The operators $\mathcal{S}^\idomlap$ for $\idomlap = \droite$ and $\idomlap = \gauche$ can also be computed using \eqref{eq:weight_s} and \eqref{op_S_lim}.
	\item[5-bis] Reconstruct for any $\varphi\in L^2_{\textit{per}}(\R)$, the solution $U^\idomlap(\varphi)$ defined by limiting absorption using \eqref{Reconstruction_formula_sseps}.
	\end{enumerate}
% \subsection{Link with a notion of group velocity}
\subsection{Quasi modes, dispersion curves, group velocity and link with the flux density}\label{sub:disp_curves}
\noindent The objective of this section is to introduce, for quasiperiodic media satisfying some assumptions, a notion of dispersion curves and the associated group velocity, using all the objects introduced previously. We want also to show the link between the group velocity and the flux density introduced in Section \ref{sub:fluxdensity}.  We suppose here that $\omega^2\in  \sigma (\mathcal{A}_\cutvec)$, $\omega$ is a regular frequency, or equivalently that Assumptions \eqref{assumption} and \eqref{ass:energy_flux} are satisfied.

\vspspe
Let $\idomlap\in\{\droite,\gauche\}$.
To each eigenfunction $\varphi^\idomlap_{k}$, for $k\in\Z$, defined in \eqref{eigenbasis}, of $\mathcal{P}^\idomlap$, we can associate a {\it 2D mode} denoted $U^\idomlap_{k}$ that is given by
\begin{equation}\label{2Dmode}
	\displaystyle
	 \yv \in \cellper^{\idomlap, 0},\;n \in \Z, \quad \big[U^\idomlap_k\big](\yv + n\, \vv{\ev}_2) = (\lambda_k^\idomlap)^n \big[ E^{\idomlap, 0} (\varphi^\idomlap_k)  + E^{\idomlap, 1} (\psi^\idomlap_k)\big] (\yv).
\end{equation}
where we recall that $\psi^\idomlap_k=\mathcal{S}^\idomlap\varphi^\idomlap_k$ is also given in \eqref{eigenbasis}.
By definition \eqref{eq:Riccati_eig} of $\lambda_k^\idomlap$ and denoting $\lambda_0^\idomlap\equiv\lambda_0^\idomlap(\omega)=e^{2\icplx\pi k_0^\idomlap(\omega)}$ with $k_0^\idomlap(\omega)\in [-1/2,1/2[$, the mode rewrites\footnote{To simplify the notation, we explicitly indicate the dependence on $\omega$ only for $\lambda_0^\idomlap$ and $k_0^\idomlap$, although all the functions introduced in this section are implicitly $\omega-$dependent.} 
\begin{equation}\label{2Dmode_bis}
	\displaystyle
	 \yv \in \cellper^{\idomlap, 0},\;n \in \Z, \quad \big[U^\idomlap_k\big](\yv + n\, \vv{\ev}_2) = \euler^{2\icplx\pi n (k\nu^\idomlap\delta+k_0^\idomlap(\omega))}\big[ E^{\idomlap, 0} (\varphi^\idomlap_k)  + E^{\idomlap, 1} (\psi^\idomlap_k)\big] (\yv).
\end{equation}
or equivalently 
\begin{equation}\label{2Dmode_ter}
	 \yv \in (0,1)\times\R, \quad U^\idomlap_k(\yv) = \euler^{2\icplx\pi y_2 (k\nu^\idomlap\delta+k_0^\idomlap(\omega))} V^\idomlap_k(\yv),\\
\end{equation}
where  
\begin{equation}\label{2Dmode_ter_2}
	\yv \in \cellper^{\idomlap, 0}, \quad V^\idomlap_k(\yv)\coloneqq \euler^{-2\icplx\pi y_2 (k\nu^\idomlap\delta+k_0^\idomlap(\omega))} \big[ E^{\idomlap, 0} (\varphi^\idomlap_k)  + E^{\idomlap, 1} (\psi^\idomlap_k)\big] (\yv)
\end{equation}
 defines a periodic function w.r.t. $\yv$ (see the first equation of the Riccati system \eqref{Riccati_sseps} applied to $\varphi^\idomlap_k$). 

\vspspe
Each function $U^\idomlap_{k}, \,k\in\Z$ is called a mode since it is a particular solution, of the form $U(\yv)=e^{2\icplx\pi\xi y_2}V(\yv)$ with $\xi\in\R$ and $V$ periodic, of 
\begin{equation} \label{modes_problem}
	\left|
	 \begin{array}{r@{\ =\ }l} 
			\displaystyle - \Dt{} \big( \mu_p \, \Dt{} U \big) - \rho_p \, \omega^2 \, U & 0 \quad \textnormal{in }\,(0,1)\times\R, 
			\ret
			\multicolumn{2}{c}{\textnormal{$U$ is periodic w.r.t. $\yvi_1$.}}
	\end{array}
	\right.
\end{equation} 
% \varphi(s)\, u^\droite_{s, \veps}(\xvi + n/\cuti_2) = \varphi^\veps_n (s + n \delta)\, e^0_{s + n\delta, \veps} (\xvi) + \psi^\veps_n (s + (n + 1) \delta)\, e^1_{s + n\delta, \veps} (\xvi).
It is a Floquet mode associated with this non elliptic equation with periodic coefficients.
\begin{rmk}\label{rmk:2dmodes_temps}
	The time-dependent function associated with $U^\idomlap_{k}$ defined by 
\begin{equation} \label{quasi-mode-temps}
	\widetilde{U}^\idomlap_{k}(\yv,t) \coloneqq U^\idomlap_{k}(\yv) \, e ^{-\icplx \omega t}  = \euler^{\icplx(2\pi y_2 (k\nu^\idomlap\delta+k_0^\idomlap(\omega))-\omega t)} V^\idomlap_k(\yv)
	\end{equation}
can be interpreted as a harmonic wave with frequency $\omega$ and wavenumber $2\pi (k\nu^\idomlap\delta+k_0^\idomlap(\omega))$ and whose amplitude is periodic w.r.t. $\yv$.
\end{rmk}

\vspspe Let us now define, for $\idomlap\in\{\droite,\gauche\}$ and $k\in\Z$, the one-dimensional mode defined by
\begin{equation}\label{1Dmode}
	\displaystyle
	 \xvi\in\R,\;s\in(0,1), \quad u^\idomlap_{k,s}(\xvi)=\big[U^\idomlap_k\big](s +\xvi\cuti_1, x\cuti_2).
\end{equation}
that rewrites using \eqref{2Dmode_bis} and Lemma \ref{lem_equivalence_ssabs} (recall that $\nu^\gauche=-\nu^\droite=1$):
for all $x\in (0,1/\cuti_2)$, for all $n\in\Z$
	\begin{multline}\label{eq:uks}
		u^\idomlap_{k,s}(\xvi + n/\cuti_2) = \euler^{2\icplx\pi n (k\nu^\idomlap\delta+k_0^\idomlap(\omega))}\, \left[\varphi^\idomlap_{k} (s - n\nu^\idomlap \delta)\, e^{\idomlap, 0}_{s - n \nu^\idomlap\delta} (\xvi) + \psi^\idomlap_{k} (s -  (n + 1)\nu^\idomlap \delta)\, e^{\idomlap, 1}_{s - n\nu^\idomlap\delta} (\xvi) \right].
	\end{multline} 
By \eqref{eigenbasis}, it is easy to see that 
\[ \forall s\in[0,1),\quad u^\idomlap_{k,s}(\xvi)=e^{2\icplx \pi ks}u^\idomlap_{0,s}(\xvi), \]
which means that for a fixed $s$, the $u^\idomlap_{k,s}$'s define the same function up to a multiplicative factor. We denote in what follows $u^\idomlap_s\coloneqq u^\idomlap_{0,s}$. Using (\ref{2Dmode_ter}-\ref{2Dmode_ter_2}), $u^\idomlap_s$ rewrites 
\begin{equation}\label{1Dmode_bis}
	\displaystyle
	 x\in\R,\;s\in(0,1), \quad u^\idomlap_{s}(x)=\euler^{2\icplx\pi x\theta_2 k_0^\idomlap(\omega)} v^\idomlap_{0,s}(x).
\end{equation}
where $v^\idomlap_{0,s}(x)\coloneqq V^\idomlap_0(s +x\theta_1, x\theta_2 )$ is a quasiperiodic function of order 2. 

\vspspe For all $s\in[0,1)$, $u^\gauche_{s}$ and $u^\droite_{s}$ are called one-dimensional modes since they are particular solutions, of the form $u(x)=e^{2\icplx\pi\xi x}v(x)$ with $\xi\in\R$ and $v$ quasiperiodic, of the one-dimensional Helmholtz equation
\begin{equation} \label{modes1D_problem}
	\displaystyle - \frac{d}{dx} \Big(\mu_{s, \cutvec} \,  \frac{d u}{dx}\Big) - \rho_{s, \cutvec} \, \omega^2  \, u = 0 \quad \mbox{in}\ \ \R.
\end{equation} 
\begin{rmk}
		Note that for constant coefficients, the modes are particular solutions of the Helmholtz equation of the form $e^{2\icplx\pi\xi x}$ with $\xi\in\R$ whereas for periodic coefficients the modes are of the form $e^{2\icplx\pi\xi x}v(x)$ with $\xi\in\R$ and $v$ periodic.
	\end{rmk}
\begin{rmk}
	The time-dependent function associated with $u^\idomlap_{s}$ defined by 
\begin{equation} \label{quasi-mode-temps-1D}
	\widetilde{u}^\idomlap_{s}(x,t) \coloneqq u^\idomlap_{s}(x) \, e ^{-\icplx \omega t}  = \euler^{\icplx(2\pi x\theta_2 k_0^\idomlap(\omega)-\omega t)}\, v^\idomlap_{0,s}(x) 
	\end{equation}
can be interpreted as a harmonic wave with frequency $\omega$ and wavenumber $2\pi \theta_2 k_0^\idomlap(\omega)$ and whose amplitude is spatially modulated in a quasiperiodic fashion through the function $v^\idomlap_{0,s}(x)$. 
Note that the wavenumber does not depend on $s$. The so-called {\it dispersion relation} that we define here as $\omega \mapsto k_0^\idomlap(\omega)$ is therefore the same for all $s$.
\end{rmk}
 \begin{rmk} \label{rmk:basedesol}
 	Note that for all $s$, $\overline{u^\idomlap_{s}}$ is also a one-dimensional mode whose wavenumber is $-2\pi \theta_2 k_0^\idomlap(\omega)$. Since $\omega$ is a regular frequency and using the proof of Proposition \ref{prop:flux_density_Riccati_eigs} (see Step 1 and the discussion on the wronskian), one can deduce that $(u^\idomlap_s,\overline{u^\idomlap_{s}})$ forms a basis of the vector space (of dimension 2) of solutions of \eqref{modes1D_problem}. In particular, since $[\lambda_0^\droite]^{-1}\equiv e^{-2\icplx\pi\theta_2 k_0^\droite(\omega)}$ is an eigenvalue of $\mathcal{P}^\gauche$ by Proposition \ref{thm_spectreaugmented_2_sseps}, there exists $p\in\Z$ such that $-k_0^\droite(\omega)=p\delta+k_0^\gauche(\omega)$ and $\overline{u^\droite_{s}}={u^\gauche_{s}}e^{2\icplx\pi p s}$.
% 	\\\\
% 	Let us denote $Q_0^\idomlap(\omega)\coloneqq Q_0^\idomlap$ where $Q_0^\idomlap$ is defined in Corollary \ref{prop:flux_density_constant_wrt_s}. A direct computation yields 
% 	\[ 
% 	Q^\idomlap_0(\omega)\coloneqq -4z\nu^\idomlap \Imag \Big(\mu_{\cutvec}\, \frac{d u^\idomlap_{s}}{d\xvi}\, \overline{u^\idomlap_{s}}\Big) (x) \quad (\text{that is independent of $x$ and $s$}).
% 	\]
% 	It is easy to deduce now that $u^\gauche_s$ is proportional to $\overline{u^\droite_s}$. Indeed, since $u^\gauche_s$ satisfies \eqref{quasi-mode-temps-1D}, there exists $(\alpha,\beta)\in\C^2$ such that $u^\gauche_s=\alpha u^\droite_s+\beta u^droite_s$. One easily shows that 
% 	\[ Q^\gauche_0(\omega) = - (|\alpha|^2-|\beta|^2) Q^\droite_0(\omega).\]
 \end{rmk}
\noindent Suppose now that there exists an interval $I$ of regular frequencies such that $\omega \in I \mapsto k_0^\idomlap(\omega)$ is regular and monotonic (and therefore bijective) on $I$. Let us denote $J = k_0^\idomlap(I)$. This function has an inverse $\omega_0^\idomlap(k) : J \rightarrow I$ that is also regular and monotonic. To define the direction of propagation of the associated mode, it is classical to introduce the so-called {\it group velocity} that is given by 
	$$
	[\omega_0^\idomlap]'(k) = \frac{1}{[k_0^\idomlap]'(\omega)} \; \mbox{ for }\quad \omega = \omega_0^\idomlap(k).
	$$
	Note that as shown by \eqref{quasi-mode-temps-1D}, $\omega_0^\idomlap$ scales inversely with time, whereas $k$ scales inversely with position; consequently, the derivative $[\omega_0^\idomlap]'(k)$ has the dimension of a velocity.
We prove now that this quantity is directly linked to the density flux of the mode $Q_0^\idomlap(\omega)\coloneqq Q_0^\idomlap$ where $Q_0^\idomlap\coloneqq$ is defined in Corollary \ref{prop:flux_density_constant_wrt_s}. Using Lemma \ref{lem:link_q_s_w_s} and the expression \eqref{eq:uks} of $u^\idomlap_{s}$, a direct computation yields 
$$
Q^\idomlap_0(\omega)\coloneqq -4z\nu^\idomlap\Imag \Big(\mu_{\cutvec}\, \frac{d u^\idomlap_{s}}{d\xvi}\, \overline{u^\idomlap_{s}}\Big) (x) \quad (\text{that is independent of $x$ and $s$}),
$$
where the invariance w.r.t. $x$ is due to the fact that $u^\idomlap_{s}$ satisfies \eqref{modes1D_problem} and the one w.r.t. $s$ is given by Corollary \ref{prop:flux_density_constant_wrt_s}.
\begin{prop}
We have 
\[ [(k_0^\idomlap)'(\omega)]\,Q_0^\idomlap(\omega)=\omega\,\int_{\cellper^{\idomlap, 0}}\rho_p\,|U_0^\idomlap|^2. \] 
\end{prop}

\begin{dem}
Let us focus on $\idomlap=\droite$ to simplify the notations. Since $U_0^\droite$ satisfies \eqref{modes_problem}, one can show easily by \eqref{2Dmode_ter} that $V_0^\droite(\omega)\coloneqq V_0^\droite$ satisfies
\begin{equation}\label{eq:V_0_omega} 
	\left(\mathcal{A}_p (k_0^\droite(\omega))-\omega^2\right)V_0^\droite(\omega)=0. 
\end{equation}
where the operator $\xi\mapsto\mathcal{A}_p (\xi)$ is defined in \eqref{eq:op_per_k} and where we need to write explicitly the dependence of $V_0^\droite$ w.r.t. $\omega$ for the rest of the proof.
By differentiating this expression with respect to $\omega$, we obtain
\[ \left(\mathcal{A}_p (k_0^\droite(\omega))-\omega^2\right)\partial_\omega V_0^\droite(\omega)=-\left([(k_0^\idomlap)'(\omega)]\,\partial_\xi \mathcal{A}_p (k_0^\droite(\omega))-2\omega\right)V_0^\droite(\omega). \]
Since the operator $\mathcal{A}_p (k_0^\droite(\omega))$ is self-adjoint and by \eqref{eq:V_0_omega}, one has
\[ \int_{\cellper^{\droite, 0}}\left([(k_0^\idomlap)'(\omega)]\partial_\xi \mathcal{A}_p (k_0^\droite(\omega))-2\omega\right)\,\rho_p\,V_0^\droite(\omega)\, \overline{V_0^\droite(\omega)}=0, \]
which rewrites, using the expression of $\partial_\xi \mathcal{A}_p$
\begin{equation}\label{eq:res_fluxvol}
	  [(k_0^\droite)'(\omega)]\,\left[\cuti_2\, \text{Im}\int_{\cellper^{\droite, 0}}\mu_p\left(\Dt{}+\icplx k_0^\droite(\omega)\,\cuti_2 \right) V_0^\droite(\omega)\, \overline{V_0^\droite(\omega)}\right]=\omega\int_{\cellper^{\droite, 0}}\rho_p|V_0^\droite(\omega)|^2.
\end{equation}
It remains to show that in the left-hand side, the quantity inside the second bracket is nothing else but $Q_0^\droite$. For this, the trick is to multiply \eqref{eq:V_0_omega} by $y_2\,\rho_p\,\overline{V_0^\droite(\omega)}$ and integrate by parts to obtain
\begin{multline} \label{eq:ipp}
	\int_{\cellper^{\droite, 0}}\mu_p\left(\Dt{}+\icplx k_0^\droite(\omega)\cuti_2\right)V_0^\droite(\omega)\,\overline{\left(\Dt{}+\icplx k_0^\droite(\omega)\cuti_2\right)[y_2\,V_0^\droite(\omega)]}-\omega^2\,\int_{\cellper^{\droite, 0}}\rho_p\, y_2|V_0^\droite(\omega)|^2 \\ =\int_{\Sigmaper^{\droite,1}}\mu_p \left(\Dt{}+\icplx k_0^\droite(\omega)\cuti_2\right)V_0^\droite(\omega)\,\overline{V_0^\droite(\omega)}.
	\end{multline}
This implies in particular,  using that 
\[ 
\left(\Dt{}+\icplx k_0^\droite(\omega)\cuti_2\right)[y_2\,V_0^\droite(\omega)]=y_2\left(\Dt{}+\icplx k_0^\droite(\omega)\cuti_2\right)[y_2\,V_0^\droite(\omega)]+\theta_2 V_0^\droite(\omega)
\] 
by taking the imaginary part of \eqref{eq:ipp}
\begin{equation*}
\cuti_2\,\left[ \text{Im}\int_{\cellper^{\droite, 0}}\mu_p\left(\Dt{}+\icplx k_0^\droite(\omega)\,\cuti_2 \right) V_0^\droite(\omega)\, \overline{V_0^\droite(\omega)}\right]=\text{Im}\int_{\Sigmaper^{\droite,1}}\mu_p \left(\Dt{}+\icplx k_0^\droite(\omega)\cuti_2\right)V_0^\droite(\omega)\,\overline{V_0^\droite(\omega)},
\end{equation*}
where the right-hand side is nothing else but $Q_0^\droite(\omega)$, by definition of $V_0^\droite(\omega)$.
\end{dem}
\begin{rmk}
	\noindent Let us go back to the solutions of \eqref{eq:whole_line_problem}. On each interval $I^\idomlap$ with $\idomlap\in\{\gauche,\droite\}$, any solution of \eqref{eq:whole_line_problem} is necessarily a linear combination of $(u^\idomlap_{s^\idomlap}, \overline{u^\idomlap_{s^\idomlap}})$ with $s^\idomlap=a^\idomlap \theta_1$. With the limiting absorption principle, we select the solution that is proportional to $u^\idomlap_{s^\idomlap}$ on $I^\idomlap$ (see \eqref{eq:sol}-\eqref{eq:recons_1d0_sseps}). Moreover, it has a positive group velocity in $I^\idomlap$ (in the sense that $[k_0^\idomlap]'(\omega)>0$).
\end{rmk}
\section{Numerical results}\label{H1Dlap:sec:results}
\begin{figure}[ht!]%[H]
	\centering
	\includegraphics[page=6]{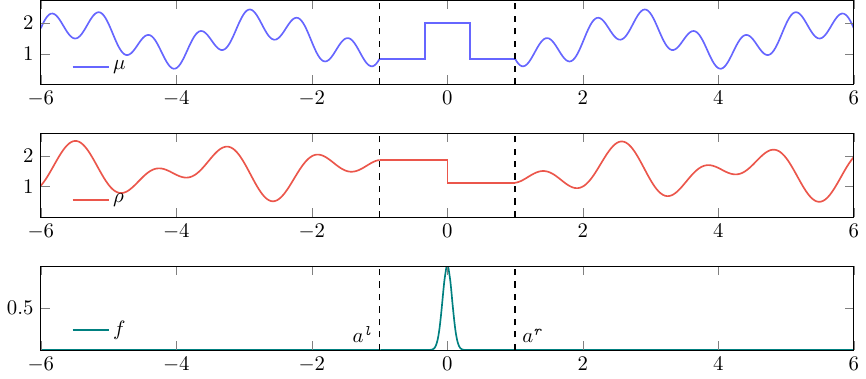}
	\caption{The locally perturbed quasiperiodic coefficients $\mu$ and $\rho$, and the source term $f$. \label{H1Dlap:fig:coefficients_mu_rho}}
\end{figure}

The procedure developed in the previous sections is illustrated through a series of numerical results. Our goal is to compute the physical solution of Problem \eqref{eq:whole_line_problem}, where the coefficients $\mu$ and $\rho$ coincide with quasiperiodic functions outside an interval $(I^\interieur) \coloneqq  (-1, 1)$, with periodic lifts
\[
  \displaystyle
  \spforall \yv \in \R^2, \quad \mu_p(\yv) = 1.5 + \cos (2\pi\yvi_1)\, \cos (2\pi\yvi_2) \quad \textnormal{and} \quad \rho_p(\yv) = 1.5 + 0.5\, \sin (2\pi\yvi_1) + 0.5\, \sin (2\pi\yvi_2).
\]
The cut vector is $\cutvec = (\cos\pi/3, \sin\pi/3)$, so that the ratio $\delta \coloneqq  \cuti_1 / \cuti_2 = 1/\sqrt{3} \in \R \setminus \Q$ is not a Liouville number. Indeed, since $\delta$ is algebraic, its irrationality measure is $m(\delta) = 2$ (see \cite{roth1955rational}). Inside $(I^\interieur)$, the local perturbations are piecewise constant functions represented in Figure \ref{H1Dlap:fig:coefficients_mu_rho}. The source $f$ is the cut-off function
\[\displaystyle
\spforall \xvi \in \R, \quad	f(\xvi) = \exp \left( 100\, \big(1 - 1 / (1 - \xvi^2)\big) \right)\; \chi_{(-1, 1)},
\]
also represented in Figure \ref{H1Dlap:fig:coefficients_mu_rho}. Finally, the impedance is set to $z = \omega$.

\vspspe
Given a frequency $\omega$, we solve the problems (\ref{cellpb0_sseps}, \ref{cellpb1_sseps}) posed in the cell $\cellper^{\droite, 0} \coloneqq  (0, 1)^2$ using Lagrange finite elements of order $1$ with $h = 5 \times 10^{-2}$, see \cite{amenoagbadji2023wave} for more details on the discretization aspects. % 
We then compute discrete local RtR operators defined by \eqref{eq:def_RtR_loc_lim} and deduce the Riccati point spectrum \eqref{eq:Ric_spec_sseps}. %
We recall from Section \ref{sec:algo_sseps} that if some element of the Riccati point spectrum is inside the unit circle, then $\omega$ is an \emph{evanescent frequency}. %
In this case, the fundamental eigenpair is the pair $(\lambda, {(\varphi_\lambda, \psi_\lambda)})$ such that $|\lambda| < 1$ and $\varphi_\lambda$ has a zero winding number. %
On the other hand, if the Riccati point spectrum is included in the unit circle, then $\omega$ is a \emph{propagative frequency}. %
In this case, we look for the eigenpairs $(\lambda,{(\varphi_\lambda, \psi_\lambda)})$ with a positive energy flux (see Proposition \ref{thm:sigmapm}), and amongst them, the fundamental eigenpair is the one such that $\varphi_\lambda$ has a zero winding number. It is worth noting that all the eigenpairs may have a zero energy flux, in which case we could not isolate the fundamental eigenpair. We refer to such a value of $\omega$ as a \emph{zero flux frequency}. If $\omega$ is not a zero flux frequency, then the fundamental eigenpair allows to construct the propagation operator via (\ref{eq:weight_p}, \ref{op_P_lim}), the scattering operator via (\ref{eq:weight_s}, \ref{op_S_lim}), and finally the RtR coefficient associated with the half-line problem set on $\R_+$. Repeating this procedure allows to deduce the RtR coefficient associated with the half-line problem set on $\R_-$. Finally, we can solve the problem \eqref{eq:int_pb_RtR} set in the interior domain $(I^\interieur)$ if it is well-posed, and deduce the solution $u$ of \eqref{eq:whole_line_problem}. 

\begin{figure}[ht!]
	\centering
	\includegraphics[page=7]{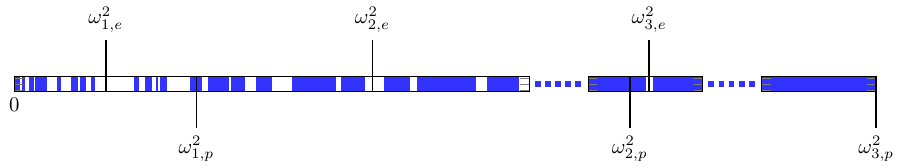}
  \caption{Position of the evanescent frequencies and of the propagative frequencies inside the spectrum of the quasiperiodic differential operator\label{fig:spec_diff_operator}}
\end{figure}

\vspspe
% La procédure décrite dans le paragraphe précédent permet de distinguer trois classes de fréquences, leading to different behaviour that we shall illustrate. We use three 
In what follows, we consider three evanescent frequencies: $\omega_{1, e} = 4$, $\omega_{2, e} = 7.912$, $\omega_{3, e} = 11.647$ and three propagative frequencies: $\omega_{1, p} = 5.642$, $\omega_{2, p} = 11.5$, $\omega_{3, p} = 20$. These frequencies have been obtained from preliminary numerical tests, and are represented in Figure \ref{fig:spec_diff_operator} with respect to the spectrum of the quasiperiodic differential operator $\mathcal{A}_\cutvec$ (see \eqref{eq:op_quasiper}), which is also obtained numerically. Note that the computation of $\sigma(\mathcal{A}_\cutvec)$ is not necessary for the method, and is only performed for illustration purposes. 

\begin{figure}[ht!]% [H]
  \centering
	\includegraphics[page=8]{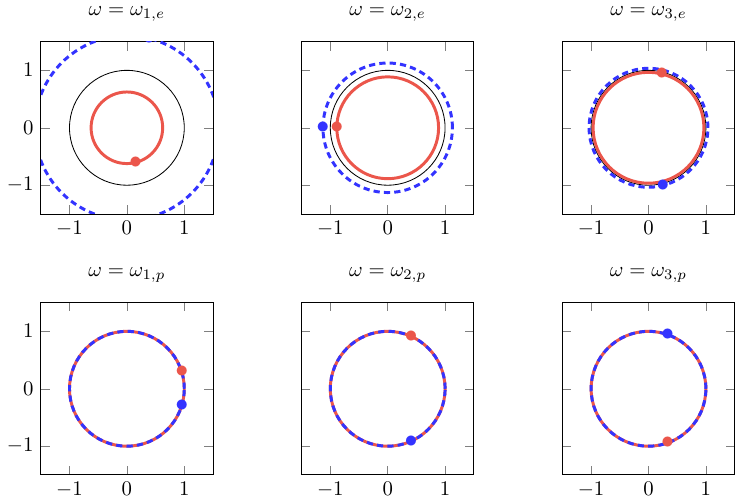}
	\caption{Riccati point spectrum reconstructed analytically from the fundamental eigenvalue (red filled circle) and its inverse (blue filled circle) for evanescent (first row) and propagative (second row) frequencies. The red circle corresponds to the spectrum of the propagation operator.\label{fig:riccati_spectrum}}
\end{figure}

\vspace{0\baselineskip} \noindent
In Figure \ref{fig:riccati_spectrum}, for each of the propagative and evanescent frequencies introduced above, we represent the fundamental eigenvalue, from which we reconstruct the Riccati point spectrum using Propositions \ref{thm_spectreaugmented} and \ref{thm_spectreaugmented_2_sseps}. %

% \vspspe

\begin{figure}[htp!]
	\centering
	\includegraphics[page=9]{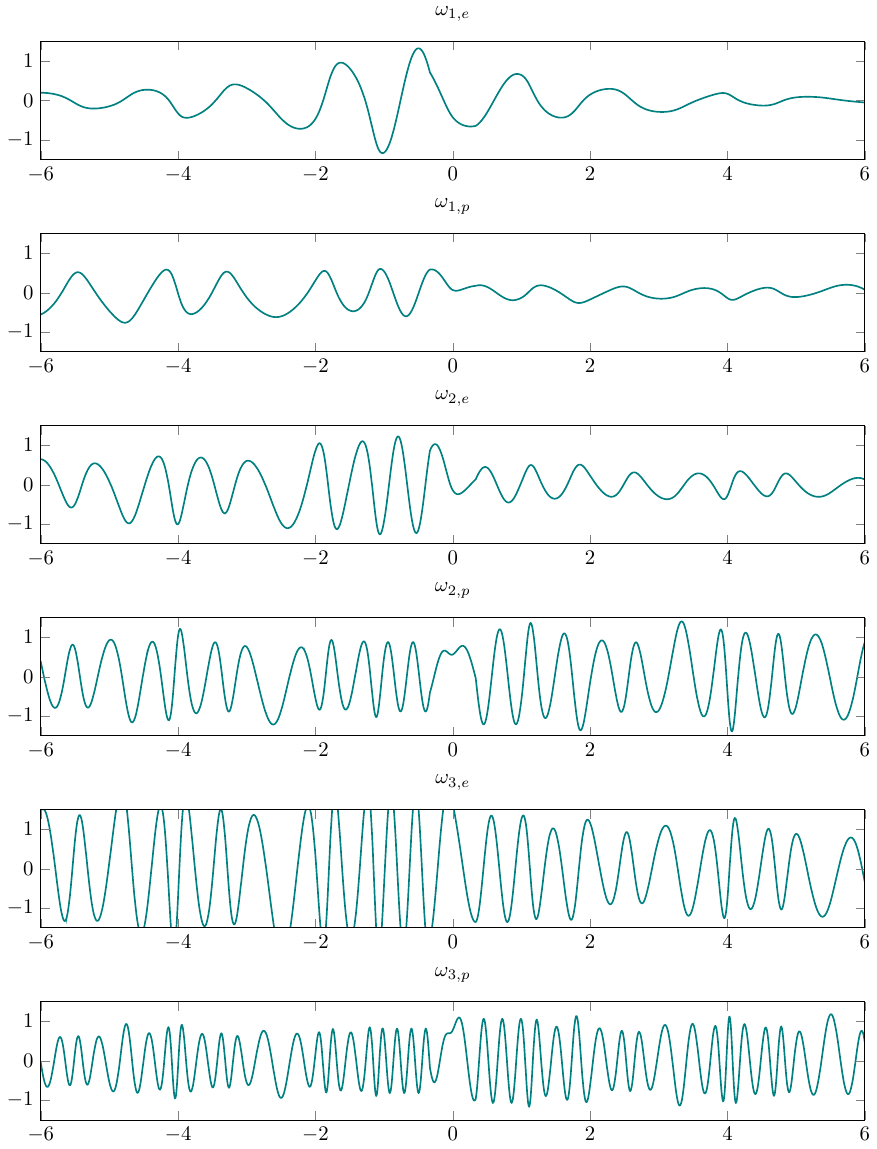}
	\caption{Real part of the solution of \eqref{eq:whole_line_problem} computed using the quasi-1D approach, with $\mathbb{P}^1$ Lagrange finite elements and $h = 5 \times 10^{-2}$, and for different evanescent frequencies and propagative frequencies.\label{H1DLAP:fig:whole_line_omega}}
\end{figure}

\begin{figure}[ht!]% [H]%[H]
	\centering
	\includegraphics[page=10]{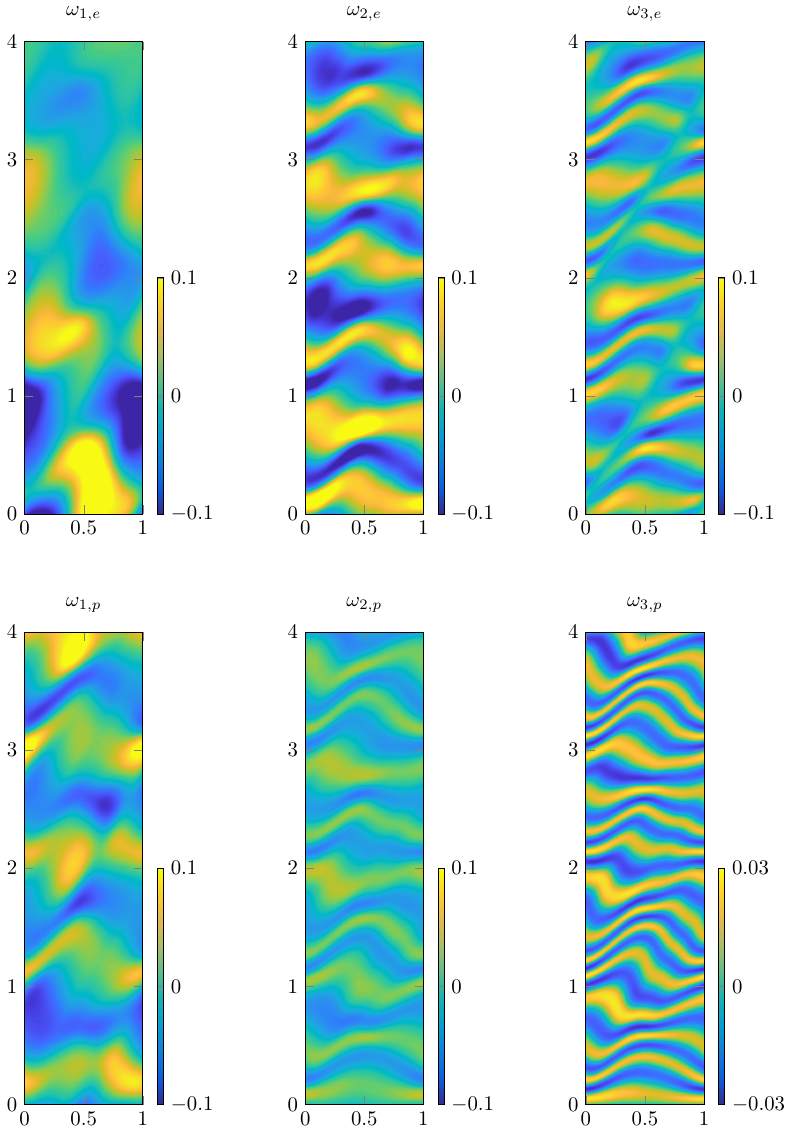}
	\caption{Real part of the half-guide solution computed using the quasi-1D approach, with $\mathbb{P}^1$ Lagrange finite elements and $h = 5 \times 10^{-2}$, and for different evanescent frequencies (first row) and propagative frequencies (second row) (Note the change of color scale for $\omega_{3, p}$).\label{H1DLAP:fig:half_guide_omega}}
\end{figure}

\paragraph*{The half-guide solution and the solution of \eqref{eq:whole_line_problem}.}
For the evanescent frequencies $(\omega_{1, e}, \omega_{2, e}, \omega_{3, e})$ and propagative frequencies $(\omega_{1, p}, \omega_{2, p}, \omega_{3, p})$, we represent the solution $U^\droite$ of the limit 2D periodic half-guide problem \eqref{RobinHalfguideproblems_sseps} in Figure \ref{H1DLAP:fig:half_guide_omega}, and the solution $u$ of \eqref{eq:whole_line_problem} in Figure \ref{H1DLAP:fig:whole_line_omega}. To compute these solutions, we choose a boundary data $\varphi = 1$. A first seemingly surprising phenomenon is that $U^\droite$ becomes smaller as $\omega$ increases (note that in Figure \ref{H1DLAP:fig:half_guide_omega}, the color scale has been adjusted for $\omega_{3, p}$). This is due to our choice of impendance $z = \omega$, which increases with $\omega$. In fact, by formally making $z$ tend to $\infty$ in the Robin boundary condition $\mathcal{R}^\droite_+ U^\droite = 1$, we see that the trace of $U^\droite$ on $\Sigma^{\droite, 0}$ becomes smaller as $z$ increases. It can be observed that for the evanescent frequencies, $U^\droite$ (\resp $u$) decays when $\yvi_2 \to +\infty$ (\resp when $\xvi \to \pm \infty$). This corresponds to an exponential decay, whose rate is linked to the modulus of the fundamental eigenvalue $|\lambda^\idomlap_0|$ (see Proposition \ref{prop:exp_decay_evanescent}). In particular, this decay can be seen clearly for $\omega_{1, e}$, whereas it is more difficult to observe for $\omega_{2, e}$ and $\omega_{3, e}$, since in these cases, $|\lambda^\idomlap_0|$ is close to $1$ as Figure \ref{fig:riccati_spectrum} shows. Moreover, the decay is absent for propagative frequencies. Note also that as expected, the solutions oscillate more as $\omega$ increase.

% \clearpage

\begin{figure}[H]
	\centering
	\includegraphics[page=11]{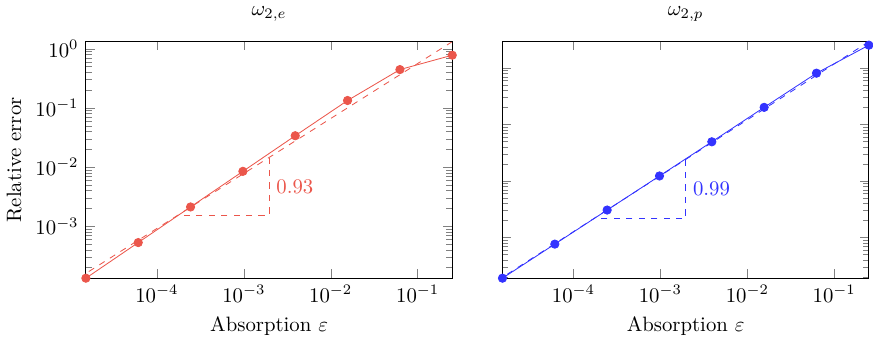}
	\caption{Relative error in $H^1$ norm between the solutions of \eqref{eq:whole_line_problem} and\eqref{eq:whole_line_problem_eps} for different values of $\omega$.\label{H1DLAP:fig:errors_veps}}
\end{figure}

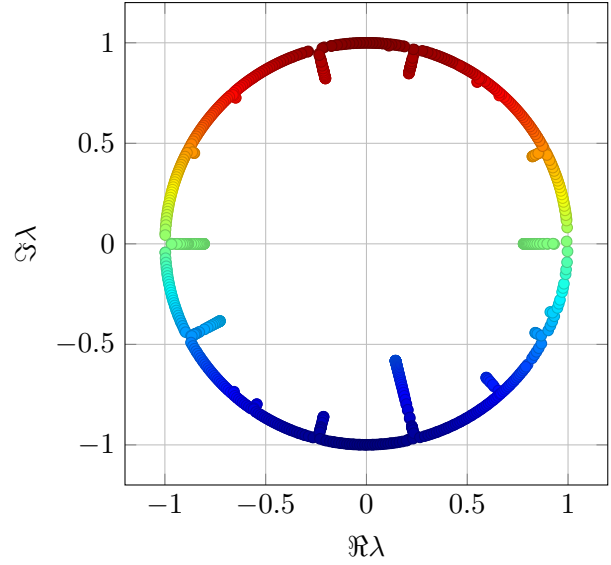
\begin{figure}[ht!]%[H]
	\begin{center} 
	\begin{tikzpicture}
	\begin{axis}[
      xlabel={$\Re \lambda$},
      ylabel={$\Im\lambda $},
      grid=both,
      width=8cm, height=8cm,
      colormap/jet,          % rainbow-like colormap
      point meta=y           % color depends on y value
  ]
	\addplot[
      only marks,
      scatter,
      scatter src=y          % use y-coordinate for color
    ] table [col sep=comma] {./eigval.txt};
	\end{axis}
\end{tikzpicture}
	\caption{Representation of the fundamental eigenvalue in the complex plane with respect to $\omega\in[0,10]$, each color corresponds to a different value of the frequency. \label{fig:lambda_omega}}
\end{center}
\end{figure}

\begin{figure}[ht!]%[H]
	\begin{center} 
		\begin{subfigure}{0.45\textwidth}
	\hspace{-1cm}\includegraphics[page=12]{tikz-picture.pdf}
\caption{The fundamental dispersion curve}
  \end{subfigure}
  \hfill\hfill
 \begin{subfigure}{0.45\textwidth}
     \centering
\includegraphics{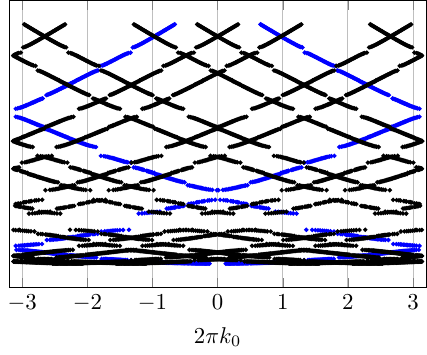}
\caption{A set of 5 dispersion curves}
  \end{subfigure}
\caption{Representation of the dispersion curves of the problem. \label{fig:disp_curve}}\end{center}
\end{figure}

\paragraph*{Convergence with respect to the absorption.}
In order to illustrate the estimates shown in Section \ref{sec:LAP}, we represent the relative error
\begin{equation*}
	\veps \mapsto \frac{\|u_\veps - u\|_{H^1(I^\interieur)}}{\|u\|_{H^1(I^\interieur)}}
\end{equation*}
where $u$ is the solution of \eqref{eq:whole_line_problem}, and where $u_\veps$ denotes the solution of \eqref{eq:whole_line_problem_eps}. This relative error is represented in Figure \ref{H1DLAP:fig:errors_veps} for $\omega_{2, e}$ and $\omega_{2, p}$. As expected, a linear convergence is observed. 
{\paragraph*{Dispersion curves and group velocity} We study finally the behavior of the fundamental eigenvalue as a function of the frequency $\omega$. In Figure \ref{fig:lambda_omega}, we represent the location of the fundamental eigenvalue $\lambda^\droite_0(\omega)$ in the complex plane for different values of $\omega$, the dependence with respect to $\omega$ being represented thanks to the color bar. We notice first that the map $\omega\mapsto\lambda^\droite_0(\omega)$ is continuous and piecewise regular. Moreover, an interesting behavior occurs when considering an interval $(\omega^2_a,\omega^2_b)\not\subset\sigma(A)$ with endpoints $\omega^2_a,\omega^2_b\in\sigma(A)$: the fundamental eigenvalue enters the unit circle and sweeps through a segment of a straight line passing through the origin. More precisely, it seems that
\[\lambda_0^\droite(\omega)=e^{-\alpha_0(\omega)}\lambda_0^\droite(\omega_a),\quad \forall\omega\in(\omega_a,\omega_b)\quad\text{with}\;\alpha_0(\omega)\in\R^+.\]
Furthermore, numerical experiments suggest that $\lambda_0^\droite$ may enter the unit circle from any location on the unit circle. For comparison not that for $q$-periodic coefficients, one can show that the equivalent of $\lambda_0^\droite$ can enter the unit circle only from the points of the form $e^{\icplx p\pi/q}$ for $p\in\{0,\ldots,2q-1\}$.  This was not part of our theoretical results and are only numerical observations.
\\\\
Next, in Figure \ref{fig:disp_curve}-(left), we represent the dispersion curves as defined in Section \ref{sub:disp_curves} in a neighborhood of all regular frequencies. More precisely, we represent $\xi\mapsto \omega_0^\droite(\xi)$ where $\omega_0^\droite$ is the inverse of the function $k_0^\droite(\omega)$ with $\lambda_0^\droite=e^{2\icplx\pi \theta_2 k_0^\droite(\omega)}$ whenever $k_0^\droite(\omega)$ is well-defined, regular and monotonic. This corresponds to the dispersion curve of the rightgoing modes $u^\droite_s$. We represent also $\xi\mapsto \omega_0^\droite(-\xi)$. By Remark \ref{rmk:basedesol}, this corresponds to the dispersion curve of the leftgoing modes $u^\gauche_s$ up to a translation $p\delta$, the (left) dispersion curve. 
\\\\
Finally, in Figure \ref{fig:disp_curve}-(right), we represent several dispersion curves corresponding to several 2D rightgoing modes whose wavenumbers have the form $2\pi (k_0^\droite(\omega)- k \delta),\;k\in\Z$ (see Remark \ref{rmk:2dmodes_temps}). The associated dispersion curves write as $\xi\mapsto \omega_k^\droite(\xi):=\omega_0^\droite(\xi)-2\pi k \delta.$ This also illustrates Proposition \ref{thm:spectrunqp_k}.
}

	%

% \newpage

%: the values $\omega_{1, p} = \sqrt{16}$, $\omega_{2, e} = \sqrt{81}$ and $\omega_{3, e} = \sqrt{338}
\section{Conclusions and extensions}\label{sec:conclusions}
	{This paper contributes to the analysis of the limiting absorption process for the Helmholtz equation in quasiperiodic media and proposes a numerical method for computing the corresponding physical solution.
	The present study can be extended without difficulty to the case where 
	\begin{itemize}
	\item $\mu_\cutvec$ and $\rho_\cutvec$ are quasiperiodic functions of order $\ord > 2$: the periodic extensions $\mu_p, \rho_p: \R^\ord \to \R$ would depend on $\ord$ variables and the cut direction $\cutvec$ would be a vector of dimension $\ord$.
	\item $\mu$ (\resp $\rho$) coincides with a different quasiperiodic function in each of the half-lines $I^\gauche \coloneqq  (-\infty, a^\gauche)$ and $I^\droite\coloneqq (a^\droite, +\infty)$:
	\[
		\spforall \idomlap \in \{\gauche, \droite\}, \quad \spforall \xvi \in I^\idomlap, \quad \mu(\xvi) = \mu_p^\idomlap(\cutvec^\idomlap\, \xvi) \quad \text{and} \quad \rho(x) = \rho_p^\idomlap(\cutvec^\idomlap\, \xvi),
	\]
	where $\mu_p^\idomlap, \rho_p^\idomlap \in \mathscr{C}^0_{\textit{per}}(\R^{\ord^\idomlap})$ with $\ord^\idomlap \geq 2$, and $\cutvec^\idomlap \in \R^{\ord^\idomlap}$ for $\idomlap \in \{\gauche, \droite\}$.% are irrational vectors.
		\end{itemize}
		There still remain theoretical open questions that concern the different assumptions we have made to achieve the study. Indeed, it would be interesting to characterize more precisely the regular frequencies. It is clear, from \cite[Section XIII-6]{reed1978iv}, that an interval of regular frequencies (if it exists) is included in the absolute continuous part of the spectrum: this is since we have shown that the limiting absorption principle holds for such frequencies. Can we find frequencies of the absolute continuous part of the frequencies that are not regular? Which assumption is not satisfied? What happens for the singular part? It is worth mentioning that for now, we have encountered in our numerical study, only frequencies that are either evanescent or propagative for which Assumption \ref{assumption} seems to be satisfied. The energy flux may vanish but only for the frequencies that are in one end of an interval of propagative frequencies. }

\section*{Aknowledgements}
		The authors would like to thank the Isaac Newton Institute for Mathematical Sciences, Cambridge, for support and hospitality during the programme Mathematical theory and applications of multiple wave scattering where work on this paper was undertaken. This work was supported by EPSRC grant EP/R014604/1. P. A. is supported by Simons Foundation Math + X Investigator Award\#376319 (Michael I. Weinstein).
\appendix

\section{The DtN approach is not adapted to limiting absorption} \label{sec_DtN}
{We show in this section that the DtN method cannot be used for the limiting absorption principle since we shall see that we cannot pass to the limit when $\veps$ tends to $0$  in the cell problems \eqref{eq:local_cell_problem}-\eqref{eq:local_cell_BC}, for frequencies above a threshold.

\vspspe
Let us introduce the unbounded self-adjoint and positive operator
\[ \begin{array}{l}
\dsp D(\mathcal{A}^{\dirbc}_\cutvec)= \Big\{ E \in H^1_\cutvec(\cellper^0), \ \mu_p \,\Dt{}E \in H^1_\cutvec(\cellper^{0}),\;E\mbox{ is periodic w.r.t. }y_1, \; E|_{\Sigmaper^{0}} =E|_{\Sigmaper^{1}} = 0 \Big\}, \\[5pt]
\dsp \mathcal{A}^\dirbc_\cutvec E= -\frac{1}{\rho_p}\Dt{} \big( \mu_p \, \Dt{} E{}\big). 
\end{array}
\]
We prove in the rest of this section, that, if $\mu_p,\rho_p$ are $\mathscr{C}^1$ functions, the spectrum of $\mathcal{A}^{\dirbc}_\cutvec$ contains a semi-infinite interval: 
\begin{equation}\label{eq:semi_int_spec}
    \spexists\omega_*>0,\quad [\omega^2_*,+\infty)\;\subset \sigma(\mathcal{A}^{\dirbc}_\cutvec).
\end{equation} 
This proves that for $\omega^2$ lying in this interval, the resolvent $[\mathcal{A}^{\dirbc}_\cutvec-(\omega^2+\icplx\,\eps)]^{-1}$ does not have a limit when $\veps$ tends to $0$, which implies that the cell solutions $\Edir{, \droite}_{\veps, \ell},\;\ell\in\{0,1\}$, satisfying \eqref{eq:local_cell_problem}-\eqref{eq:local_cell_BC}, seen as operators from $L^2(\Sigmaper^\droite)$ to  $H^1_\cutvec(\cellper^0)$, do not have a limit neither.}
% \begin{rmk}
% A natural question could be: what happens if  $\mu_p$ and $\rho_p$ are not $\mathscr{C}^1$ functions? We did not investigate this question since our objective was to justify why the DtN approach is not relevant in general and why another approach (developped in Section ...) is necessary.
% \end{rmk}

\begin{figure}[ht!]
	\centering
		\centering
		\includegraphics[page=3]{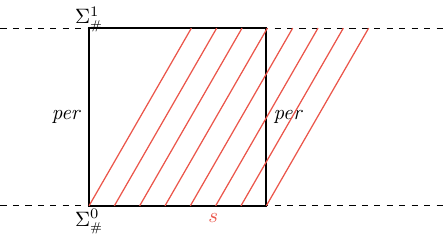}
	\caption{Illustration of \eqref{eq:op_quasiper2Dtheta_D}}\label{fig:correspondence1D2D}
\end{figure}

\vspspe
By definition of the differential operator $\Dt{}$ and by periodicity, this operator can be seen, roughly speaking, as a \textquote{concatenation} of 1D differential operators along the lines 
\begin{equation} \label{deflines}
	\big\{ (s + \cuti_1\vts \xvi, \cuti_2\vts \xvi),\ \ \xvi \in [0, L_\cutvec] \big\} \quad \textnormal{with} \quad L_\theta= 1/\cuti_2.
\end{equation}
A quite intuitive illustration of the result is provided by Figure \ref{fig:correspondence1D2D}. 
More precisely, let us define the unbounded self-adjoint and positive operator $\mathcal{A}^{\dirbc}_s$ defined by% in $L^2\big((0, L_\cutvec), \rho_{s, \cutvec} \, dx \big)$ defined by
\[
D(\mathcal{A}^{\dirbc}_s)= \Big\{ u \in H^1_0(0,L_\cutvec) \ / \ \mu_{s, \cutvec} \, \frac{d u}{d \xvi} \in H^1(0,L_\cutvec) \Big\}, \qquad \mathcal{A}^{\dirbc}_s u= -\frac{1}{\rho_{s, \cutvec}}\frac{d}{d \xvi} \Big( \mu_{s, \cutvec} \, \frac{d u}{d \xvi} \Big), %\quad \spforall u \in D(\mathcal{A}^{\emph{cell}}_s).
\]
where $\mu_{s,\cutvec}$ and $\rho_{s,\cutvec}$ are defined in \eqref{defrhosmus}.
From the family $\mathcal{A}^{\dirbc}_s, s \in [0,1)$, we can build a self-adjoint operator in $L^2_{\textit{per}}\big(\R, L^2(0,L_\cutvec)\big)$, as
\begin{equation}
\label{eq:op_quasiper2Dtheta_D}
\widetilde{\mathcal{A}}^{\dirbc}_\cutvec = \int^\oplus_{(0,1)}	\mathcal{A}^\dirbc_{s} \, ds
\end{equation}
meaning that (see \cite[Section XIII.16]{reed1978iv} for instance for direct integrals of operators), 
\begin{equation}
\label{eq:op_quasiper2Dtheta2_D}
U \in D(\widetilde{\mathcal{A}}^{\dirbc}_\cutvec) \quad \Longleftrightarrow \quad \spforall s\in[0,1], \quad U(s,\cdot) \in D(\mathcal{A}_{s}^\dirbc) \quad \textnormal{and} \quad  [\widetilde{\mathcal{A}}^{\dirbc}_\cutvec\, U] (s,\cdot)  = \mathcal{A}_{s}^\dirbc [U(s,\cdot)].
\end{equation} 
It is known, see \cite[Theorem XIII.85]{reed1978iv}, that the spectra of $	\widetilde{\mathcal{A}}^{\dirbc}_\cutvec $ and $\mathcal{A}_{s}^\dirbc$ are related by 
\begin{equation}\label{floquet_D}
\sigma\big(\widetilde{\mathcal{A}}^{\dirbc}_\cutvec\big)=\bigcup_{s \in [0,1]} \sigma\big(\mathcal{A}_{s}^\dirbc\big).
\end{equation}
Moreover, the operators $\widetilde{\mathcal{A}}^{\dirbc}_\cutvec$ and ${\mathcal{A}}^{\dirbc}_\cutvec$ are equivalent in $L^2_{\textit{per}}(\R^2)$:
$$
\widetilde{\mathcal{A}}^{\dirbc}_\cutvec = \widetilde{\mathcal{U}}_\cutvec \, {\mathcal{A}}^{\dirbc}_\cutvec\, \widetilde{\mathcal{U}}_\cutvec^{-1},
$$
where $\widetilde{\mathcal{U}}_\cutvec$ is the invertible operator from $L^2_{\textit{per}}(\R\times(0,1))$ into $L^2_{\textit{per}}(\R\times (0,L_\cutvec))$ defined by 
\begin{equation*}% \label{defStheta}
 \spforall U \in L^2_{\textit{per}}(\R\times(0,1)), \quad (\widetilde{\mathcal{U}}_\cutvec\, U) (s,x) = U (s + \theta_1 \, x, \theta_2 \, x).
\end{equation*} 
As a consequence, we have 
\begin{equation}\label{lem_inclusion_spectra}
  \sigma \big({\mathcal{A}}^{\dirbc}_\cutvec\big)=\bigcup_{s \in [0,1]} \; \sigma(\mathcal{A}^{\dirbc}_s).
\end{equation}
Since for each $s \in [0, 1]$, $\mathcal{A}^{\dirbc}_s$ is positive, self-adjoint and with compact resolvent,  its spectrum $\sigma(\mathcal{A}^{\dirbc}_s)$ consists of distinct simple eigenvalues (by Sturm-Liouville theory \cite{titchmarsh2011}) 
$$
\sigma(\mathcal{A}^{\dirbc}_s) = \big\{ 0 < \lambda_1(s) < \lambda_2(s) < \cdots \lambda_n(s)<  \cdots\; \}, \quad \lim_{n \to +\infty} \lambda_n(s) \rightarrow + \infty.% \; ( n \rightarrow + \infty)
$$
Moreover, from the regularity of $\rho_p$ and $\mu _p$, we deduce that the maps $s \in [0,1] \mapsto \rho_{s, \cutvec} \in \mathscr{C}^0(0, L_\cutvec)$ are continuous. As a consequence, the functions $s \mapsto\lambda_n(s)$ are continuous and we have
\begin{equation} \label{forbidden-set}
	\bigcup_{s \in [0,1]} \; \sigma(\mathcal{A}^{\dirbc}_s)= \bigcup_{n \geq 1} [a_n ,b_n] \quad \textnormal{where} \quad [a_n ,b_n] \coloneqq  \lambda_n( [0,1]), \quad a_n > 0, \quad a_n \rightarrow + \infty. 
	\end{equation}
The intervals $[a_n ,b_n]$ may overlap and will generically do so for $n$ large enough, as explained in the following result. We introduce the
function $\xvi \mapsto c_s(\xvi)$, namely the wave velocity along the line $[0, L_\cutvec]$, and its harmonic mean 
\begin{equation}
	\spforall \xvi \in [0,L_\theta], \quad c_s(x) = \left( \frac{\mu_{s, \cutvec}(x)}{\rho_{s, \cutvec}(x)}\right)^{1/2}\quad \textnormal{and} \quad \overline{c}_s=\left[ \frac{1}{L_\cutvec} \, \int_0^{L_\cutvec} 	\frac{dx}{c_s(x)}\right]^{-1}.
\end{equation}
\begin{lem}\label{lem:sp_overlap}
	As soon as $s\mapsto \overline{c}_s $ is not constant, $b_n \geq a_{n+1}$ for $n$ large enough.
\end{lem}
\begin{dem} 
    This is essentially a matter of applying known results from the spectral theory of Sturm-Liouville operators 
	with smooth coefficients (see \cite{titchmarsh2011} or \cite[Chapter 4]{eastham1973spectral} for instance). Indeed, for all $s\in[0,1]$,  $\lambda_n(s)$ behaves asymptotically as
	\begin{equation} \label{asymptotics0_ok} 
	\lambda_n(s) = \Big( \frac{n \, \pi\, \overline{c}_s }{L_{\cutvec}} \Big)^2 + \mathcal{O}(n) \quad \mbox{uniformly in $s$} .
	\end{equation} 
	As a consequence, with $\overline{c}_+= \sup_{s\in(0,1)} \overline{c}_s$ and $\overline{c}_-= \inf_{s\in(0,1)} \overline{c}_s$, the bounds $a_n \leq b_n$ of the interval $\lambda_n( [0,1])$ admit the asymptotic expansion
	\begin{equation*} 
	a_n = \Big( \frac{n \, \pi\, \overline{c}_-}{L_{\cutvec}} \Big)^2 + \mathcal{O}(n) \quad \textnormal{and} \quad b_n= \Big(\frac{n \, \pi \, \overline{c}_+}{L_{\cutvec}} \Big)^2 + \mathcal{O}(n), 
	\end{equation*} 
so that if $\overline{c}_-<\overline{c}_+$, we can conclude. We give some hints of the proof of \eqref{asymptotics0_ok}, for the ease of the reader. 
  We can rewrite the eigenvalue 
	equation for $\mathcal{A}^{{\dirbc}}_s,\,s\in[0,1]$ as follows
\begin{equation} \label{regEVP_ok} 
	-\frac{1}{\rho_{s, \cutvec}}\frac{d}{d \xvi} \Big( \mu_{s, \cutvec} \, \frac{d u}{d \xvi} \Big) = \lambda \, u \quad \Longrightarrow \quad -c_s^2 \, \frac{d^2 u}{d \xvi^2}- \beta_s \, \frac{d u}{d \xvi} = \lambda \, u \quad \textnormal{with} \quad \beta_s \coloneqq  \frac{\mu_{s, \cutvec}'}{\rho_{s, \cutvec}}
\end{equation} 
where we have used the $\mathscr{C}^1$--regularity of $\mu_{s, \cutvec}$.
	Then, we apply a change of variables (known as the \emph{Prüfer transform}) to get rid of the variable coefficient in factor of the second order derivative 
	$$
	\spforall x \in [0, L_\cutvec], \quad t = {\cal T}(x)= \int_0^x \frac{d \xi}{c_s(\xi)} \in [0, L_\cutvec/ \overline{c}_s] \quad (\mbox{travel time variable}). 
	$$
	Setting $U(t) \coloneqq  u\big({\cal T}^{-1}(t) \big)\ \Longleftrightarrow\ u(x) = U\big({\cal T}(x)\big)$, we observe that
	$$
	\frac{d u}{d\xvi} (x) = c_s(x)^{-1} \, \frac{d U}{d t}\big({\cal T}(x)\big) \qquad \textnormal{and} \qquad \frac{d^2 u}{d\xvi^2} (\xvi) = c_s(x)^{-2}\, \big[ \, \frac{d^2 U}{d t^2}\big({\cal T}(x)\big) - c_s'(x) \, \frac{d U}{d t}\big({\cal T}(x)\big) \, \big].
	$$
		Substituting this into \eqref{regEVP_ok}, we get an eigenvalue equation for $U(t)$
	\[
		\left|
			\begin{array}{l@{\ =\ }l}
				\displaystyle- \frac{d^2 U}{d t^2} + a_s(t) \, \frac{d U}{d t} & \lambda \, U \quad \textnormal{ in } [0, L_\cutvec/ \overline{c}_s], 
				\ret
				U(0) = U\big( L_\cutvec/ \overline{c}_s\big) & 0,				
			\end{array}
		\right.
		\quad \textnormal{with} \quad a_s(t) \coloneqq  \Big[c_s'-\frac{\beta_s}{c_s}\Big] \big( {\cal T}^{-1}(t) \big),
	\]
	which can be viewed as a (compact) perturbation of the Dirichlet eigenvalue problem for the operator $U \mapsto -\, U''$ in the interval $[0, L_\cutvec/ \overline{c}_s]$. For this unperturbed operator, the eigenvalues are 
	$$
\{{(n\, \pi\,\overline{c}_s}/ {L_\cutvec})^2, \quad n\geq 1\}
	$$
	which provide the dominant term in the asymptotics \eqref{asymptotics0_ok}. One then concludes by a perturbation argument.
	\end{dem}
\noindent It follows from \eqref{forbidden-set} and the previous lemma that setting $\omega_*^2= {a_N}$ with $N= \min \{ n ,\;a_{k+1} \leq b_{k} \ \ \spforall k \geq n \big \}$, we have
\begin{equation} \label{lowerboundOmegaforb}
	 [ \, \omega_*^2, + \infty ) \; \subset\bigcup_{s \in [0,1]} \; \sigma(\mathcal{A}^{\dirbc}_s).
	\end{equation}
Relation \eqref{lem_inclusion_spectra} allows finally to deduce \eqref{eq:semi_int_spec}.

\begin{figure}[ht!]
	\centering
	\includegraphics[page=4]{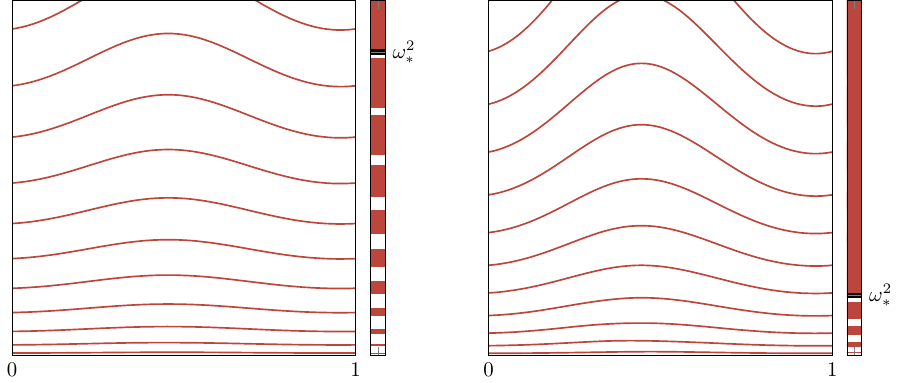}
	\caption{Curves $s \mapsto \lambda_n(s)$ for the medium given in \eqref{medium1_ok} for $\alpha = 1/2$ (left) and $1$ (right).}\label{fig:forbidden}
\end{figure}

\vspspe
Let us illustrate Lemma \eqref{lem:sp_overlap}, numerically, with the medium defined by
\begin{equation} \label{medium1_ok}
	\mu_p = 1, \quad \rho_p(\yv)	= 1.5 + \alpha \, (\sin 2\pi y_1 + \sin 2\pi y_2), \quad \cutvec = (\cos(\pi/3), \sin(\pi/3)).
\end{equation}
Note that a larger $\alpha$ gives a larger contrast $\overline{c}_+ / \, \overline{c}_-$ and thus a lower $\omega_*$ (see \eqref{lowerboundOmegaforb}). In Figure \ref{fig:forbidden}, we represent for $\alpha = 1/2$ and $\alpha =1$, the curves $s \mapsto \lambda_n(s)$ for $1 \leq n \leq 10$ (these eigenvalues have been computed numerically) as well as $\omega_*$. One clearly sees that $\omega_*$ decreases as $\alpha$ increases. 

\begin{rmk}% [Comparison with the DtN method for Helmholtz in waveguides]
	When applying the DtN method to the elliptic Helmholtz equation $-\dive \mu_p \nabla u - \rho_p\, \omega^2\, u = f$ with periodic coefficients (as done by \cite{fliss2009analyse}) in the non-absorbing case, one also has to exclude the spectrum of the associated Dirichlet cell operator. However, using compactness arguments, the spectrum is merely discrete in that case, and therefore, is not a limiting factor in practice. 
\end{rmk} 

\begin{rmk}
	In \cite[Section 5.5.1]{amenoagbadji2023wave}, we mentioned that the quality of the numerical results obtained with the DtN method was clearly deteriorated when  $\veps$ is smaller and smaller. This is linked to the issue presented in this section. 
\end{rmk}

\section{Invariance of spectra of differential operators}\label{app:spectrum}
The goal of this appendix is to prove Proposition \ref{thm:spectrunqp_s} (\resp Proposition \ref{thm:spectrunqp_k}) which states that the spectrum of $\mathcal{A}_{s, \cutvec}$ (\resp $\mathcal{A}_p (\fbvar)$) is independent of $s$ (\resp $\fbvar$). We first begin with the family of operators $\mathcal{A}_{s, \cutvec}$ with $\mathcal{A}_{0, \cutvec}={\cal A}_\theta$. Since $\mathcal{A}_{s+1, \cutvec}=\mathcal{A}_{s, \cutvec}$ for all $s\in\R$, it suffices to study $\mathcal{A}_{s, \cutvec}$ for $s\in[0,1)$. Furthermore, they are linked to one another as follows.
\begin{lem}\label{lem:unitary_equivalence}
   Let $\tau_a$ be the translation operator defined by $\tau_a \varphi = \varphi(\cdot + a)$. Then we have the identity
   \begin{equation}
     \displaystyle
     \tau^{-1}_{1/\theta_2}\, \mathcal{A}_{s+\delta, \cutvec}\, \tau_{1/\theta_2} = \mathcal{A}_{s, \cutvec}.
   \end{equation}
 \end{lem}
 \begin{dem}
   Since $\rho_p$ is $1$--periodic with respect to its second variable, for any $(s, x) \in \R^2$, we have
   \begin{align}
     \rho_{s + \delta, \cutvec}(x) &= \rho_p\Big(s + \frac{\theta_1}{\theta_2} + x\theta_1, x\theta_2\Big) \nonumber
     \\
     &= \rho_p\left(s + \Big(x + \frac{1}{\theta_2}\Big)\, \theta_1, \Big(x + \frac{1}{\theta_2}\Big)\, \theta_2\right) = \rho_{s, \cutvec}\Big(x + \frac{1}{\theta_2}\Big).
   \end{align}
   and we have similarly $\mu_{s + \delta, \cutvec}=\tau_{1/\theta_2}\, \mu_{s, \cutvec}$. We deduce the following
   \begin{equation}
     \spforall u \in D(\mathcal{A}_{s, \cutvec}), \quad \tau_{1/\theta_2}u \in D(\mathcal{A}_{s+\delta, \cutvec}) \quad \text{and} \quad  \mathcal{A}_{s+\delta, \cutvec}\, u\Big(\cdot + \frac{1}{\theta_2}\Big) = \big[\mathcal{A}_{s, \cutvec} u \big](\cdot + \frac{1}{\theta_2}\Big),
   \end{equation}
   which is the expected identity.
 \end{dem}

 \noindent
 Proposition \ref{thm:spectrunqp_s} which we recall and prove below is known to hold in a more general context, namely for Schrödinger operators with almost-periodic \cite{simon1982almost} and random \cite{pastur1980spectral} coefficients. In these references, the proof relies on the notion of spectral projections. Here for the sake of completeness, we propose an alternative proof which is specific to quasiperiodic functions, and which rely on more elementary objects.
 \begin{prop}\label{prop:spec_1D_2D_ligne_indt_s}
   If $\delta$ is irrational, then the spectrum of $\mathcal{A}_{s, \cutvec}$ does not depend on $s$.
 \end{prop}

 \noindent
 The proof of this result uses a perturbation result whose proof can be found in \cite[Theorem V.4.10]{kato2013perturbation}. The Hausdorff distance between two sets $\Lambda_1, \Lambda_2 \in \R$ is defined by
 \begin{equation}
   \text{dist}_{\mathscr{H}}(\Lambda_1, \Lambda_2) \coloneqq  \max \left\{\sup_{\lambda_1 \in \Lambda_1} \text{dist}(\lambda_1, \Lambda_2), \sup_{\lambda_2 \in \Lambda_2} \text{dist}(\lambda_2, \Lambda_1) \right\},
 \end{equation}
 where $\displaystyle\smash{\text{dist}(\xi, \Lambda) \coloneqq  \inf_{\lambda \in \Lambda} |\lambda - \xi|}$ represents the distance from $\xi \in \R$ to $\Lambda \subset \R$.
 
 \begin{lem}\label{lem:perturbation_spectrum}
   Consider a self-adjoint operator $B$ and a bounded symetric operator $T$. Then $B + T$ is also self-adjoint and $\text{dist}_{\mathscr{H}}(\sigma(B + T), \sigma(B)) \leq \|T\|$.
 \end{lem}
 
 \noindent
 Now, let us prove Proposition \ref{prop:spec_1D_2D_ligne_indt_s}.
 \begin{dem}[of Proposition \ref{prop:spec_1D_2D_ligne_indt_s}] %
  Since $\mathcal{A}_{s, \cutvec}$ is self-adjoint and positive, we can introduce the resolvent operator ${\cal R}_{s, \cutvec} \coloneqq  (\mathcal{A}_{s, \cutvec} + 1)^{-1} \in \mathscr{L}(L^2(\R))$ for any $s \in \R$. Contrary to $\mathcal{A}_{s, \cutvec}$, ${\cal R}_{s, \cutvec} $ has the advantage to be \textit{bounded}, hence allowing us to use the perturbation theorem stated in Lemma \ref{lem:perturbation_spectrum}. Furthermore, the spectra of ${\cal R}_{s, \cutvec} $ and $\mathcal{A}_{s, \cutvec}$ are related by the following characterization:
   \begin{equation}
   \displaystyle
   \label{eq:spec_1D_2D_ligne_indt_s:preuve_3}
   \lambda \in \sigma(\mathcal{A}_{s, \cutvec}) \quad \Longleftrightarrow \quad \frac{1}{\lambda + 1} \in \sigma({\cal R}_{s, \cutvec}).
   \end{equation}
   Therefore, proving Proposition \ref{prop:spec_1D_2D_ligne_indt_s} is equivalent to proving that $\sigma({\cal R}_{s, \cutvec})$ does not depend on $s$. The idea to do so will be to show that $s \mapsto \sigma({\cal R}_{s, \cutvec})$ defines a (uniformly) continuous mapping (with respect to the Hausdorff distance) which is both $1$--periodic and $\delta$--periodic. Since $\delta$ is irrational, it will then follow from Kronecker's approximation theorem for instance that $s \mapsto \sigma({\cal R}_{s, \cutvec})$ is constant.
 
   \vspspe
   \textbf{Step 1: Continuity of $s \mapsto \sigma({\cal R}_{s, \cutvec})$} --- Fix $s, t \in \R$, and let $f \in L^2(\R)$. 
   % It is clear that $f$ also belongs to both $L^2_{\rho_{s, \cutangle}}(\R)$ and $L^2_{\rho_{t, \cutangle}}(\R)$, given that $\rho_p$ is bounded from above and below. Hence, 
   We can introduce the functions $u_s, u_t \in D(\mathcal{A}_{s, \cutvec})\equiv D(\mathcal{A}_{t, \cutvec})$ defined by
   \begin{equation*}
     \displaystyle
     u_s \coloneqq  {\cal R}_{s, \cutvec}\, f \quad \textnormal{and} \quad u_t \coloneqq  {\cal R}_{t, \cutvec}\, f.
   \end{equation*}
   By substracting the variational formulations satisfied by $u_s$ and $u_t$, we obtain
   \begin{multline*}
     \spforall v \in H^1(\R), \quad \int_{\R} \mu_{s, \cutvec}\,\frac{d}{dx}(u_s - u_t)\, \frac{d \overline{v}}{dx} + \rho_{s, \cutvec}\, (u_s - u_t)\, \overline{v} \; dx \\[3pt]= -\int_{\R} (\mu_{s, \cutvec}-\mu_t)\,\frac{d u_t}{dx}\, \frac{d \overline{v}}{dx} -\int_{\R} (\rho_{s, \cutvec} - \rho_{t, \cutvec})\, u_t\, \overline{v} + \int_{\R} (\rho_{s, \cutvec} - \rho_{t, \cutvec})\, f\, \overline{v}.
   \end{multline*}
   % \begin{equation}
   %   \int_{\R} \Big| \frac{d}{dx}(u_s - u_t) \Big|^2 + \rho_{s,\cutangle}\, |u_s - u_t|^2 \; dx = -\int_{\R} (\rho_{s,\cutangle} - \rho_{t,\cutangle})\; u_t\, \overline{(u_s - u_t)} + \int_{\R} (\rho_{s,\cutangle} - \rho_{t,\cutangle})\; f\, \overline{(u_s - u_t)}.
   % \end{equation}
   Now let us choose $v = u_s - u_t$ in this equality. By using the boundedness of $\rho_{s, \cutvec}$ and $\mu_{s, \cutvec}$ from below on the left side as well as the Cauchy-Schwarz inequality on the right side, and by dividing both sides by $\|u_s - u_t\|_{H^1(\R)}$, we obtain the existence of a constant $c_1 > 0$ such that
   \begin{equation*}
     \displaystyle
     c_1\, \|u_s - u_t\|_{H^1(\R)} \leq \left( \|\mu_{s} - \mu_{t, \cutvec}\|_{L^\infty(\R)}+ \|\rho_{s, \cutvec} - \rho_{t, \cutvec}\|_{L^\infty(\R)}\right) \; \|u_t\|_{H^1(\R)} + \|\rho_{s, \cutvec} - \rho_{t, \cutvec}\|_{L^\infty(\R)}\,\|f\|_{L^2(\R)}.
   \end{equation*}
   Furthermore, the fact that $-1 \notin \sigma(\mathcal{A}_{t, \cutvec})$ (or equivalently the application of Lax-Milgram's theorem) leads to the estimate $\|u_t\|_{H^1(\R)} \leq c_2 \|f\|_{L^2(\R)}$, where $c_2$ depends only on $\rho_\pm$. Consequently, there exists a constant $c > 0$ such that% with $c \coloneqq  c_1\, c_2 + 1$, one has
   \begin{equation*}
     \displaystyle
     \big\|({\cal R}_{s, \cutvec} - {\cal R}_{t, \cutvec})\,f\big\|_{H^1(\R)} = \|u_s - u_t\|_{H^1(\R)} \leq c\; \|f\|_{L^2(\R)}\;  \left( \|\mu_{s} - \mu_{t, \cutvec}\|_{L^\infty(\R)}+ \|\rho_{s, \cutvec} - \rho_{t, \cutvec}\|_{L^\infty(\R)}\right).
     \label{eq:spec_1D_2D_ligne_indt_s:preuve_4}
   \end{equation*}
   The advantage of working with the resolvent operator lies in the fact that ${\cal R}_{s, \cutvec} - {\cal R}_{t, \cutvec}$ is \textit{bounded}. This allows one to apply directly Lemma \ref{lem:perturbation_spectrum} to $B \coloneqq  {\cal R}_{t, \cutvec}$ and $T \coloneqq  {\cal R}_{s, \cutvec} - {\cal R}_{t, \cutvec}$, and to derive from the above that
   \[
     \displaystyle
     \spforall s, t \in \R, \; \text{dist}_{\mathscr{H}}\big[\sigma({\cal R}_{s, \cutvec}), \sigma({\cal R}_{t, \cutvec})\big] \leq \|{\cal R}_{s, \cutvec} - {\cal R}_{t, \cutvec}\| \leq c\,  \left( \|\mu_{s} - \mu_{t, \cutvec}\|_{L^\infty(\R)}+ \|\rho_{s, \cutvec} - \rho_{t, \cutvec}\|_{L^\infty(\R)}\right).
   \]
   But, as the functions $\rho_p$ and $\mu_p$ are continuous and $1$--periodic in each direction, it follows from Heine's theorem that they are uniformly continuous. Therefore, the previous estimate implies in particular the (uniform) continuity of the mapping $s \mapsto \sigma({\cal R}_{s, \cutvec})$, which can be expressed as follows:
   \begin{equation}
     \label{eq:seq_continuity_spectrum_resolvent}
     \spforall (s_n)_n \in \R^{\N},\  \spforall s \in \R, \quad |s_n - s| \to 0 \; \Longrightarrow \; \text{dist}_{\mathscr{H}}\big[\sigma({\cal R}_{\theta,s_n}), \sigma({\cal R}_{s, \cutvec})\big] \to 0, \; n \to +\infty.
   \end{equation}
 
   \vspspe
   \textbf{Step 2} --- Since $\mathcal{A}_{s+1, \cutvec} = \mathcal{A}_{s, \cutvec}$, it is obvious that $\sigma(\mathcal{A}_{s+1, \cutvec}) = \sigma(\mathcal{A}_{s, \cutvec})$. Furthermore, according to Lemma \ref{lem:unitary_equivalence}, $\mathcal{A}_{s, \cutvec}$ and $\mathcal{A}_{s+\delta, \cutvec}$ are equivalent for any $s \in \R$. Thus, they have the same spectrum, i.e. $\sigma(\mathcal{A}_{s, \cutvec}) = \sigma(\mathcal{A}_{s+\delta, \cutvec})$. These observations, combined with the link \eqref{eq:spec_1D_2D_ligne_indt_s:preuve_3} between $\sigma(\mathcal{A}_{s, \cutvec})$ and $\sigma({\cal R}_{s, \cutvec})$, implies that $s \mapsto \sigma({\cal R}_{s, \cutvec})$ is both $1$ and $\delta$--periodic, that is,
   \begin{equation}
     \label{eq:spec_1D_2D_ligne_indt_s:preuve_1_2}
     \spforall s \in \R, \quad \spforall (k, \ell) \in \Z \times \N, \quad \sigma({\cal R}_{\theta,s+\ell \delta+k}) = \sigma({\cal R}_{s, \cutvec}).
   \end{equation}
   But since $\delta$ is irrational, Kronecker's theorem states that $\N \delta + \Z$ is dense in $\R$. In particular,
   \begin{equation}
     \displaystyle
     \label{eq:spec_1D_2D_ligne_indt_s:preuve_2}
     \spforall s, t \in \R, \quad \spexists (k_n, \ell_n) \in \big(\Z \times \N\big)^{\N}, \quad \big|\ell_n \delta + k_n + s - t\big| \to 0, \quad n \to +\infty.
   \end{equation}
   Now let $s, t \in \R$, and pick a sequence $(k_n, \ell_n) \in \big(\Z \times \N\big)^{\N}$ such that \eqref{eq:spec_1D_2D_ligne_indt_s:preuve_2} is satisfied. From \eqref{eq:spec_1D_2D_ligne_indt_s:preuve_1_2} and from the continuity \eqref{eq:seq_continuity_spectrum_resolvent} of $s \mapsto \sigma({\cal R}_{s, \cutvec})$, we deduce the following
   \begin{equation*}
     \displaystyle
     \text{dist}_{\mathscr{H}}\big[\sigma({\cal R}_{s, \cutvec}), \sigma({\cal R}_{t, \cutvec})\big] = \text{dist}_{\mathscr{H}}\big[\sigma({\cal R}_{\theta,s + \ell_n \delta + k_n}), \sigma({\cal R}_{t, \cutvec})\big] \to 0, \quad n \to +\infty,
   \end{equation*}
   which implies that $\text{dist}_{\mathscr{H}}\big[\sigma({\cal R}_{s, \cutvec}), \sigma({\cal R}_{t, \cutvec})\big] = 0$, or equivalently that $\sigma({\cal R}_{s, \cutvec}) = \sigma({\cal R}_{t, \cutvec})$.
 \end{dem}

\vspspe
Let us now  study the family of operators ${\cal A}_p(\xi)$ with $\xi\in]-\pi,\pi]$. They are linked to one another as follows.
\begin{lem}\label{lem:unitary_equivalence2}
  Let $\xi\in(-\pi,\pi]$ and $T_\kv$ be the multiplication operator defined by 
  \begin{equation*}
    \spforall \kv \coloneqq  (k_1,k_2)\in\Z^2, \quad  \spforall u\in L^2((0,1)^2),\quad T_\kv\; u(\yv) \coloneqq  \euler^{2\icplx\pi \kv\cdot\yv} u(\yv),\quad \aeforall \yv \in (0,1)^2.  
  \end{equation*}
  Then we have the identity
  \begin{equation}
    \displaystyle
   (T_\kv)^{-1}\, {\cal A}_p(\xi)\, T_\kv = {\cal A}_p(\xi+2\pi(k_1\delta+k_2)).
  \end{equation}
\end{lem}

\begin{dem}
   Let us remark that if $u\in D({\cal A}_p(\xi))$ then $T_\kv u\in D({\cal A}_p(\xi))$ where we recall that $D({\cal A}_p(\xi))$ is independent of $\xi$. Moreover, for any $u\in D({\cal A}_p(\xi))$
   \[
       (\Dt{}+\icplx \vts \xi\vts \theta_2)(T_\kv u)(\yv)=\left((2\icplx\pi\kv\cdot\cutvec+\icplx \vts \xi\vts \theta_2)\, u(\yv)+\Dt{} u(\yv)\right)\, \euler^{2\icplx\pi\kv\cdot \yv},\quad \aeforall \yv \in (0,1)^2
   \]
    which implies that % $u\in D({\cal A}_p(\xi))$
   \[
        {\cal A}_p(\xi)\, (T_\kv u) = T_\kv\, \big({\cal A}_p(\xi+2\pi(k_1\delta+k_2))\, u\big),
   \]
  which is the expected identity.
\end{dem}

\vspspe
We can now show the following result, namely Proposition \ref{thm:spectrunqp_k}.
\begin{prop}\label{prop:spec_1D_2D_ligne_indt_xi}
  If $\delta$ is irrational, then the spectrum of ${\cal A}_p(\xi)$ does not depend on $\xi$.
\end{prop}

\begin{dem}
   The proof is really similar to the one of Proposition \ref{prop:spec_1D_2D_ligne_indt_xi}. We first show that
   \[
       \xi\mapsto {\cal R}_p(\xi)\coloneqq ({\cal A}_p(\xi)+1)^{-1} \text{ is continuous in operator norm,}
   \]
   which implies that $\xi \mapsto \sigma({\cal A}_p(\xi))$ is continuous for the Hausdorff distance. Lemma \ref{lem:unitary_equivalence2} implies that 
   \[
       \spforall \xi \in (-\pi,\pi],\quad \spforall (k_1,k_2)\in\Z^2,\quad \sigma({\cal A}_p(\xi)) = \sigma\big({\cal A}_p(\xi+2\pi(k_1\delta+k_2))\big).
   \]
   Since $\delta\in\R\setminus\Q$, it then follows that $\xi\mapsto  \sigma\big({\cal A}_p(\xi)\big)$ is independent of $\xi$.
 \end{dem}
%  Finally, we show the following result.
%  \begin{proposition}{}{}
%  \[
%    \sigma({\cal A}_\theta)=\sigma({\cal A}_p)
%  \]
%  \end{proposition}
%  \begin{dem}
%   By adapting the proof of Floquet theorem, one shows that $\sigma({\cal A}_p) = \bigcup{s \in [0, 1]} \sigma(\mathcal{A}_{\cutvec, s})$. However, since $\sigma(\mathcal{A}_{\cutvec, s})$ is independent of $s$ the result follows. 
%  \end{dem}

% % ====================================================================
% \stopcontents[chapters]                               % END OF CHAPTER
% % ====================================================================
%%%%%%%%%%%%%%%%%%%%%%%%%%%%%%%%%%%%%%%%%%%%%%%%%%%%%%%%%
%%      Bibliographie
\printbibliography%[title={Bibliographie}]
\end{document}